# SPECIAL FUZZY MATRICES FOR SOCIAL SCIENTISTS

**W. B. Vasantha Kandasamy**
**Florentin Smarandache**
**K. Ilanthenral**

# SPECIAL FUZZY MATRICES FOR SOCIAL SCIENTISTS

**W. B. Vasantha Kandasamy**
e-mail: vasanthakandasamy@gmail.com
web: http://mat.iitm.ac.in/~wbv
www.vasantha.net

**Florentin Smarandache**
e-mail: smarand@unm.edu

**K. Ilanthenral**
e-mail: ilanthenral@gmail.com

# CONTENTS

# PREFACE

This book is a continuation of the book, "Elementary fuzzy matrix and fuzzy models for socio-scientists" by the same authors. This book is a little advanced because we introduce a multi-expert fuzzy and neutrosophic models. It mainly tries to help social scientists to analyze any problem in which they need multi-expert systems with multi-models.

To cater to this need, we have introduced new classes of fuzzy and neutrosophic special matrices. The first chapter is essentially spent on introducing the new notion of different types of special fuzzy and neutrosophic matrices, and the simple operations on them which are needed in the working of these multi expert models.

In the second chapter, new set of multi expert models are introduced; these special fuzzy models and special fuzzy neutrosophic models that can cater to adopt any number of experts. The working of the model is also explained by illustrative examples.

However, these special fuzzy models can also be used by applied mathematicians to study social and psychological problems. These models can also be used by doctors, engineers, scientists and statisticians. The SFCM, SMFCM, SNCM, SMNCM, SFRM, SNRM, SMFRM, SMNRM, SFNCMs, SFNRMs, etc. can give the special hidden pattern for any given special input vector.

The working of these SFREs, SMFREs and their neutrosophic analogues depends heavily upon the problems and the experts' expectation. The authors have given a long list for further reading which may help the socio scientists to know more about SFRE and SMFREs.

We thank Dr. K. Kandasamy and Meena, without their unflinching support, this book would have never been possible.

W.B.VASANTHA KANDASAMY
FLORENTIN SMARANDACHE
ILANTHENRAL. K

**Chapter One**

# A NEW CLASS OF SPECIAL FUZZY MATRICES AND SPECIAL NEUTROSOPHIC MATRICES

In this chapter for the first time we introduce some new classes of special fuzzy matrices and illustrate them with examples. Also we give the main type of operations carried out on them. All these special fuzzy matrices will be used in the special fuzzy models which will be introduced in chapter two of this book. This chapter has three sections. In sections one we introduce the notion of fuzzy matrices and give the operations used on them like min max operations or max min operations. In section two we introduce the new classes of special fuzzy matrices define special operations on them and illustrate them with examples. In section three neutrosophic matrices, special operations on them are introduced and described. Several illustrative examples are given to make the operations explicit.

## 1.1 Introduction to Fuzzy Matrices

Here we just recall the definition of fuzzy matrices for more about these concepts please refer [106]. Throughout this book the unit interval [0, 1] denotes the fuzzy interval. However in certain fuzzy matrices we also include the interval [–1, 1] to be the fuzzy interval. So any element $a_{ij} \in [-1, 1]$ can be positive or negative. If $a_{ij}$ is positive then $0 < a_{ij} \leq 1$, if $a_{ij}$ is negative then $-1 \leq a_{ij} \leq 0$; $a_{ij} = 0$ can also occur. So [0, 1] or [–1, 1] will be known as fuzzy interval.

Thus if $A = (a_{ij})$ is a matrix and if in particular $a_{ij} \in [0, 1]$ (or [–1, 1]) we call A to be a fuzzy matrix. So all fuzzy matrices are matrices but every matrix in general need not be a fuzzy matrix. Hence fuzzy matrices forms a subclass of matrices.

Now we give some examples of fuzzy matrices.

***Example 1.1.1:*** Let

$$A = \begin{bmatrix} 0.3 & 0.1 & 0.4 & 1 \\ 0.2 & 1 & 0.7 & 0 \\ 0.9 & 0.8 & 0.5 & 0.8 \end{bmatrix}$$

be a matrix. Every element in A is in the unit interval [0, 1]. Thus A is a fuzzy matrix. We can say A is a 3 × 4 rectangular fuzzy matrix.

***Example 1.1.2:*** Let

$$B = \begin{bmatrix} 0.1 & 1 & 0 & 0.3 & 0.6 & 0.2 \\ 1 & 0.5 & 1 & 0.8 & 0.9 & 1 \\ 0 & 1 & 0.3 & 0.9 & 0.7 & 0.5 \\ 0.4 & 0 & 1 & 0.6 & 0.3 & 1 \\ 1 & 0.8 & 0.9 & 0.3 & 0.8 & 0 \end{bmatrix}$$

be a fuzzy matrix B is a 5 × 5 square fuzzy matrix. It is clear all the entries of B are from the unit interval [0, 1].

***Example 1.1.3:*** Consider the matrix

$$C = \begin{bmatrix} 1 & 0 & 1 \\ -1 & 1 & -1 \\ 0 & -1 & 1 \end{bmatrix}$$

C is a matrix, its entries are from the set {–1, 0,1}. C is a 3 × 3 square fuzzy matrix.

***Example 1.1.4:*** Let A = [0 0.3 0.1 0.5 1 0.8 0.9 1 0]. A is a 1 × 9 fuzzy matrix will also be known as the fuzzy row vector or row fuzzy matrix.

***Example 1.1.5:*** Let

$$T = \begin{bmatrix} 1 \\ 0.3 \\ 0.7 \\ 0.5 \\ 0.1 \\ 0.9 \end{bmatrix}$$

be 6 × 1 fuzzy matrix. T is also know as the fuzzy column vector or fuzzy column matrix. Thus if A = [$a_1$ $a_2$ … $a_n$] where $a_i$ ∈ [0, 1]; 1 ≤ i ≤ n, A will be known as the fuzzy row matrix or the fuzzy row vector.

Let

$$B = \begin{bmatrix} b_1 \\ b_2 \\ b_3 \\ \vdots \\ b_m \end{bmatrix}$$

where $b_j$ ∈ [0, 1], 1 ≤ j ≤ m, B will be known as the fuzzy column vector or the fuzzy column matrix.

Let A = ($a_{ij}$) with $a_{ij}$ ∈ [0, 1], 1 ≤ i ≤ n and 1 ≤ j ≤ n. A will be known as the n × n fuzzy square matrix. Suppose C = ($c_{ij}$) with $c_{ij}$ ∈ [0, 1]; 1 ≤ i ≤ n and 1 ≤ j ≤ m then C is known as the n × m rectangular fuzzy matrix. We have seen examples of these types of fuzzy matrices. A = [0 0 … 0] will be know as the zero fuzzy row vector, B = $\begin{bmatrix} 0 \\ 0 \\ \vdots \\ 0 \end{bmatrix}$ will be known as the zero fuzzy column vector and X = [1 1 … 1] is the fuzzy row unit vector and Y = $\begin{bmatrix} 1 \\ 1 \\ \vdots \\ 1 \end{bmatrix}$ will be known as the fuzzy column unit vector. Thus the unit fuzzy row vector and unit row vector are one and the same so is the zero fuzzy row vector and zero fuzzy column vector they are identical with the zero row vector and the zero column vector respectively.

Now we see the usual matrix addition of fuzzy matrices in general does not give a fuzzy matrix.

This is clearly evident from the following example.

***Example 1.1.6:*** Consider the 3 × 3 fuzzy matrices A and B, where

$$A = \begin{bmatrix} 0.3 & 0.8 & 1 \\ 1 & 0.8 & 0 \\ 0.9 & 0.5 & 0.7 \end{bmatrix}$$

and

$$B = \begin{bmatrix} 0.5 & 0.7 & 0.3 \\ 0.2 & 0.1 & 0.6 \\ 0.3 & 0.8 & 0.9 \end{bmatrix}.$$

Now under usual matrix addition

$$A + B = \begin{bmatrix} 0.3 & 0.8 & 1 \\ 1 & 0.8 & 0 \\ 0.9 & 0.5 & 0.7 \end{bmatrix} + \begin{bmatrix} 0.5 & 0.7 & 0.3 \\ 0.2 & 0.1 & 0.6 \\ 0.3 & 0.8 & 0.9 \end{bmatrix}$$

$$= \begin{bmatrix} 0.8 & 1.5 & 1.3 \\ 1.2 & 0.9 & 0.6 \\ 1.2 & 1.3 & 1.6 \end{bmatrix}.$$

Clearly all entries in A + B are not in [0, 1]. Thus A + B is only a 3 × 3 matrix and is not a 3 × 3 fuzzy matrix.

On similar lines we see the product of two fuzzy matrices under usual matrix multiplication in general does not lead to a fuzzy matrix. This is evident from the following example.

***Example 1.1.7:*** Let

$$A = \begin{bmatrix} 0.8 & 0.9 \\ 1 & 0.3 \end{bmatrix}$$

and

$$B = \begin{bmatrix} 1 & 0.8 \\ 0 & 0.9 \end{bmatrix}$$

be two 2 × 2 fuzzy matrices.
Now the product of matrices A with B is given by

$$A \times B = \begin{bmatrix} 0.8 & 0.9 \\ 1 & 0.3 \end{bmatrix} \begin{bmatrix} 1 & 0.8 \\ 0 & 0.9 \end{bmatrix}$$

$$= \begin{bmatrix} 0.8+0 & 0.64+0.81 \\ 1+0 & 0.8+0.27 \end{bmatrix}$$

$$= \begin{bmatrix} 0.8 & 1.45 \\ 1 & 1.07 \end{bmatrix}.$$

We see all the entries in A × B which will also be denoted by AB are not in [0, 1] so AB is not a fuzzy matrix. Thus under the usual multiplication the product of two fuzzy matrices in general need not yield a fuzzy matrix. So we are forced to find different operations on fuzzy matrices so that under those operations we get the resultant to be also a fuzzy matrix.

Let A = ($a_{ij}$) and B = ($b_{ij}$) be any two m × n fuzzy matrices; define Max (A, B) = (Max ($a_{ij}$, $b_{ij}$)), $a_{ij} \in A$ and $b_{ij} \in B$, $1 \leq i \leq m$ and $1 \leq j \leq n$. Then Max (A, B) is a fuzzy matrix. This operation will be known as Max operation.

***Example 1.1.8:*** Let

$$A = \begin{bmatrix} 0.8 & 1 & 0 & 0.3 \\ 0.3 & 0.2 & 0.4 & 1 \\ 0.1 & 0 & 0.7 & 0.8 \end{bmatrix}$$

and

$$B = \begin{bmatrix} 0.9 & 0.8 & 0.7 & 0 \\ 0.1 & 1 & 0 & 0.3 \\ 0.2 & 0.5 & 0.5 & 0.8 \end{bmatrix}$$

be any two 3 × 4 fuzzy matrices Max {A, B} =

$$\begin{bmatrix} \text{Max}(0.8, 0.9) & \text{Max}(1, 0.8) & \text{Max}(0, 0.7) & \text{Max}(0.3, 0) \\ \text{Max}(0.3, 0.1) & \text{Max}(0.2, 1) & \text{Max}(0.4, 0) & \text{Max}(1, 0.3) \\ \text{Max}(0.1, 0.2) & \text{Max}(0, 0.5) & \text{Max}(0.7, 0.5) & \text{Max}(0.8, 0.8) \end{bmatrix}$$

$$= \begin{bmatrix} 0.9 & 1 & 0.7 & 0.3 \\ 0.3 & 1 & 0.4 & 1 \\ 0.2 & 0.5 & 0.7 & 0.8 \end{bmatrix};$$

clearly Max (A, B) is again a fuzzy matrix as every entry given by Max (A, B) belongs to the interval [0, 1].

It is interesting to note Max (A, A) = A and Max ((0), A) = A where (0) is the zero matrix of the same order as that of A.

Now we proceed on to define yet another operation on fuzzy matrices. Let A and B be any two m × n fuzzy matrices

Min (A, B) = (Min ($a_{ij}$, $b_{ij}$)) where A = ($a_{ij}$) and B = ($b_{ij}$), $1 \leq i < m$ and $1 \leq j \leq n$.

We illustrate this Min operation on two fuzzy matrices by the following example.

***Example 1.1.9:*** Let

$$A = \begin{bmatrix} 0.3 & 1 & 0.8 \\ 1 & 0.3 & 0.9 \\ 0 & 0.8 & 0.3 \\ 0.7 & 0.2 & 1 \\ 1 & 0 & 0.8 \end{bmatrix}$$

$$B = \begin{bmatrix} 1 & 0.8 & 0 \\ 0.3 & 0.2 & 0.5 \\ 0.1 & 1 & 0.5 \\ 1 & 0.3 & 0 \\ 0.5 & 0 & 1 \end{bmatrix}$$

be any two 5 × 3 fuzzy matrices.
Now

$$\text{Min (A, B)} = \begin{bmatrix} \text{Min}(0.3,1) & \text{Min}(1,0.8) & \text{Min}(0.8,0) \\ \text{Min}(1,0.3) & \text{Min}(0.3,0.2) & \text{Min}(0.9,0.5) \\ \text{Min}(0,0.1) & \text{Min}(0.8,1) & \text{Min}(0.3,0.5) \\ \text{Min}(0.7,1) & \text{Min}(0.2,0.3) & \text{Min}(1,0) \\ \text{Min}(1,0.5) & \text{Min}(0,0) & \text{Min}(0.8,1) \end{bmatrix}$$

$$= \begin{bmatrix} 0.3 & 0.8 & 0 \\ 0.3 & 0.2 & 0.5 \\ 0 & 0.8 & 0.3 \\ 0.7 & 0.2 & 0 \\ 0.5 & 0 & 0.8 \end{bmatrix}.$$

Clearly Min(A, B) is a fuzzy matrix as all entries in Min(A, B) belong to the unit interval [0, 1].

Now it is interesting to note that Min(A, A) =A where as Min (A, (0)) = (0).

Further we see Min (A, B) = Min(B, A) and Max(A, B) = Max(B, A).

We can have other types of operations on the matrices A and B called max min operation or min max operations. Let P and Q be two fuzzy matrices where P = $(p_{ik})$ be a m × n matrix; $1 \leq i \leq m$ and $1 \leq k \leq n$, Q = $(q_{kj})$ be a n × t matrix where $1 \leq k \leq n$ and $1 \leq j \leq t$; then max min operations of P and Q is given by R = $(r_{ij})$ = max min $(p_{ik}, q_{kj})$, where $1 \leq i \leq m$ and $1 \leq j \leq t$ and R is a m × t matrix. Clearly R = $(r_{ij})$ is a fuzzy matrix.

We illustrate this by the following example.

***Example 1.1.10:*** Let

$$A = \begin{bmatrix} 0.3 & 0.1 & 0.6 \\ 0 & 0.7 & 1 \\ 0.4 & 0.2 & 0.3 \end{bmatrix}$$

and

$$B = \begin{bmatrix} 0.6 & 0.2 & 0 & 0.7 \\ 0.3 & 0.8 & 0.2 & 0 \\ 1 & 0.1 & 0.4 & 1 \end{bmatrix}$$

be any two fuzzy matrices where A is a 3 × 3 fuzzy matrix and B is a 3 × 4 fuzzy matrix, max min (A, B) is given by R:

$$R = \max\min\left\{\begin{bmatrix} 0.3 & 0.1 & 0.6 \\ 0 & 0.7 & 1 \\ 0.4 & 0.2 & 0.3 \end{bmatrix}, \begin{bmatrix} 0.6 & 0.2 & 0 & 0.7 \\ 0.3 & 0.8 & 0.2 & 0 \\ 1 & 0.1 & 0.4 & 1 \end{bmatrix}\right\}$$

$$R = \begin{bmatrix} r_{11} & r_{12} & r_{13} & r_{14} \\ r_{21} & r_{22} & r_{23} & r_{24} \\ r_{31} & r_{32} & r_{33} & r_{34} \end{bmatrix};$$

where

$$
\begin{aligned}
r_{11} &= \max \{\min(0.3, 0.6), \min (0.1, 0.3), \min (0.6, 1)\} \\
&= \max \{0.3, 0.1, 0.6\} \\
&= 0.6. \\
r_{12} &= \max \{\min (0.3, 0.2), \min (0.1, 0.8), \min (0.6, 0.1)\} \\
&= \max \{0.2, 0.1, 0.1\} \\
&= 0.2. \\
r_{13} &= \max \{\min (0.3, 0), \min (0.1, 0.2), \min (0.6, 0.4)\} \\
&= \max \{0, 0.1, 0.4\} \\
&= 0.4. \\
r_{14} &= \max \{\min (0.3, 0.7), \min (0.1, 0), \min (0.6, 1)\} \\
&= \max \{0.3, 0, 0.6\} \\
&= 0.6. \\
r_{21} &= \max \{\min (0, 0.6), \min (0.7, 0.3), \min (1, 1)\} \\
&= \max \{0, 0.3, 1\} \\
&= 1;
\end{aligned}
$$

and so on and

$$
\begin{aligned}
r_{34} &= \max \{\min (0.4, 0.7), \min (0.2, 0), \min (0.3, 1)\} \\
&= \max \{0.4, 0, 0.3\} \\
&= 0.4.
\end{aligned}
$$

$$
R = \begin{bmatrix} 0.6 & 0.2 & 0.4 & 0.6 \\ 1 & 0.7 & 0.4 & 1 \\ 0.4 & 0.2 & 0.3 & 0.4 \end{bmatrix}
$$

$$= (r_{ij}).$$

Likewise we can also define the notion of min max operator on A, B as follows; if $A = (a_{ik})$ be a $m \times n$ matrix with $1 \leq i \leq m$ and $1 \leq k \leq n$ and $B = (b_{kj})$ be a $n \times t$ matrix where $1 \leq k \leq n$ and $1 \leq j \leq t$, then $C = (c_{ij}) = \min \max \{a_{ik}, b_{kj}\}$, where $1 \leq i \leq m$ and $1 \leq j \leq t$,

We illustrate this by the following example. Further it is pertinent to mention here that in general

$$\max_k \min \{a_{ik}, b_{kj}\} \neq \min_k \max \{a_{ik}, b_{kj}\}.$$

We find the min max (A, B) for the same A and B given in example 1.1.10.

***Example 1.1.11:*** Let

$$A = \begin{bmatrix} 0.3 & 0.1 & 0.6 \\ 0 & 0.7 & 1 \\ 0.4 & 0.2 & 0.3 \end{bmatrix}$$

and

$$B = \begin{bmatrix} 0.6 & 0.2 & 0 & 0.7 \\ 0.3 & 0.8 & 0.2 & 0 \\ 1 & 0.1 & 0.4 & 1 \end{bmatrix}$$

be the fuzzy matrices given in example 1.1.10.
Now

min max (A, B)

$$= \min\{\max \left\{ \begin{bmatrix} 0.3 & 0.1 & 0.6 \\ 0 & 0.7 & 1 \\ 0.4 & 0.2 & 0.3 \end{bmatrix}, \begin{bmatrix} 0.6 & 0.2 & 0 & 0.7 \\ 0.3 & 0.8 & 0.2 & 0 \\ 1 & 0.1 & 0.4 & 1 \end{bmatrix} \right\}$$

$$= \begin{bmatrix} p_{11} & p_{12} & p_{13} & p_{14} \\ p_{21} & p_{22} & p_{23} & p_{24} \\ p_{31} & p_{32} & p_{33} & p_{34} \end{bmatrix} = P$$

$p_{11}$ = min {max (0.3, 0.6), max (0.1, 0.3), max (0.6, 1)}
= min {0.6, 0.3, 1}
= 0.3.
$p_{12}$ = min {max (0.3, 0.2), max (0.1, 0.8), max (0.6, 1)}
= min {0.3, 0.8, 1}
= 0.3.
$p_{13}$ = min {max (0.3, 0), max (0.1, 0.2), max (0.6, 0.4)}

$= \min \{0.3, 0.2, 0.6\}$
$= 0.2.$
$p_{14} = \min \{\max (0.3, 0.7), \max (0.1, 0), \max (0.6, 1)\}$
$= \min \{0.7, 0.1, 1\}$
$= 0.1.$
$p_{21} = \min \{\max (0, 0.6), \max (0.7, 0.3), \max (1, 1)\}$
$= \min \{0.6, 0.7, 1\}$
$= 0.6$
and so on and
$p_{34} = \min \{\max (0.4, 0.7), \max (0.2, 0), \max (0.3, 1)\}$
$= \min \{0.7, 0.2, 1\}$
$= 0.2.$

Thus we get

$$P = \begin{bmatrix} 0.3 & 0.3 & 0.2 & 0.1 \\ 0.6 & 0.2 & 0 & 0.7 \\ 0.3 & 0.3 & 0.2 & 0.2 \end{bmatrix}.$$

We see $\max_k \{\min(a_{ik}, b_{kj})\} \neq \min_k \{\max (a_{ik}, b_{kj})\}$ from the examples 1.1.10 and 1.1.11.

**Note:** It is important to note that the expert who works with the fuzzy models may choose to have the max min operator or the min max operator. For we see in general the resultant given by max min operator will be always greater than or equal to the min max operator. That is in other words we can say that min max operator in general will yield a value which will always be less than or equal to the max min operator.

Now we have to answer the question will max min operator be always defined for any two fuzzy matrices. We see the max min operator or the min max operator will not in general be defined for any fuzzy matrices A and B. The max min operator {A, B} or the min max operator {A, B} will be defined if and only if the number of columns in A equal to the number of rows of B otherwise max min or min max of A, B will not be defined. Further max min operator of A, B in general will not be equal to the max min operator of B, A. In general even if max min operator of A, B be defined than max min operator of B, A may not be even defined. This is true of min max operator also we see max min of A, B is defined in the example 1.1.10 but max min {B, A} cannot be found at all i.e., max min {A B} is not defined.

Similarly in example 1.1.11 we see min max of A, B is defined where as min max {B, A} is not defined. Thus we can say that for any two fuzzy matrices A and B, in general max min (A, B) ≠ max min (B, A) (may not be even compatible). Like wise min max {A, B} ≠ min max (B, A) (may not be defined or compatible at all).

Now we see even if both max min (A, B) and max min (B, A) are defined they may not be equal in general.

For we will give an example to this effect.

***Example 1.1.12:*** Let

$$A = \begin{bmatrix} 0.3 & 0.1 & 1 \\ 0.6 & 0.3 & 0.8 \\ 0 & 0.4 & 0.5 \end{bmatrix}$$

and

$$B = \begin{bmatrix} 1 & 0.4 & 0.9 \\ 0.8 & 0.6 & 0.2 \\ 0.7 & 0.4 & 1 \end{bmatrix}$$

be any two 3 × 3 fuzzy matrices. Now

max min (A, B)

$$= \max \left\{ \min \left\{ \begin{bmatrix} 0.3 & 0.1 & 1 \\ 0.6 & 0.3 & 0.8 \\ 0 & 0.4 & 0.5 \end{bmatrix}, \begin{bmatrix} 1 & 0.4 & 0.9 \\ 0.8 & 0.6 & 0.2 \\ 0.7 & 0.4 & 1 \end{bmatrix} \right\} \right\}$$

$$= \begin{bmatrix} 0.7 & 0.4 & 1 \\ 0.7 & 0.4 & 0.8 \\ 0.5 & 0.4 & 0.5 \end{bmatrix}.$$

max min (B, A)

$$= \max\left\{\min\left\{\begin{bmatrix} 1 & 0.4 & 0.9 \\ 0.8 & 0.6 & 0.2 \\ 0.7 & 0.4 & 1 \end{bmatrix}, \begin{bmatrix} 0.3 & 0.1 & 1 \\ 0.6 & 0.3 & 0.8 \\ 0 & 0.4 & 0.5 \end{bmatrix}\right\}\right\}$$

$$= \begin{bmatrix} 0.4 & 0.4 & 1 \\ 0.6 & 0.3 & 0.8 \\ 0.4 & 0.4 & 0.7 \end{bmatrix}.$$

We see max min (A, B) $\neq$ max (min) (B, A).

$$\begin{bmatrix} 0.7 & 0.4 & 1 \\ 0.7 & 0.4 & 0.8 \\ 0.5 & 0.4 & 0.5 \end{bmatrix} \neq \begin{bmatrix} 0.4 & 0.4 & 1 \\ 0.6 & 0.3 & 0.8 \\ 0.4 & 0.4 & 0.7 \end{bmatrix}.$$

We see in general min max (A, B) $\neq$ min max (A, B). We illustrate this by the following.

***Example 1.1.13:*** Let

$$A = \begin{bmatrix} 1 & 0.3 & 0.2 \\ 0.4 & 1 & 0.5 \\ 0.7 & 0.3 & 1 \end{bmatrix}$$

and

$$B = \begin{bmatrix} 0.3 & 1 & 0.8 \\ 0.7 & 0.7 & 1 \\ 1 & 0.6 & 0.3 \end{bmatrix}$$

be any two fuzzy matrices.
min max (A, B)

$$= \min\left\{\max\left\{\begin{bmatrix} 1 & 0.3 & 0.2 \\ 0.4 & 1 & 0.5 \\ 0.7 & 0.3 & 1 \end{bmatrix}, \begin{bmatrix} 0.3 & 1 & 0.8 \\ 0.7 & 0.7 & 1 \\ 1 & 0.6 & 0.3 \end{bmatrix}\right\}\right\}$$

$$= \begin{bmatrix} 0.7 & 0.6 & 0.3 \\ 0.4 & 0.6 & 0.5 \\ 0.7 & 0.7 & 0.8 \end{bmatrix}.$$

Now we find
min max (B, A)

$$= \min \left\{ \max \left\{ \begin{bmatrix} 0.3 & 1 & 0.8 \\ 0.7 & 0.7 & 1 \\ 1 & 0.6 & 0.3 \end{bmatrix}, \begin{bmatrix} 1 & 0.3 & 0.2 \\ 0.4 & 1 & 0.5 \\ 0.7 & 0.3 & 1 \end{bmatrix} \right\} \right\}$$

$$= \begin{bmatrix} 0.8 & 0.3 & 0.3 \\ 0.7 & 0.7 & 0.7 \\ 0.6 & 0.3 & 0.6 \end{bmatrix}.$$

We see clearly min max (A, B) ≠ min max (B, A).

$$\text{i.e. } \begin{bmatrix} 0.7 & 0.6 & 0.3 \\ 0.4 & 0.6 & 0.5 \\ 0.7 & 0.7 & 0.8 \end{bmatrix} \neq \begin{bmatrix} 0.8 & 0.3 & 0.3 \\ 0.7 & 0.7 & 0.7 \\ 0.6 & 0.3 & 0.6 \end{bmatrix}.$$

Now we also work with special type of products using fuzzy matrices. Suppose we have a n × n square fuzzy matrix $A = (a_{ij})$; $a_{ij} \in \{-1, 0, 1\}$ and $1 \leq i, j \leq n$.

Let $X = [x_1\ x_2 \ldots x_n]$ be a fuzzy row vector where $x_i \in \{0, 1\}$, $1 \leq i \leq n$. Now we find the product $X \text{ o } A = [y_1\ y_2 \ldots y_n]$ where the product X o A is the usual matrix product. Now $y_i$'s $1 \leq i \leq n$ need not in general belong to {0, 1}, we make the following operations so that $y_i$'s belong to {0, 1} which will be known as the updating operation and thresholding operation.

We define if $y_i > 0$ replace them by 1, $1 \leq i \leq n$ and if $y_i$'s are $\leq 0$ then replace $y_i$ by 0. Now this operation is known as the thresholding operation. It is important to mention here that in general we can replace $y_i > 0$ by a value say $a > 0$, $a \in [0, 1]$ and $y_i \leq 0$ by 0. It is left to the choice of the expert who works with it. Now the updating operation is dependent on the X = [$x_1$ … $x_n$] with which we started to work with. If $x_k$ in X was 1 in 1 ≤ k ≤ n; then we demand in the resultant Y = X o A = [$y_1$ $y_2$ … $y_n$] the $y_k$ must be 1, 1 ≤ k ≤ n, this operation on the resultant Y is known as the updating operation.

Now using Y' which is the updated and thresholded fuzzy row vector of Y = X o A we find Y' o A = Z (say) Z = [$z_1$ $z_2$ … $z_n$]. We now threshold and update Z to get Z' this process is repeated till we arrive at a fixed point or a limit cycle i.e. if X o A = Y, Y' after updating and thresholding Y we find Y' o A = Z, if Z' is got after updating and thresholding Z we proceed on till some R' o A = S and S' is Y' or Z' and this process will repeat, then S' will be known as the limit cycle given by the fuzzy row vector X using the fuzzy matrix A. [108, 112]

Now we illustrate this by the following example.

***Example 1.1.14:*** Let A= ($a_{ij}$) be a 5 × 5 fuzzy matrix with $a_{ij}$ ∈ {–1, 0, 1},

$$\text{i.e., } A = \begin{bmatrix} 0 & 1 & 0 & -1 & 0 \\ -1 & 0 & -1 & 0 & 1 \\ 0 & -1 & 0 & 1 & -1 \\ 0 & 0 & 0 & 0 & 1 \\ 1 & -1 & 1 & 0 & 0 \end{bmatrix}.$$

Let X = [0 1 0 0 1] be the fuzzy row vector (given), to find the resultant of X using A.

$$X \text{ o } A = \begin{bmatrix} 0 & 1 & 0 & 0 & 1 \end{bmatrix} \text{ o } \begin{bmatrix} 0 & 1 & 0 & -1 & 0 \\ -1 & 0 & -1 & 0 & 1 \\ 0 & -1 & 0 & 1 & -1 \\ 0 & 0 & 0 & 0 & 1 \\ 1 & -1 & 1 & 0 & 0 \end{bmatrix}$$

= [0 –1 0 0 1]
= Y.

Now

$Y' = [0\ 1\ 0\ 0\ 1]$

which is got after updating and thresholding Y. Now

$$\begin{aligned} Y' \text{ o } A &= [0\ -1\ 0\ 0\ 1] \\ &= Z. \end{aligned}$$

Let Z' be the vector got after thresholding and updating Z and we see Z' =Y'. Thus it gives us the fixed point. Suppose X = [1 0 1 0 1] to find the effect of X on A.

$$X \text{ o } A = \begin{bmatrix} 1 & 0 & 1 & 0 & 1 \end{bmatrix} \text{ o } \begin{bmatrix} 0 & 1 & 0 & -1 & 0 \\ -1 & 0 & -1 & 0 & 1 \\ 0 & -1 & 0 & 1 & -1 \\ 0 & 0 & 0 & 0 & 1 \\ 1 & -1 & 1 & 0 & 0 \end{bmatrix}$$

$$\begin{aligned} &= [1\ -1\ 1\ 0\ -1] \\ &= Y. \end{aligned}$$

Y' got by updating and thresholding Y is got as [1 0 1 0 1], which shows X remains unaffected by this product.

Let X = [0 0 1 0 0] to find the effect of X on A.

$$X \text{ o } A = \begin{bmatrix} 0 & 0 & 1 & 0 & 0 \end{bmatrix} \text{ o } \begin{bmatrix} 0 & 1 & 0 & -1 & 0 \\ -1 & 0 & -1 & 0 & 1 \\ 0 & -1 & 0 & 1 & -1 \\ 0 & 0 & 0 & 0 & 1 \\ 1 & -1 & 1 & 0 & 0 \end{bmatrix}$$

$$\begin{aligned} &= [0\ -1\ 0\ 1\ -1] \\ &= Y; \end{aligned}$$

after updating and thresholding Y we get Y' = [0 0 1 1 0].

Now

$$Y' \text{ o } A = \begin{bmatrix} 0 & 0 & 1 & 1 & 0 \end{bmatrix} \text{ o } \begin{bmatrix} 0 & 1 & 0 & -1 & 0 \\ -1 & 0 & -1 & 0 & 1 \\ 0 & -1 & 0 & 1 & -1 \\ 0 & 0 & 0 & 0 & 1 \\ 1 & -1 & 1 & 0 & 0 \end{bmatrix}$$

= [0 –1 0 1 0]
= Z.

After updating and thresholding Z we get Z' = [0 0 1 1 0].

Thus the special hidden pattern is a special fixed point given by [0 0 1 1 0].

It is very usual we get back the same X after a sequence of operations such cases one says the special hidden pattern is a special limit cycle. Now we use after set of operations when the given fuzzy matrix is not a square matrix. Let B = ($b_{ij}$) where $b_{ij} \in \{-1, 0, 1\}$ and B is a n × m matrix m ≠ n with $1 \leq i \leq n$ and $1 \leq j \leq m$. Suppose X = [$x_1 \ldots x_n$] be a fuzzy row vector, $x_i \in \{0, 1\}$; $1 \leq i \leq n$. Then X o B = Y = [$y_1 \ldots y_m$] where $y_i \notin \{0, 1\}$ we update and threshold Y to Y'; Y' is a 1 × m fuzzy row vector.

We find Y' o $B^T$ = Z where Z = [$z_1$ $z_2$ … $z_n$]. Let Z' be the fuzzy row vector after thresholding and updating Z. Now we calculate Z' o A = P where P is a 1 × m row matrix and P' is the updated and thresholded resultant of P. Find P' o $B^T$ and so on is continued till one arrives at a fixed binary pair or a limit bicycle.

We illustrate this by the following example.

***Example 1.1.15:*** Let B = ($b_{ij}$) be a 6 × 4 fuzzy matrix, $b_{ij} \in \{-1, 0, 1\}$ where

$$B = \begin{bmatrix} 1 & -1 & 0 & 1 \\ 0 & 1 & 0 & 0 \\ -1 & 0 & 1 & 0 \\ 0 & 0 & 0 & -1 \\ 0 & 1 & 1 & 0 \\ 1 & 1 & -1 & 1 \end{bmatrix}.$$

Suppose X = [1 0 1 0 1 1] be given to find the resultant of X and B.

$$X \text{ o } B = \begin{bmatrix}1 & 0 & 1 & 0 & 1 & 1\end{bmatrix} \text{ o } \begin{bmatrix} 1 & -1 & 0 & 1 \\ 0 & 1 & 0 & 0 \\ -1 & 0 & 1 & 0 \\ 0 & 0 & 0 & -1 \\ 0 & 1 & 1 & 0 \\ 1 & 1 & -1 & 1 \end{bmatrix}$$

= [1 1 1 2]
= Y ;

after thresholding Y we get Y' = [1 1 1 1]. Clearly as Y is not a given fuzzy vector we need not update Y. Now we find

$$Y \text{ o } B^T = \begin{bmatrix}1 & 1 & 1 & 1\end{bmatrix} \text{ o } \begin{bmatrix} 1 & 0 & -1 & 0 & 0 & 1 \\ -1 & 1 & 0 & 0 & 1 & 1 \\ 0 & 0 & 1 & 0 & 1 & -1 \\ 1 & 0 & 0 & -1 & 0 & 1 \end{bmatrix}$$

= [1 1 0 –1 2 2]
= Z.

Now Z' got after updating and thresholding Z is given by Z' = [1 1 1 0 1 1].
We find

$$Z' \text{ o } B = \begin{bmatrix}1 & 1 & 1 & 0 & 1 & 1\end{bmatrix} \text{ o } \begin{bmatrix} 1 & -1 & 0 & 1 \\ 0 & 1 & 0 & 0 \\ -1 & 0 & 1 & 0 \\ 0 & 0 & 0 & -1 \\ 0 & 1 & 1 & 0 \\ 1 & 1 & -1 & 1 \end{bmatrix}$$

= [1 2 1 2]
= P.

P after thresholding we get P' = [1 1 1 1]. We see P' o $B^T$ gives Q and Q' after thresholding and updating is [1 1 1 0 1 1] = Z'. Thus we get the resultant as a fixed binary pair given by {[1 1 1 0 1 1], [1 1 1 1]}.

Now we study the effect of a 1 × 4 row vector T = [1 0 0 1] on B

$$T \text{ o } B^T = \begin{bmatrix} 1 & 0 & 0 & 1 \end{bmatrix} \text{ o } \begin{bmatrix} 1 & 0 & -1 & 0 & 0 & 1 \\ -1 & 1 & 0 & 0 & 1 & 1 \\ 0 & 0 & 1 & 0 & 1 & -1 \\ 1 & 0 & 0 & -1 & 0 & 1 \end{bmatrix}$$

= [2 0 –1 –1 0 2]
= X.

X’ after thresholding X is given by X' = [1 0 0 0 0 1].

Now effect of X' on B is given by

$$X' \text{ o } B = \begin{bmatrix} 1 & 0 & 0 & 0 & 0 & 1 \end{bmatrix} \begin{bmatrix} 1 & -1 & 0 & 1 \\ 0 & 1 & 0 & 0 \\ -1 & 0 & 1 & 0 \\ 0 & 0 & 0 & -1 \\ 0 & 1 & 1 & 0 \\ 1 & 1 & -1 & 1 \end{bmatrix}$$

= [2 0 –1 2]
= Y.

Y' = [1 0 0 1] the updated and thresholded vector which is T.

Thus we get the fixed binary pair in case of T = [1 0 0 1] is given by {[1 0 0 1], [1 0 0 0 0 1]}.

We now proceed on to give yet another operation viz max min operation on fuzzy matrices.

Let A = ($a_{ij}$) be a m × n fuzzy matrix, m = [$m_1$ $m_2$ … $m_n$] where $a_{ij} \in [0, 1]$ and $m_k \in \{0, 1\}$ $1 \leq i \leq m$, $1 \leq j \leq n$ and $1 \leq k \leq n$, we calculate X = A o M that is $x_i = \max \min (a_{ij}, m_j)$ where $1 \leq j \leq n$ and i = 1, 2, …, m.
X o A = Y is calculated using this X.

We explicitly show this by the following example.

***Example 1.1.16:*** Let us consider the fuzzy matrix A = ($a_{ij}$) where A is a 5 × 7 matrix with $a_{ij} \in [0, 1]$.

$$\text{i.e. } A = \begin{bmatrix} 0.1 & 0.3 & 1 & 0.2 & 1 & 0 & 0.8 \\ 0.5 & 1 & 0.5 & 0.6 & 1 & 0.7 & 0.2 \\ 1 & 0.4 & 0.5 & 1 & 0.7 & 0.7 & 1 \\ 0.7 & 0 & 1 & 0.2 & 0.6 & 1 & 0 \\ 0.6 & 0.8 & 0.6 & 0.3 & 1 & 0.2 & 0.3 \end{bmatrix}$$

Let M = [1 0 0 1 0 1 1] be the given fuzzy vector.

A o M = max min {$a_{ij}$, $m_j$}
= max ({0.1 0 0 0.2 0 0 0.8}, {0.5 0 0 0.6 0 0.7 0.2}, {1 0 0 1 0 0.7 1}, {0.7 0 0 0.2 0 1 0}, {0.6 0 0 0.3 0 0.2 0.3})
= [0.8 0.7 1 1 0.6]
= B
B o A = max min {$b_j$ $a_{ji}$}
= [1 0.7 1 1 0.8 1 1]
= $M_1$.

Now we can find $M_1$ o A and so on.

This type of operator will be used in the fuzzy models.

## 1.2 Special Classes of Fuzzy Matrices

In this section we introduce a special classes of fuzzy matrices. We illustrate them with examples. Certain special type of operations are defined on them which will be used in the fuzzy special models which will be constructed in chapter two of this book.

**DEFINITION 1.2.1:** *Let $M = M_1 \cup M_2 \cup \ldots \cup M_n$; ($n \geq 2$) where each $M_i$ is a t × t fuzzy matrix. We call M to be a special fuzzy t* $\times$ *t square matrix i.e.* $M = M_1 \cup M_2 \cup \ldots \cup M_n = (m^1_{ij}) \cup (m^2_{ij}) \cup \ldots \cup (m^n_{ij})$, $1 \le i \le t$ *and* $1 \le j \le t$.

We illustrate this by the following example.

***Example 1.2.1:*** Let $M = M_1 \cup M_2 \cup M_3 \cup M_4$ be a $4 \times 4$ square fuzzy matrix where

$$M_1 = \begin{bmatrix} 1 & 0.2 & 0.7 & 0.9 \\ 0.3 & 1 & 0.2 & 0.5 \\ 0.8 & 0.9 & 1 & 1 \\ 0 & 0.6 & 0.4 & 0 \end{bmatrix}$$

$$M_2 = \begin{bmatrix} 0.3 & 0.4 & 0.5 & 0.6 \\ 1 & 0.7 & 0.9 & 0.5 \\ 0.2 & 0.3 & 0.4 & 1 \\ 0.5 & 1 & 0 & 0.8 \end{bmatrix}$$

$$M_3 = \begin{bmatrix} 0.8 & 1 & 0 & 0.9 \\ 0.3 & 0.2 & 1 & 0.2 \\ 0.6 & 0.7 & 0.5 & 1 \\ 1 & 0 & 0.8 & 0 \end{bmatrix}$$

and

$$M_4 = \begin{bmatrix} 1 & 0 & 1 & 0 \\ -1 & 0.7 & 0.2 & 0.5 \\ 0.3 & 0.8 & 1 & 0.9 \\ 0.1 & 0.6 & 0.7 & 1 \end{bmatrix}.$$

M is a special fuzzy $4 \times 4$ square matrix.

Now we proceed on to define the notion of special fuzzy row matrix / vector and special fuzzy column vector / matrix.

**DEFINITION 1.2.2:** *Let $X = X_1 \cup X_2 \cup \ldots \cup X_M$ ($M \geq 2$) where each $X_i$ is a $1 \times s$ fuzzy row vector / matrix then we define X to be a special fuzzy row vector / matrix ($i = 1, 2, \ldots, M$). If in particular $X = X_1 \cup X_2 \cup \ldots \cup X_M$ ($M \geq 2$) where each $X_i$ is a $1 \times s_i$ fuzzy row vector/matrix where for atleast one $s_i \neq s_j$ with $i \neq j$, $1 \leq i, j \leq M$ then we define X to be a special fuzzy mixed row vector / matrix.*

We illustrate first these two concepts by some examples.

***Example 1.2.2:*** Let

$$\begin{aligned} X &= X_1 \cup X_2 \cup X_3 \cup X_4 \cup X_5 \\ &= [0.1\ 0.2\ 0.1\ 0.5\ 0.7] \cup [1\ 0\ 1\ 0.9\ 1] \cup [0\ 0.1\ 0.8\ 0.9\ 0.4] \\ &\quad \cup [0.8\ 0.6\ 0.4\ 0.5\ 0.7] \cup [0.3\ 0.1\ 0.5\ 0.3\ 0.2] \end{aligned}$$

where each $X_i$ is a $1 \times 6$ fuzzy row vector; i = 1, 2, 3, 4, 5. X is a special fuzzy row vector / matrix.

Now we proceed on to give an example of a special fuzzy mixed row vector / matrix.

***Example 1.2.3:*** Let

$$\begin{aligned} Y &= Y_1 \cup Y_2 \cup Y_3 \cup Y_4 \cup Y_5 \\ &= [1\ 0.3\ 0.2\ 0.9\ 1\ 0.3] \cup [1\ 1\ 1\ 0.3] \cup [0.3\ 1\ 0\ 0.2\ 0.5] \cup \\ &\quad [1\ 0.3\ 0.4] \cup [1\ 0.8\ 0.9\ 0.7\ 0.6\ 0.5] \end{aligned}$$

where $Y_1$ is a $1 \times 6$ fuzzy row vector, $Y_2$ is a $1 \times 4$ fuzzy row vector, $Y_3$ is a $1 \times 5$ fuzzy row vector, $Y_4$ is a $1 \times 3$ fuzzy row vector and $Y_5$ is a $l \times 6$ fuzzy row vector. We see Y is a special fuzzy mixed row vector/matrix.

Now we proceed on to define the notion of special fuzzy column vector/matrix and the notion of special fuzzy mixed column vector/matrix.

**DEFINITION 1.2.3:** *Let $Y = Y_1 \cup Y_2 \cup \ldots \cup Y_m$ ($m \geq 2$) we have each $Y_i$ to be a $t \times 1$ fuzzy column vector/ matrix then we define*

*Y to be a special fuzzy column vector / matrix. If in particular in $Y = Y_1 \cup Y_2 \cup \ldots \cup Y_m$ ($m \geq 2$) we have each $Y_i$ to be $t_i \times 1$ fuzzy column vector where atleast for one or some $t_i \neq t_j$ for $i \neq j$, $1 \leq i, j \leq m$ then we define Y to be a special fuzzy mixed column vector/matrix.*

Now we proceed on to describe these two concepts with examples.

***Example 1.2.4:*** Let us consider

$$Z = Z_1 \cup Z_2 \cup Z_3 \cup Z_4 \cup Z_5$$

$$= \begin{bmatrix} 0.3 \\ 0.8 \\ 0.1 \\ 0.5 \\ 1 \\ 0.7 \end{bmatrix} \cup \begin{bmatrix} 1 \\ 0 \\ 1 \\ 1 \\ 0.8 \\ 0.1 \end{bmatrix} \cup \begin{bmatrix} 0.5 \\ 0.7 \\ 0.6 \\ 0.9 \\ 1 \\ 0 \end{bmatrix} \cup \begin{bmatrix} 0.8 \\ 0.1 \\ 0.3 \\ 0.1 \\ 0.8 \\ 0.4 \end{bmatrix} \cup \begin{bmatrix} 0.5 \\ 0.7 \\ 0.8 \\ 0.6 \\ 0.3 \\ 0.1 \end{bmatrix}$$

where each $Z_i$ is a $6 \times 1$ fuzzy column vector. We see Z is a special fuzzy column vector/matrix, $1 \leq i \leq 5$.

***Example 1.2.5:*** Let

$$T = T_1 \cup T_2 \cup T_3 \cup T_4 \cup T_5 \cup T_6 \cup T_7$$

$$= \begin{bmatrix} 0.3 \\ 1 \\ 0 \\ 0.8 \\ 0.7 \end{bmatrix} \cup \begin{bmatrix} 0.8 \\ 0.71 \\ 0.512 \\ 0.3 \\ 0.031 \\ 0.54 \\ 0.111 \\ 0.14 \end{bmatrix} \cup \begin{bmatrix} 0.91 \\ 0.7 \\ 0.11 \\ 0.44 \\ 0.13 \\ 0.71 \\ 0.8 \end{bmatrix} \cup \begin{bmatrix} 1 \\ 0 \\ 0.31 \end{bmatrix}$$

$$\cup \begin{bmatrix} 1 \\ 0.81 \\ 0.6 \\ 1 \\ 0.9 \end{bmatrix} \cup \begin{bmatrix} 0.9 \\ 0.3 \\ 0.81 \\ 0.116 \\ 0.08 \\ 1 \\ 0 \\ 0.4 \end{bmatrix} \cup \begin{bmatrix} 0.9 \\ 0.8 \\ 0.15 \end{bmatrix};$$

clearly T is a special fuzzy mixed column vector/matrix.

The following facts are important and is to be noted. We see in a special fuzzy mixed column (row) matrix we can have two or more fuzzy column (row) vector to have same number of rows. But in the special fuzzy column (row) matrix we must have each fuzzy column (row) vector should have the same number of rows. Further in case of special fuzzy row (column) vector/matrix we can also have the same fuzzy row (column) vector to repeat itself.

For instance if

$$X = X_1 \cup X_2 \cup X_3 \cup X_4$$

$$= \begin{bmatrix} 1 \\ 0 \\ 1 \\ 1 \end{bmatrix} \cup \begin{bmatrix} 1 \\ 0 \\ 1 \\ 1 \end{bmatrix} \cup \begin{bmatrix} 1 \\ 0 \\ 1 \\ 1 \end{bmatrix} \cup \begin{bmatrix} 0 \\ 1 \\ 1 \\ 1 \end{bmatrix}$$

then also X is a special fuzzy column matrix.
Like wise if

$$\begin{aligned} Y &= Y_1 \cup Y_2 \cup Y_3 \cup Y_4 \cup X_5 \\ &= [1\ 0\ 0\ 0\ 1\ 0.7\ 0.5] \cup [1\ 0\ 0\ 0\ 1\ 0.7\ 0.5] \cup [1\ 0\ 0\ 0\ 1\ 0.7\ 0.5] \cup [1\ 0\ 0\ 0\ 1\ 0.7\ 0.5] \cup [1\ 0\ 0\ 0\ 1\ 0.7\ 0.5] \cup [1\ 0\ 0\ 0\ 1\ 0.7\ 0.5] \cup [1\ 0\ 0\ 0\ 1\ 0.7\ 0.5], \end{aligned}$$

we say Y is a special fuzzy row matrix, each $Y_i$ is $1 \times 7$ fuzzy row vector and we have

$$Y_1 = Y_2 = Y_3 = Y_4 = Y_5 = [1\ 0\ 0\ 0\ 1\ 0.7\ 0.5].$$

Now we proceed on to define the notion of special fuzzy mixed square matrix.

**DEFINITION 1.2.4:** *Let $V = V_1 \cup V_2 \cup \ldots \cup V_n$ $(n \geq 2)$ where each $V_i$ is a $n_i \times n_i$ square fuzzy matrix where atleast one $n_i \neq n_j$, $i \neq j$, $(1 \leq i, j \leq n)$. Then we define V to be a special fuzzy mixed square matrix.*

We illustrate this by the following example.

***Example 1.2.6:*** Let

$$V = V_1 \cup V_2 \cup V_3 \cup V_4 \cup V_5 \cup V_6$$

$$= \begin{bmatrix} 0.3 & 0.8 & 1 \\ 0.9 & 1 & 0.5 \\ 0.2 & 0.4 & 0 \end{bmatrix} \cup \begin{bmatrix} 0.9 & 1 & 0 & 0.9 \\ 1 & 0.7 & 1 & 0.5 \\ 0.8 & 0.1 & 0.5 & 1 \\ 0.6 & 1 & 0.7 & 0 \end{bmatrix}$$

$$\cup \begin{bmatrix} 0.91 & 0.82 \\ 0.45 & 1 \end{bmatrix} \cup \begin{bmatrix} 0.91 & 0 & 1 & 0.8 & 0.7 \\ 0 & 1 & 0 & 1 & 0 \\ 1 & 0.7 & 0.7 & 0.6 & 0.2 \\ 0.8 & 0.6 & 0.3 & 0.1 & 1 \\ 0.3 & 0.1 & 0.11 & 0 & 0.3 \end{bmatrix}$$

$$\cup \begin{bmatrix} 1 & 0 & 0.3 \\ 0.9 & 1 & 0.8 \\ 1 & 0.4 & 1 \end{bmatrix} \cup \begin{bmatrix} 0.8 & 0.1 & 0.4 & 0.9 & 1 \\ 0.3 & 1 & 0.2 & 0.8 & 0.3 \\ 0.4 & 0.2 & 1 & 0.2 & 0.5 \\ 0.3 & 1 & 0.3 & 0.6 & 1 \\ 1 & 0.7 & 1 & 0.2 & 0.6 \end{bmatrix}.$$

Clearly V is a special fuzzy square mixed matrix.

Now we proceed on to define the notion of special fuzzy rectangular matrix.

**DEFINITION 1.2.5:** *Let $S = S_1 \cup S_2 \cup \ldots \cup S_m$ $(m \geq 2)$ where each $S_i$ is a $t \times s$ rectangular fuzzy matrix $t \neq s$; $1 \leq i \leq m$ then we define S to be a special fuzzy rectangular matrix.*

We illustrate this notion by an example.

***Example 1.2.7:*** Let $S = S_1 \cup S_2 \cup \ldots \cup S_5$

$$= \begin{bmatrix} 0.2 & 1 & 0 \\ 0.1 & 0.1 & 1 \\ 0.3 & 1 & 0 \\ 0.4 & 0 & 0.7 \end{bmatrix} \cup \begin{bmatrix} 0.3 & 1 & 0.31 \\ 0.2 & 0 & 1 \\ 0.7 & 0.3 & 0.9 \\ 0.2 & 1 & 0.1 \end{bmatrix} \cup \begin{bmatrix} 1 & 0.5 & 1 \\ 0 & 1 & 0.2 \\ 0.1 & 0 & 0.3 \\ 0.4 & 1 & 1 \end{bmatrix}$$

$$\cup \begin{bmatrix} 1 & 0.3 & 0.9 \\ 0.6 & 0.2 & 0.1 \\ 0.7 & 0.8 & 0.5 \\ 0.1 & 1 & 0.4 \end{bmatrix} \cup \begin{bmatrix} 0.3 & 0.2 & 0.9 \\ 0.8 & 1 & 0.4 \\ 1 & 0.5 & 1 \\ 0.7 & 0.4 & 0.3 \end{bmatrix}.$$

We see each $S_i$ is a $4 \times 3$ rectangular fuzzy matrix, $1 \leq i \leq 5$; hence S is a special fuzzy rectangular matrix.

Now we proceed on to define the notion of special fuzzy mixed rectangular matrix.

**DEFINITION 1.2.6:** *Let $P = P_1 \cup P_2 \cup \ldots \cup P_m$ where each $P_i$ is a $s_i \times t_i$ rectangular fuzzy matrix $(s_i \neq t_i)$, $1 \leq i \leq m$ and atleast one $P_i \neq P_j$ for $i \neq j$ i.e., $s_i \neq s_j$ or $t_i \neq t_j$, $1 \leq i, j \leq m$, then we define P to be a special fuzzy mixed rectangular matrix.*

We now illustrate this by the following example.

***Example 1.2.8:*** Let

$$S = S_1 \cup S_2 \cup S_3 \cup S_4 \cup S_5 \cup S_6$$

$$= \begin{bmatrix} 1 & 0.3 \\ 0.8 & 1 \\ 0.7 & 0.4 \\ 0.5 & 0.8 \\ 1 & 0.9 \end{bmatrix} \cup \begin{bmatrix} 0.3 & 0.8 & 1 & 0.7 & 0.5 & 0.3 & 0.2 \\ 0.1 & 1 & 0.9 & 0 & 1 & 0.9 & 0.14 \end{bmatrix}$$

$$\cup \begin{bmatrix} 0.3 & 1 \\ 0.8 & 0.9 \\ 1 & 0 \\ 0.9 & 0.4 \end{bmatrix} \cup \begin{bmatrix} 1 & 0.1 \\ 0.7 & 1 \\ 0.8 & 0.7 \\ 0.4 & 0 \\ 1 & 0.2 \end{bmatrix}$$

$$\cup \begin{bmatrix} 0.3 & 0.2 & 1 \\ 0.1 & 1 & 0.3 \\ 0.9 & 0 & 0.9 \\ 0.1 & 1 & 0.4 \end{bmatrix} \cup \begin{bmatrix} 0.3 & 1 & 0.4 & 1 & 0.5 \\ 0.9 & 0.4 & 1 & 0.8 & 1 \\ 1 & 0.6 & 0.7 & 0.6 & 0.3 \\ 0.5 & 1 & 0.1 & 1 & 0.7 \\ 0.8 & 0.7 & 0 & 0.8 & 0.4 \end{bmatrix}.$$

We see S is a special fuzzy mixed rectangular matrix.

***Example 1.2.9:*** Let

$$X = X_1 \cup X_2 \cup X_3 \cup X_4$$

$$= \begin{bmatrix} 0.1 & 0.8 \\ 0.9 & 1 \\ 0.3 & 0.4 \\ 0.7 & 0 \end{bmatrix} \cup \begin{bmatrix} 0.1 & 0.8 \\ 0.9 & 1 \\ 0.3 & 0.4 \\ 0.7 & 0 \end{bmatrix} \cup \begin{bmatrix} 1 & 0 \\ 0 & 1 \\ 1 & 1 \\ 0 & 0 \end{bmatrix} \cup \begin{bmatrix} 0.9 & 0.3 \\ 0.8 & 1 \\ 0 & 0.1 \\ 1 & 0 \end{bmatrix}$$

is the special fuzzy rectangular matrix.

We see

$$X_1 = X_2 = \begin{bmatrix} 0.1 & 0.8 \\ 0.9 & 1 \\ 0.3 & 0.4 \\ 0.7 & 0 \end{bmatrix}$$

a 4 × 2 fuzzy matrix.

Now we proceed on to define the notion of special fuzzy mixed matrix.

**DEFINITION 1.2.7:** *Let $X = X_1 \cup X_2 \cup X_3 \cup \ldots \cup X_n$ ($n \geq 2$) where $X_i$ is a $t_i \times t_i$ fuzzy square matrix and some $X_j$ is a $p_j \times q_j$ ($p_j \neq q_j$) fuzzy rectangular matrix. Then we define X to be a special fuzzy mixed matrix.*

We now illustrate this by the following examples.

***Example 1.2.10:*** Let

$$T = T_1 \cup T_2 \cup T_3 \cup T_4 \cup T_5 =$$

$$\begin{bmatrix} 0.3 & 1 & 0.8 \\ 1 & 0.9 & 0.5 \\ 0.6 & 0.7 & 1 \end{bmatrix} \cup \begin{bmatrix} 0.8 & 1 \\ 0.7 & 0.6 \\ 0.5 & 0.4 \\ 0.6 & 0.3 \\ 0.2 & 0.1 \\ 1 & 0 \end{bmatrix} \cup \begin{bmatrix} 1 & 0.3 & 0.4 & 0.7 & 0.5 \\ 0.3 & 1 & 0.3 & 0.5 & 1 \\ 1 & 0 & 1 & 0 & 0.7 \end{bmatrix}$$

$$\cup \begin{bmatrix} 0.4 & 1 & 0.3 & 0.9 \\ 0.8 & 0.6 & 1 & 0.5 \\ 0.1 & 1 & 0.8 & 1 \\ 0.7 & 0.5 & 0.7 & 0.3 \end{bmatrix} \cup \begin{bmatrix} 0.8 & 1 \\ 0.9 & 0.2 \end{bmatrix},$$

T is a special fuzzy mixed matrix.

***Example 1.2.11:*** Let

$$S = S_1 \cup S_2 \cup S_3 \cup S_4$$

$$= \begin{bmatrix} 0.3 & 0.2 & 1 \\ 0.9 & 1 & 0.3 \\ 0.2 & 0.7 & 0.5 \end{bmatrix} \cup \begin{bmatrix} 0.1 & 0.2 \\ 0.4 & 0.5 \\ 1 & 0.9 \\ 0.7 & 0.4 \\ 1 & 0.8 \end{bmatrix} \cup$$

$$\begin{bmatrix} 1 & 0 & 0.8 \\ 0 & 1 & 0.3 \\ 0.7 & 0 & 1 \end{bmatrix} \cup \begin{bmatrix} 0.8 & 1 \\ 1 & 0.3 \end{bmatrix};$$

S is a special fuzzy mixed matrix.

***Example 1.2.12:*** Let

$$T = T_1 \cup T_2 \cup T_3 \cup T_4 \cup T_5 =$$

$$[1\ 0\ 0.3\ 1\ 0.5] \cup \begin{bmatrix} 1 \\ 0 \\ 0.8 \\ 0.9 \\ 0.6 \end{bmatrix} \cup \begin{bmatrix} 0.3 & 1 & 0.8 & 0.5 & 0.31 & 0 & 1 \end{bmatrix}$$

$$\cup \begin{bmatrix} 1 \\ 0 \\ 0.8 \\ 0.9 \\ 1 \\ 0 \\ 0.1 \end{bmatrix} \cup \begin{bmatrix} 0.1 & 0.3 & 0.2 & 1 & 0.7 & 0.8 & 0.9 & 1 \end{bmatrix}$$

T is a special fuzzy rectangular mixed matrix. We see T is not a special fuzzy mixed matrix.

Now we proceed on to define how in the special fuzzy matrices we get its transpose.

**DEFINITION 1.2.8:** *Let $P = P_1 \cup P_2 \cup \ldots \cup P_n$ ($n \geq 2$) be a special fuzzy square (rectangular) matrix. Now the transpose of P denoted by $P^T = (P_1 \cup P_2 \cup \ldots \cup P_n)^T = P_1^T \cup P_2^T \cup \ldots \cup P_n^T$, $P^T$ is also a special fuzzy square (rectangular) matrix.*

**DEFINITION 1.2.9:** *Let $X = X_1 \cup X_2 \cup \ldots \cup X_n$ be a special fuzzy column matrix, then $X^T$ the transpose of X is $(X_1 \cup X_2 \cup \ldots \cup X_n)^T = X_1^T \cup X_2^T \cup \ldots \cup X_n^T = X^T$. Clearly $X^T$ is a special fuzzy row matrix. Thus we see if $Y = Y_1 \cup Y_2 \cup \ldots \cup Y_m$ is a special fuzzy row matrix then $Y^T = (Y_1 \cup Y_2 \cup \ldots \cup Y_m)^T = Y_1^T \cup Y_2^T \cup \ldots \cup Y_m^T$ is a special fuzzy column matrix. Like wise if $P = P_1 \cup P_2 \cup \ldots \cup P_m$ is a special fuzzy mixed column matrix then the transpose of P, $P^T = P^T_1 \cup P^T_2 \cup \ldots \cup P^T_m$ is a special fuzzy mixed row matrix and vice versa.*

We just illustrate these by the following examples.

***Example 1.2.13:*** Let $T = T_1 \cup T_2 \cup T_3 \cup T_4$ = [1 0 0.3 0.7 1 0.8] ∪ [1 1 1 0 0.8 1] ∪ [0.8 0.1 1 0 0.7 0.9] ∪ [1 1 1 0 1 1] be a special fuzzy mixed row matrix.

$$T^T = (T_1 \cup T_2 \cup T_3 \cup T_4)^T = T_1^T \cup T_2^T \cup T_3^T \cup T_4^T =$$

$$\begin{bmatrix} 1 \\ 0 \\ 0.3 \\ 0.7 \\ 1 \\ 0.8 \end{bmatrix} \cup \begin{bmatrix} 1 \\ 1 \\ 1 \\ 0 \\ 0.8 \\ 1 \end{bmatrix} \cup \begin{bmatrix} 0.8 \\ 0.1 \\ 1 \\ 0 \\ 0.7 \\ 0.9 \end{bmatrix} \cup \begin{bmatrix} 1 \\ 1 \\ 1 \\ 0 \\ 1 \\ 1 \end{bmatrix}.$$

Clearly $T^T$ is a special fuzzy mixed column matrix.

***Example 1.2.14:*** Consider the special fuzzy mixed row matrix

$$P = P_1 \cup P_2 \cup P_3 \cup P_4 \cup P_5$$
$$= [0\ 1\ 0.3\ 0.9] \cup [1\ 0\ 1\ 1\ 0.9\ 0.3\ 0.7] \cup [0\ 1\ 0.8] \cup [0.9\ 0.8\ 0.7\ 0.6\ 0.5\ 0.4\ 0.3\ 0.2\ 0.1] \cup [1\ 0.9\ 0.6].$$

Now the transpose of P is

$$P^T = \begin{bmatrix} 0 \\ 1 \\ 0.3 \\ 0.9 \end{bmatrix} \cup \begin{bmatrix} 1 \\ 0 \\ 1 \\ 1 \\ 0.9 \\ 0.3 \\ 0.7 \end{bmatrix} \cup \begin{bmatrix} 0 \\ 1 \\ 0.8 \end{bmatrix} \cup \begin{bmatrix} 0.9 \\ 0.8 \\ 0.7 \\ 0.6 \\ 0.5 \\ 0.4 \\ 0.3 \\ 0.2 \\ 0.1 \end{bmatrix} \cup \begin{bmatrix} 1 \\ 0.9 \\ 0.6 \end{bmatrix}.$$

We see $P^T$ is a special fuzzy mixed column matrix.

It can be easily verified that the transpose of a special fuzzy square (mixed) matrix will once again be a special fuzzy square (mixed) matrix. Likewise the transpose of a special fuzzy rectangular (mixed) matrix will once again will be a special fuzzy rectangular (mixed).

We are forced to discuss about this because we will be using these concepts in the special fuzzy models which we will be constructing in chapter two.

Now we give an example of how the transpose of a special fuzzy mixed matrix looks like.

***Example 1.2.15:*** Let

$$S = S_1 \cup S_2 \cup S_3 \cup S_4 \cup S_5$$

$$= [1\ 0\ 0.9\ 0.2\ 1\ 0.7\ 0.8] \cup \begin{bmatrix} 0.9 \\ 1 \\ 0.2 \\ 0.1 \end{bmatrix} \cup \begin{bmatrix} 0.3 & 0.8 & 1 & 0.7 \\ 0.1 & 0.3 & 0 & 1 \\ 1 & 0.4 & 0.3 & 0.7 \\ 0.9 & 1 & 0.4 & 1 \end{bmatrix}$$

$$\cup \begin{bmatrix} 0.3 & 0.2 \\ 1 & 0.7 \\ 0.9 & 0.5 \\ 0.8 & 0.3 \\ 0.6 & 0.1 \\ 1 & 0 \end{bmatrix} \cup \begin{bmatrix} 0.3 & 0.2 & 0.7 & 0.9 & 0.2 & 0.5 \\ 1 & 0.5 & 0.8 & 1 & 0.3 & 0.6 \\ 0.9 & 0.6 & 0 & 0.1 & 0.4 & 0.7 \end{bmatrix}$$

be a special fuzzy mixed matrix.

Now the transpose of S denoted by

$$S^T = S_1^T \cup S_2^T \cup S_3^T \cup S_4^T \cup S_5^T$$

$$= \begin{bmatrix} 1 \\ 0 \\ 0.9 \\ 0.2 \\ 1 \\ 0.7 \\ 0.8 \end{bmatrix} \cup \begin{bmatrix} 0.9 & 1 & 0.2 & 0.1 \end{bmatrix} \cup \begin{bmatrix} 0.3 & 0.1 & 1 & 0.9 \\ 0.8 & 0.3 & 0.4 & 1 \\ 1 & 0 & 0.3 & 0.4 \\ 0.7 & 1 & 0.7 & 1 \end{bmatrix}$$

$$\cup \begin{bmatrix} 0.3 & 1 & 0.9 & 0.8 & 0.6 & 1 \\ 0.2 & 0.7 & 0.5 & 0.3 & 0.1 & 0 \end{bmatrix} \cup \begin{bmatrix} 0.3 & 1 & 0.9 \\ 0.2 & 0.5 & 0.6 \\ 0.7 & 0.8 & 0 \\ 0.9 & 1 & 0.1 \\ 0.2 & 0.3 & 0.4 \\ 0.5 & 0.6 & 0.7 \end{bmatrix}.$$

We see $S^T$ is also a special fuzzy mixed matrix.

Now we proceed on to introduce some special fuzzy operations on special fuzzy square matrix which will be used on these special fuzzy models. Thus at this juncture we do not promise to give all types of operations that can be carried out on these class of matrices we give here only those relevant operations on these new class of special fuzzy matrices which will be described in chapter two of this book. Operations on special fuzzy square matrices with a special fuzzy operator 'o' which yields a fixed point or a limit cycle.

Let

$T = T_1 \cup T_2 \cup \ldots \cup T_m$

be a special fuzzy square matrix where each $T_i$ is a $n \times n$ fuzzy matrix with entries from the set $\{-1, 0, 1\}$, $1 \leq i \leq m$.

Now suppose we have a special fuzzy row vector say

$X = X_1 \cup X_2 \cup \ldots \cup X_m$

where $X_i = \left[x_1^i \; x_2^i \; \ldots \; x_n^i\right]$, $1 \leq i \leq m$. To find the effect of X on the special fuzzy square matrix T. It is to be noted that the $X_i$'s need not all be distinct, we can have $X_i = X_j$ without i = j. Now the entries in each $X_i = \left[x_1^i \; x_2^i \; \ldots \; x_n^i\right]$ is such that $x^i_k \in \{0, 1\}$; i = 1, 2, …, m; k = 1, 2, …, n.

Now the operation

$$\begin{aligned} X \text{ o } T &= (X_1 \cup X_2 \cup \ldots \cup X_m) \text{ o } (T_1 \cup T_2 \cup \ldots \cup T_n) \\ &= X_1 \text{ o } T_1 \cup X_2 \text{ o } T_2 \cup \ldots \cup X_m \text{ o } T_m. \end{aligned}$$

The value of $X_i$ o $T_i$ is calculated as mentioned in pages 20-1 of this book, for i = 1, 2, …, m.

We illustrate this by the following example.

***Example 1.2.16:*** Let

$$T = T_1 \cup T_2 \cup T_3 \cup T_4 \cup T_5$$

$$= \begin{bmatrix} 0 & 1 & 0 & 0 & -1 \\ 1 & 0 & 1 & 1 & 0 \\ 0 & 1 & 0 & -1 & 0 \\ 0 & 0 & 0 & 0 & 1 \\ 1 & 0 & -1 & 0 & 0 \end{bmatrix} \cup$$

$$\begin{bmatrix} 0 & 0 & 1 & 0 & 1 \\ -1 & 0 & 0 & -1 & 0 \\ 0 & 0 & 0 & 1 & 0 \\ 1 & 0 & 0 & 0 & 1 \\ 0 & 1 & 0 & 0 & 0 \end{bmatrix} \cup \begin{bmatrix} 0 & 1 & 0 & 1 & 0 \\ 0 & 0 & -1 & 0 & 1 \\ 1 & 0 & 0 & 1 & 0 \\ 0 & 1 & 0 & 0 & 1 \\ 1 & 0 & 0 & 1 & 0 \end{bmatrix} \cup$$

$$\begin{bmatrix} 0 & 0 & 0 & 1 & 1 \\ -1 & 0 & 1 & 0 & 0 \\ 0 & 0 & 0 & 1 & 0 \\ 0 & 1 & 0 & 0 & -1 \\ 1 & -1 & 0 & 1 & 0 \end{bmatrix} \cup \begin{bmatrix} 0 & 1 & 0 & -1 & 0 \\ 1 & 0 & -1 & 0 & 0 \\ 0 & 1 & 0 & -1 & 0 \\ 0 & 0 & -1 & 0 & 1 \\ -1 & 0 & 0 & 1 & 0 \end{bmatrix}$$

be a special fuzzy square matrix. Let

$$\begin{aligned} X &= X_1 \cup X_2 \cup X_3 \cup X_4 \cup X_5 \\ &= [0\ 1\ 0\ 0\ 0] \cup [1\ 0\ 0\ 0\ 1] \cup [0\ 0\ 1\ 0\ 0] \cup [0\ 1\ 0\ 0\ 0] \\ &\quad \cup [0\ 0\ 1\ 0\ 0] \end{aligned}$$

be a special fuzzy row matrix.

To find X o T, to be in more technical terms the resultant of X on T.

Now

$$\begin{aligned} X \text{ o } T &= \{[X_1 \cup X_2 \cup X_3 \cup X_4 \cup X_5] \text{ o } [T_1 \cup T_2 \cup T_3 \cup T_4 \\ &\quad \cup T_5]\} \\ &= X_1 \text{ o } T_1 \cup X_2 \text{ o } T_2 \cup X_3 \text{ o } T_3 \cup X_4 \text{ o } T_4 \cup X_5 \text{ o } T_5 \end{aligned}$$

$$= \begin{bmatrix} 0 & 1 & 0 & 0 & 0 \end{bmatrix} \circ \begin{bmatrix} 0 & 1 & 0 & 0 & -1 \\ 1 & 0 & 1 & 1 & 0 \\ 0 & 1 & 0 & -1 & 0 \\ 0 & 0 & 0 & 0 & 1 \\ 1 & 0 & -1 & 0 & 0 \end{bmatrix} \cup$$

$$\begin{bmatrix} 1 & 0 & 0 & 0 & 1 \end{bmatrix} \circ \begin{bmatrix} 0 & 0 & 1 & 0 & 1 \\ -1 & 0 & 0 & -1 & 0 \\ 0 & 0 & 0 & 1 & 0 \\ 1 & 0 & 0 & 0 & 1 \\ 0 & 1 & 0 & 0 & 0 \end{bmatrix} \cup$$

$$\begin{bmatrix} 0 & 0 & 1 & 0 & 0 \end{bmatrix} \circ \begin{bmatrix} 0 & 1 & 0 & 1 & 0 \\ 0 & 0 & -1 & 0 & 1 \\ 1 & 0 & 0 & 1 & 0 \\ 0 & 1 & 0 & 0 & 1 \\ 1 & 0 & 0 & 1 & 0 \end{bmatrix} \cup$$

$$\begin{bmatrix} 0 & 1 & 0 & 0 & 0 \end{bmatrix} \circ \begin{bmatrix} 0 & 0 & 0 & 1 & 1 \\ -1 & 0 & 1 & 0 & 0 \\ 0 & 0 & 0 & 1 & 0 \\ 0 & 1 & 0 & 0 & -1 \\ 1 & -1 & 0 & 1 & 0 \end{bmatrix} \cup$$

$$\begin{bmatrix} 0 & 0 & 1 & 0 & 0 \end{bmatrix} \circ \begin{bmatrix} 0 & 1 & 0 & -1 & 0 \\ 1 & 0 & -1 & 0 & 0 \\ 0 & 1 & 0 & -1 & 0 \\ 0 & 0 & -1 & 0 & 1 \\ -1 & 0 & 0 & 1 & 0 \end{bmatrix}$$

$= [1\ 0\ 1\ 1\ 0] \cup [0\ 1\ 1\ 0\ 1] \cup [1\ 0\ 0\ 1\ 0] \cup [-1\ 0\ 1\ 0\ 0] \cup [0\ 1\ 0\ -1\ 0]$

$= Y'$.

After updating and thresholding mentioned in pages 20-1 we get

$Y = [1\ 1\ 1\ 1\ 0] \cup [1\ 1\ 1\ 0\ 1] \cup [1\ 0\ 1\ 1\ 0] \cup [0\ 1\ 1\ 0\ 0] \cup [0\ 1\ 1\ 0\ 0]$

$= Y_1 \cup Y_2 \cup Y_3 \cup Y_4 \cup Y_5$.

Now we find the effect of Yon T.

$Y \text{ o } T = [Y_1 \cup Y_2 \cup Y_3 \cup Y_4 \cup Y_5] \text{ o } [T_1 \cup T_2 \cup T_3 \cup T_4 \cup T_5]$

$= Y_1 \text{ o } T_1 \cup Y_2 \text{ o } T_2 \cup Y_3 \text{ o } T_3 \cup Y_4 \text{ o } T_4 \cup Y_5 \text{ o } T_5$

$$= \begin{bmatrix} 1 & 1 & 1 & 1 & 0 \end{bmatrix} \text{ o } \begin{bmatrix} 0 & 1 & 0 & 0 & -1 \\ 1 & 0 & 1 & 1 & 0 \\ 0 & 1 & 0 & -1 & 0 \\ 0 & 0 & 0 & 0 & 1 \\ 1 & 0 & -1 & 0 & 0 \end{bmatrix} \cup$$

$$\begin{bmatrix} 1 & 1 & 1 & 0 & 1 \end{bmatrix} \text{ o } \begin{bmatrix} 0 & 0 & 1 & 0 & 1 \\ -1 & 0 & 0 & -1 & 0 \\ 0 & 0 & 0 & 1 & 0 \\ 1 & 0 & 0 & 0 & 1 \\ 0 & 1 & 0 & 0 & 0 \end{bmatrix} \cup$$

$$\begin{bmatrix} 1 & 0 & 1 & 1 & 0 \end{bmatrix} \text{ o } \begin{bmatrix} 0 & 1 & 0 & 1 & 0 \\ 0 & 0 & -1 & 0 & 1 \\ 1 & 0 & 0 & 1 & 0 \\ 0 & 1 & 0 & 0 & 1 \\ 1 & 0 & 0 & 1 & 0 \end{bmatrix} \cup$$

$$[0 \;\; 1 \;\; 1 \;\; 0 \;\; 0] \text{ o } \begin{bmatrix} 0 & 0 & 0 & 1 & 1 \\ -1 & 0 & 1 & 0 & 0 \\ 0 & 0 & 0 & 1 & 0 \\ 0 & 1 & 0 & 0 & -1 \\ 1 & -1 & 0 & 1 & 0 \end{bmatrix} \cup$$

$$[0 \;\; 1 \;\; 1 \;\; 0 \;\; 0] \text{ o } \begin{bmatrix} 0 & 1 & 0 & -1 & 0 \\ 1 & 0 & -1 & 0 & 0 \\ 0 & 1 & 0 & -1 & 0 \\ 0 & 0 & -1 & 0 & 1 \\ -1 & 0 & 0 & 1 & 0 \end{bmatrix}.$$

= [1 2 1 0 0] ∪ [–1 1 1 0 1] ∪ [1 2 0 2 1] ∪ [–1 0 1 1 0] ∪ [1 1 –1 –1 0]
= Z',

after updating and thresholding Z' we get

Z = [1 1 1 0 0] ∪ [1 1 1 0 1] ∪ [1 1 1 1 1] ∪ [0 1 1 1 0] ∪ [1 1 1 0 0]
= $Z_1 \cup Z_2 \cup Z_3 \cup Z_4 \cup Z_5$.

Now we find

Z o T = $[Z_1 \cup Z_2 \cup Z_3 \cup Z_4 \cup Z_5]$ o $[T_1 \cup T_2 \cup T_3 \cup T_4 \cup T_5]$
= $Z_1$ o $T_1 \cup Z_2$ o $T_2 \cup Z_3$ o $T_3 \cup Z_4$ o $T_4 \cup Z_5$ o $T_5$

$$= [1 \;\; 1 \;\; 1 \;\; 0 \;\; 0] \text{ o } \begin{bmatrix} 0 & 1 & 0 & 0 & -1 \\ 1 & 0 & 1 & 1 & 0 \\ 0 & 1 & 0 & -1 & 0 \\ 0 & 0 & 0 & 0 & 1 \\ 1 & 0 & -1 & 0 & 0 \end{bmatrix} \cup$$

$$[1\ \ 1\ \ 1\ \ 0\ \ 1] \text{ o } \begin{bmatrix} 0 & 0 & 1 & 0 & 1 \\ -1 & 0 & 0 & -1 & 0 \\ 0 & 0 & 0 & 1 & 0 \\ 1 & 0 & 0 & 0 & 1 \\ 0 & 1 & 0 & 0 & 0 \end{bmatrix} \cup$$

$$[1\ \ 1\ \ 1\ \ 1\ \ 1] \text{ o } \begin{bmatrix} 0 & 1 & 0 & 1 & 0 \\ 0 & 0 & -1 & 0 & 1 \\ 1 & 0 & 0 & 1 & 0 \\ 0 & 1 & 0 & 0 & 1 \\ 1 & 0 & 0 & 1 & 0 \end{bmatrix} \cup$$

$$[0\ \ 1\ \ 1\ \ 1\ \ 0] \text{ o } \begin{bmatrix} 0 & 0 & 0 & 1 & 1 \\ -1 & 0 & 1 & 0 & 0 \\ 0 & 0 & 0 & 1 & 0 \\ 0 & 1 & 0 & 0 & -1 \\ 1 & -1 & 0 & 1 & 0 \end{bmatrix} \cup$$

$$[1\ \ 1\ \ 1\ \ 0\ \ 0] \text{ o } \begin{bmatrix} 0 & 1 & 0 & -1 & 0 \\ 1 & 0 & -1 & 0 & 0 \\ 0 & 1 & 0 & -1 & 0 \\ 0 & 0 & -1 & 0 & 1 \\ -1 & 0 & 0 & 1 & 0 \end{bmatrix}$$

= [1 2 1 0 –1] ∪ [–1 1 1 0 1] ∪ [2 2 –1 3 2] ∪ [–1 1 1 1 –1] ∪ [1 2 –1 –2 0]
= P' .

Let P be the special fuzzy row vector got by thresholding and updating P'.

P = [1 1 1 0 0] ∪ [1 1 1 0 1] ∪ [1 1 1 1 1] ∪ [0 1 1 1 0] ∪ [1 1 1 0 0].

Thus we get the resultant of the special fuzzy row vector X to be a special fixed point which is a special fuzzy row vector. This sort of operations will be used in the special fuzzy cognitive models which will be introduced in chapter two of this book.

Next we proceed on to give some special fuzzy operations in case of special fuzzy mixed square matrix.

Let

$$P = P_1 \cup P_2 \cup \ldots \cup P_m$$

be a special fuzzy mixed square matrix where $P_k = (p_{ij})$, k = 1, 2, …, m and $P_k$ is a $t_k \times t_k$ square fuzzy matrix with $p_{ij} \in \{-1, 0, 1\}$, $1 \leq i, j \leq t_k$.

Let $X = X_1 \cup X_2 \cup \ldots \cup X_m$ be a special fuzzy mixed row vector such that each $X_k$ is a $1 \times t_k$ fuzzy row vector k = 1, 2, …, m with

$$X_k = \left[x_1^k \; x_2^k \; \ldots \; x_{t_k}^k\right]$$

and $x_i^k \in \{0, 1\}$, $1 \leq i \leq t_k$; k = 1, 2, …, m.

Now

$$\begin{aligned} X \text{ o } P &= [X_1 \cup X_2 \cup \ldots \cup X_m] \text{ o } [P_1 \cup P_2 \cup \ldots \cup P_m] \\ &= X_1 \text{ o } P_1 \cup X_2 \text{ o } P_2 \cup \ldots \cup X_m \text{ o } P_m. \end{aligned}$$

The operation $X_k$ o $P_k$ is carried out and described in page 22.

Now we illustrate this explicitly by an example.

***Example 1.2.17:*** Let

$$P = P_1 \cup P_2 \cup P_3 \cup P_4$$

$$= \begin{bmatrix} 0 & 1 & 0 & 1 \\ 1 & 0 & -1 & 0 \\ -1 & 1 & 0 & 0 \\ 0 & 0 & 1 & 0 \end{bmatrix} \cup \begin{bmatrix} 0 & 1 & 0 & 0 & 1 \\ 1 & 0 & 1 & 0 & 0 \\ 0 & 0 & 0 & 1 & 1 \\ 1 & 1 & 0 & 0 & 0 \\ 0 & 0 & 1 & 0 & 1 \end{bmatrix} \cup$$

$$\begin{bmatrix} 0 & 1 & 1 & 0 & 0 & 0 \\ 1 & 0 & 0 & 0 & 1 & 0 \\ 0 & 0 & 0 & 1 & 0 & 1 \\ 1 & 0 & 0 & 0 & 1 & 0 \\ 0 & 1 & 0 & 1 & 0 & 0 \\ 0 & 0 & 1 & 0 & 1 & 0 \end{bmatrix} \cup \begin{bmatrix} 0 & 0 & 0 & 1 \\ 0 & 1 & 0 & 0 \\ 0 & 1 & 0 & 1 \\ 1 & 0 & 1 & 0 \end{bmatrix}$$

be a special fuzzy mixed square matrix. Suppose

$X = [1\ 0\ 0\ 0] \cup [0\ 1\ 0\ 0\ 1] \cup [0\ 0\ 1\ 1\ 0\ 0] \cup [0\ 0\ 0\ 1]$
$= X_1 \cup X_2 \cup X_3 \cup X_4,$

be the special fuzzy mixed row vector.

To find

$X \text{ o } P = [X_1 \cup X_2 \cup X_3 \cup X_4] \text{ o } [P_1 \cup P_2 \cup P_3 \cup P_4]$
$= X_1 \text{ o } P_1 \cup X_2 \text{ o } P_2 \cup X_3 \text{ o } P_3 \cup X_4 \text{ o } P_4$

$$= \begin{bmatrix} 1 & 0 & 0 & 0 \end{bmatrix} \text{ o } \begin{bmatrix} 0 & 1 & 0 & 1 \\ 1 & 0 & -1 & 0 \\ -1 & 1 & 0 & 0 \\ 0 & 0 & 1 & 0 \end{bmatrix} \cup$$

$$\begin{bmatrix} 0 & 1 & 0 & 0 & 1 \end{bmatrix} \text{ o } \begin{bmatrix} 0 & 1 & 0 & 0 & 1 \\ 1 & 0 & 1 & 0 & 0 \\ 0 & 0 & 0 & 1 & 1 \\ 1 & 1 & 0 & 0 & 0 \\ 0 & 0 & 1 & 0 & 1 \end{bmatrix} \cup$$

$$\begin{bmatrix} 0 & 0 & 1 & 1 & 0 & 0 \end{bmatrix} \text{ o } \begin{bmatrix} 0 & 1 & 1 & 0 & 0 & 0 \\ 1 & 0 & 0 & 0 & 1 & 0 \\ 0 & 0 & 0 & 1 & 0 & 1 \\ 1 & 0 & 0 & 0 & 1 & 0 \\ 0 & 1 & 0 & 1 & 0 & 0 \\ 0 & 0 & 1 & 0 & 1 & 0 \end{bmatrix} \cup$$

$$[0 \quad 0 \quad 0 \quad 1] \text{ o } \begin{bmatrix} 0 & 0 & 0 & 1 \\ 0 & 1 & 0 & 0 \\ 0 & 1 & 0 & 1 \\ 1 & 0 & 1 & 0 \end{bmatrix}$$

$= [0\ 1\ 0\ 1] \cup [1\ 0\ 2\ 0\ 1] \cup [1\ 0\ 0\ 1\ 1\ 1] \cup [1\ 0\ 1\ 0]$

$= Y'$,

after updating and thresholding Y' we get

$Y = Y_1 \cup Y_2 \cup Y_3 \cup Y_4$

$= [1\ 1\ 0\ 1] \cup [1\ 1\ 1\ 0\ 1] \cup [1\ 0\ 1\ 1\ 1\ 1] \cup [1\ 0\ 1\ 1]$.

Now we find

$Y \text{ o } P = (Y_1 \cup Y_2 \cup Y_3 \cup Y_4) \text{ o } (P_1 \cup P_2 \cup P_3 \cup P_4)$

$= Y_1 \text{ o } P_1 \cup Y_2 \text{ o } P_2 \cup Y_3 \text{ o } P_3 \cup Y_4 \text{ o } P_4$

$$= [1 \quad 1 \quad 0 \quad 1] \text{ o } \begin{bmatrix} 0 & 1 & 0 & 1 \\ 1 & 0 & -1 & 0 \\ -1 & 1 & 0 & 0 \\ 0 & 0 & 1 & 0 \end{bmatrix} \cup$$

$$[1 \quad 1 \quad 1 \quad 0 \quad 1] \text{ o } \begin{bmatrix} 0 & 1 & 0 & 0 & 1 \\ 1 & 0 & 1 & 0 & 0 \\ 0 & 0 & 0 & 1 & 1 \\ 1 & 1 & 0 & 0 & 0 \\ 0 & 0 & 1 & 0 & 1 \end{bmatrix} \cup$$

$$[1 \quad 0 \quad 1 \quad 1 \quad 1 \quad 1] \text{ o } \begin{bmatrix} 0 & 1 & 1 & 0 & 0 & 0 \\ 1 & 0 & 0 & 0 & 1 & 0 \\ 0 & 0 & 0 & 1 & 0 & 1 \\ 1 & 0 & 0 & 0 & 1 & 0 \\ 0 & 1 & 0 & 1 & 0 & 0 \\ 0 & 0 & 1 & 0 & 1 & 0 \end{bmatrix} \cup$$

$$[1 \quad 0 \quad 1 \quad 1] \text{ o } \begin{bmatrix} 0 & 0 & 0 & 1 \\ 0 & 1 & 0 & 0 \\ 0 & 1 & 0 & 1 \\ 1 & 0 & 1 & 0 \end{bmatrix}$$

$= [1\ 1\ 0\ 1] \cup [1\ 1\ 2\ 1\ 3] \cup [1\ 2\ 2\ 2\ 2\ 1] \cup [1\ 1\ 1\ 2]$
$= R'$.

Let R be the special fuzzy state vector obtained after thresholding and updating R'

$R = [1\ 1\ 0\ 1] \cup [1\ 1\ 1\ 1\ 1] \cup [1\ 1\ 1\ 1\ 1\ 1] \cup [1\ 1\ 1\ 1]$.

We see R o P = R. It is left as an exercise for the reader to calculate R o P. Thus we arrive at a special fixed point which is also a special fuzzy mixed row vector.

Next we proceed on to define special fuzzy operations on special fuzzy rectangular matrix $R = R_1 \cup R_2 \cup \ldots \cup R_t$ ($t \geq 2$) where $R_i$'s are $m \times n$ rectangular fuzzy matrices $i = 1, 2, \ldots, t$ ($m \neq n$); Let $X = X_1 \cup X_2 \cup \ldots \cup X_t$ where each $X_i$ is a $1 \times m$ fuzzy row vector with entries of $X_i$ from the set {0,1}. X is clearly a special fuzzy row vector, $1 \leq i \leq t$.

Now how to study the effect of X on R i.e., how does X operate on R.

$X \text{ o } R = (X_1 \cup X_2 \cup \ldots \cup X_t) \text{ o } (R_1 \cup R_2 \cup \ldots \cup R_t)$
$= X_1 \text{ o } R_1 \cup X_2 \text{ o } R_2 \cup \ldots \cup X_t \text{ o } R_t$

where $X_i$ o $R_i$ operation is described in pages 23-4 where i = 1, 2, …, t.

Now let

$X \text{ o } R = Y'$
$= Y'_1 \cup Y'_2 \cup \ldots \cup Y'_t$,

clearly each $Y_i$ is a $1 \times n$ row vector but now entries of $Y'_i$ need not even belong to the unit fuzzy interval [0, 1]. So we threshold and update $Y'_i$ as described in pages 20-1 for i = 1, 2, 3, …, t.

Now let the updated and thresholded $Y'_i$ be denoted by $Y_i = [Y^i_1, Y^i_2, \ldots, Y^i_t]$ we see $Y^i_j \in \{0, 1\}$ for j = 1, 2, …, t, $1 \le i \le t$. Let $Y = Y_1 \cup Y_2 \cup \ldots \cup Y_t$, Y is a special fuzzy row vector.

Now using Y we can work only on the transpose of R i.e., $R^T$, otherwise we do not have the compatibility of 'o' operation with R.

We find $Y \text{ o } R^T = (Y_1 \cup Y_2 \cup \ldots \cup Y_t) \text{ o } (R_1 \cup R_2 \cup \ldots \cup R_t)^T = Y_1 \text{ o } R^T_1 \cup Y_2 \text{ o } R^T_2 \cup \ldots \cup Y_t \text{ o } R^T_t = Z'_1 \cup Z'_2 \cup \ldots \cup Z'_t = Z'$ we see Z' in general need not be a special fuzzy row vector. We update and threshold each $Z'_i$ to $Z_i$ $(1 \le i \le t)$ to get $Z = Z_1 \cup Z_2 \cup \ldots \cup Z_t$. Z is a special fuzzy row vector. We can work with Z as Z o R and so on until we arrive at an equilibrium i.e., a special fixed point or a limit cycle.

We illustrate this by the following example.

***Example 1.2.18:*** Let

$$T = T_1 \cup T_2 \cup T_3 \cup T_4 \cup T_5$$

be a special fuzzy rectangular matrix where

$$T = \begin{bmatrix} 1 & 0 & -1 \\ 0 & 1 & 0 \\ 0 & -1 & 1 \\ 1 & 0 & 1 \\ -1 & 1 & 0 \\ 0 & 0 & 0 \end{bmatrix} \cup \begin{bmatrix} 0 & 1 & 1 \\ 1 & 0 & 1 \\ 1 & 1 & 0 \\ 0 & 1 & 0 \\ 1 & 0 & 0 \\ 0 & 0 & 1 \end{bmatrix} \cup \begin{bmatrix} 0 & 1 & 0 \\ 1 & 0 & 0 \\ 0 & 0 & 1 \\ 1 & -1 & 0 \\ -1 & 0 & 1 \\ 1 & 1 & -1 \end{bmatrix}$$

$$\cup \begin{bmatrix} 1 & 0 & 0 \\ 0 & 1 & 0 \\ 0 & 0 & 1 \\ 1 & -1 & 0 \\ 0 & 1 & -1 \\ -1 & 0 & 1 \end{bmatrix} \cup \begin{bmatrix} 1 & 1 & -1 \\ 0 & 0 & 1 \\ 1 & 0 & 0 \\ 0 & 1 & 0 \\ -1 & 0 & 1 \\ 1 & -1 & 0 \end{bmatrix};$$

T is clearly a 6 × 3 rectangular special fuzzy matrix with entries from {–1, 0, 1}.

Let X = [1 0 0 1 0 0] ∪ [0 0 0 1 0 0] ∪ [0 0 0 0 0 1] ∪ [1 0 0 0 0 0] ∪ [0 1 0 0 0 0] = $X_1 \cup X_2 \cup X_3 \cup X_4 \cup X_5$ be a special fuzzy row vector with entries from the set {0, 1}.

The effect of X on T is given by

$$X \text{ o } T = (X_1 \cup X_2 \cup X_3 \cup X_4 \cup X_5) \text{ o } (T_1 \cup T_2 \cup T_3 \cup T_4 \cup T_5)$$

$$= X_1 \text{ o } T_1 \cup X_2 \text{ o } T_2 \cup X_3 \text{ o } T_3 \cup X_4 \text{ o } T_4 \cup X_5 \text{ o } T_5$$

$$= \begin{bmatrix} 1 & 0 & 0 & 1 & 0 & 0 \end{bmatrix} \text{ o } \begin{bmatrix} 1 & 0 & -1 \\ 0 & 1 & 0 \\ 0 & -1 & 1 \\ 1 & 0 & 1 \\ -1 & 1 & 0 \\ 0 & 0 & 0 \end{bmatrix} \cup$$

$$\begin{bmatrix} 0 & 0 & 0 & 1 & 0 & 0 \end{bmatrix} \text{ o } \begin{bmatrix} 0 & 1 & 1 \\ 1 & 0 & 1 \\ 1 & 1 & 0 \\ 0 & 1 & 0 \\ 1 & 0 & 0 \\ 0 & 0 & 1 \end{bmatrix} \cup$$

$$\begin{bmatrix} 0 & 0 & 0 & 0 & 0 & 1 \end{bmatrix} \text{ o } \begin{bmatrix} 0 & 1 & 0 \\ 1 & 0 & 0 \\ 0 & 0 & 1 \\ 1 & -1 & 0 \\ -1 & 0 & 1 \\ 1 & 1 & -1 \end{bmatrix} \cup$$

$$[1 \quad 0 \quad 0 \quad 0 \quad 0 \quad 0] \text{ o } \begin{bmatrix} 1 & 0 & 0 \\ 0 & 1 & 0 \\ 0 & 0 & 1 \\ 1 & -1 & 0 \\ 0 & 1 & -1 \\ -1 & 0 & 1 \end{bmatrix} \cup$$

$$[0 \quad 1 \quad 0 \quad 0 \quad 0 \quad 0] \text{ o } \begin{bmatrix} 1 & 1 & -1 \\ 0 & 0 & 1 \\ 1 & 0 & 0 \\ 0 & 1 & 0 \\ -1 & 0 & 1 \\ 1 & -1 & 0 \end{bmatrix}$$

$= [2\ 0\ 0] \cup [0\ 1\ 0] \cup [1\ 1\ -1] \cup [1\ 0\ 0] \cup [0\ 0\ 1]$
$= Z'$
$= Z'_1 \cup Z'_2 \cup Z'_3 \cup Z'_4 \cup Z'_5$

clearly Z' is not a special fuzzy row vector so we threshold it and obtain

$Z = Z_1 \cup Z_2 \cup \ldots \cup Z_5$
$= [1\ 0\ 0] \cup [0\ 1\ 0] \cup [1\ 1\ 0] \cup [1\ 0\ 0] \cup [0\ 0\ 1].$

Now

$Z \text{ o } T^T = (Z_1 \cup Z_2 \cup Z_3 \cup Z_4 \cup Z_5) \text{ o } (T^T{}_1 \cup T^T{}_2 \cup \ldots \cup T^T{}_5)$
$= Z_1 \text{ o } T^T{}_1 \cup Z_2 \text{ o } T^T{}_2 \cup \ldots \cup Z_5 \text{ o } T^T{}_5$

$$= [1 \quad 0 \quad 0] \text{ o } \begin{bmatrix} 1 & 0 & 0 & 1 & -1 & 0 \\ 0 & 1 & -1 & 0 & 1 & 0 \\ -1 & 0 & 1 & 1 & 0 & 0 \end{bmatrix} \cup$$

$$[0 \quad 1 \quad 0] \text{ o } \begin{bmatrix} 0 & 1 & 1 & 0 & 1 & 0 \\ 1 & 0 & 1 & 1 & 0 & 0 \\ 1 & 1 & 0 & 0 & 0 & 1 \end{bmatrix} \cup$$

$$[1 \quad 1 \quad 0] \text{ o } \begin{bmatrix} 0 & 1 & 0 & 1 & -1 & 1 \\ 1 & 0 & 0 & -1 & 0 & 1 \\ 0 & 0 & 1 & 0 & 1 & -1 \end{bmatrix} \cup$$

$$[1 \quad 0 \quad 0] \text{ o } \begin{bmatrix} 1 & 0 & 0 & 1 & 0 & -1 \\ 0 & 1 & 0 & -1 & 1 & 0 \\ 0 & 0 & 1 & 0 & -1 & 1 \end{bmatrix} \cup$$

$$[0 \quad 0 \quad 1] \text{ o } \begin{bmatrix} 1 & 0 & 1 & 0 & -1 & 1 \\ 1 & 0 & 0 & 1 & 0 & -1 \\ -1 & 1 & 0 & 0 & 1 & 0 \end{bmatrix}$$

$= [1\ 0\ 0\ 1\ {-1}\ 0] \cup [1\ 0\ 1\ 1\ 0\ 0] \cup [1\ 1\ 0\ 0\ {-1}\ 2] \cup [1\ 0\ 0\ 1\ 0\ {-1}] \cup [-1\ 1\ 0\ 0\ 1\ 0]$
$= P'$
$= P'_1 \cup P'_2 \cup \ldots \cup P'_5$

clearly P' is a not a special fuzzy row matrix.

We update and threshold P' to obtain

$P = P_1 \cup P_2 \cup P_3 \cup P_4 \cup P_5$
$= [1\ 0\ 0\ 1\ 0\ 0] \cup [1\ 0\ 1\ 1\ 0\ 0] \cup [1\ 1\ 0\ 0\ 0\ 1] \cup [1\ 0\ 0\ 1\ 0\ 0] \cup [0\ 1\ 0\ 0\ 1\ 0]$.

P is a special fuzzy row vector. Now we can find P o T and so on.

We now define special fuzzy operator on special fuzzy mixed rectangular matrix

Suppose

$P = P_1 \cup P_2 \cup \ldots \cup P_n$

($n \geq 2$) be a special fuzzy mixed rectangular matrix i.e., each $P_i$ is a fuzzy rectangular $s_i \times t_i$ ($s_i \neq t_i$) matrix and $i = 1, 2, \ldots, n$.

Let $X = X_1 \cup X_2 \cup \ldots \cup X_n$ be a special fuzzy mixed row vector were each $X_i$ is a $l \times s_i$ fuzzy row vector with entries from the set {0, 1}; i = 1, 2, …, n.

Now we find

$$X \text{ o } P = [X_1 \cup X_2 \cup \ldots \cup X_n] \text{ o } [P_1 \cup P_2 \cup \ldots \cup P_n]$$
$$= X_1 \text{ o } P_1 \cup X_2 \text{ o } P_2 \cup \ldots \cup X_n \text{ o } P_n.$$

Each $X_i$ o $P_i$ is obtained as described in pages 23-4; i = 1, 2, …, n now let

$$X \text{ o } P = Y'_1 \cup Y'_2 \cup \ldots \cup Y'_n$$
$$= Y'.$$

Clearly Y' is not a special fuzzy row vector. We threshold Y' and obtain

$$Y = Y_1 \cup Y_2 \cup \ldots \cup Y_n;$$

Y is a special fuzzy mixed row vector with each $Y_i$ a $1 \times t_i$ fuzzy row vector with entries from {0, 1}; i = 1, 2, …, n.

We find

$$Y \text{ o } P^T = (Y_1 \cup Y_2 \cup \ldots \cup Y_n) \text{ o } (P^T{}_1 \cup P^T{}_2 \cup \ldots \cup P^T{}_n)$$
$$= Y_1 \text{ o } P^T{}_1 \cup Y_2 \circ P_2^T \cup \ldots \cup Y_n \circ P_n^T.$$

Let

$$Y \text{ o } P^T = Z'_1 \cup Z'_2 \cup \ldots \cup Z'_n$$
$$= Z';$$

Z' need not be a special fuzzy row vector we update and threshold Z' and obtain $Z = Z_1 \cup Z_2 \cup \ldots \cup Z_n$ which is a special fuzzy row vector and each $Z_i$ is a $1 \times s_i$ fuzzy row vector for i = 1, 2, …, n. Thus $Z = Z_1 \cup Z_2 \cup \ldots \cup Z_n$ is a special fuzzy mixed row vector. Using Z we can find Z o P and so on.

Now we illustrate this by the following example.

***Example 1.2.19:*** Let T be a special fuzzy mixed rectangular matrix, where

$$T = \begin{bmatrix} 1 & 0 & 1 & -1 \\ 0 & 1 & 0 & 0 \\ 0 & 0 & 1 & 0 \\ 1 & 0 & 0 & -1 \\ 0 & 1 & 0 & 1 \\ -1 & 0 & 1 & 0 \end{bmatrix} \cup \begin{bmatrix} 0 & 1 & 0 & 1 & 0 & -1 & 0 \\ 1 & 0 & -1 & 0 & 0 & 0 & 0 \\ 0 & 0 & 1 & 0 & 1 & 0 & -1 \end{bmatrix} \cup$$

$$\begin{bmatrix} 0 & 1 & 0 & 0 & 0 & 1 \\ 1 & 0 & -1 & 0 & 0 & 0 \\ 0 & 0 & 0 & 1 & -1 & 0 \\ 0 & 1 & 0 & 0 & 0 & -1 \\ 1 & 0 & -1 & 0 & 1 & 0 \\ -1 & 1 & -1 & 0 & 0 & 0 \\ 0 & 0 & 0 & -1 & 1 & 1 \end{bmatrix} \cup \begin{bmatrix} -1 & 0 & 0 \\ 0 & 1 & 0 \\ 1 & -1 & 1 \\ 0 & 0 & 0 \\ 1 & 0 & -1 \end{bmatrix}$$

$$\cup \begin{bmatrix} 1 & 0 & 0 & 0 & -1 \\ 0 & 1 & -1 & 0 & 0 \\ 0 & 1 & 0 & 1 & 0 \\ 1 & 0 & 0 & 0 & -1 \\ 0 & -1 & 0 & 1 & 0 \\ 0 & 0 & 1 & 0 & 1 \\ -1 & 0 & 1 & 0 & 0 \\ 0 & 0 & 1 & -1 & 1 \\ 0 & 1 & 0 & 0 & 1 \end{bmatrix}$$

we see $T_1$ is a 6 × 4 matrix fuzzy matrix, $T_2$ a 3 × 7 matrix, $T_3$ a 7 × 6 fuzzy matrix, $T_4$ a 5 × 3 fuzzy matrix and a $T_5$ a 9 × 5 fuzzy matrix. Thus T is a special fuzzy mixed rectangular matrix.
Suppose

$$\begin{aligned} X &= X_1 \cup X_2 \cup X_3 \cup X_4 \cup X_5 \\ &= [1\ 0\ 0\ 1\ 0\ 0] \cup [0\ 1\ 0] \cup [0\ 0\ 1\ 1\ 0\ 0\ 0] \cup [1\ 0\ 0\ 0\ 0] \cup [1\ 0\ 0\ 0\ 0\ 0\ 0\ 0\ 0] \end{aligned}$$

be a special fuzzy mixed row vector to find the special product of X with T.

$$\begin{aligned} X \text{ o } T &= [X_1 \cup X_2 \cup \ldots \cup X_5] \text{ o } [T_1 \cup T_2 \cup \ldots \cup T_5] \\ &= X_1 \text{ o } T_1 \cup X_2 \text{ o } T_2 \cup X_3 \text{ o } T_3 \cup X_4 \text{ o } T_4 \cup X_5 \text{ o } T_5 \end{aligned}$$

$$= \begin{bmatrix} 1 & 0 & 0 & 1 & 0 & 0 \end{bmatrix} \text{ o } \begin{bmatrix} 1 & 0 & 1 & -1 \\ 0 & 1 & 0 & 0 \\ 0 & 0 & 1 & 0 \\ 1 & 0 & 0 & -1 \\ 0 & 1 & 0 & 1 \\ -1 & 0 & 1 & 0 \end{bmatrix} \cup$$

$$\begin{bmatrix} 0 & 1 & 0 \end{bmatrix} \text{ o } \begin{bmatrix} 0 & 1 & 0 & 1 & 0 & -1 & 0 \\ 1 & 0 & -1 & 0 & 0 & 0 & 0 \\ 0 & 0 & 1 & 0 & 1 & 0 & -1 \end{bmatrix} \cup$$

$$\begin{bmatrix} 0 & 0 & 1 & 1 & 0 & 0 & 0 \end{bmatrix} \text{ o } \begin{bmatrix} 0 & 1 & 0 & 0 & 0 & 1 \\ 1 & 0 & -1 & 0 & 0 & 0 \\ 0 & 0 & 0 & 1 & -1 & 0 \\ 0 & 1 & 0 & 0 & 0 & -1 \\ 1 & 0 & -1 & 0 & 1 & 0 \\ -1 & 1 & -1 & 0 & 0 & 0 \\ 0 & 0 & 0 & -1 & 1 & 1 \end{bmatrix} \cup$$

$$[1 \;\; 0 \;\; 0 \;\; 0 \;\; 0] \text{ o } \begin{bmatrix} -1 & 0 & 0 \\ 0 & 1 & 0 \\ 1 & -1 & 1 \\ 0 & 0 & 0 \\ 1 & 0 & -1 \end{bmatrix} \cup$$

$$[1 \;\; 0 \;\; 0 \;\; 0 \;\; 0 \;\; 0 \;\; 0 \;\; 0 \;\; 0 \;\; 0] \text{ o } \begin{bmatrix} 1 & 0 & 0 & 0 & -1 \\ 0 & 1 & -1 & 0 & 0 \\ 0 & 1 & 0 & 1 & 0 \\ 1 & 0 & 0 & 0 & -1 \\ 0 & -1 & 0 & 1 & 0 \\ 0 & 0 & 1 & 0 & 1 \\ -1 & 0 & 1 & 0 & 0 \\ 0 & 0 & 1 & -1 & 1 \\ 0 & 1 & 0 & 0 & 1 \end{bmatrix}$$

$$\begin{aligned} &= [2\ 0\ 1\ –2] \cup [1\ 0\ –1\ 0\ 0\ 0\ 0] \cup [0\ 1\ 0\ 1\ –1\ –1] \cup \\ &\quad [–1\ 0\ 0] \cup [1\ 0\ 0\ 0\ –1] \\ &= Y'_1 \cup Y'_2 \cup Y'_3 \cup Y'_4 \cup Y'_5 \\ &= Y'. \end{aligned}$$

Clearly Y' is not a special fuzzy mixed row vector. After thresholding Y' we get

$$\begin{aligned} Y &= [1\ 0\ 1\ 0] \cup [1\ 0\ 0\ 0\ 0\ 0\ 0] \cup [0\ 1\ 0\ 1\ 0\ 0] \cup [0\ 0\ 0] \\ &\quad \cup [1\ 0\ 0\ 0\ 0]. \end{aligned}$$

Y is a special fuzzy mixed row vector. We find

$$\begin{aligned} Y \text{ o } T^T &= [Y_1 \cup Y_2 \cup \ldots \cup Y_5] \text{ o } [T_1^T \cup T_2^T \cup T_3^T \cup T_4^T \cup T_5^T] \\ &= Y_1 \text{ o } T_1^T \cup Y_2 \text{ o } T_2^T \cup Y_3 \text{ o } T_3^T \cup Y_4 \text{ o } T_4^T \cup Y_5 \text{ o } T_5^T \end{aligned}$$

$$= \begin{bmatrix} 1 & 0 & 1 & 0 \end{bmatrix} \text{ o } \begin{bmatrix} 1 & 0 & 0 & 1 & 0 & -1 \\ 0 & 1 & 0 & 0 & 1 & 0 \\ 1 & 0 & 1 & 0 & 0 & 1 \\ -1 & 0 & 0 & -1 & 1 & 0 \end{bmatrix} \cup$$

$$\begin{bmatrix} 1 & 0 & 0 & 0 & 0 & 0 & 0 \end{bmatrix} \text{ o } \begin{bmatrix} 0 & 1 & 0 \\ 1 & 0 & 0 \\ 0 & -1 & 1 \\ 1 & 0 & 0 \\ 0 & 0 & 1 \\ -1 & 0 & 0 \\ 0 & 0 & -1 \end{bmatrix} \cup$$

$$\begin{bmatrix} 0 & 1 & 0 & 1 & 0 & 0 \end{bmatrix} \text{ o } \begin{bmatrix} 0 & 1 & 0 & 0 & 1 & -1 & 0 \\ 1 & 0 & 0 & 1 & 0 & 1 & 0 \\ 0 & -1 & 0 & 0 & -1 & -1 & 0 \\ 0 & 0 & 1 & 0 & 0 & 0 & -1 \\ 0 & 0 & -1 & 0 & 1 & 0 & 1 \\ 1 & 0 & 0 & -1 & 0 & 0 & 1 \end{bmatrix} \cup$$

$$\begin{bmatrix} 0 & 0 & 0 \end{bmatrix} \text{ o } \begin{bmatrix} -1 & 0 & 1 & 0 & 1 \\ 0 & 1 & -1 & 0 & 0 \\ 0 & 0 & 1 & 0 & -1 \end{bmatrix} \cup$$

$$\begin{bmatrix} 1 & 0 & 0 & 0 & 0 \end{bmatrix} \text{ o } \begin{bmatrix} 1 & 0 & 0 & 1 & 0 & 0 & -1 & 0 & 0 \\ 0 & 1 & 1 & 0 & -1 & 0 & 0 & 0 & 1 \\ 0 & -1 & 0 & 0 & 0 & 1 & 1 & 1 & 0 \\ 0 & 0 & 1 & 0 & 1 & 0 & 0 & -1 & 0 \\ -1 & 0 & 0 & -1 & 0 & 1 & 0 & 1 & 1 \end{bmatrix}.$$

$= [2\ 0\ 1\ 1\ 0\ 0] \cup [0\ 1\ 0] \cup [1\ 0\ 1\ 1\ 0\ 1\ –1] \cup [0\ 0\ 0\ 0\ 0] \cup [1\ 0\ 0\ 1\ 0\ 0\ –1\ 0\ 0]$
$= Z'$
$= Z'_1 \cup Z'_2 \cup Z'_3 \cup Z'_4 \cup Z'_5.$

We see Z' is not a special fuzzy mixed row vector. Now we update and threshold Z' to

$$
\begin{aligned}
Z &= Z_1 \cup Z_2 \cup Z_3 \cup Z_4 \cup Z_5 \\
&= [1\ 0\ 1\ 1\ 0\ 0] \cup [0\ 1\ 0] \cup [1\ 0\ 1\ 1\ 0\ 1\ 0] \cup [1\ 0\ 0\ 0\ 0] \cup [1\ 0\ 0\ 1\ 0\ 0\ 0\ 0\ 0]
\end{aligned}
$$

which is a special fuzzy mixed row vector.

Now we can work with

$$
\begin{aligned}
Z \text{ o } T &= [Z_1 \cup Z_2 \cup \ldots \cup Z_5] \text{ o } [T_1 \cup T_2 \cup \ldots \cup T_5] \\
&= Z_1 \text{ o } T_1 \cup Z_2 \text{ o } T_2 \cup \ldots \cup Z_5 \text{ o } T_5
\end{aligned}
$$

and so on.

Thus we have seen a special type of operation on the special fuzzy mixed row vector with special fuzzy mixed rectangular matrix.

Now we proceed on to find yet another special type of operations on the four types of special fuzzy matrices viz. the special min max operator on the special fuzzy square matrix.

Suppose $W = (W_1 \cup W_2 \cup \ldots \cup W_m)$ $(m \geq 2)$ be a special fuzzy square matrix with entries from [0, 1]. Further let us assume each $W_i$ is a $n \times n$ square fuzzy matrix, $1 \leq i \leq m$.

Let

$$X = X_1 \cup X_2 \cup \ldots \cup X_m$$

with $X_j = \left[x_1^j\ x_2^j\ \ldots\ x_m^j\right]$ be a special fuzzy row vector, where $x_t^j \in [0, 1]$ for j = 1, 2, …, m and $1 \leq t \leq n$. To use special min max operator, to find the special product X and W.

Now

min max {X, W}

$=$ min max $\{(X_1 \cup X_2 \cup \ldots \cup X_m), (W_1 \cup W_2 \cup \ldots \cup W_m)\}$
$=$ min max$\{X_1, W_1\} \cup$ min max$\{X_2, W_2\} \cup \ldots \cup$ min max$\{X_m, W_m\}$.

The min max operator has been defined between a fuzzy row vector and a fuzzy square matrix in pages 14-21.

Now we illustrate this by the following example.

***Example 1.2.20:*** Let $W = W_1 \cup W_2 \cup W_3 \cup W_4 \cup W_5$ be a special fuzzy square matrix where W

$$= \begin{bmatrix} 0 & 0.3 & 0.5 & 0.1 & 0.2 \\ 1 & 0.4 & 0.6 & 0.8 & 0 \\ 0.2 & 0.2 & 0 & 0.9 & 0.4 \\ 0.7 & 0 & 0.5 & 0 & 0.8 \\ 0.9 & 0.1 & 0 & 0.7 & 0 \end{bmatrix}$$

$$\cup \begin{bmatrix} 1 & 0.2 & 0.4 & 0.6 & 0.8 \\ 0 & 0.3 & 0.5 & 0.7 & 0.9 \\ 1 & 0 & 0.2 & 1 & 0.7 \\ 0 & 0.5 & 0.4 & 0 & 0.3 \\ 0.1 & 1 & 0.7 & 0.8 & 1 \end{bmatrix}$$

$$\cup \begin{bmatrix} 0.3 & 1 & 0.2 & 1 & 0.4 \\ 1 & 0.2 & 0 & 0.5 & 0.7 \\ 0 & 0.8 & 1 & 0.4 & 0.9 \\ 0.8 & 1 & 0 & 0.7 & 0 \\ 0.4 & 0 & 0.3 & 1 & 0.1 \end{bmatrix}$$

$$\cup \begin{bmatrix} 1 & 0.2 & 0.4 & 1 & 0.7 \\ 0 & 0.2 & 1 & 0.7 & 0.5 \\ 0.7 & 1 & 0.5 & 1 & 0.4 \\ 0.9 & 0.7 & 0.3 & 0 & 1 \\ 1 & 0.2 & 0.4 & 0.7 & 0.4 \end{bmatrix}$$

$$\cup \begin{bmatrix} 1 & 0 & 0.3 & 0.5 & 0.7 \\ 0 & 0.7 & 1 & 0.9 & 0.2 \\ 0.4 & 0.6 & 0.8 & 1 & 0 \\ 0.8 & 0.3 & 1 & 0.3 & 0.7 \\ 1 & 0.4 & 0 & 0.2 & 1 \end{bmatrix}$$

Clearly entries of each $W_i$ are in [0, 1], i = 1, 2, 3, 4, 5.

Let

$$\begin{aligned} X &= X_1 \cup X_2 \cup X_3 \cup X_4 \cup X_5 \\ &= [1\ 0\ 0\ 0\ 1] \cup [0\ 1\ 0\ 1\ 0] \cup [0\ 0\ 1\ 1\ 0] \cup [1\ 0\ 1\ 0\ 0] \\ &\quad \cup [0\ 1\ 0\ 1\ 0] \end{aligned}$$

be a special fuzzy row vector with entries form the set {0, 1}. To find the special product of X with W using the min max operator

min{max(X, W)}

$= \min\{\max(X_1 \cup X_2 \cup \ldots \cup X_5), (W_1 \cup W_2 \cup \ldots \cup W_5)\}$

$= \min \max\{X_1, W_1\} \cup \min \max\{X_2, W_2\} \cup \min \max\{X_3, W_3\} \cup \min \max\{X_4, W_4\} \cup \min \max\{X_5, W_5\}$

$$= \min\max\left\{ \begin{bmatrix} 1 & 0 & 0 & 0 & 1 \end{bmatrix}, \begin{bmatrix} 0 & 0.3 & 0.5 & 0.1 & 0.2 \\ 1 & 0.4 & 0.6 & 0.8 & 0 \\ 0.2 & 0.2 & 0 & 0.9 & 0.4 \\ 0.7 & 0 & 0.5 & 0 & 0.8 \\ 0.9 & 0.1 & 0 & 0.7 & 0 \end{bmatrix} \right\}$$

$$\cup \text{ min max } \left\{ \begin{bmatrix} 0 & 1 & 0 & 1 & 0 \end{bmatrix}, \begin{bmatrix} 1 & 0.2 & 0.4 & 0.6 & 0.8 \\ 0 & 0.3 & 0.5 & 0.7 & 0.9 \\ 1 & 0 & 0.2 & 1 & 0.7 \\ 0 & 0.5 & 0.4 & 0 & 0.3 \\ 0.1 & 1 & 0.7 & 0.8 & 1 \end{bmatrix} \right\}$$

$$\cup \text{ min max } \left\{ \begin{bmatrix} 0 & 0 & 1 & 1 & 0 \end{bmatrix}, \begin{bmatrix} 0.3 & 1 & 0.2 & 1 & 0.4 \\ 1 & 0.2 & 0 & 0.5 & 0.7 \\ 0 & 0.8 & 1 & 0.4 & 0.9 \\ 0.8 & 1 & 0 & 0.7 & 0 \\ 0.4 & 0 & 0.3 & 1 & 0.1 \end{bmatrix} \right\}$$

$$\cup \text{ min max } \left\{ \begin{bmatrix} 1 & 0 & 1 & 0 & 0 \end{bmatrix}, \begin{bmatrix} 1 & 0.2 & 0.4 & 1 & 0.7 \\ 0 & 0.2 & 1 & 0.7 & 0.5 \\ 0.7 & 1 & 0.5 & 1 & 0.4 \\ 0.9 & 0.7 & 0.3 & 0 & 1 \\ 1 & 0.2 & 0.4 & 0.7 & 0.4 \end{bmatrix} \right\}$$

$$\cup \text{ min max } \left\{ \begin{bmatrix} 0 & 1 & 0 & 1 & 0 \end{bmatrix}, \begin{bmatrix} 1 & 0 & 0.3 & 0.5 & 0.7 \\ 0 & 0.7 & 1 & 0.9 & 0.2 \\ 0.4 & 0.6 & 0.8 & 1 & 0 \\ 0.8 & 0.3 & 1 & 0.3 & 0.7 \\ 1 & 0.4 & 0 & 0.2 & 1 \end{bmatrix} \right\}$$

= [0.2 0 0 0 0] $\cup$ [0.1 0 0.2 0.6 0.7] $\cup$ [0.3 0 0 0.5 0.4] $\cup$ [0 0.2 0.3 0 0.4] $\cup$ [0.4 0 0 0.2 0]

= Y

= $Y_1 \cup Y_2 \cup Y_3 \cup Y_4 \cup Y_5$.

Now using the same min max operation, we can find min max {Y o W} = Z and so on.

Here we define a special max min operator of special fuzzy mixed square matrix. Let $M = M_1 \cup M_2 \cup \ldots \cup M_n$ be a special fuzzy mixed square matrix i.e., each $M_i$ is a $t_i \times t_i$ fuzzy matrix; i = 1, 2, …, n. Suppose $X = X_1 \cup X_2 \cup \ldots \cup X_n$ be a special fuzzy mixed row vector i.e., each $X_i$ is a $\left[x_1^i, \ldots, x_{t_i}^i\right]$ fuzzy row vector $x_k^i \in \{0, 1\}$, i = 1, 2, …, n; $1 \leq k \leq t_i$.

$$\begin{aligned} \text{max min}\{X, M\} &= \text{max min}\{X_1, M_1\} \cup \text{max min}\{X_2, M_2\} \cup \ldots \cup \\ &\quad \text{max min}\{X_n, M_n\} \end{aligned}$$

where max min $\{X_i, M_i\}$ is found in pages 25 and 26 of this book, i = 1, 2, …, n.

Now we illustrate this situation by the following example.

***Example 1.2.21:*** Let $M = M_1 \cup M_2 \cup M_3 \cup M_4 \cup M_5 \cup M_6$ be a special fuzzy mixed square matrix; that is entries in each $M_i$ are from the unit interval [0, 1]; i = 1, 2, …, 6. Suppose $X = X_1 \cup X_2 \cup \ldots \cup X_6$ be a special fuzzy mixed row vector. We will illustrate how max min {X, M} is obtained. Given

$$M = M_1 \cup M_2 \cup M_3 \cup M_4 \cup M_5 \cup M_6$$

$$= \begin{bmatrix} 0.4 & 0.6 & 1 & 0 & 0.5 \\ 0.2 & 0 & 0.6 & 1 & 0.7 \\ 1 & 0.3 & 1 & 0 & 0.5 \\ 0.6 & 1 & 0.3 & 1 & 0.2 \\ 0 & 0.4 & 0.2 & 0 & 1 \\ 0.5 & 1 & 0.6 & 0.7 & 0.5 \end{bmatrix} \cup \begin{bmatrix} 0.3 & 1 & 0.8 \\ 1 & 0.4 & 0.6 \\ 0.5 & 0 & 0.7 \end{bmatrix}$$

$$\cup \begin{bmatrix} 0.3 & 1 & 0.2 & 0.9 \\ 1 & 0.8 & 1 & 0.6 \\ 0.9 & 0.6 & 0.1 & 0.5 \\ 0.7 & 0.8 & 0.7 & 0.3 \end{bmatrix} \cup \begin{bmatrix} 1 & 0.6 & 0.3 & 1 & 0.8 \\ 0.3 & 1 & 0.6 & 0.7 & 1 \\ 0.6 & 0.5 & 0.4 & 0.5 & 0.5 \\ 0.8 & 0.2 & 0.1 & 0.2 & 0.1 \\ 0.5 & 0 & 0.3 & 0 & 0.4 \end{bmatrix}$$

$$\cup \begin{bmatrix} 0.3 & 1 & 0.8 \\ 1 & 0.3 & 0.4 \\ 0.7 & 0.6 & 0.2 \end{bmatrix} \cup \begin{bmatrix} 1 & 0.3 \\ 0.7 & 0.2 \end{bmatrix}$$

is a special fuzzy mixed square matrix.
Given

$$X = X_1 \cup X_2 \cup X_3 \cup X_4 \cup X_5 \cup X_6$$
$$= [1\ 0\ 0\ 0\ 0] \cup [0\ 1\ 0] \cup [0\ 0\ 1\ 0] \cup [0\ 0\ 0\ 0\ 1] \cup [0\ 0\ 1] \cup [1\ 0]$$

be the special fuzzy mixed row vector. To find X o M using the max min operation.

max min {X o M}

$$= \text{max min } \{(X_1 \cup X_2 \cup X_3 \cup X_4 \cup X_5 \cup X_6), (M_1 \cup M_2 \cup M_3 \cup M_4 \cup M_5 \cup M_6)\}$$

$$= \text{max min}(X_1, M_1) \cup \text{max min}(X_2, M_2) \cup \ldots \cup \text{max min}(X_6, M_6)$$

$$= \max\min\left\{ \begin{bmatrix} 1 & 0 & 0 & 0 & 0 \end{bmatrix}, \begin{bmatrix} 0.4 & 0.6 & 1 & 0 & 0.5 \\ 0.2 & 0 & 0.6 & 1 & 0.7 \\ 1 & 0.3 & 1 & 0 & 0.5 \\ 0.6 & 1 & 0.3 & 1 & 0.2 \\ 0 & 0.4 & 0.2 & 0 & 1 \\ 0.5 & 1 & 0.6 & 0.7 & 0.5 \end{bmatrix} \right\}$$

$$\cup \max \min \left\{ \begin{bmatrix} 0 & 1 & 0 \end{bmatrix}, \begin{bmatrix} 0.3 & 1 & 0.8 \\ 1 & 0.4 & 0.6 \\ 0.5 & 0 & 0.7 \end{bmatrix} \right\}$$

$$\cup \max \min \left\{ \begin{bmatrix} 0 & 0 & 1 & 0 \end{bmatrix}, \begin{bmatrix} 0.3 & 1 & 0.2 & 0.9 \\ 1 & 0.8 & 1 & 0.6 \\ 0.9 & 0.6 & 0.1 & 0.5 \\ 0.7 & 0.8 & 0.7 & 0.3 \end{bmatrix} \right\}$$

$$\cup \max \min \left\{ \begin{bmatrix} 0 & 0 & 0 & 0 & 1 \end{bmatrix}, \begin{bmatrix} 1 & 0.6 & 0.3 & 1 & 0.8 \\ 0.3 & 1 & 0.6 & 0.7 & 1 \\ 0.6 & 0.5 & 0.4 & 0.5 & 0.5 \\ 0.8 & 0.2 & 0.1 & 0.2 & 0.1 \\ 0.5 & 0 & 0.3 & 0 & 0.4 \end{bmatrix} \right\}$$

$$= \max \min \left\{ \begin{bmatrix} 0 & 0 & 1 \end{bmatrix}, \begin{bmatrix} 0.3 & 1 & 0.8 \\ 1 & 0.3 & 0.4 \\ 0.7 & 0.6 & 0.2 \end{bmatrix} \right\}$$

$$\cup \max \min \left\{ \begin{bmatrix} 1 & 0 \end{bmatrix}, \begin{bmatrix} 1 & 0.3 \\ 0.7 & 0.2 \end{bmatrix} \right\}$$

$$= [0.4\ 0.6\ 1\ 0\ 0.5] \cup [1\ 0.4\ 0.6] \cup [0.9\ 0.6\ 0.1\ 0.5] \cup [0.5\ 0\ 0.3\ 0\ 0.4] \cup [0.7\ 0.6\ 0.2] \cup [1\ 0.3]$$

$$= Y$$

$$= Y_1 \cup Y_2 \cup Y_3 \cup Y_4 \cup Y_5 \cup Y_6;$$

which is a special fuzzy mixed row vector. Now one can work with max min of Y o M.

This work of finding max min Y o M is left as an exercise for the reader.

Now we proceed on to show how we use the special max min operator and work with the special fuzzy rectangular matrix.

Let $P = P_1 \cup P_2 \cup \ldots \cup P_n$ $(n \geq 2)$ be a special fuzzy rectangular matrix where each $P_i$ is a $t \times s$ rectangular matrix $t \neq s$, $i = 1, 2, \ldots, n$.

Suppose

$X = X_1 \cup X_2 \cup \ldots \cup X_n$

be a special fuzzy row matrix where each $X_i$ is a $1 \times t$ fuzzy row vector $i = 1, 2, \ldots, n$. To find using the special max min operator

$$\begin{aligned} \text{max min}\{X \text{ o } P\} &= \text{max min } (X, P) \\ &= \text{max min}\{(X_1 \cup X_2 \cup \ldots \cup X_n), (P_1 \cup P_2 \cup \ldots \cup P_n)\} \end{aligned}$$

(say) now Y is again a special fuzzy row vector with each $Y_j$ a $1 \times s$ fuzzy row vector. $j = 1, 2, \ldots, n$.

Now using Y one can find

$$\begin{aligned} \text{max min } \{Y, P^T\} &= \text{max min } \{(Y_1 \cup Y_2 \cup \ldots \cup Y_n), (P_1^T \cup P_2^T \cup \ldots \cup P_n^T)\} \\ &= \text{max min } \{Y_1, P_1^T\} \cup \text{max min } \{Y_2, P_2^T\} \cup \ldots \cup \text{max min } \{Y_n, P_n^T\} \\ &= (Z_1 \cup Z_2 \cup \ldots \cup Z_n) \\ &= Z \text{ (say)}. \end{aligned}$$

Now once again Z is a special fuzzy row vector and each $Z_k$ is a $1 \times t$ fuzzy row vector for $k = 1, 2, \ldots, n$. Using $Z$ one can find max min Z, P and so on.

Now we illustrate the situation by an example.

***Example 1.2.22:*** Let $P = P_1 \cup P_2 \cup P_3 \cup P_4 \cup P_5$ be a special fuzzy rectangular matrix. Suppose $X = X_1 \cup X_2 \cup X_3 \cup X_4 \cup X_5$ be a special fuzzy row vector. Let each $P_i$ be a $6 \times 5$ rectangular fuzzy matrix $i = 1, 2, \ldots, 5$. Clearly each $X_i$ is a $1 \times 6$ fuzzy row vector $i = 1, 2, \ldots, 5$.

Given

$$P = P_1 \cup P_2 \cup P_3 \cup P_4 \cup P_5$$

$$= \begin{bmatrix} 0.3 & 1 & 0.2 & 1 & 0.8 \\ 1 & 0.7 & 0.4 & 0 & 0.1 \\ 0.5 & 1 & 0.8 & 0.9 & 1 \\ 0.7 & 0.9 & 1 & 0.6 & 0 \\ 0 & 1 & 0.7 & 0 & 0.9 \\ 1 & 0.8 & 1 & 0.9 & 0 \end{bmatrix} \cup$$

$$\begin{bmatrix} 1 & 0.8 & 0.6 & 0.4 & 0.2 \\ 0.3 & 0.5 & 1 & 0.7 & 0.9 \\ 0 & 0.1 & 0.3 & 0 & 0.5 \\ 1 & 0 & 0.4 & 0.6 & 1 \\ 0 & 0.7 & 0 & 0.9 & 0.1 \\ 0.2 & 0.8 & 1 & 0.4 & 0 \end{bmatrix} \cup$$

$$\begin{bmatrix} 0 & 1 & 0.9 & 0.3 & 1 \\ 0.1 & 0.2 & 1 & 0.5 & 0.8 \\ 0.7 & 1 & 0.7 & 1 & 0.1 \\ 1 & 0.6 & 0.4 & 0.9 & 0.5 \\ 0.5 & 1 & 0 & 1 & 0.4 \\ 0 & 0.7 & 0.8 & 0 & 1 \end{bmatrix} \cup$$

$$\begin{bmatrix} 0 & 0.8 & 0 & 1 & 0.7 \\ 0.9 & 1 & 0.2 & 0 & 0.4 \\ 0.1 & 0.2 & 1 & 0.3 & 0 \\ 0.5 & 0.6 & 0.7 & 1 & 0.8 \\ 1 & 0 & 0.9 & 0.4 & 1 \\ 0.2 & 0.4 & 1 & 0 & 0.6 \end{bmatrix} \cup$$

$$\begin{bmatrix} 0.5 & 0.7 & 1 & 0.9 & 0 \\ 1 & 0.9 & 0 & 0.1 & 0.2 \\ 0 & 1 & 0.2 & 0.4 & 0.6 \\ 0.8 & 0 & 1 & 0.5 & 0.4 \\ 1 & 0.9 & 0 & 0.4 & 1 \\ 0.7 & 0.6 & 0.4 & 1 & 0 \end{bmatrix}.$$

We see that each $P_i$ is a 6 × 5 fuzzy matrix. Given

$$X = [0\ 0\ 0\ 0\ 0\ 1] \cup [1\ 0\ 0\ 0\ 0\ 0] \cup [0\ 1\ 0\ 0\ 0\ 0] \cup [0\ 0\ 1\ 0\ 0\ 0] \cup [0\ 0\ 0\ 0\ 1\ 0]$$

to be a 1 × 6 special fuzzy row vector. Using special max min operator we find
max min {X, P}

$$= \max\min \{(X_1 \cup X_2 \cup X_3 \cup X_4 \cup X_5) , (P_1 \cup P_2 \cup P_3 \cup P_4 \cup P_5)\}$$

$$= \max\min(X_1, P_1) \cup \max\min(X_2, P_2) \cup \ldots \cup \max\min(X_5, P_5)$$

$$= \max\min \left\{ \begin{bmatrix} 0 & 0 & 0 & 0 & 0 & 1 \end{bmatrix}, \begin{bmatrix} 0.3 & 1 & 0.2 & 1 & 0.8 \\ 1 & 0.7 & 0.4 & 0 & 0.1 \\ 0.5 & 1 & 0.8 & 0.9 & 1 \\ 0.7 & 0.9 & 1 & 0.6 & 0 \\ 0 & 1 & 0.7 & 0 & 0.9 \\ 1 & 0.8 & 1 & 0.9 & 0 \end{bmatrix} \right\}$$

$$\cup \max\min \left\{ \begin{bmatrix} 1 & 0 & 0 & 0 & 0 & 0 \end{bmatrix}, \begin{bmatrix} 1 & 0.8 & 0.6 & 0.4 & 0.2 \\ 0.3 & 0.5 & 1 & 0.7 & 0.9 \\ 0 & 0.1 & 0.3 & 0 & 0.5 \\ 1 & 0 & 0.4 & 0.6 & 1 \\ 0 & 0.7 & 0 & 0.9 & 0.1 \\ 0.2 & 0.8 & 1 & 0.4 & 0 \end{bmatrix} \right\}$$

$$\cup \max\min\left\{\begin{bmatrix}0 & 1 & 0 & 0 & 0 & 0\end{bmatrix}, \begin{bmatrix} 0 & 1 & 0.9 & 0.3 & 1 \\ 0.1 & 0.2 & 1 & 0.5 & 0.8 \\ 0.7 & 1 & 0.7 & 1 & 0.1 \\ 1 & 0.6 & 0.4 & 0.9 & 0.5 \\ 0.5 & 1 & 0 & 1 & 0.4 \\ 0 & 0.7 & 0.8 & 0 & 1 \end{bmatrix}\right\}$$

$$\cup \max\min\left\{\begin{bmatrix}0 & 0 & 1 & 0 & 0 & 0\end{bmatrix}, \begin{bmatrix} 0 & 0.8 & 0 & 1 & 0.7 \\ 0.9 & 1 & 0.2 & 0 & 0.4 \\ 0.1 & 0.2 & 1 & 0.3 & 0 \\ 0.5 & 0.6 & 0.7 & 1 & 0.8 \\ 1 & 0 & 0.9 & 0.4 & 1 \\ 0.2 & 0.4 & 1 & 0 & 0.6 \end{bmatrix}\right\}$$

$$\cup \max\min\left\{\begin{bmatrix}0 & 0 & 0 & 0 & 1 & 0\end{bmatrix}, \begin{bmatrix} 0.5 & 0.7 & 1 & 0.9 & 0 \\ 1 & 0.9 & 0 & 0.1 & 0.2 \\ 0 & 1 & 0.2 & 0.4 & 0.6 \\ 0.8 & 0 & 1 & 0.5 & 0.4 \\ 1 & 0.9 & 0 & 0.4 & 1 \\ 0.7 & 0.6 & 0.4 & 1 & 0 \end{bmatrix}\right\}$$

$$= [1\ 0.8\ 1\ 0.9\ 0] \cup [1\ 0.8\ 0.6\ 0.4\ 0.2] \cup [0.1\ 0.2\ 1\ 0.5\ 0.8] \cup [0.1\ 0.2\ 1\ 0.3\ 0] \cup [1\ 0.9\ 0\ 0.4\ 1]$$

$$= Y_1 \cup Y_2 \cup Y_3 \cup Y_4 \cup Y_5$$

$$= Y.$$

Y is again a special fuzzy row vector, each $Y_i$ is $1 \times 5$ fuzzy vector; $i = 1, 2, \ldots, 5$.

Now using the special max min operator we find

$$Y \text{ o } P^T = \max\min\{(Y_1 \cup Y_2 \cup \ldots \cup Y_5), \left(P_1^T \cup P_2^T \cup \ldots \cup P_5^T\right)\}$$

$$= \max\min\{Y_1, P_1^T\} \cup \max\min\{Y_2, P_2^T\} \cup \ldots \cup \max\min\{Y_5, P_5^T\}$$

=

$$\max\min\left\{\begin{bmatrix}1 & 0.8 & 1 & 0.9 & 0\end{bmatrix}, \begin{bmatrix}0.3 & 1 & 0.5 & 0.7 & 0 & 1 \\ 1 & 0.7 & 1 & 0.9 & 1 & 0.8 \\ 0.2 & 0.4 & 0.8 & 1 & 0.7 & 1 \\ 1 & 0 & 0.9 & 0.6 & 0 & 0.9 \\ 0.8 & 0.1 & 1 & 0 & 0.9 & 0\end{bmatrix}\right\}$$

$\cup$

$$\max\min\left\{\begin{bmatrix}1 & 0.8 & 0.6 & 0.4 & 0.2\end{bmatrix}, \begin{bmatrix}1 & 0.3 & 0 & 1 & 0 & 0.2 \\ 0.8 & 0.5 & 0.1 & 0 & 0.7 & 0.8 \\ 0.6 & 1 & 0.3 & 0.4 & 0 & 1 \\ 0.4 & 0.7 & 0 & 0.6 & 0.9 & 0.4 \\ 0.2 & 0.9 & 0.5 & 1 & 0.1 & 0\end{bmatrix}\right\}$$

$\cup$

$$\max\min\left\{\begin{bmatrix}0.1 & 0.2 & 1 & 0.5 & 0.8\end{bmatrix}, \begin{bmatrix}0 & 0.1 & 0.7 & 1 & 0.5 & 0 \\ 1 & 0.2 & 1 & 0.6 & 1 & 0.7 \\ 0.9 & 1 & 0.7 & 0.4 & 0 & 0.8 \\ 0.3 & 0.5 & 1 & 0.9 & 1 & 0 \\ 1 & 0.8 & 0.1 & 0.5 & 0.4 & 1\end{bmatrix}\right\}$$

$\cup$

$$\max\min\left\{\begin{bmatrix}0.1 & 0.2 & 1 & 0.3 & 0\end{bmatrix}, \begin{bmatrix}0 & 0.9 & 0.1 & 0.5 & 1 & 0.2 \\ 0.8 & 1 & 0.2 & 0.6 & 0 & 0.4 \\ 0 & 0.2 & 1 & 0.7 & 0.9 & 1 \\ 1 & 0 & 0.3 & 1 & 0.4 & 0 \\ 0.7 & 0.4 & 0 & 0.8 & 1 & 0.6\end{bmatrix}\right\}$$

$\cup$

$$\max\min\left\{\begin{bmatrix}1 & 0.9 & 0 & 0.4 & 1\end{bmatrix}, \begin{bmatrix}0.5 & 1 & 0 & 0.8 & 1 & 0.7\\ 0.7 & 0.9 & 1 & 0 & 0.9 & 0.6\\ 1 & 0 & 0.2 & 1 & 0 & 0.4\\ 0.9 & 0.1 & 0.4 & 0.5 & 0.4 & 1\\ 0 & 0.2 & 0.6 & 0.4 & 1 & 0\end{bmatrix}\right\}$$

$= [0.9\ 1\ 0.9\ 1\ 0.8\ 1] \cup [1\ 0.6\ 0.3\ 1\ 0.7\ 0.8] \cup [0.9\ 1\ 0.7\ 0.5\ 0.5\ 0.8] \cup [0.3\ 0.2\ 1\ 0.7\ 0.9\ 1] \cup [0.7\ 1\ 0.9\ 0.8\ 1\ 0.7]$

$= Z_1 \cup Z_2 \cup Z_3 \cup Z_4 \cup Z_5$

$= Z$;

is once again a special fuzzy row vector and each $Z_i$ in this case is a $1 \times 6$ fuzzy vector. Now using Z we can calculate Z o P using the special max min operator.

Now we define special max min operator when we have the special fuzzy matrix to be a special fuzzy mixed rectangular matrix and the special fuzzy row vector which works on it is a special fuzzy mixed row vector.

Let $V = V_1 \cup V_2 \cup V_3 \cup V_4 \cup \ldots \cup V_n$ where V is a special fuzzy mixed rectangular matrix where each $V_i$ is a $s_i \times t_i$ fuzzy matrix $s_i \neq t_i$, i = 1, 2, …, n. Let $X = X_1 \cup X_2 \cup X_3 \cup X_4 \cup \ldots \cup X_n$ be a special fuzzy mixed row vector where each $X_i$ is a $1 \times s_i$ fuzzy row vector. To use the special max min operator and find

max min {X, V}

$= \max\min \{X, V\}$

$= \max\min \{(X_1 \cup X_2 \cup \ldots \cup X_n), (V_1 \cup V_2 \cup \ldots \cup V_n)\}$

$= \max\min \{X_1, V_1\} \cup \max\min \{X_2, V_2\} \cup \ldots \cup \max\min \{X_n, V_n\}$

$= Y_1 \cup Y_2 \cup \ldots \cup Y_n$

$= Y$

where Y is a special fuzzy mixed row vector with each $Y_j$ a 1 × $t_i$ fuzzy vector for j = 1, 2, …, n.

We now find max min {Y, $V^T$} and so on, which will once again give a special fuzzy mixed row vector and so on.

Now we illustrate this situation by the following example.

***Example 1.2.23:*** Let us consider the special fuzzy mixed rectangular matrix

$$V = V_1 \cup V_2 \cup V_3 \cup V_4 \cup V_5$$

$$= \begin{bmatrix} 0.4 & 1 & 0.3 \\ 1 & 0.2 & 0.5 \\ 0.6 & 1 & 0.2 \\ 0 & 0.3 & 1 \\ 0.7 & 1 & 0.3 \end{bmatrix} \cup \begin{bmatrix} 0.3 & 0 & 1 & 0.8 \\ 1 & 0.9 & 0 & 0.1 \\ 0.2 & 1 & 0.3 & 0.7 \\ 0.5 & 0 & 0.1 & 0 \\ 1 & 0.5 & 0.8 & 0 \\ 0.3 & 0 & 0.9 & 1 \\ 0.2 & 0.5 & 1 & 0.7 \end{bmatrix}$$

$$\cup \begin{bmatrix} 1 & 0.2 & 0.4 & 0.6 & 0.8 & 0 \\ 0 & 1 & 0.6 & 1 & 0.2 & 1 \\ 0.2 & 0.3 & 0.5 & 0 & 0.9 & 0.3 \\ 0.9 & 0.8 & 1 & 0.5 & 0 & 0.2 \end{bmatrix}$$

$$\cup \begin{bmatrix} 0.3 & 0.8 & 1 & 0.5 & 0.3 \\ 0.1 & 1 & 0 & 1 & 0.7 \\ 0.4 & 0 & 0.2 & 0.2 & 1 \\ 0.8 & 0.2 & 1 & 0 & 1 \\ 0.9 & 0 & 0.3 & 1 & 0.4 \\ 0.7 & 1 & 0 & 0.2 & 0.9 \end{bmatrix} \cup \begin{bmatrix} 0.2 & 0.3 & 0.1 \\ 1 & 0 & 1 \\ 0.2 & 1 & 0.3 \\ 1 & 0.2 & 1 \\ 0.6 & 0.7 & 0 \\ 0.7 & 1 & 0.4 \end{bmatrix}.$$

Let

X = [1 0 0 0 1] ∪ [0 0 0 0 0 0 1] ∪ [0 1 0 0] ∪ [0 0 0 1 0 0] ∪ [0 0 0 0 1 0]

be the given special fuzzy mixed row vector. To find using the special max min operator the value of X, V

max min (X, V)

= max min $\{(X_1 \cup X_2 \cup X_3 \cup X_4 \cup X_5), (V_1 \cup V_2 \cup V_3 \cup V_4 \cup V_5)\}$

= max min$(X_1, V_1) \cup$ max min$(X_2, V_2) \cup \ldots \cup$ max min$(X_5, V_5)$

$$= \max \min \left\{ \begin{bmatrix} 1 & 0 & 0 & 0 & 1 \end{bmatrix}, \begin{bmatrix} 0.4 & 1 & 0.3 \\ 1 & 0.2 & 0.5 \\ 0.6 & 1 & 0.2 \\ 0 & 0.3 & 1 \\ 0.7 & 1 & 0.3 \end{bmatrix} \right\}$$

$$\cup \max \min \left\{ \begin{bmatrix} 0 & 0 & 0 & 0 & 0 & 0 & 1 \end{bmatrix}, \begin{bmatrix} 0.3 & 0 & 1 & 0.8 \\ 1 & 0.9 & 0 & 0.1 \\ 0.2 & 1 & 0.3 & 0.7 \\ 0.5 & 0 & 0.1 & 0 \\ 1 & 0.5 & 0.8 & 0 \\ 0.3 & 0 & 0.9 & 1 \\ 0.2 & 0.5 & 1 & 0.7 \end{bmatrix} \right\}$$

$$\cup \max \min \left\{ \begin{bmatrix} 0 & 1 & 0 & 0 \end{bmatrix}, \begin{bmatrix} 1 & 0.2 & 0.4 & 0.6 & 0.8 & 0 \\ 0 & 1 & 0.6 & 1 & 0.2 & 1 \\ 0.2 & 0.3 & 0.5 & 0 & 0.9 & 0.3 \\ 0.9 & 0.8 & 1 & 0.5 & 0 & 0.2 \end{bmatrix} \right\}$$

$$\cup \max\min\left\{\begin{bmatrix}0 & 0 & 0 & 1 & 0 & 0\end{bmatrix}, \begin{bmatrix}0.3 & 0.8 & 1 & 0.5 & 0.3\\0.1 & 1 & 0 & 1 & 0.7\\0.4 & 0 & 0.2 & 0.2 & 1\\0.8 & 0.2 & 1 & 0 & 1\\0.9 & 0 & 0.3 & 1 & 0.4\\0.7 & 1 & 0 & 0.2 & 0.9\end{bmatrix}\right\}$$

$$\cup \max\min\left\{\begin{bmatrix}0 & 0 & 0 & 0 & 1 & 0\end{bmatrix}, \begin{bmatrix}0.2 & 0.3 & 0.1\\1 & 0 & 1\\0.2 & 1 & 0.3\\1 & 0.2 & 1\\0.6 & 0.7 & 0\\0.7 & 1 & 0.4\end{bmatrix}\right\}$$

= [0.7 1 0.3] ∪ [0.2 0.5 1 0.7] ∪ [0 1 0.6 1 0.2 1] ∪ [0.8 0.2 1 0 1] ∪ [0.6 0.7 0]

= $Y_1 \cup Y_2 \cup Y_3 \cup Y_4 \cup Y_5$

= Y.

We see Y is once again a special fuzzy mixed row vector. Now we see Y o V is not defined, only max min{Y, $V^T$} is defined so we find Y o $V^T$ using the special max min operator.

max min{Y, $V^T$}

= max min $\{(Y_1 \cup Y_2 \cup Y_3 \cup Y_4 \cup Y_5), (V_1^T \cup V_2^T \cup V_3^T \cup V_4^T \cup V_5^T)\}$

= max min $\{Y_1, V_1^T\}$ ∪ max min $\{Y_2, V_2^T\}$ ∪ … ∪ max min $\{Y_5, V_5^T\}$

$$= \max\min\left\{\begin{bmatrix}0.7 & 1 & 0.3\end{bmatrix}, \begin{bmatrix}0.4 & 1 & 0.6 & 0 & 0.7\\1 & 0.2 & 1 & 0.3 & 1\\0.3 & 0.5 & 0.2 & 1 & 0.3\end{bmatrix}\right\}$$

$$\cup \max\min \left\{ \begin{bmatrix} 0.2 & 0.5 & 1 & 0.7 \end{bmatrix}, \begin{bmatrix} 0.3 & 1 & 0.2 & 0.5 & 1 & 0.3 & 0.2 \\ 0 & 0.9 & 1 & 0 & 0.5 & 0 & 0.5 \\ 1 & 0 & 0.3 & 0.1 & 0.8 & 0.9 & 1 \\ 0.8 & 0.1 & 0.7 & 0 & 0 & 1 & 0.7 \end{bmatrix} \right\}$$

$$\cup \max\min \left\{ \begin{bmatrix} 0 & 1 & 0.6 & 1 & 0.2 & 1 \end{bmatrix}, \begin{bmatrix} 1 & 0 & 0.2 & 0.9 \\ 0.2 & 1 & 0.3 & 0.8 \\ 0.4 & 0.6 & 0.5 & 1 \\ 0.6 & 1 & 0 & 0.5 \\ 0.8 & 0.2 & 0.9 & 0 \\ 0 & 1 & 0.3 & 0.2 \end{bmatrix} \right\}$$

$$\cup \max\min \left\{ \begin{bmatrix} 0.8 & 0.2 & 1 & 0 & 1 \end{bmatrix}, \begin{bmatrix} 0.3 & 0.1 & 0.4 & 0.8 & 0.9 & 0.7 \\ 0.8 & 1 & 0 & 0.2 & 0 & 1 \\ 1 & 0 & 0.2 & 1 & 0.3 & 0 \\ 0.5 & 1 & 0.2 & 0 & 1 & 0.2 \\ 0.3 & 0.7 & 1 & 1 & 0.4 & 0.9 \end{bmatrix} \right\}$$

$$\cup \max\min \left\{ \begin{bmatrix} 0.6 & 0.7 & 0 \end{bmatrix}, \begin{bmatrix} 0.2 & 1 & 0.2 & 1 & 0.6 & 0.7 \\ 0.3 & 0 & 1 & 0.2 & 0.7 & 1 \\ 0.1 & 1 & 0.3 & 1 & 0 & 0.4 \end{bmatrix} \right\}$$

= [1 0.7 1 0.3 1] $\cup$ [1 0.5 0.7 0.2 0.8 0.9 1] $\cup$ [0.6 1 0.5 0.8] $\cup$ [1 0.7 1 1 0.8 0.9] $\cup$ [0.3 0.6 0.7 0.6 0.7 0.7]

= $T_1 \cup T_2 \cup T_3 \cup T_4 \cup T_5$

= T;

is a special fuzzy mixed row vector. Now using this T and V one can calculate the special max min value of {T, V} and so on.

Now we proceed on to define the special min max operator on the special class of fuzzy rectangular matrices and the special class of fuzzy mixed rectangular matrices.

Let $S = S_1 \cup S_2 \cup \ldots \cup S_n$ be the set of special fuzzy rectangular matrices with each $S_i$ a $t \times s$ rectangular matrix,

$$X = X_1 \cup X_2 \cup \ldots \cup X_n;$$

where $X_i$ is a $1 \times t$ fuzzy row vector $i = 1, 2, \ldots, n$, the special fuzzy row vector. To define the special min max operator using X and S.

min max (X, S)

$$= \text{min max } \{(X_1 \cup X_2 \cup \ldots \cup X_n), (S_1 \cup S_2 \cup \ldots \cup S_n)\}$$
$$= \text{min max } \{X_1, S_1\} \cup \text{min max } \{X_2, S_2\} \cup \ldots \cup \text{min max } \{X_n, S_n\}.$$

Each min max $\{X_i, S_i\}$ is calculated in pages 14-21 of this book.

Now
min max (X, S)
$$= Y_1 \cup Y_2 \cup \ldots \cup Y_n.$$
$$= Y$$

is a special fuzzy row vector; each $Y_i$ is a $1 \times s$ fuzzy row vector $i = 1, 2, \ldots, n$. Special min max operator is obtained for $\{Y, S^T\}$ and so on.

Now we illustrate this by the following example.

***Example 1.2.24:*** Let $W = W_1 \cup W_2 \cup \ldots \cup W_6$ be a special fuzzy rectangular matrix where

$$W = \begin{bmatrix} 0.1 & 1 & 0 & 0.5 & 0.6 & 1 \\ 0.5 & 0 & 0.7 & 1 & 1 & 0.8 \\ 0.4 & 1 & 0.2 & 0.3 & 0.4 & 1 \\ 1 & 0.6 & 1 & 0.7 & 1 & 0 \end{bmatrix} \cup$$

$$\begin{bmatrix} 1 & 0.6 & 0.9 & 0 & 0 & 0.9 \\ 0 & 0.8 & 1 & 0.6 & 0.9 & 0.2 \\ 0.2 & 1 & 0 & 1 & 0 & 0.3 \\ 0.4 & 0 & 0.5 & 0.2 & 0.4 & 1 \end{bmatrix} \cup$$

$$\begin{bmatrix} 1 & 0.2 & 0.6 & 0.4 & 0.8 & 0.9 \\ 0 & 1 & 0.7 & 0.9 & 0.5 & 0.3 \\ 0.5 & 0 & 0.2 & 0 & 1 & 0.1 \\ 0.6 & 0.8 & 1 & 0.5 & 0.3 & 0 \end{bmatrix} \cup$$

$$\begin{bmatrix} 0.8 & 0.1 & 0 & 1 & 0.9 & 0.5 \\ 0 & 1 & 0.7 & 0.5 & 0.7 & 0.3 \\ 0.6 & 0.5 & 1 & 0 & 0.9 & 0.4 \\ 0.7 & 0 & 0.6 & 0.3 & 1 & 0 \end{bmatrix} \cup$$

$$\begin{bmatrix} 1 & 0.8 & 0 & 1 & 0 & 0.8 \\ 0 & 0 & 1 & 1 & 0.1 & 1 \\ 0 & 1 & 0.1 & 0.3 & 0.5 & 0.7 \\ 0.3 & 0.2 & 0.4 & 0.2 & 1 & 0 \end{bmatrix} \cup$$

$$\begin{bmatrix} 0.3 & 1 & 0.2 & 0 & 0.1 & 1 \\ 0.5 & 0 & 0.4 & 1 & 0.2 & 0 \\ 0.7 & 1 & 0.6 & 0 & 0.3 & 0 \\ 0.9 & 0 & 0.8 & 1 & 0.4 & 1 \end{bmatrix}$$

Given

$$X = X_1 \cup X_2 \cup \ldots \cup X_6$$
$$= [1\ 0\ 0\ 0] \cup [0\ 1\ 0\ 0] \cup [0\ 0\ 1\ 0] \cup [0\ 0\ 0\ 1] \cup [0\ 1\ 1\ 0] \cup [1\ 0\ 0\ 1]$$

a special fuzzy row vector. To find the special min max value of X o W.

min max (X, W)

$$= \text{min max } \{(X_1 \cup X_2 \cup \ldots \cup X_6), (W_1 \cup W_2 \cup \ldots \cup W_6)\}$$
$$= \text{min max } \{X_1, W_1\} \cup \text{min max } \{X_2, W_2\} \cup \ldots \cup \text{min max } \{X_6, W_6\}$$

$$= \text{min max } \left\{ \begin{bmatrix} 1 & 0 & 0 & 0 \end{bmatrix}, \begin{bmatrix} 0.1 & 1 & 0 & 0.5 & 0.6 & 1 \\ 0.5 & 0 & 0.7 & 1 & 1 & 0.8 \\ 0.4 & 1 & 0.2 & 0.3 & 0.4 & 1 \\ 1 & 0.6 & 1 & 0.7 & 1 & 0 \end{bmatrix} \right\}$$

$$\cup \text{ min max } \left\{ \begin{bmatrix} 0 & 1 & 0 & 0 \end{bmatrix}, \begin{bmatrix} 1 & 0.6 & 0.9 & 0 & 0 & 0.9 \\ 0 & 0.8 & 1 & 0.6 & 0.9 & 0.2 \\ 0.2 & 1 & 0 & 1 & 0 & 0.3 \\ 0.4 & 0 & 0.5 & 0.2 & 0.4 & 0.1 \end{bmatrix} \right\}$$

$$\cup \text{ min max } \left\{ \begin{bmatrix} 0 & 0 & 1 & 0 \end{bmatrix}, \begin{bmatrix} 1 & 0.2 & 0.6 & 0.4 & 0.8 & 0.9 \\ 0 & 1 & 0.7 & 0.9 & 0.5 & 0.3 \\ 0.5 & 0 & 0.2 & 0 & 1 & 0.1 \\ 0.6 & 0.8 & 1 & 0.5 & 0.3 & 0 \end{bmatrix} \right\}$$

$$\cup \text{ min max } \left\{ \begin{bmatrix} 0 & 0 & 0 & 1 \end{bmatrix}, \begin{bmatrix} 0.8 & 0.1 & 0 & 1 & 0.9 & 0.5 \\ 0 & 1 & 0.7 & 0.5 & 0.7 & 0.3 \\ 0.6 & 0.5 & 1 & 0 & 0.9 & 0.4 \\ 0.7 & 0 & 0.6 & 0.3 & 1 & 0 \end{bmatrix} \right\}$$

$$\cup \text{ min max } \left\{ \begin{bmatrix} 0 & 1 & 1 & 0 \end{bmatrix}, \begin{bmatrix} 1 & 0.8 & 0 & 1 & 0 & 0.8 \\ 0 & 0 & 1 & 1 & 0.1 & 1 \\ 0 & 1 & 0.1 & 0.3 & 0.5 & 0.7 \\ 0.3 & 0.2 & 0.4 & 0.2 & 1 & 0 \end{bmatrix} \right\}$$

$$\cup \text{ min max } \left\{ \begin{bmatrix} 1 & 0 & 0 & 1 \end{bmatrix}, \begin{bmatrix} 0.3 & 1 & 0.2 & 0 & 0.1 & 1 \\ 0.5 & 0 & 0.4 & 1 & 0.2 & 0 \\ 0.7 & 1 & 0.6 & 0 & 0.3 & 0 \\ 0.9 & 0 & 0.8 & 1 & 0.4 & 1 \end{bmatrix} \right\}$$

= [0.4 0 0.2 0.3 0.4 0] ∪ [0.2 0 0 0 0 0.1] ∪ [0 0.2 0.6 0.4 0.3 0] ∪ [0 0.1 0 0 0.7 0.3] ∪ [0.3 0.2 0 0.2 0 0] ∪ [0.5 0 0.4 0 0.2 0]

= $Z_1 \cup Z_2 \cup Z_3 \cup Z_4 \cup Z_5 \cup Z_6$

= Z;

is a special fuzzy row vector and each $Z_i$ is a 1 × 6 fuzzy row vector; i = 1, 2, …, 6.

Now using Z we find min max of

min max $\{Z, W^T\}$

= min max $\{(Z_1 \cup Z_2 \cup \ldots \cup Z_6), (W_1^T \cup W_2^T \cup \ldots \cup W_6^T)$

= min max $\{Z_1, W_1^T\} \cup$ max mix $\{Z_2, W_2^T\} \cup \ldots \cup$ max mix $\{Z_6, W_6^T\}$

=

$$\text{min max} \left\{ \begin{bmatrix} 0.4 & 0 & 0.2 & 0.3 & 0.4 & 0 \end{bmatrix}, \begin{bmatrix} 0.1 & 0.5 & 0.4 & 1 \\ 1 & 0 & 1 & 0.6 \\ 0 & 0.7 & 0.2 & 1 \\ 0.5 & 1 & 0.3 & 0.7 \\ 0.6 & 1 & 0.4 & 1 \\ 1 & 0.8 & 1 & 0 \end{bmatrix} \right\}$$

$$\cup \min\max\left\{\begin{bmatrix}0.2 & 0 & 0 & 0 & 0 & 0.1\end{bmatrix}, \begin{bmatrix}1 & 0 & 0.2 & 0.4\\ 0.6 & 0.8 & 1 & 0\\ 0.9 & 1 & 0 & 0.5\\ 0 & 0.6 & 1 & 0.2\\ 0 & 0.9 & 0 & 0.4\\ 0.9 & 0.2 & 0.3 & 0.1\end{bmatrix}\right\}$$

$$\cup \min\max\left\{\begin{bmatrix}0 & 0.2 & 0.6 & 0.4 & 0.3 & 0\end{bmatrix}, \begin{bmatrix}1 & 0 & 0.5 & 0.6\\ 0.2 & 1 & 0 & 0.8\\ 0.6 & 0.7 & 0.2 & 1\\ 0.4 & 0.9 & 0 & 0.5\\ 0.8 & 0.5 & 1 & 0.3\\ 0.9 & 0.3 & 0.1 & 0\end{bmatrix}\right\}$$

$$\cup \min\max\left\{\begin{bmatrix}0 & 0.1 & 0 & 0 & 0.7 & 0.3\end{bmatrix}, \begin{bmatrix}0.8 & 0 & 0.6 & 0.7\\ 0.1 & 1 & 0.5 & 0\\ 0 & 0.7 & 1 & 0.6\\ 1 & 0.5 & 0 & 0.3\\ 0.9 & 0.7 & 0.9 & 1\\ 0.5 & 0.3 & 0.4 & 0\end{bmatrix}\right\}$$

$$\cup \min\max\left\{\begin{bmatrix}0.3 & 0.2 & 0 & 0.2 & 0 & 0\end{bmatrix}, \begin{bmatrix}1 & 0 & 0 & 0.3\\ 0.8 & 0 & 1 & 0.2\\ 0 & 1 & 0.1 & 0.4\\ 1 & 1 & 0.3 & 0.2\\ 0 & 0.1 & 0.5 & 1\\ 0.8 & 1 & 0.7 & 0\end{bmatrix}\right\}$$

$$\cup \min\max\left\{\begin{bmatrix}0.5 & 0 & 0.4 & 0 & 0.2 & 0\end{bmatrix}, \begin{bmatrix}0.3 & 0.5 & 0.7 & 0.9 \\ 1 & 0 & 1 & 0 \\ 0.2 & 0.4 & 0.6 & 0.8 \\ 0 & 1 & 0 & 1 \\ 0.1 & 0.2 & 0.3 & 0.4 \\ 1 & 0 & 0 & 1\end{bmatrix}\right\}$$

= [0.2 0 0.2 0] ∪ [0 0.2 0 0] ∪ [0.2 0 0.1 0] ∪ [0 0 0 0.1] ∪ [0 0.2 0.1 0] ∪ [0 0 0 0]

= $T_1 \cup T_2 \cup T_3 \cup T_4 \cup T_5 \cup T_6$

= T;

where T is again a special fuzzy row vector and one can find min max of {T, W} and so on.

Now we proceed on to work with special min max product on a special fuzzy mixed rectangular matrix $N = N_1 \cup N_2 \cup \ldots \cup N_t$ ($t \geq 2$) where each $N_i$ is a $p_i \times q_i$ rectangular fuzzy matrix i = 1, 2, …, t. ($p_i \neq q_i$).

Let $X = X_1 \cup X_2 \cup \ldots \cup X_t$ be a special mixed row vector where each $X_i$ is $1 \times p_i$ fuzzy row vector i = 1, 2, …, t.

To find the special min max product of X with N ,

min max (X, N)

= min max $\{(X_1 \cup X_2 \cup \ldots \cup X_t), (N_1 \cup N_2 \cup \ldots \cup N_t)\}$

= min max $\{X_1, N_1\}$ ∪ min max $\{X_2, N_2\}$ ∪ … ∪ min max $\{X_t, N_t\}$

= $Y_1 \cup Y_2 \cup \ldots \cup Y_t$.

= Y

where Y is a special fuzzy mixed row vector where each $Y_i$ is a $1 \times q_i$ fuzzy row vector i = 1, 2, …, t.

Now we will find special min max value of

min max $\{Y, N^T\}$

$$= \min\max \{(Y_1 \cup Y_2 \cup \ldots \cup Y_t), (N_1^T \cup N_2^T \cup \ldots \cup N_t^T)\}$$
$$= \min\max \{Y_1, N_1^T\} \cup \min\max \{Y_2, N_2^T\} \cup \ldots \cup \min\max \{Y_t, N_t^T\}$$
$$= Z_1 \cup Z_2 \cup \ldots \cup Z_t$$
$$= Z,$$

which is again a special fuzzy mixed row vector where each $Z_i$ is a $1 \times p_i$ fuzzy row vector; $i = 1, 2, \ldots, t$.

Now using Z and N we can find min max (Z, N) and so on.

We shall illustrate this by the following example.

***Example 1.2.25:*** Let $N = N_1 \cup N_2 \cup \ldots \cup N_5$ be a special fuzzy mixed rectangular matrix where

$$N = \begin{bmatrix} 0.3 & 0.2 & 0.1 & 0.5 & 0.2 & 0 & 1 \\ 0.1 & 0.7 & 0.8 & 1 & 0 & 1 & 0.6 \\ 0 & 1 & 0.9 & 0.2 & 0.7 & 0.9 & 0.8 \end{bmatrix} \cup$$

$$\begin{bmatrix} 0.9 & 0.3 & 0.5 & 0.9 & 0.7 \\ 0.6 & 0.7 & 1 & 0.1 & 0 \\ 1 & 0.8 & 0 & 1 & 0.2 \\ 0.7 & 1 & 0.5 & 0.4 & 0.5 \\ 0.4 & 0 & 1 & 0.8 & 0.1 \\ 0 & 0.7 & 0.2 & 1 & 0.9 \end{bmatrix} \cup$$

$$\begin{bmatrix} 0.2 & 0.4 \\ 1 & 0 \\ 0.7 & 0.3 \\ 0.2 & 1 \end{bmatrix} \cup$$

$$\begin{bmatrix} 0.3 & 0.2 & 0.8 & 1 \\ 0.2 & 1 & 0 & 0 \\ 0.7 & 0.9 & 1 & 1 \\ 1 & 0.8 & 0.4 & 0 \\ 0.5 & 1 & 0.2 & 0.7 \end{bmatrix} \cup$$

$$\begin{bmatrix} 0.2 & 0.3 & 0.5 \\ 1 & 0.1 & 0.8 \\ 0.4 & 1 & 0 \\ 0.5 & 0.7 & 1 \\ 0.8 & 0 & 0.5 \\ 1 & 0.6 & 0.7 \end{bmatrix}.$$

Suppose $X = X_1 \cup X_2 \cup \ldots \cup X_5$ be a special fuzzy mixed row vector where

$$X = [0\ 1\ 0] \cup [0\ 0\ 0\ 0\ 0\ 1] \cup [0\ 0\ 1\ 0] \cup [0\ 0\ 1\ 0\ 0] \cup [1\ 0\ 0\ 0\ 0\ 0].$$

Now using the special min max operator we calculate

min max (X, N)

$$= \text{min max } \{(X_1 \cup X_2 \cup \ldots \cup X_5), (N_1 \cup N_2 \cup \ldots \cup N_5)\}$$
$$= \text{min max } \{X_1, N_1\} \cup \text{min max } \{X_2, N_2\} \cup \ldots \cup \text{min max } \{X_5, N_5\}$$

$$= \text{min max} \left\{ \begin{bmatrix} 0 & 1 & 0 \end{bmatrix}, \begin{bmatrix} 0.3 & 0.2 & 0.1 & 0.5 & 0.2 & 0 & 1 \\ 0.1 & 0.7 & 0.8 & 1 & 0 & 1 & 0.6 \\ 0 & 1 & 0.9 & 0.2 & 0.7 & 0.9 & 0.8 \end{bmatrix} \right\}$$

$$\cup \min\max\left\{\begin{bmatrix}0 & 0 & 0 & 0 & 0 & 1\end{bmatrix}, \begin{bmatrix}0.9 & 0.3 & 0.5 & 0.9 & 0.7\\ 0.6 & 0.7 & 1 & 0.1 & 0\\ 1 & 0.8 & 0 & 1 & 0.2\\ 0.7 & 1 & 0.5 & 0.4 & 0.5\\ 0.4 & 0 & 1 & 0.8 & 0.1\\ 0 & 0.7 & 0.2 & 1 & 0.9\end{bmatrix}\right\}$$

$$\cup \min\max\left\{\begin{bmatrix}0 & 0 & 1 & 0\end{bmatrix}, \begin{bmatrix}0.2 & 0.4\\ 1 & 0\\ 0.7 & 0.3\\ 0.2 & 1\end{bmatrix}\right\}$$

$$\cup \min\max\left\{\begin{bmatrix}0 & 0 & 1 & 0 & 0\end{bmatrix}, \begin{bmatrix}0.3 & 0.2 & 0.8 & 1\\ 0.2 & 1 & 0 & 0\\ 0.7 & 0.9 & 1 & 1\\ 1 & 0.8 & 0.4 & 0\\ 0.5 & 1 & 0.2 & 0.7\end{bmatrix}\right\}$$

$$\cup \min\max\left\{\begin{bmatrix}1 & 0 & 0 & 0 & 0 & 0\end{bmatrix}, \begin{bmatrix}0.2 & 0.3 & 0.5\\ 1 & 0.1 & 0.8\\ 0.4 & 1 & 0\\ 0.5 & 0.7 & 1\\ 0.8 & 0 & 0.5\\ 1 & 0.6 & 0.7\end{bmatrix}\right\}$$

= [0 0.2 0.1 0.2 0.2 0 0.8] $\cup$ [0.4 0 0 0.1 0] $\cup$ [0.2 0]
$\cup$ [0.2 0.2 0 0] $\cup$ [0.2 0 0]

= $Z_1 \cup Z_2 \cup \ldots \cup Z_5$

= Z;

is the special fuzzy mixed row vector. Clearly min max {Z, N} is not defined so we have to define only min max $\{Z, N^T\}$ and find the special min max of $\{Z, N^T\}$.

min max $\{(Z, N^T)\}$

$= \text{min max } \{(Z_1 \cup Z_2 \cup \ldots \cup Z_5), (N_1^T \cup N_2^T \cup \ldots \cup N_5^T)\}$

$= \text{min max } \{Z_1, N_1^T\} \cup \text{min max } \{Z_2, N_2^T\} \cup \ldots \cup \text{min max } \{Z_5, N_5^T\}$

=

$$\text{min max}\left\{\begin{bmatrix}0 & 0.2 & 0.1 & 0.2 & 0.2 & 0 & 0.8\end{bmatrix}, \begin{bmatrix}0.3 & 0.1 & 0 \\ 0.2 & 0.7 & 1 \\ 0.1 & 0.8 & 0.9 \\ 0.5 & 1 & 0.2 \\ 0.2 & 0 & 0.7 \\ 0 & 1 & 0.9 \\ 1 & 0.6 & 0.8\end{bmatrix}\right\} \cup$$

$$\text{min max}\left\{\begin{bmatrix}0.4 & 0 & 0 & 0.1 & 0\end{bmatrix}, \begin{bmatrix}0.9 & 0.6 & 1 & 0.7 & 0.4 & 0 \\ 0.3 & 0.7 & 0.8 & 1 & 0 & 0.7 \\ 0.5 & 1 & 0 & 0.5 & 1 & 0.2 \\ 0.9 & 0.1 & 1 & 0.4 & 0.8 & 1 \\ 0.7 & 0 & 0.2 & 0.5 & 0.1 & 0.9\end{bmatrix}\right\}$$

$$\cup \text{ min max } \left\{\begin{bmatrix}0.2 & 0\end{bmatrix}, \begin{bmatrix}0.2 & 1 & 0.7 & 0.2 \\ 0.4 & 0 & 0.3 & 1\end{bmatrix}\right\}$$

$$\cup \text{ min max } \left\{\begin{bmatrix}0.2 & 0.2 & 0 & 0\end{bmatrix}, \begin{bmatrix}0.3 & 0.2 & 0.7 & 1 & 0.5 \\ 0.2 & 1 & 0.9 & 0.8 & 1 \\ 0.8 & 0 & 1 & 0.4 & 0.2 \\ 1 & 0 & 1 & 0 & 0.7\end{bmatrix}\right\}$$

$$\cup \min\max \left\{ \begin{bmatrix} 0.2 & 0 & 0 \end{bmatrix}, \begin{bmatrix} 0.2 & 1 & 0.4 & 0.5 & 0.8 & 1 \\ 0.3 & 0.1 & 1 & 0.7 & 0 & 0.6 \\ 0.5 & 0.8 & 0 & 1 & 0.5 & 0.7 \end{bmatrix} \right\}$$

= [0 0.1 0] $\cup$ [0.3 0 0 0.4 0 0.2] $\cup$ [0.2 0 0.3 0.2] $\cup$ [0.2 0 0.7 0 0.2] $\cup$ [0.2 0.1 0 0.5 0 0.6]
= $Z_1 \cup Z_2 \cup \ldots \cup Z_5$
= Z.

We see Z is a special fuzzy mixed row vector. Now we can find min max {Z, N} and so on.

Now we proceed on to work with the special min max operator on special fuzzy mixed matrix. Let $M = M_1 \cup M_2 \cup \ldots \cup M_s$ ($s \geq 2$) be the special fuzzy mixed matrix, i.e., some $M_i$ are $n_i \times n_i$ fuzzy square matrices and some of the $M_j$' s are $m_j \times t_j$ fuzzy rectangular matrices; $1 \leq i, j \leq s$ ($m_j \neq t_j$).
Suppose

$$X = X_1 \cup X_2 \cup \ldots \cup X_s$$

be a special fuzzy mixed row vector some $X_i$' s are $1 \times n_i$ fuzzy row vectors and other $X_j$'s are $1 \times m_j$ fuzzy row vectors. $1 \leq i, j \leq s$. To use the special min max operator and find min max {X, M}; i.e.,

min max (X, M)

= min max $\{(X_1 \cup X_2 \cup \ldots \cup X_s), (M_1 \cup M_2 \cup \ldots \cup M_s)\}$
= min max $\{X_1, M_1\}$ $\cup$ min max $\{X_2, M_2\}$ $\cup \ldots \cup$ min max $\{X_s, M_s\}$
= $Y_1 \cup Y_2 \cup \ldots \cup Y_s$.
= Y

is once again a special fuzzy mixed row vector. Now we see min max {Y, M} is not compatible under special min max operator so we define special transpose of M since M is a very peculiar fuzzy mixed matrix containing both square and rectangular matrices.

Suppose $M_1 \cup M_2 \cup \ldots \cup M_s$; be a special mixed matrix (s $\geq 2$) and if some of the $M_i$'s are $n_i \times n_i$ square fuzzy matrices and some of the $M_j$'s are $m_j \times t_j$ ($m_j \neq t_j$) rectangular fuzzy matrices $1 \leq i, j \leq s$ then the special transpose of M is defined and denoted by $M^{ST}$;

$$M^{ST} = M_1^T \cup \ldots \cup M_i \cup \ldots \cup M_j^T \cup \ldots \cup M_s$$

where if $M_i$ is a square fuzzy matrix we do not take its transpose in $M^{ST}$ only for the rectangular fuzzy matrices we take the transpose. This transpose is defined in a special way is called the special transpose.

We just for the sake of better understanding exhibit the notion of special fuzzy transpose by the following example.

Suppose

$$M = M_1 \cup M_2 \cup \ldots \cup M_5$$

$$= \begin{bmatrix} 0.3 & 1 & 0 \\ 0.8 & 0.9 & 1 \\ 0.7 & 0.4 & 0.5 \end{bmatrix} \cup \begin{bmatrix} 0.7 & 0.8 \\ 1 & 0 \\ 0 & 0.9 \\ 0.2 & 0.7 \\ 0.3 & 0.2 \\ 0.4 & 0.1 \\ 0.9 & 1 \end{bmatrix}$$

$$\cup \begin{bmatrix} 0.3 & 0.5 & 0.8 & 0.9 \\ 1 & 0 & 0.2 & 1 \\ 0.4 & 0.7 & 1 & 0.4 \\ 0.9 & 0.8 & 0.5 & 1 \end{bmatrix} \cup \begin{bmatrix} 0.9 & 1 & 0.2 & 0.7 & 0.4 & 0.5 \\ 1 & 0 & 0.9 & 1 & 0.7 & 0 \\ 0.7 & 1 & 0.2 & 0 & 0.6 & 1 \end{bmatrix}$$

$$\cup \begin{bmatrix} 0.8 & 1 & 0.7 & 0.5 & 0.6 \\ 1 & 0 & 0.6 & 1 & 0.9 \\ 0.9 & 1 & 0.2 & 0.5 & 1 \\ 0.6 & 0.3 & 1 & 0.7 & 0.4 \\ 0.4 & 0.2 & 0.5 & 1 & 0.8 \\ 0.7 & 1 & 1 & 0 & 0.1 \end{bmatrix}$$

be a special fuzzy mixed matrix. The special transpose of M denoted by

$$M^{ST} = (M_1 \cup M_2 \cup \ldots \cup M_5)^{ST}$$
$$= M_1 \cup M^T{}_2 \cup M_3 \cup M^T{}_4 \cup M^T{}_5.$$

$$= \begin{bmatrix} 0.3 & 1 & 0 \\ 0.8 & 0.9 & 1 \\ 0.7 & 0.4 & 0.5 \end{bmatrix} \cup$$

$$\begin{bmatrix} 0.7 & 1 & 0 & 0.2 & 0.3 & 0.4 & 0.9 \\ 0.8 & 0 & 0.9 & 0.7 & 0.2 & 0.1 & 1 \end{bmatrix} \cup$$

$$\begin{bmatrix} 0.3 & 0.5 & 0.8 & 0.9 \\ 1 & 0 & 0.2 & 1 \\ 0.4 & 0.7 & 1 & 0.4 \\ 0.9 & 0.8 & 0.5 & 1 \end{bmatrix} \cup$$

$$\begin{bmatrix} 0.9 & 1 & 0.7 \\ 1 & 0 & 1 \\ 0.2 & 0.9 & 0.2 \\ 0.7 & 1 & 0 \\ 0.4 & 0.7 & 0.6 \\ 0.5 & 0 & 1 \end{bmatrix} \cup$$

$$\begin{bmatrix} 0.8 & 1 & 0.9 & 0.6 & 0.4 & 0.7 \\ 1 & 0 & 1 & 0.3 & 0.2 & 1 \\ 0.7 & 0.6 & 0.2 & 1 & 0.5 & 1 \\ 0.5 & 1 & 0.5 & 0.7 & 1 & 0 \\ 0.6 & 0.9 & 1 & 0.4 & 0.8 & 0.1 \end{bmatrix}$$

$$= \quad M_1 \cup M_2^T \cup M_3 \cup M_4^T \cup M_5^T \; .$$

Thus one can see that special transpose of a special fuzzy mixed matrix is different form the usual transpose of any special mixed square or rectangular matrix.

Now we calculate the special min max operator X on with M where $M = M_1 \cup M_2 \cup \ldots \cup M_s$ is a special fuzzy mixed matrix and $X_1 \cup X_2 \cup \ldots \cup X_s$ is the special fuzzy mixed row vector with conditions described just above in the example.

If min max $\{X, M\} = Y_1 \cup Y_2 \cup \ldots \cup Y_S = Y$ be the special fuzzy mixed row vector then we find the special value of Y on $M^{ST}$ using the special min max operator as follows.

min max $\{Y, M^{ST}\}$

$$\begin{aligned} &= \text{min max } \{(Y_1 \cup Y_2 \cup \ldots \cup Y_s), (M_1 \cup M_2^T \cup \ldots \\ &\quad \cup M_i \cup \ldots \cup M_j^T \cup \ldots \cup \quad M_s)\} \\ &= \text{min max } \{Y_1, M_1\} \cup \text{min max } \{Y_2, M_2^T\} \cup \ldots \cup \\ &\quad \text{min max } \{Y_i, M_i\} \cup \ldots \cup \text{min max } \{Y_j, M_j^T\} \cup \ldots \\ &\quad \cup \text{min max } \{Y_s, M_s\} \\ &= Z_1 \cup Z_2 \cup \ldots \cup Z_s \\ &= Z; \end{aligned}$$

where Z is a special fuzzy mixed row vector. Now we find min max {Z, M} and so on.

Now we shall illustrate this situation by the following example.

***Example 1.2.26:*** Let $M = M_1 \cup M_2 \cup \ldots \cup M_6$ be a special fuzzy mixed matrix where

$$M = \begin{bmatrix} 0.4 & 1 & 0 & 0.2 \\ 1 & 0.7 & 1 & 0.3 \\ 0.8 & 1 & 0.4 & 1 \\ 1 & 0.7 & 0 & 0.2 \end{bmatrix} \cup \begin{bmatrix} 0.8 & 1 & 0.7 \\ 1 & 0.5 & 0 \\ 0.9 & 1 & 0.9 \\ 0.1 & 0.6 & 1 \\ 0.5 & 1 & 0.2 \\ 0.3 & 0.1 & 0.7 \\ 0.5 & 1 & 0.4 \end{bmatrix} \cup$$

$$\begin{bmatrix} 0.7 & 0.6 & 1 & 0.7 & 0.9 & 0.3 & 0.8 \\ 0 & 1 & 0 & 0.5 & 1 & 0.6 & 1 \\ 1 & 0.7 & 0.9 & 0.2 & 0.7 & 1 & 0.9 \\ 0.8 & 0.5 & 1 & 0.3 & 1 & 0.8 & 0.4 \end{bmatrix} \cup$$

$$\begin{bmatrix} 0.6 & 1 & 0.8 \\ 1 & 0.7 & 1 \\ 0.2 & 0.8 & 0 \end{bmatrix} \cup \begin{bmatrix} 0.8 & 1 & 0.9 & 1 \\ 1 & 0.7 & 0.1 & 0.4 \\ 0.5 & 1 & 0.5 & 1 \\ 0.9 & 0 & 1 & 0.7 \\ 0.2 & 1 & 0.8 & 0.3 \\ 0.7 & 0 & 1 & 0.2 \\ 1 & 0.7 & 0.6 & 0.5 \\ 0.3 & 1 & 0.4 & 1 \end{bmatrix} \cup$$

$$\begin{bmatrix} 0.3 & 1 & 0.4 & 0.9 \\ 1 & 0 & 0.2 & 1 \\ 0.5 & 0.7 & 1 & 0 \\ 1 & 0.9 & 0 & 1 \\ 0.8 & 1 & 0.7 & 0.4 \end{bmatrix}.$$

Suppose the special fuzzy mixed row vector

$$X = X_1 \cup X_2 \cup \ldots \cup X_6$$

$= [1\ 0\ 0\ 0] \cup [0\ 1\ 0\ 0\ 0\ 0\ 0] \cup [0\ 0\ 1\ 0] \cup [0\ 0\ 1] \cup [0\ 0\ 0\ 0\ 1\ 0\ 0\ 0] \cup [1\ 0\ 0\ 0\ 0]$

operates on M under special min max product to find

min max {X, M}

$= \min\max \{(X_1 \cup X_2 \cup \ldots \cup X_6), (M_1 \cup M_2 \cup \ldots \cup M_6)\}$

$= \min\max \{X_1, M_1\} \cup \min\max \{X_2, M_2\} \cup \ldots \cup \min\max \{X_6, M_6\}$

$$= \min\max \left\{ \begin{bmatrix} 1 & 0 & 0 & 0 \end{bmatrix}, \begin{bmatrix} 0.4 & 1 & 0 & 0.2 \\ 1 & 0.7 & 1 & 0.3 \\ 0.8 & 1 & 0.4 & 1 \\ 1 & 0.7 & 0 & 0.2 \end{bmatrix} \right\}$$

$$\cup \min\max \left\{ \begin{bmatrix} 0 & 1 & 0 & 0 & 0 & 0 & 0 \end{bmatrix}, \begin{bmatrix} 0.8 & 1 & 0.7 \\ 1 & 0.5 & 0 \\ 0.9 & 1 & 0.9 \\ 0.1 & 0.6 & 1 \\ 0.5 & 1 & 0.2 \\ 0.3 & 0.1 & 0.7 \\ 0.5 & 1 & 0.4 \end{bmatrix} \right\}$$

$$\cup \min\max \left\{ \begin{bmatrix} 0 & 0 & 1 & 0 \end{bmatrix}, \begin{bmatrix} 0.7 & 0.6 & 1 & 0.7 & 0.9 & 0.3 & 0.8 \\ 0 & 1 & 0 & 0.5 & 1 & 0.6 & 1 \\ 1 & 0.7 & 0.9 & 0.2 & 0.7 & 1 & 0.9 \\ 0.8 & 0.5 & 1 & 0.3 & 1 & 0.8 & 0.4 \end{bmatrix} \right\}$$

$$\cup \min\max \left\{ \begin{bmatrix} 0 & 0 & 1 \end{bmatrix}, \begin{bmatrix} 0.6 & 1 & 0.8 \\ 1 & 0.7 & 1 \\ 0.2 & 0.8 & 0 \end{bmatrix} \right\}$$

$$\cup \min\max\left\{\begin{bmatrix}0 & 0 & 0 & 0 & 1 & 0 & 0 & 0\end{bmatrix}, \begin{bmatrix}0.8 & 1 & 0.9 & 1 \\ 1 & 0.7 & 0.1 & 0.4 \\ 0.5 & 1 & 0.5 & 1 \\ 0.9 & 0 & 1 & 0.7 \\ 0.2 & 1 & 0.8 & 0.3 \\ 0.7 & 0 & 1 & 0.2 \\ 1 & 0.7 & 0.6 & 0.5 \\ 0.3 & 1 & 0.4 & 1\end{bmatrix}\right\}$$

$$\cup \min\max\left\{\begin{bmatrix}1 & 0 & 0 & 0 & 0\end{bmatrix}, \begin{bmatrix}0.3 & 1 & 0.4 & 0.9 \\ 1 & 0 & 0.2 & 1 \\ 0.5 & 0.7 & 1 & 0 \\ 1 & 0.9 & 0 & 1 \\ 0.8 & 1 & 0.7 & 0.4\end{bmatrix}\right\}$$

= {[0.8 0.7 0 0.2] ∪ [0.1 0.1 0.2] ∪ [0 0.5 0 0.3 0.9 0.3 0.4] ∪ [0.6 0.7 0.8] ∪ [0.3 0 0.1 0.2] ∪ [0.5 0 0 0]}
= $Y_1 \cup Y_2 \cup \ldots \cup Y_6$.
= Y

is again a special fuzzy mixed row vector.
Now we use special min max function to find the value of

Y on $M^{ST}$

= min max $(Y, M^{ST})$.
= min max $\{(Y_1 \cup Y_2 \cup \ldots \cup Y_6), (M_1 \cup M_2^T \cup M_3^T \cup M_4 \cup M_5^T \cup M_6^T)$.
= min max $(Y_1, M_1) \cup$ min max $(Y_2, M_2^T) \cup$ min max $(Y_3, M_3^T) \cup$ min max $(Y_4, M_4) \cup$ min max $(Y_5, M_5^T) \cup$ min max $(Y_6, M_6^T)$.

$$= \min\max\left\{[0.8 \quad 0.7 \quad 0 \quad 0.2],\ \begin{bmatrix} 0.4 & 1 & 0 & 0.2 \\ 1 & 0.7 & 1 & 0.3 \\ 0.8 & 1 & 0.4 & 1 \\ 1 & 0.7 & 0 & 0.2 \end{bmatrix}\right\}$$

$$\cup$$

$$\min\max\left\{[0.1 \quad 0.1 \quad 0.2],\ \begin{bmatrix} 0.8 & 1 & 0.9 & 0.1 & 0.5 & 0.3 & 0.5 \\ 1 & 0.5 & 1 & 0.6 & 1 & 0.1 & 1 \\ 0.7 & 0 & 0.9 & 1 & 0.2 & 0.7 & 0.4 \end{bmatrix}\right\}$$

$$\cup$$

$$\min\max\left\{[0 \quad 0.5 \quad 0 \quad 0.3 \quad 0.9 \quad 0.3 \quad 0.4],\ \begin{bmatrix} 0.7 & 0 & 1 & 0.8 \\ 0.6 & 1 & 0.7 & 0.5 \\ 1 & 0 & 0.9 & 1 \\ 0.7 & 0.5 & 0.2 & 0.3 \\ 0.9 & 1 & 0.7 & 1 \\ 0.3 & 0.6 & 1 & 0.8 \\ 0.8 & 1 & 0.9 & 0.4 \end{bmatrix}\right\}$$

$$\cup \min\max\left\{[0.6 \quad 0.7 \quad 0.8],\ \begin{bmatrix} 0.6 & 1 & 0.8 \\ 1 & 0.7 & 1 \\ 0.2 & 0.8 & 0 \end{bmatrix}\right\}$$

$$\cup \min\max$$

$$\left\{[0.3 \quad 0 \quad 0.1 \quad 0.2],\ \begin{bmatrix} 0.8 & 1 & 0.5 & 0.9 & 0.2 & 0.7 & 1 & 0.3 \\ 1 & 0.7 & 1 & 0 & 1 & 0 & 0.7 & 1 \\ 0.9 & 0.1 & 0.5 & 1 & 0.8 & 1 & 0.6 & 0.4 \\ 1 & 0.4 & 1 & 0.7 & 0.3 & 0.2 & 0.5 & 1 \end{bmatrix}\right\}$$

$$\cup \min\max\left\{[0.5 \quad 0 \quad 0 \quad 0],\ \begin{bmatrix} 0.3 & 1 & 0.5 & 1 & 0.8 \\ 1 & 0 & 0.7 & 0.9 & 1 \\ 0.4 & 0.2 & 1 & 0 & 0.7 \\ 0.9 & 1 & 0 & 1 & 0.4 \end{bmatrix}\right\}$$

$$= [0.8\ 0.7\ 0.2\ 0.2] \cup [0.7\ 0.2\ 0.9\ 0.1\ 0.2\ 0.1\ 0.4] \cup [0.3\ 0\ 0.3\ 0.3] \cup [0.6\ 0.7\ 0.8] \cup [0.8\ 0.1\ 0.5\ 0\ 0.3\ 0\ 0.5\ 0.3] \cup [0.4\ 0\ 0\ 0\ 0.4]$$
$$= P_1 \cup P_2 \cup P_3 \cup P_4 \cup P_5 \cup P_6$$
$$= P$$

is a special fuzzy mixed row vector. Now we can find min max {P, M} and so on.

On similar lines we can work with special max min operator also. Now having defined some of the major operations we now proceed on to define mixed operations on these special fuzzy matrices.

**DEFINITION 1.2.10:** *Let $S = S_1 \cup S_2 \cup \ldots \cup S_n$ ($n \geq 2$) be a special fuzzy mixed matrix $X = X_1 \cup X_2 \cup \ldots \cup X_n$ ($n \geq 2$) be a special fuzzy mixed row vector. We define a new operation of X on S called the special mixed operation denoted by $o_m$.*

$$X\ o_m\ S = (X_1 \cup X_2 \cup \ldots \cup X_n)\ o_m\ (S_1 \cup S_2 \cup \ldots \cup S_n)$$
$$= X_1\ o_m\ S_1 \cup X_2\ o_m\ S_2 \cup \ldots \cup X_n\ o_m\ S_n.$$

*is defined as follows.*

If $S_i$ is a square fuzzy matrix with entries from the set {–1, 0, 1} then $X_i\ o_m\ S_i$ is the operation described in pages 20-1 of this book. If $S_j$ is a rectangular fuzzy matrix with entries from the set {–1, 0, 1} then $X_j\ o_m\ S_j$ is the operation described in page 23 of this book. Suppose $S_k$ is a square fuzzy matrix with entries from [0, 1] then $X_k\ o_m\ S_k$ is the min max or max min operation defined in pages 14-5 of this book.

If $S_t$ is a rectangular fuzzy matrix with entries from [0, 1] then $X_t\ o_m\ S_t$ is the min max operation described in page 15 max min operation defined in page 14 of this book for $1 \leq i, j \leq k, t \leq n$. Thus '$o_m$' defined between X and S will be known as the special mixed operation and denoted by $o_m$.

We will illustrate this by the following example.

***Example 1.2.27:*** Let S = $S_1 \cup S_2 \cup \ldots \cup S_5$ be a special fuzzy mixed matrix where

$$S = \begin{bmatrix} 0 & 1 & -1 & 0 & 0 \\ 1 & 0 & 0 & 0 & 1 \\ -1 & 0 & 0 & 1 & 0 \\ 0 & 0 & 1 & 0 & 0 \\ 0 & -1 & 0 & 1 & 0 \end{bmatrix} \cup$$

$$\begin{bmatrix} 0.3 & 0.8 & 1 & 0.5 & 0 & 0.7 \\ 0.5 & 0.7 & 0 & 1 & 0.6 & 1 \\ 1 & 0.5 & 0.2 & 0.7 & 0.2 & 0 \end{bmatrix} \cup$$

$$\begin{bmatrix} 1 & -1 & 0 & 1 \\ 0 & 0 & 1 & 0 \\ 0 & 1 & 0 & 1 \\ 1 & 0 & 1 & -1 \\ 0 & 1 & 0 & 0 \\ 0 & 1 & 1 & 0 \end{bmatrix} \cup \begin{bmatrix} 0.9 & 0.2 & 1 & 0 & 0.5 \\ 0.3 & 0 & 0.5 & 1 & 0.7 \\ 1 & 0.2 & 1 & 0 & 0.3 \\ 0 & 1 & 0.3 & 0.2 & 1 \\ 0.9 & 0.7 & 0.6 & 0.7 & 0 \end{bmatrix} \cup$$

$$\begin{bmatrix} 1 & 0 & 1 & 1 & 1 & 0 & 0 & 1 & -1 \\ 0 & 1 & 0 & 0 & 0 & 1 & 0 & 0 & 0 \\ 1 & -1 & 0 & 0 & 0 & 0 & 0 & -1 & 1 \\ 0 & 0 & 1 & 0 & 1 & 0 & 1 & 0 & 0 \\ -1 & 0 & 0 & -1 & 1 & 0 & 0 & 0 & 0 \\ 0 & 0 & -1 & 0 & 0 & 1 & 1 & 0 & 0 \end{bmatrix}.$$

Let

$$\begin{aligned} X &= X_1 \cup X_2 \cup X_3 \cup X_4 \cup X_5 \\ &= [1\ 0\ 0\ 0\ 0] \cup [0\ 0\ 1] \cup [0\ 0\ 0\ 0\ 0\ 1] \cup [0\ 1\ 0\ 0\ 1] \cup \\ &\quad [0\ 0\ 0\ 1\ 0\ 0] \end{aligned}$$

be special fuzzy mixed row vector. We find X $o_m$ S using the special mixed operation.

$$X\ o_m\ S = (X_1 \cup X_2 \cup X_3 \cup X_4 \cup X_5)\ o_m\ (S_1 \cup S_2 \cup \ldots \cup S_5)$$
$$= X_1\, o_m^1\, S_1 \cup X_2\, o_m^2\, S_2 \cup X_3\, o_m^3\, S_3 \cup X_4\, o_m^4\, S_4 \cup X_5\, o_m^5\, S_5$$

where $o_m^1$ is the thresholding and updating resultant row vectors after usual matrix multiplication. $o_m^2$ is the min max operator between $X_2$ and $S_2$. $o_m^3$ and $o_m^4$ are usual matrix multiplication which are updated and thresholded and are found sequentially by finding $X_i$ o $S_i$ and next $Y_i$ o $S_i^T$ and so on. $o_m^5$ is the max min operation.

Now we explicitly show how the special mixed operation $o_m$ function works:

$$X\ o_m\ S = X_1\, o_m^1\, S_1 \cup \ldots \cup X_5\, o_m^5\, S_5$$

$$= \begin{bmatrix} 1 & 0 & 0 & 0 & 0 \end{bmatrix}\ o_m^1 \begin{bmatrix} 0 & 1 & -1 & 0 & 0 \\ 1 & 0 & 0 & 0 & 1 \\ -1 & 0 & 0 & 1 & 0 \\ 0 & 0 & 1 & 0 & 0 \\ 0 & -1 & 0 & 1 & 0 \end{bmatrix} \quad (o_m^1 = \text{'o'})$$

$$\text{min max} \left\{ \begin{bmatrix} 0 & 0 & 1 \end{bmatrix}\ o_m^2 \begin{bmatrix} 0.3 & 0.8 & 1 & 0.5 & 0 & 0.7 \\ 0.5 & 0.7 & 0 & 1 & 0.6 & 1 \\ 1 & 0.5 & 0.2 & 0.7 & 0.2 & 0.1 \end{bmatrix} \right\}$$

$(o_m^2 = \text{','})$

$$\left\{ \begin{bmatrix} 0 & 0 & 0 & 0 & 0 & 1 \end{bmatrix} o_m^3 \begin{bmatrix} 1 & -1 & 0 & 1 \\ 0 & 0 & 1 & 0 \\ 0 & 1 & 0 & 1 \\ 1 & 0 & 1 & -1 \\ 0 & 1 & 0 & 0 \\ 0 & 1 & 1 & 0 \end{bmatrix} \right\} (\, o_m^3 = \text{'o'}$$

operations)

$\cup$ max min

$$\left\{ \begin{bmatrix} 0 & 1 & 0 & 0 & 1 \end{bmatrix} o_m^4 \begin{bmatrix} 0.9 & 0.2 & 1 & 0 & 0.5 \\ 0.3 & 0 & 0.5 & 1 & 0.7 \\ 1 & 0.2 & 1 & 0 & 0.3 \\ 0 & 1 & 0.3 & 0.2 & 1 \\ 0.9 & 0.7 & 0.6 & 0.7 & 0 \end{bmatrix} \right\} \; (\, o_m^4 = \text{','})$$

$$\left\{ \begin{bmatrix} 0 & 0 & 0 & 1 & 0 & 0 \end{bmatrix} o_m^5 \begin{bmatrix} 1 & 0 & 1 & 1 & 1 & 0 & 0 & 1 & -1 \\ 0 & 1 & 0 & 0 & 0 & 1 & 0 & 0 & 0 \\ 1 & -1 & 0 & 0 & 0 & 0 & 0 & -1 & 1 \\ 0 & 0 & 1 & 0 & 1 & 0 & 1 & 0 & 0 \\ -1 & 0 & 0 & -1 & 1 & 0 & 0 & 0 & 0 \\ 0 & 0 & -1 & 0 & 0 & 1 & 1 & 0 & 0 \end{bmatrix} \right\}$$

$(\, o_m^5 = \text{'o'})$

$$\begin{aligned} Y &= [0\ 1\ {-1}\ 0\ 0] \cup [0.3\ 0.7\ 0\ 0.5\ 0\ 0.7] \cup [0\ 1\ 1\ 0] \cup \\ &\quad [0.9\ 0.7\ 0.6\ 1\ 0.7] \cup [0\ 0\ 1\ 0\ 1\ 0\ 1\ 0\ 0] \\ &= Y'_1 \cup Y'_2 \cup Y'_3 \cup Y'_4 \cup Y'_5 \\ &= Y'. \end{aligned}$$

We update and threshold the resultant wherever applicable since $Y' = Y'_1 \cup Y'_2 \cup Y'_3 \cup Y'_4 \cup Y'_5$ is not a special fuzzy mixed row vector and find

$$Y = Y_1 \cup Y_2 \cup Y_3 \cup Y_4 \cup Y_5.$$

= [1 1 0 0 0] ∪ [0.3 0.7 0 0.5 0 0.7] ∪ [0 1 1 0] ∪ [0.9 0.7 0.6 1 0.7] ∪ [0 0 1 0 1 0 1 0 0].

Now we find

$$Y \, o_m \, S^{ST} = Y_1 \, o_m^1 \, S_1 \cup Y_2 \, o_m^2 \, S_2^T \cup Y_3 \, o_m^3 \, S_3^T \cup Y_4 \, o_m^4 \, S_4 \cup Y_5 \, o_m^5 \, S_5^T$$

$$= \begin{bmatrix} 1 & 1 & 0 & 0 & 0 \end{bmatrix} \begin{bmatrix} 0 & 1 & -1 & 0 & 0 \\ 1 & 0 & 0 & 0 & 1 \\ -1 & 0 & 0 & 1 & 0 \\ 0 & 0 & 1 & 0 & 0 \\ 0 & -1 & 0 & 1 & 0 \end{bmatrix}$$

$$\cup \text{ min max } \left\{ \begin{bmatrix} 0.3 & 0.7 & 0 & 0.5 & 0 & 0.7 \end{bmatrix}, \begin{bmatrix} 0.3 & 0.5 & 1 \\ 0.8 & 0.7 & 0.5 \\ 1 & 0 & 0.2 \\ 0.5 & 1 & 0.7 \\ 0 & 0.6 & 0.2 \\ 0.7 & 1 & 0.1 \end{bmatrix} \right\}$$

$$\cup \begin{bmatrix} 0 & 1 & 1 & 0 \end{bmatrix} \, o \begin{bmatrix} 1 & 0 & 0 & 1 & 0 & 0 \\ -1 & 0 & 1 & 0 & 1 & 1 \\ 0 & 1 & 0 & 1 & 0 & 1 \\ 1 & 0 & 1 & -1 & 0 & 0 \end{bmatrix}$$

$$\cup \text{ max min } \left\{ \begin{bmatrix} 0.9 & 0.7 & 0.6 & 1 & 0.7 \end{bmatrix}, \begin{bmatrix} 0.9 & 0.2 & 1 & 0 & 0.5 \\ 0.3 & 0 & 0.5 & 1 & 0.7 \\ 1 & 0.2 & 1 & 0 & 0.3 \\ 0 & 1 & 0.3 & 0.2 & 1 \\ 0.9 & 0.7 & 0.6 & 0.7 & 0 \end{bmatrix} \right\}$$

$$[0\ \ 0\ \ 1\ \ 0\ \ 1\ \ 0\ \ 1\ \ 0\ \ 0] \text{ o } \begin{bmatrix} 1 & 0 & 1 & 0 & -1 & 0 \\ 0 & 1 & -1 & 0 & 0 & 0 \\ 1 & 0 & 0 & 1 & 0 & -1 \\ 1 & 0 & 0 & 0 & -1 & 0 \\ 1 & 0 & 0 & 1 & 1 & 0 \\ 0 & 1 & 0 & 0 & 0 & 1 \\ 0 & 0 & 0 & 1 & 0 & 1 \\ 1 & 0 & -1 & 0 & 0 & 0 \\ -1 & 0 & 1 & 0 & 0 & 0 \end{bmatrix}$$

$$\begin{aligned} &= [1\ 1\ {-1}\ 0\ 1] \cup [0\ 0\ 0.2] \cup [-1\ 1\ 1\ 1\ 1\ 2] \cup [0.9\ 1\ 0.9\ 0.7\ 1] \cup [2\ 0\ 0\ 3\ 1\ 0] \\ &= Z'_1 \cup Z'_2 \cup Z'_3 \cup Z'_4 \cup Z'_5 \\ &= Z'; \end{aligned}$$

after updating and thresholding Z' we get Z

$$\begin{aligned} &= Z_1 \cup Z_2 \cup Z_3 \cup Z_4 \cup Z_5 \\ &= [1\ 1\ 0\ 0\ 1] \cup [0\ 0\ 0.2] \cup [0\ 1\ 1\ 1\ 1\ 1] \cup [0.9\ 1\ 0.9\ 0.7\ 1] \cup [1\ 0\ 0\ 1\ 1\ 0]. \end{aligned}$$

Using Z and S we find Z $o_m$ S using the special mixed operator and so on. This sort of special mixed operator will be used when we use special mixed fuzzy models which will be described in Chapter two.

## 1.3 Special Neutrosophic matrices and fuzzy neutrosophic matrices and some essential operators using them

In this section we just introduce the special neutrosophic matrices and special fuzzy neutrosophic matrices and give some of the essential operators using them. We just give some information about neutrosophy. For more about these concepts please refer [187-200].

We denote the indeterminate by I and $I^2 = I$; further

$$\underbrace{I + I + \ldots + I}_{n-\text{times}} = nI, \; n \geq 2.$$

We call [0 1] ∪ [0 I] to be the fuzzy neutrosophic interval. If we take the set generated by ⟨Z ∪ I⟩ we call this as the neutrosophic ring of integers. Likewise ⟨Q ∪ I⟩ denotes the neutrosophic ring of rationals. ⟨R ∪ I⟩ the set generated by R and I is called as the neutrosophic ring of reals.

For more about these concepts please refer [187-190].

Thus a pure neutrosophic number is nI; nI ∈ R and a mixed neutrosophic integer is m + nI; n, m ∈ R. Now a matrix is called a neutrosophic matrix if it takes its entries from ⟨Z ∪ I⟩ or ⟨Q ∪ I⟩ or ⟨R ∪ I⟩.

We now illustrate different types of neutrosophic matrices.

***Example 1.3.1:*** Consider the neutrosophic matrix

$$M = \begin{bmatrix} I & 0 & 2I+1 & 0.8 \\ 7 & 8I & 3I-5 & -2I \\ 9-3I & 6I+3 & 3+0.9I & 12I \end{bmatrix}.$$

We call M to be 3 × 4 rectangular neutrosophic matrix with entries from the neutrosophic ring of reals.

***Example 1.3.2:*** Let N be a neutrosophic matrix where

$$N = \begin{bmatrix} 3I & 2-5I & 8 \\ 0 & 1 & 9+I \\ 2+7I & 8I+1 & 0 \end{bmatrix}.$$

N is a neutrosophic square matrix.

***Example 1.3.3:*** Consider

$$X = \begin{bmatrix} 0.9I \\ 8+2I \\ 5I-7 \\ I \\ 0 \\ 4I \\ -8I \\ 1 \end{bmatrix}.$$

X is a neutrosophic matrix which is known as the neutrosophic column vector/matrix.

***Example 1.3.4:*** Consider the neutrosophic matrix

$$Y = \begin{bmatrix} 0.9 & I+8 & 7I & 8I-1 & 0.2 & I+0.3 & 0.1 \end{bmatrix}$$

is called the neutrosophic row vector/matrix.

***Example 1.3.5:*** Let

$$T = \begin{bmatrix} 0.2I & 0.3 & 0.4I & 0.2I-0.7 & 1 & 0 & 0.9I \end{bmatrix}$$

be a neutrosophic matrix, T will be called the fuzzy neutrosophic row vector/matrix.

***Example 1.3.6:*** Let

$$V = \begin{bmatrix} 0.9 \\ 0.2+I \\ 0.3I-1 \\ 1 \\ 0.8I \\ I-0.7 \\ I \\ 0 \end{bmatrix}$$

neutrosophic matrix. V is known as the fuzzy neutrosophic column vector/matrix.

***Example 1.3.7:*** Let us consider the neutrosophic matrix

$$S = \begin{bmatrix} 4 & 0 & 0 & 0 & 0 \\ 0 & 2I & 0 & 0 & 0 \\ 0 & 0 & 1+9I & 0 & 0 \\ 0 & 0 & 0 & I+1 & 0 \\ 0 & 0 & 0 & 0 & 0.1 \end{bmatrix}$$

S is called the neutrosophic diagonal matrix.

Now we will recall the definition of neutrosophic matrix addition and multiplication.

The neutrosophic zero matrix is the usual zero matrix i.e., a matrix in which all elements in it are zero and denoted by

$$(0) = \begin{bmatrix} 0 & 0 & \dots & 0 \\ 0 & 0 & \dots & 0 \\ \vdots & \vdots & & \vdots \\ 0 & 0 & \dots & 0 \end{bmatrix}.$$

***Example 1.3.8:*** Let S and T be two 3 × 5 rectangular matrices where

$$S = \begin{bmatrix} 0 & 4 & 7I-1 & 2 & -I \\ 2I & 0 & 4I+1 & 0 & 9I \\ 7+3I & 0.9-2I & 3 & 9I & 2-I \end{bmatrix}$$

and

$$T = \begin{bmatrix} 5+I & 7+I & 2-I & 9+I & 0 \\ 7I & 2-I & 3+8I & 1 & 7 \\ 3-0.8I & 0.9 & 9 & 0 & 4+I \end{bmatrix}$$

The neutrosophic matrix addition of S and T

$$S + T = \begin{bmatrix} 5+I & 11+I & 6I+1 & 11+I & -I \\ 9I & 2-I & 4+12I & 1 & 7+9I \\ 10+2.2I & 1.8-2I & 12 & 9I & 6 \end{bmatrix}.$$

We see S + T is also a neutrosophic matrix. We will state in general sum of two neutrosophic matrices need not yield back a neutrosophic matrix. This is shown by the following example.

***Example 1.3.9:*** Let

$$P = \begin{bmatrix} 3I+1 & 2I \\ 17-I & 4 \end{bmatrix}$$

and

$$Q = \begin{bmatrix} 7-3I & 8-2I \\ 3+I & 2 \end{bmatrix}$$

be any two 2 × 2 square matrices. Now

$$P + Q = \begin{bmatrix} 8 & 8 \\ 20 & 6 \end{bmatrix}.$$

Clearly P + Q is not a neutrosophic matrix. Likewise in general the product of two neutrosophic matrices under matrix product need not yield a neutrosophic matrix.

***Example 1.3.10:*** Let

$$A = \begin{bmatrix} 7+I & I \\ I & -6I \end{bmatrix} \text{ and } B = \begin{bmatrix} 7-I & 0 \\ I & 0 \end{bmatrix}$$

be any two neutrosophic matrices AB the matrix product of neutrosophic matrices

$$AB = \begin{bmatrix} 7+I & I \\ I & -6I \end{bmatrix} \begin{bmatrix} 7-I & 0 \\ I & 0 \end{bmatrix}$$

$$= \begin{bmatrix} 49 & 0 \\ 0 & 0 \end{bmatrix}$$

which is not a neutrosophic matrix.

However we illustrate by the following example that even if the neutrosophic product AB is defined the product BA need not be defined. If A is a m × t matrix and B is a t × s matrix then AB is defined but BA is not defined (m ≠ s).

***Example 1.3.11:*** Let

$$A = \begin{bmatrix} 0 & I & 2-I \\ 4-I & 0 & 7 \\ 8I & -1 & 0 \end{bmatrix}$$

and

$$B = \begin{bmatrix} 7I-1 & 2+I & 3-I & 5-I & 0 \\ 0 & 7I & 2 & 0 & 3 \\ 8+I & 3I & -I & 1 & 0 \end{bmatrix}$$

Now

$$AB = \begin{bmatrix} 16-7I & 10I & I & 2-I & 3I \\ 29I+52 & 8+22I & 12-13I & 27-8I & 0 \\ 48I & 17I & 16I-2 & 32I & -3 \end{bmatrix}$$

which is a neutrosophic matrix. But BA is not defined as B is a 3 × 5 matrix and A is a 3 × 3 matrix, so BA is not defined.

Now we proceed on to define a new class of special neutrosophic matrices.

**DEFINITION 1.3.1:** *Let M = $M_1 \cup M_2 \cup \ldots \cup M_n$ be a collection of neutrosophic matrices where each $M_i$ is a t × t neutrosophic matrix; i = 1, 2, …, n. We call M to be a special neutrosophic square matrix.*

We illustrate this by the following example.

***Example 1.3.12:*** Let $M = M_1 \cup M_2 \cup \ldots \cup M_5$

$$= \begin{bmatrix} 3I & I+7 & 8-I & 0 \\ 1 & 7I & 5I-1 & 8I \\ 2I+4 & 11I-1 & 1 & 0 \\ 0 & 1 & 9-5I & 2-I \end{bmatrix} \cup$$

$$\begin{bmatrix} 0 & 7-I & 5+I & 2I \\ 7I & 0 & 8 & 9I-1 \\ 2I-1 & 8I & 1 & 21I \\ 9+3I & 3 & 2-I & 5I \end{bmatrix} \cup$$

$$\begin{bmatrix} 3I-1 & 0 & 1 & 3+2I \\ 2+4I & 3I & 4I-1 & 0 \\ 3-I & 4 & 3I & 1 \\ 5I & 9I & 6-3I & 0 \end{bmatrix} \cup$$

$$\begin{bmatrix} 0 & 3I & 7I-1 & 2-3I \\ 4I & 0 & 6+5I & 9I \\ 2-I & 4-2I & 1 & 17I \\ 6I+1 & 12I & 5-I & 3-I \end{bmatrix} \cup$$

$$\begin{bmatrix} 1 & 0 & I & -2+I \\ I & 7-I & 8I & 5I-1 \\ 6I+1 & 2I & 7-5I & 8I \\ 9I-1 & 6I-1 & 0 & 1 \end{bmatrix}.$$

be a special neutrosophic square matrix.

*Note:* Even if $M_i = M_j$ ($i \neq j$) still M is a special neutrosophic square matrix.

***Example 1.3.13:*** Let G = $G_1 \cup G_2 \cup G_3$

$$= \begin{bmatrix} 2I & 6-4I \\ 5I & 7+3I \end{bmatrix} \cup \begin{bmatrix} 0 & I \\ 5I-1 & 2I+1 \end{bmatrix} \cup \begin{bmatrix} 0 & I \\ 5I-1 & 2I+1 \end{bmatrix}$$

G is a special neutrosophic square matrix.

**DEFINITION 1.3.2**: *Let $T = T_1 \cup T_2 \cup \ldots \cup T_s$ if the $T_i$'s are $n_i \times n_i$ square neutrosophic matrices then we call T to be a special neutrosophic mixed square matrix (if $s \geq 2$) $n_i \neq n_j$ ($i \neq j$ at least for some i and j) $1 \leq i, j \leq s$.*

***Example 1.3.14:*** Let T = $T_1 \cup T_2 \cup T_3 \cup T_4$

$$= \begin{bmatrix} 0.3I & 5I-1 & 3I+2 & 4I \\ 0 & 3I+1 & 3-8I & 9 \\ I & 2I-1 & 4I-1 & 21 \\ 2I+7 & 9-3I & 9I-1 & 9-7I \end{bmatrix} \cup$$

$$\begin{bmatrix} 1 & 0.3+I & 7-2I \\ 0 & 4I-1 & 12I-1 \\ 6I+1 & 0 & 34 \end{bmatrix} \cup \begin{bmatrix} 0.I & 2 \\ 7I & 4-I \end{bmatrix} \cup$$

$$\begin{bmatrix} 0.1 & I-7 & 8+I & 5I-1 & 8I \\ 9 & 12+I & 0 & 17I-1 & 6I \\ 8+2I & 0 & I & 2-5I & 7 \\ I & 9I-1 & 9+3I & 8I+1 & 1-5I \\ 5I-1 & 12-I & 1 & 3 & 9+I \end{bmatrix}.$$

T is a special neutrosophic mixed square matrix.

**DEFINITION 1.3.3:** *Let $W = W_1 \cup W_2 \cup \ldots \cup W_n$ ($n \geq 2$) be the collection of $p \times q$ neutrosophic ($p \neq q$) rectangular matrices. We define W to be a special neutrosophic rectangular matrix.*

We illustrate this by the following example.

***Example 1.3.15:*** Let W = $W_1 \cup W_2 \cup W_2 \cup W_4$ =

$$\begin{bmatrix} 0.3I & 7I-1 & 8I+1 \\ 8 & 5I+1 & 9I \\ 21I-4 & 9+6I & 4+6I \\ 9-5I & 7 & 0.5 \\ 2+6 & 5I & 21 \end{bmatrix} \cup$$

$$\begin{bmatrix} 9I & 2I-1 & 5 \\ 0 & 7+2I & 6+I \\ 21 & 6I-1 & 12I-1 \\ 7I+3 & 15+I & 8I-1 \\ 9I-5I & 2 & 3I \end{bmatrix} \cup$$

$$\begin{bmatrix} 6 & 8I-1 & 9I \\ 5I-1 & 0 & 8 \\ 7I & 9I+8 & 2I+5 \\ 3-I & 21 & 0 \\ 1 & 15I & 9-6I \end{bmatrix} \cup$$

$$\begin{bmatrix} 21-I & 2+5I & 8 \\ 9 & 0 & 20I \\ 4+5I & 7I & 5+3I \\ 12+I & 5 & 0 \\ 9I & 7-5I & 12+I \end{bmatrix},$$

W is a special neutrosophic rectangular matrix.

Next we define the notion of special neutrosophic mixed rectangular matrix.

**DEFINITION 1.3.4:** *Let V = $V_1 \cup V_2 \cup \ldots \cup V_n$ (n ≥ 2) be a collection of neutrosophic rectangular matrices. Each $V_i$ is a $t_i \times s_i$ ($t_i \neq s_i$) neutrosophic rectangular matrix; i = 1, 2, …, n. We define V to be a special neutrosophic mixed rectangular matrix.*

We illustrate this by the following example.

***Example 1.3.16:*** Let V = $V_1 \cup V_2 \cup V_3 \cup V_4 \cup V_5$ =

$$\begin{bmatrix} 0.3I & I+7 \\ 2I- & 8I \\ 5 & 9I-4 \\ 0 & 21I+1 \\ 12I & 3 \\ 17 & 6-I \end{bmatrix} \cup$$

$$\begin{bmatrix} 0.5I & 7-5I & 15I+4 & 8I+7 & 21 \\ 7I & 0 & 21-7I & 6I+1 & 8I+4 \\ 1 & 15I+6 & 9I & 12 & 0 \end{bmatrix} \cup$$

$$\begin{bmatrix} 0.8 & 9I & 5I+1 \\ 3I+1 & 7I-1 & 16I \\ 7 & 8I & 1 \\ 6I-7 & 2+I & 0 \\ 10I & 8+5I & 9I-1 \\ 7-8I & 5 & 3+5I \end{bmatrix} \cup$$

$$\begin{bmatrix} 2I & 9 & 5-I & 5I-1 & 9-2I \\ 0 & 5I+1 & 25I & 4-I & 2I \\ 7 & 21-7I & 4 & 9+3I & 0 \\ 15I- & 0 & 1+I & 2II & 19 \end{bmatrix} \cup$$

$$\begin{bmatrix} 8I & -9 & 25-I & 7+8I & 0 \\ 0 & I+1 & 4+2I & 8 & 5I \\ 17-I & 1 & 0 & 9I & 12 \\ 21 & 12+5I & 8I+1 & 0 & 7I-5 \\ 4+5I & 25I & 16I-1 & 6-7I & 8+I \\ 9I-1 & 8I-1 & 25I+1 & 6I & 75 \end{bmatrix}.$$

V is a special neutrosophic mixed rectangular matrix.

**DEFINITION 1.3.5:** *Let $X = X_1 \cup X_2 \cup \ldots \cup X_n$ ($n \geq 2$) be a collection of neutrosophic row matrices. Here each $X_i$ is a $1 \times t$ row neutrosophic matrix. We call X to be a special neutrosophic row vector / matrix. If some $X_j$ is a $1 \times t_j$ neutrosophic row matrix and some $X_k$ is a $1 \times t_k$ neutrosophic row vector and $t_j \neq t_k$ for some $k \neq j$. Then we define X to be a special neutrosophic mixed row vector.*

We illustrate them by the following example.

***Example 1.3.17:*** Let

$$\begin{aligned} X &= X_1 \cup X_2 \cup X_3 \cup X_4 \cup X_5 \cup X_6 \\ &= [I\ 0\ 0.8\ 0.8\ I\ 0.9\ I\ 0.2] \cup [0\ I\ I\ 6\ 1\ I\ 8I\ 0.7] \cup [0.5\ 0.6\ I\ 0.7\ 4\ 0.9I\ 0.6\ 9] \cup [9\ 1\ 4\ 1\ I\ 0\ 0\ I] \cup [0.2\ I\ 0.2\ 0.2\ 0.3I\ 0.4I\ 8\ 0] \cup [0\ I\ 0\ I\ 0\ 0\ I\ 8], \end{aligned}$$

we see X is a special neutrosophic row vector for we see each $X_i$ is a $1 \times 8$ neutrosophic row vector / matrix, $1 \leq i \leq 6$.

Next we give an example of a special neutrosophic mixed neutrosophic row vector.

***Example 1.3.18:*** Let

$$\begin{aligned} X &= X_1 \cup X_2 \cup X_3 \cup X_4 \\ &= [0\ I\ 0.7\ 0.8] \cup [I\ 1\ I\ 0\ 0\ 0\ 1\ 0\ 0\ 5] \cup [0.2\ 0.9\ I\ 0.6\ 0.8\ 1\ 0\ 0] \cup [0\ I\ 0\ 0\ 7\ 1\ I\ 0], \end{aligned}$$

X is a special neutrosophic mixed row vector / matrix.

We define the notion of special neutrosophic column vector and special neutrosophic mixed column vector / matrix.

**DEFINITION 1.3.6:** *Let $Y = Y_1 \cup Y_2 \cup \ldots \cup Y_m$ ($m \geq 2$) where $Y_i$ is a $t \times 1$ column neutrosophic vector/matrix for each i = 1, 2, …, m. Then we call Y to be a special neutrosophic column vector/matrix. Suppose if $Y = Y_1 \cup Y_2 \cup \ldots \cup Y_m$ be such that some $Y_j$ is a $t_j \times 1$ neutrosophic column vector/matrix for j = 1, 2, …, m and $t_j \neq t_i$ for at least one $i \neq j$, $1 \leq i, j \leq m$. Then we call Y to be a special neutrosophic mixed column vector / matrix.*

Now we illustrate these concepts by the following examples.

***Example 1.3.19:*** Let $Y = Y_1 \cup Y_2 \cup \ldots \cup Y_6 =$

$$\begin{bmatrix} 0.2I \\ 1 \\ 0.7 \\ I \\ 9.3 \\ 0 \\ 0.6I \end{bmatrix} \cup \begin{bmatrix} I \\ 9I \\ 0.8 \\ I \\ 0.6I \\ 1 \\ I \end{bmatrix} \cup \begin{bmatrix} 6 \\ 9I \\ 0.8I \\ 0.7 \\ I \\ 9 \\ 0 \end{bmatrix}$$

$$\cup \begin{bmatrix} I \\ 14 \\ 4I \\ I \\ 1 \\ 6 \\ I \end{bmatrix} \cup \begin{bmatrix} 8 \\ 0 \\ I \\ 0 \\ 6 \\ 6.8I \\ 9.9 \end{bmatrix} \cup \begin{bmatrix} 1 \\ 6.5 \\ I \\ 0 \\ I \\ 6I \\ I \end{bmatrix}.$$

Y is a special neutrosophic column vector / matrix and each of the $Y_i$ is a 7 × 1 fuzzy column vector, $1 \leq i \leq 6$.

***Example 1.3.20:*** Let $Y = Y_1 \cup Y_2 \cup Y_3 \cup Y_4 \cup Y_5 =$

$$\begin{bmatrix} I \\ 6 \\ 0.7I \\ 0.8 \end{bmatrix} \cup \begin{bmatrix} 0.9I \\ 0.6I \\ 6.3 \end{bmatrix} \cup \begin{bmatrix} 0.4I \\ I \\ 1 \\ I \\ 6.7 \end{bmatrix} \cup \begin{bmatrix} 0.8I \\ 6.9 \\ 4I \\ 9.6I \\ 0.I \\ 12 \\ 8I \\ 0.15 \end{bmatrix} \cup \begin{bmatrix} 6.16I \\ 8.02 \\ 6.7I \\ 6I \\ 9 \\ 0.7I \\ 8.9 \\ 0.8 \\ 7I \end{bmatrix}$$

be the special neutrosophic mixed column vector / matrix. We see each neutrosophic column vector/matrix has different order.

Now we proceed onto define some basic operations called the special $m\underline{a}x\ m\underline{i}n$ operator on special neutrosophic square matrix and illustrate them with examples.

Let $S = S_1 \cup S_2 \cup \ldots \cup S_n$ be a special neutrosophic square matrix where each $S_i$ is a m × m neutrosophic matrix, i = 1, 2, …, n. Let $X = X_1 \cup X_2 \cup \ldots \cup X_n$ be a special neutrosophic row vector where each $X_i$ is a 1 × m neutrosophic row vector. We define the special neutrosophic operation using X and S.

Suppose we define a $m\underline{a}x\ m\underline{i}n$ operator, to find

$$\begin{aligned} & m\underline{a}x\ m\underline{i}n\ (X \text{ o } S) \\ &= m\underline{a}x\ m\underline{i}n\ \{(X_1 \cup X_2 \cup \ldots \cup X_n), (S_1 \cup S_2 \cup \ldots \cup S_n)\} \\ &= m\underline{a}x\ m\underline{i}n\{(X_1, S_1)\} \cup m\underline{a}x\ m\underline{i}n\ \{(X_2, S_2)\} \cup \ldots \cup m\underline{a}x\ m\underline{i}n\ \{(X_n, S_n)\} \\ &= Y'_1 \cup Y'_2 \cup \ldots \cup Y'_n \\ &= Y' \end{aligned}$$

now Y' may be a special neutrosophic row vector. Now we calculate using Y' on S we find

$\underline{max}\ \underline{min}$ (Y', S)

$= \underline{max}\ \underline{min}\ \{(Y_1 \cup Y_2 \cup \ldots \cup Y_n), (S_1 \cup S_2 \cup \ldots \cup S_n)\}$

$= \underline{max}\ \underline{min}\ \{Y_1, S_1\} \cup \underline{max}\ \underline{min}\ (Y_2, S_2) \cup \ldots \cup \underline{max}\ \underline{min}\ \{Y_n, S_n\}$

$= T_1 \cup T_2 \cup \ldots \cup T_n$

and now we find max min {T, S} and so on.

We illustrate this by the following example.

***Example 1.3.21:*** Let $S = S_1 \cup S_2 \cup S_3 \cup S_4 \cup S_5$ be a special neutrosophic square matrix and $X = X_1 \cup X_2 \cup \ldots \cup X_5$ be a special neutrosophic row vector; where

$$S = \begin{bmatrix} 0.7 & I & 3 & 2I \\ 4I & 0 & I & 0 \\ I & 1 & 0 & 0 \\ 7 & 2I & 0 & 5I \end{bmatrix} \cup \begin{bmatrix} I & 2I & 0 & 1 \\ 0 & 7 & I & 0 \\ 1 & 0 & 0 & I \\ 0 & I & 0 & 2 \end{bmatrix} \cup$$

$$\begin{bmatrix} 3 & I & 0 & 0 \\ 0 & 4 & 2I & 0 \\ 5 & 0 & 3I & 0 \\ 0 & 2 & 0 & 4I \end{bmatrix} \cup \begin{bmatrix} 0 & 2 & 0 & 5I \\ I & 0 & 7 & 0 \\ 3I & 0 & 0 & 2 \\ 0 & 4 & I & 0 \end{bmatrix} \cup \begin{bmatrix} 0 & 2I & 1 & 0 \\ 2I & 0 & 0 & 1 \\ I & 0 & 1 & 0 \\ 0 & 4I & 0 & I \end{bmatrix}$$

is a special neutrosophic square matrix with each $S_i$ a $4 \times 4$ neutrosophic square matrix for i = 1, 2, 3, 4, 5 and

$X = X_1 \cup X_2 \cup \ldots \cup X_5$

$= [I\ 0\ 0\ 0] \cup [0\ I\ 0\ 0] \cup [0\ 0\ I\ 0] \cup [0\ 0\ 0\ I] \cup [0\ 1\ 0\ I\ 0]$

be the special neutrosophic row vector each $X_i$ is a $1 \times 4$ neutrosophic row vector. To find special max min of X, S, i.e. $\underline{max}\ \underline{min}$ {X, S} $= \underline{\max} \underline{\min} \{(X_1 \cup X_2 \cup X_3 \cup X_4 \cup X_5), (S_1 \cup S_2 \cup \ldots \cup S_5)\}$

$= \underline{\max} \underline{\min} \{(X_1, S_1)\} \cup \underline{\max} \underline{\min} \{(X_2, S_2)\} \cup \underline{\max} \underline{\min} \{(X_3, S_3)\} \cup \ldots \cup \underline{\max} \underline{\min} \{(X_5, S_5)\}$

$$= \underline{\max} \underline{\min} \left\{ \begin{bmatrix} I & 0 & 0 & 0 \end{bmatrix}, \begin{bmatrix} 0.7 & I & 3 & 2I \\ 4I & 0 & I & 0 \\ I & 1 & 0 & 0 \\ 0.7 & 2I & 0 & 5I \end{bmatrix} \right\} \cup$$

$$\underline{\max} \underline{\min} \left\{ \begin{bmatrix} 0 & I & 0 & 0 \end{bmatrix}, \begin{bmatrix} I & 2I & 0 & 1 \\ 0 & 7 & I & 0 \\ 1 & 0 & 0 & I \\ 0 & I & 0 & 2 \end{bmatrix} \right\} \cup$$

$$\underline{\max} \underline{\min} \left\{ \begin{bmatrix} 0 & 0 & I & 0 \end{bmatrix}, \begin{bmatrix} 3 & I & 0 & 0 \\ 0 & 4 & 2I & 0 \\ 5 & 0 & 3I & 0 \\ 0 & 2 & 0 & 4I \end{bmatrix} \right\} \cup$$

$$\underline{\max} \underline{\min} \left\{ \begin{bmatrix} 0 & 0 & 0 & I \end{bmatrix}, \begin{bmatrix} 0 & 2 & 0 & 5I \\ I & 0 & 7 & 0 \\ 3I & 0 & 0 & 2 \\ 0 & 4 & I & 0 \end{bmatrix} \right\} \cup$$

$$\underline{\max} \underline{\min} \left\{ \begin{bmatrix} 0 & 1 & 0 & I \end{bmatrix}, \begin{bmatrix} 0 & 2I & 1 & 0 \\ 2I & 0 & 0 & 1 \\ I & 0 & 1 & 0 \\ 0 & 4I & 0 & I \end{bmatrix} \right\}$$

How do we define m$\underline{i}$n(m, n), for $m < n$ and $n < m$ cannot be defined in a natural way as comparison of a neutrosophic number and the usual reals is not possible. So we have defined depending on the wishes of the expert according as the expert wants to give what type of weightage to the indeterminacy. Suppose n is a neutrosophic number like n = 5I and m = 8 then m$\underline{i}$n (5I, 8) = 5I and m$\underline{a}$x (5I, 8) is 8; m$\underline{i}$n(2, 7I) is 2 and m$\underline{a}$x(2, 7I) is 7I and like wise when m and n are of the form xI and yI, x and y real numbers m$\underline{a}$x (xI, yI) = xI or yI according as $x > y$ or $y > x$ and m$\underline{i}$n (xI, yI) = xI or yI according as $x < y$ or $y < x$ and m$\underline{a}$x (nI, n) is nI and m$\underline{i}$n (nI, n) = nI.

It is pertinent and important to mention here that if the expert feels the presence indeterminacy is important he/she defines m$\underline{i}$n (n, nI) = nI and m$\underline{a}$x (n, nI) = nI this does not affect any logic for our definition of min max is different. Clearly one can not have the usual max or min using reals and indeterminacy.

Now using this mode of calculation we find

m$\underline{a}$x m$\underline{i}$n (S, T)

$$= [0.7\ I\ I\ I] \cup [0\ I\ I\ 0] \cup [5\ 0\ I\ 0] \cup [0\ I\ I\ 0] \cup [1\ I\ 0\ I]$$
$$= Y_1 \cup Y_2 \cup Y_3 \cup Y_4 \cup Y_5$$
$$= Y,$$

Y is again a special neutrosophic row vector/matrix. Now we can find max min (Y, S) and so on.

Now we proceed on to define special m$\underline{i}$n m$\underline{a}$x operator on special neutrosophic mixed square matrix.

Now suppose

$$T = T_1 \cup T_2 \cup T_3 \cup T_4 \cup T_5 \cup \ldots \cup T_p$$

($p \geq 2$) be a special neutrosophic mixed square matrix, where $T_i$ is a $n_i \times n_i$ square neutrosophic matrix i = 1, 2, …, p with $n_i \neq n_j$ $i \neq j$ for atleast some i and j. $1 \leq i, j \leq p$. Suppose $X = X_1 \cup X_2 \cup \ldots \cup X_p$ be a special neutrosophic mixed row vector where each $X_i$ is a $1 \times p_i$ neutrosophic row vector i = 1, 2, …, p. To find

$\underline{\max}$ $\underline{\min}$ (X, T)

$$= \underline{\max}\ \underline{\min}\ \{(X_1 \cup X_2 \cup \ldots \cup X_p), (T_1 \cup T_2 \cup \ldots \cup T_p)\}$$
$$= \underline{\max}\ \underline{\min}\ (X_1, T_1) \cup \underline{\max}\ \underline{\min}\ (X_2, T_2) \cup \ldots \cup \underline{\max}\ \underline{\min}\ (X_p, T_p)$$
$$= Y_1 \cup Y_2 \cup \ldots \cup Y_p$$
$$= Y$$

is once again a special neutrosophic mixed row vector.

We can find $\underline{\max}$ $\underline{\min}$ (Y, S) and so on.

We illustrate this by the following example.

***Example 1.3.22:*** Let $V = V_1 \cup V_2 \cup \ldots \cup V_5$

$$= \begin{bmatrix} I & 0 & 2 & 3I & 1 \\ 0 & 3I & 0 & 0 & 2 \\ 1 & 0 & 1 & 5 & 0 \\ 0 & 0 & 2I & 0 & 1 \\ 7 & 3 & 0 & 0 & I \end{bmatrix} \cup \begin{bmatrix} 7 & 0 & 0 & I & 0 \\ 0 & 4I & 0 & 0 & 1 \\ 3I & 0 & 1 & 0 & 0 \\ 0 & 0 & 0 & 4 & 3I \\ 0 & 0 & 5I & 0 & 0 \end{bmatrix} \cup$$

$$\begin{bmatrix} 3 & 0 & 4I & 0 & 0 \\ 0 & 0 & 0 & 8 & I \\ 3I & 0 & 0 & 1 & 0 \\ 0 & 2 & 0 & 0 & 0 \\ 0 & 5I & 0 & 1 & 0 \end{bmatrix} \cup \begin{bmatrix} I & 0 & 0 & 0 & 1 \\ 2 & 0 & 0 & 0 & 5I \\ 0 & I & 3 & 0 & 0 \\ 0 & 0 & 2I & 0 & 0 \\ 0 & 5 & 0 & 7I & 0 \end{bmatrix} \cup$$

$$\begin{bmatrix} 0 & 0 & I & 7 & 0 \\ I & 0 & 0 & 0 & 8 \\ 0 & 8 & 0 & I & 0 \\ 2I & 0 & 4 & 0 & I \\ 0 & 0 & 0 & 0 & 6 \end{bmatrix}$$

be a special neutrosophic square matrix and

$$X = X_1 \cup X_2 \cup X_3 \cup X_4 \cup X_5$$
$$= [1\ 0\ 0\ 0\ I] \cup [0\ I\ 0\ 0\ 0] \cup [0\ 1\ I\ 0\ 0] \cup [0\ 0\ 0\ I\ 1] \cup [0\ 0\ I\ 0\ 1]$$

be a special neutrosophic row vector. Now using the special operator viz. min max (X, V) we get the following resultant vector

$\min_{_}\max_{_}(X, V)$

$$= \text{min max } \{(X_1 \cup X_2 \cup \ldots \cup X_5), (V_1 \cup V_2 \cup \ldots \cup V_5)\}$$
$$= \text{min max } (X_1, V_1) \cup \ldots \text{ min max } (X_5, V_5)$$

$$= \text{min max } \left\{ \begin{bmatrix} 1 & 0 & 0 & 0 & I \end{bmatrix}, \begin{bmatrix} I & 0 & 2 & 3I & 1 \\ 0 & 3I & 0 & 0 & 2 \\ 1 & 0 & 1 & 5 & 0 \\ 0 & 0 & 2I & 0 & 1 \\ 7 & 3 & 0 & 0 & I \end{bmatrix} \right\}$$

$$\cup \text{ min max } \left\{ \begin{bmatrix} 0 & I & 0 & 0 & 0 \end{bmatrix}, \begin{bmatrix} 7 & 0 & 0 & I & 0 \\ 0 & 4I & 0 & 0 & I \\ 3I & 0 & 1 & 0 & 0 \\ 0 & 0 & 0 & 4 & 3I \\ 0 & 0 & 5I & 0 & 0 \end{bmatrix} \right\}$$

$$\cup \text{ min max } \left\{ \begin{bmatrix} 0 & 1 & I & 0 & 0 \end{bmatrix}, \begin{bmatrix} 3 & 0 & 4I & 0 & 0 \\ 0 & 0 & 0 & 8 & I \\ 3I & 0 & 0 & 1 & 0 \\ 0 & 2 & 0 & 0 & 0 \\ 0 & 5I & 0 & 1 & 0 \end{bmatrix} \right\}$$

$$\cup \ \underline{\min}\ \underline{\max} \left\{ \begin{bmatrix} 0 & 0 & 0 & I & 1 \end{bmatrix}, \begin{bmatrix} I & 0 & 0 & 0 & 1 \\ 2 & 0 & 0 & 0 & 5I \\ 0 & I & 3 & 0 & 0 \\ 0 & 0 & 2I & 0 & 0 \\ 0 & 5 & 0 & 7I & 0 \end{bmatrix} \right\}$$

$$\cup \ \underline{\min}\ \underline{\max} \left\{ \begin{bmatrix} 0 & 0 & I & 0 & 1 \end{bmatrix}, \begin{bmatrix} 0 & 0 & I & 7 & 0 \\ I & 0 & 0 & 0 & 8 \\ 0 & 8 & 0 & I & 0 \\ 2I & 0 & 4 & 0 & I \\ 0 & 0 & 0 & 0 & 6 \end{bmatrix} \right\}$$

$$= [0\ 0\ 0\ 0\ 0] \cup [0\ 0\ 0\ 0\ 0] \cup [0\ 0\ 0\ 0\ 0] \cup [0\ 0\ 0\ 0\ 0] \cup [0\ 0\ 0\ 0\ 0].$$

Thus the resultant of X on V is zero.

Now we proceed on to describe by an example the special min max and max min operations using on special neutrosophic mixed square matrices.

***Example 1.3.23:*** Let $X = X_1 \cup X_2 \cup X_3 \cup X_4 \cup X_5$

$$= \begin{bmatrix} 0 & 3 & I & 0 & 1 \\ 1 & 0 & 0 & I & 0 \\ 0 & I & 0 & 0 & 1 \\ I & 0 & 0 & 0 & 0 \\ 0 & 0 & 1 & I & I \end{bmatrix} \cup$$

$$\begin{bmatrix} 0 & 4 & I & 0 \\ 1 & 0 & 0 & I \\ I & 0 & 4 & 0 \\ 0 & I & 0 & 0 \end{bmatrix} \cup$$

$$\begin{bmatrix} 0 & 2 & 0 & I & 0 & 1 \\ 1 & 0 & 0 & 0 & 2I & 0 \\ I & 0 & 3 & 0 & 0 & 0 \\ 0 & 0 & 0 & 0 & 0 & 0 \\ 0 & 3I & 0 & 4 & 0 & 3I \\ 0 & 0 & 2I & 0 & 7 & 0 \end{bmatrix} \cup$$

$$\begin{bmatrix} 0 & 7 & I \\ I & 0 & 1 \\ 0 & I & 0 \end{bmatrix} \cup \begin{bmatrix} 0 & 1 & I & 0 \\ 2I & 0 & 0 & 4 \\ 0 & 3I & 0 & 1 \\ 5 & 0 & 6 & 0 \end{bmatrix}$$

be the special neutrosophic mixed square matrix. Suppose

$$\begin{aligned} T &= T_1 \cup T_2 \cup \ldots \cup T_5 \\ &= [I\ 0\ 0\ 0\ 1] \cup [0\ 1\ 0\ 0] \cup [I\ 0\ 1\ 0\ 0\ 1] \cup [0\ 1\ 0] \cup [I\ 0\ 1\ 0] \end{aligned}$$

be the mixed special neutrosophic mixed row matrix. We find using the $\min \max$ operator on T and X and obtain the resultant.

$\min \max (T, X)$

$$\begin{aligned} &= \min \max \{(T_1 \cup T_2 \ldots \cup T_5), (X_1 \cup X_2 \cup \ldots \cup X_5)\} \\ &= \min \max \{T_1, X_1\} \cup \min \max \{(T_2, X_2)\} \cup \min \max \{T_3, X_3\} \cup \min \max\{T_4, X_4\} \cup \min \max\{T_5, X_5\} \end{aligned}$$

$$= \min \max \left\{ \begin{bmatrix} I & 0 & 0 & 0 & 1 \end{bmatrix}, \begin{bmatrix} 0 & 3 & I & 0 & 1 \\ 1 & 0 & 0 & I & 0 \\ 0 & I & 0 & 0 & 1 \\ I & 0 & 0 & 0 & 0 \\ 0 & 0 & 1 & I & I \end{bmatrix} \right\} \cup$$

$$\underline{\min}\ \underline{\max}\left\{\begin{bmatrix}0 & 1 & 0 & 0\end{bmatrix},\ \begin{bmatrix}0 & 4 & I & 0\\ 1 & 0 & 0 & I\\ I & 0 & 4 & 0\\ 0 & I & 0 & 0\end{bmatrix}\right\} \cup$$

$$\underline{\min}\ \underline{\max}\left\{\begin{bmatrix}I & 0 & 1 & 0 & 0 & 1\end{bmatrix},\ \begin{bmatrix}0 & 2 & 0 & I & 0 & 1\\ 1 & 0 & 0 & 0 & 2I & 0\\ I & 0 & 3 & 0 & 0 & 0\\ 0 & 0 & 0 & 0 & 0 & 0\\ 0 & 3I & 0 & 4 & 0 & 3I\\ 0 & 0 & 2I & 0 & 7 & 0\end{bmatrix}\right\} \cup$$

$$\underline{\min}\ \underline{\max}\left\{\begin{bmatrix}0 & 1 & 0\end{bmatrix},\ \begin{bmatrix}0 & 7 & I\\ I & 0 & 1\\ 0 & I & 0\end{bmatrix}\right\} \cup$$

$$\underline{\min}\ \underline{\max}\left\{\begin{bmatrix}I & 0 & 1 & 0\end{bmatrix},\ \begin{bmatrix}0 & 1 & I & 0\\ 2I & 0 & 0 & 4\\ 0 & 3I & 0 & 1\\ 5 & 0 & 6 & 0\end{bmatrix}\right\}$$

$= [0\ 0\ 0\ 0\ 0] \cup [0\ 0\ 0\ 0] \cup [0\ 0\ 0\ 0\ 0\ 0] \cup [0\ I\ 0] \cup [I\ 0\ 0\ 0]$
$= S_1 \cup S_2 \cup S_3 \cup S_4 \cup S_5$
$= S$

is again a special neutrosophic mixed row vector.

Now using S we can find $\underline{\min}\ \underline{\max}$ {S, X} and so on.

Next we proceed on to find using the special min max function the value of a special neutrosophic row vector and the special neutrosophic rectangular matrix.

Let $W = W_1 \cup W_2 \cup \ldots \cup W_n$ ($n \geq 2$) be a special neutrosophic rectangular matrix where each $W_i$ is a $t \times s$ ($t \neq s$) rectangular neutrosophic matrix ($i = 1, 2, \ldots, n$).
Let $X = X_1 \cup X_2 \cup \ldots \cup X_n$ ($n \geq 2$) be the special neutrosophic row vector/matrix where each $X_i$ is a $1 \times t$ neutrosophic vector $i = 1, 2, \ldots, n$.

To find

$$\max \min (\{X, W\})$$

$$\begin{aligned}
&= \max \min \{(X_1 \cup X_2 \cup \ldots \cup X_n), (W_1 \cup W_2 \cup \ldots \cup W_n)\} \\
&= \max \min \{X_1, W_1\} \cup \max \min \{X_2, W_2\} \cup \ldots \cup \max \min \{X_n, W_n\} \\
&= Y_1 \cup Y_2 \cup \ldots \cup Y_n \\
&= Y
\end{aligned}$$

now Y is a special neutrosophic row vector where each $Y_i$ is a $1 \times s$ neutrosophic row vector for $i = 1, 2, \ldots, n$.

Now we see $\max \min \{Y, W\}$ is not defined so we find the value of

$$\max \min \{Y, W^T\}$$

$$\begin{aligned}
&= \max \min \{(Y_1 \cup Y_2 \cup \ldots \cup Y_n), (W_1^T \cup W_2^T \cup \ldots \cup W_n^T)\} \\
&= \max \min (Y_1, W_1^T) \cup \max \min (Y_2, W_2^T) \cup \ldots \cup \max \min (Y_n, W_n^T) \\
&= T_1 \cup T_2 \cup \ldots \cup T_n \\
&= T
\end{aligned}$$

where T is a special neutrosophic row vector with each $T_i$ a $1 \times t$ neutrosophic row vector for i = 1, 2, …, n. We now find out $\max\min$ {T, W} and so on.

Now we illustrate this situation by the following example.

***Example 1.3.24:*** Let $W = W_1 \cup W_2 \cup \ldots \cup W_5$

$$= \begin{bmatrix} 0 & 4 & I & 2 & 4I & 0 & 0 & 7 & 1 \\ 0 & 0 & 0 & 0 & 0 & 8I & 5 & 0 & 1 \\ 1 & 0 & 7I & 0 & 0 & 0 & 3 & 0 & 1 \end{bmatrix} \cup$$

$$\begin{bmatrix} 1 & I & 4 & 5 & 7 & 0 & 0 & 0 & 0 \\ 0 & 1 & 0 & I & 2 & 1 & 0 & 0 & 8 \\ 1 & 0 & 0 & 0 & 1 & 0 & I & 0 & 4I \end{bmatrix} \cup$$

$$\begin{bmatrix} 9I & 0 & 2 & 0 & 0 & 0 & 1 & 4 & I \\ 0 & I & 3 & 0 & 0 & 7 & 0 & 5 & 0 \\ 2 & 0 & 1 & 0 & 8 & 0 & I & 0 & 1 \end{bmatrix} \cup$$

$$\begin{bmatrix} 0 & 3I & 0 & 0 & 0 & 2 & 7 & 1 & 0 \\ 0 & 0 & 5 & 7 & I & 3 & 0 & 1 & 0 \\ 7 & 0 & 0 & 0 & 2 & 0 & 8 & 0 & 8 \end{bmatrix} \cup$$

$$\begin{bmatrix} 0 & 0 & 7 & 0 & 0 & 0 & 2 & 5 & I \\ I & 3 & 0 & 0 & 0 & 4 & 0 & 0 & 2 \\ 0 & 0 & I & 7I & 8 & 0 & 0 & 0 & 1 \end{bmatrix}$$

be a special neutrosophic rectangular matrix where each $W_i$ is a $3 \times 9$ neutrosophic rectangular matrix; i = 1, 2, …, 5. Suppose

$$\begin{aligned} X &= X_1 \cup X_2 \cup \ldots \cup X_5 \\ &= [0\ 0\ I] \cup [I\ 3I\ 0] \cup [5\ 0\ 2I] \cup [0\ 7I\ 0] \cup [2\ I\ 6] \end{aligned}$$

be the special neutrosophic row vector where each $X_j$ is a $1 \times 3$ row vector, j = 1, 2, 3, …, 5. To find the value of

$$\begin{aligned}\underline{\max}\,\underline{\min}\,(X, W) &= \underline{\max}\,\underline{\min}\,\{(X_1 \cup X_2 \cup \ldots \cup X_5), (W_1 \cup W_2 \cup \ldots \cup W_5)\} \\ &= \underline{\max}\,\underline{\min}\,(X_1, W_1) \cup \underline{\max}\,\underline{\min}\,(X_2, W_2) \cup \ldots \cup \underline{\max}\,\underline{\min}\,(X_5, W_5)\end{aligned}$$

=

$$\underline{\max}\,\underline{\min}\left\{\begin{bmatrix}0 & 0 & I\end{bmatrix}, \begin{bmatrix}0 & 4 & I & 2 & 4I & 0 & 0 & 7 & 1\\ 0 & 0 & 0 & 0 & 0 & 8I & 5 & 0 & 1\\ 1 & 0 & 7I & 0 & 0 & 0 & 3 & 0 & 1\end{bmatrix}\right\} \cup$$

$$\underline{\max}\,\underline{\min}\left\{\begin{bmatrix}I & 3I & 0\end{bmatrix}, \begin{bmatrix}1 & I & 4 & 5 & 7 & 0 & 0 & 0 & 0\\ 0 & 1 & 0 & I & 2 & 1 & 0 & 0 & 8\\ 1 & 0 & 0 & 0 & 1 & 0 & I & 0 & 4I\end{bmatrix}\right\} \cup$$

$$\underline{\max}\,\underline{\min}\left\{\begin{bmatrix}5 & 0 & 2I\end{bmatrix}, \begin{bmatrix}9I & 0 & 2 & 0 & 0 & 0 & 1 & 4 & I\\ 0 & I & 3 & 0 & 0 & 7 & 0 & 5 & 0\\ 2 & 0 & 1 & 0 & 8 & 0 & I & 0 & 1\end{bmatrix}\right\} \cup$$

$$\underline{\max}\,\underline{\min}\left\{\begin{bmatrix}0 & 7I & 0\end{bmatrix}, \begin{bmatrix}0 & 3I & 0 & 0 & 0 & 2 & 7 & 1 & 0\\ 0 & 0 & 5 & 7 & I & 3 & 0 & 1 & 0\\ 7 & 0 & 0 & 0 & 2 & 0 & 0 & 0 & 8\end{bmatrix}\right\} \cup$$

$$\underline{\max}\,\underline{\min}\left\{\begin{bmatrix}2 & I & 6\end{bmatrix}, \begin{bmatrix}0 & 0 & 7 & 0 & 0 & 0 & 2 & 5 & I\\ I & 3 & 0 & 0 & 0 & 4 & 0 & 0 & 2\\ 0 & 0 & I & 7I & 8 & 0 & 0 & 0 & 1\end{bmatrix}\right\}$$

$$\begin{aligned}&= [I\ 0\ I\ 0\ 0\ 0\ I\ 0\ I] \cup [I\ 1\ I\ I\ 2\ 1\ 0\ 0\ 3I] \cup [5\ 0\ 2\ 0\ 2I\ 0\ I\ 4\ I] \cup [0\ 0\ 5\ 7I\ I\ 3\ 0\ 1\ 0] \cup [I\ I\ 2\ 6\ 6\ I\ 2\ 2\ I] \\ &= Y_1 \cup Y_2 \cup Y_3 \cup Y_4 \cup Y_5 \\ &= Y;\end{aligned}$$

now Y is a special neutrosophic row vector and each $Y_i$ is a $1 \times 9$ neutrosophic row vector, i = 1, 2, 3, 4, 5.

Now we calculate the using special max min find the value of $\{Y, W^T\}$ i.e.

$\max \min (Y, W^T)$

$= \max \min \{[Y_1 \cup Y_2 \cup \ldots \cup Y_5], [W_1^T \cup W_2^T \cup \ldots \cup W_5^T]\}$

$= \max \min (Y_1, W_1^T) \cup \max \min (Y_2, W_2^T) \cup \ldots \cup \max \min (Y_5, W_5^T)$

$$= \max \min \left\{ \begin{bmatrix} I & 0 & I & 0 & 0 & 0 & I & 0 & I \end{bmatrix}, \begin{bmatrix} 0 & 0 & 1 \\ 4 & 0 & 0 \\ I & 0 & 7I \\ 2 & 0 & 0 \\ 4I & 0 & 0 \\ 0 & 8I & 0 \\ 0 & 5 & 3 \\ 7 & 0 & 0 \\ 1 & 1 & 1 \end{bmatrix} \right\} \cup$$

$$\max \min \left\{ \begin{bmatrix} I & 1 & I & I & 2 & 1 & 0 & 0 & 3I \end{bmatrix}, \begin{bmatrix} 1 & 0 & 1 \\ I & 1 & 0 \\ 4 & 0 & 0 \\ 5 & I & 0 \\ 7 & 2 & 1 \\ 0 & 1 & 0 \\ 0 & 0 & I \\ 0 & 0 & 0 \\ 0 & 8 & 4I \end{bmatrix} \right\} \cup$$

$$\max\min \left\{ \begin{bmatrix} 5 & 0 & 2 & 0 & 2I & 0 & I & 4 & I \end{bmatrix}, \begin{bmatrix} 9I & 0 & 2 \\ 0 & I & 0 \\ 2 & 3 & 1 \\ 0 & 0 & 0 \\ 0 & 0 & 8 \\ 0 & 7 & 0 \\ 1 & 0 & I \\ 4 & 5 & 0 \\ I & 0 & 1 \end{bmatrix} \right\} \cup$$

$$\max\min \left\{ \begin{bmatrix} 0 & 0 & 5 & 7I & I & 3 & 0 & 1 & 0 \end{bmatrix}, \begin{bmatrix} 0 & 0 & 7 \\ 3I & 0 & 0 \\ 0 & 5 & 0 \\ 0 & 7 & 0 \\ 0 & I & 2 \\ 2 & 3 & 0 \\ 7 & 0 & 0 \\ 1 & 1 & 0 \\ 0 & 0 & 8 \end{bmatrix} \right\} \cup$$

$$\max\min \left\{ \begin{bmatrix} I & I & I & 6 & 6 & I & 2 & 2 & I \end{bmatrix}, \begin{bmatrix} 0 & I & 0 \\ 0 & 3 & 0 \\ 7 & 0 & I \\ 0 & 0 & 7I \\ 0 & 0 & 8 \\ 0 & 4 & 0 \\ 2 & 0 & 0 \\ 5 & 0 & 0 \\ I & 2 & 1 \end{bmatrix} \right\}$$

$= [I\ I\ I] \cup [2\ 3I\ 3I] \cup [5\ 4\ 2] \cup [2\ 5\ I] \cup [2\ I\ 6]$
$= P_1 \cup P_2 \cup P_3 \cup P_4 \cup P_5$
$= P,$

P is a special neutrosophic row vector, where each $P_i$ is a $1 \times 3$ neutrosophic row vector i = 1, 2, …, 5. Now we can as before calculate $m\underline{a}x\ m\underline{i}n$ (P, W) and so on.

Next we proceed on to show how $m\underline{a}x\ m\underline{i}n$ function works when we have a special neutrosophic mixed rectangular matrix. Suppose $V = V_1 \cup V_2 \cup \ldots \cup V_n$ ($n \geq 2$) be a special neutrosophic mixed rectangular matrix where each $V_i$ is a $t_i \times s_i$ ($t_i \neq s_i$) rectangular neutrosophic matrix i = 1, 2, …, n. We have atleast for one pair i and j; $i \neq j$, $t_i \neq t_j$ (or $s_i \neq s_j$), $1 \leq j, i \leq n$.

Let $X = X_1 \cup X_2 \cup \ldots \cup X_n$ ($n \geq 2$) special neutrosophic mixed row vector where each $X_i$ is a $1 \times t_i$ neutrosophic row vector, i = 1, 2, …, n. Now special $m\underline{a}x\ m\underline{i}n$ of {(X, V)}

$= m\underline{a}x\ m\underline{i}n \{(X_1 \cup X_2 \cup \ldots \cup X_n), (V_1 \cup V_2 \cup \ldots \cup V_n)\}$
$= m\underline{a}x\ m\underline{i}n \{X_1, V_1\} \cup m\underline{a}x\ m\underline{i}n \{X_2, V_2\} \cup \ldots \cup m\underline{a}x\ m\underline{i}n (X_n, V_n\}$
$= Y_1 \cup Y_2 \cup \ldots \cup Y_n$
$= Y$

where Y is a special neutrosophic mixed row vector. Now we find $m\underline{a}x\ m\underline{i}n$ of $\{Y, V^T\}$

$= m\underline{a}x\ m\underline{i}n \{(Y_1 \cup Y_2 \cup \ldots \cup Y_n), (V_1^T \cup V_2^T \cup \ldots \cup V_n^T)\}$
$= m\underline{a}x\ m\underline{i}n \{Y_1, V_1^T\} \cup m\underline{a}x\ m\underline{i}n \{Y_2, V_2^T\} \cup \ldots \cup m\underline{a}x\ m\underline{i}n \{Y_n, V_n^T\}$
$= Z_1 \cup Z_2 \cup \ldots \cup Z_n$
$= Z,$

where Z is a special neutrosophic mixed row vector. We can if need be find $m\underline{a}x\ m\underline{i}n$ {Z, V} and so on.

Now we illustrate this situation by the following example.

***Example 1.3.25:*** Let $V = V_1 \cup V_2 \cup V_3 \cup V_4 \cup V_5$ be a special neutrosophic mixed rectangular matrix given by

$$V = \begin{bmatrix} 0 & 7 & 0 & I & 5 & 3I & 1 \\ I & 0 & 0 & 0 & 6 & 5 & 2 \\ 6 & 2I & 1 & 0 & 7 & 0 & 3 \end{bmatrix} \cup$$

$$\begin{bmatrix} 5 & 3 & I & 7 \\ 0 & 2I & 0 & 0 \\ 8I & 0 & 0 & 6 \\ 1 & 7 & 0 & 0 \\ 0 & 0 & 8 & 6I \\ 1 & 0 & 0 & 2 \\ 0 & 0 & 6I & 1 \end{bmatrix} \cup$$

$$\begin{bmatrix} 0 & 3 & 1 & 0 & 9I \\ 9 & 0 & 0 & 3I & 2 \\ 1 & 4 & 6 & 0 & 7 \\ 0 & 0 & 0 & 8 & 0 \\ 0 & 2I & 0 & 0 & 0 \\ 7 & 0 & 2I & 1 & 9 \end{bmatrix} \cup$$

$$\begin{bmatrix} 5 & 0 & I & 0 & 1 & 0 & 6 & 0 \\ 0 & 9 & 0 & 7 & 2 & 4I & 7 & 2 \\ 0 & 2I & 6 & 0 & 3 & 5 & 8I & 3 \\ 2 & 0 & 0 & I & 0 & 0 & 9 & 4 \end{bmatrix} \cup$$

$$\begin{bmatrix} 8 & 0 & 0 & 7I & 0 & 0 & 0 & 9 \\ 4 & 0 & 6 & 0 & 0 & 2I & 0 & 1 \\ 0 & 0 & I & 0 & 0 & 0 & I & 3 \\ 2I & 3 & 0 & 2 & 7 & 0 & 0 & 2I \\ 1 & 2 & 3 & 4 & 5 & 6 & 7 & 8 \end{bmatrix}.$$

Suppose

$$\begin{aligned} X &= X_1 \cup X_2 \cup \ldots \cup X_5 \\ &= [6\ 0\ I] \cup [4\ 0\ 0\ 0\ 0\ I\ 0] \cup [0\ 0\ 2I\ 0\ 9\ 0] \cup [0\ 0\ 0\ 6] \\ &\quad \cup [I\ 0\ 0\ 7\ 0] \end{aligned}$$

where X is a special neutrosophic mixed row matrix.

$$\begin{aligned} &\max \min (X, V) \\ &= \max \min \{(X_1 \cup X_2 \cup \ldots \cup X_5), (V_1 \cup V_2 \cup \ldots \cup V_5)\} \\ &= \max \min \{X_1, V_1\} \cup \max \min \{X_2, V_2\} \cup \ldots \cup \max \min \{X_5, V_5\} \end{aligned}$$

$$= \max \min \left\{ \begin{bmatrix} 6 & 0 & I \end{bmatrix}, \begin{bmatrix} 0 & 7 & 0 & I & 5 & 3I & 1 \\ I & 0 & 0 & 0 & 6 & 5 & 2 \\ 6 & 2I & 1 & 0 & 7 & 0 & 3 \end{bmatrix} \right\} \cup$$

$$\max \min \left\{ \begin{bmatrix} 4 & 0 & 0 & 0 & 0 & I & 0 \end{bmatrix}, \begin{bmatrix} 5 & 3 & I & 7 \\ 0 & 2I & 0 & 0 \\ 8I & 0 & 0 & 6 \\ 1 & 7 & 0 & 0 \\ 0 & 0 & 8 & 6I \\ 1 & 0 & 0 & 2 \\ 0 & 0 & 6I & 1 \end{bmatrix} \right\} \cup$$

$$\max\min \left\{ \begin{bmatrix} 0 & 0 & 2I & 0 & 9 & 0 \end{bmatrix}, \begin{bmatrix} 0 & 3 & 1 & 0 & 9I \\ 9 & 0 & 0 & 3I & 2 \\ 1 & 4 & 6 & 0 & 7 \\ 0 & 0 & 0 & 8 & 0 \\ 0 & 2I & 0 & 0 & 0 \\ 7 & 0 & 2I & 1 & 9 \end{bmatrix} \right\} \cup$$

$$\max\min \left\{ \begin{bmatrix} 0 & 0 & 0 & 6 \end{bmatrix}, \begin{bmatrix} 5 & 0 & I & 0 & 1 & 0 & 6 & 0 \\ 0 & 9 & 0 & 7 & 2 & 4I & 7 & 2 \\ 0 & 2I & 6 & 0 & 3 & 5 & 8I & 3 \\ 2 & 0 & 0 & I & 0 & 0 & 9 & 4 \end{bmatrix} \right\} \cup$$

$$\max\min \left\{ \begin{bmatrix} I & 0 & 0 & 7 & 0 \end{bmatrix}, \begin{bmatrix} 8 & 0 & 0 & 4I & 0 & 0 & 0 & 9 \\ 4 & 0 & 6 & 0 & 0 & 2I & 0 & 1 \\ 0 & 0 & I & 0 & 0 & 0 & I & 3 \\ 2I & 3 & 0 & 2 & 7 & 0 & 0 & 2I \\ 1 & 2 & 3 & 4 & 5 & 6 & 7 & 8 \end{bmatrix} \right\}$$

= [I 6 I I 5 3I I] $\cup$ [4 3 I 4] $\cup$ [1 2I 2I 0 2I] $\cup$ [2 0 0 I 0 0 6 4] $\cup$ [2I 3 0 2 7 0 0 2I]
= $Y_1 \cup Y_2 \cup Y_3 \cup Y_4 \cup Y_5$
= Y

where Y is a special neutrosophic mixed row vector. Now we use Y and $V^T$ and find the value of

max min $(Y, V^T)$
= max min $\{(Y_1 \cup Y_2 \cup \ldots \cup Y_5), (V_1^T \cup V_2^T \cup \ldots \cup V_5^T)\}$
= max min $\{Y_1, V_1^T\}$ $\cup$ max min $\{Y_2, V_2^T\}$ $\cup$ max min $\{Y_3, V_3^T\}$ $\cup$ max min $\{Y_4, V_4^T\}$ $\cup$ max min $\{Y_5, V_5^T\}$

$$= \max_{_} \min_{_} \left\{ \begin{bmatrix} I & 6 & I & I & 5 & 3I & I \end{bmatrix}, \begin{bmatrix} 0 & I & 6 \\ 7 & 0 & 2I \\ 0 & 0 & 1 \\ I & 0 & 0 \\ 5 & 6 & 7 \\ 3I & 5 & 0 \\ 1 & 2 & 3 \end{bmatrix} \right\} \cup$$

$$\max_{_} \min_{_} \left\{ \begin{bmatrix} 4 & 3 & I & 4 \end{bmatrix}, \begin{bmatrix} 5 & 0 & 8I & 1 & 0 & 1 & 0 \\ 3 & 2I & 0 & 7 & 0 & 0 & 0 \\ I & 0 & 0 & 0 & 8 & 0 & 6I \\ 7 & 0 & 6 & 0 & 6I & 2 & 1 \end{bmatrix} \right\} \cup$$

$$\max_{_} \min_{_} \left\{ \begin{bmatrix} 1 & 2I & 2I & 0 & 2I \end{bmatrix}, \begin{bmatrix} 0 & 9 & 1 & 0 & 0 & 7 \\ 3 & 0 & 4 & 0 & 2I & 0 \\ 1 & 0 & 6 & 0 & 0 & 2I \\ 0 & 3I & 0 & 8 & 0 & 1 \\ 9I & 2 & 7 & 0 & 0 & 9 \end{bmatrix} \right\} \cup$$

$$\max_{_} \min_{_} \left\{ \begin{bmatrix} 2 & 0 & 0 & I & 0 & 0 & 6 & 4 \end{bmatrix}, \begin{bmatrix} 5 & 0 & 0 & 2 \\ 0 & 9 & 2I & 0 \\ I & 0 & 6 & 0 \\ 0 & 7 & 0 & I \\ 1 & 2 & 3 & 0 \\ 0 & 4I & 5 & 0 \\ 6 & 7 & 8I & 9 \\ 0 & 2 & 3 & 4 \end{bmatrix} \right\} \cup$$

$$\underline{\max}\ \underline{\min}\ \left\{ \begin{bmatrix} 2I & 3 & 0 & 2 & 7 & 0 & 0 & 2I \end{bmatrix}, \begin{bmatrix} 8 & 4 & 0 & 2I & 1 \\ 0 & 0 & 0 & 3 & 2 \\ 0 & 6 & I & 0 & 3 \\ 4I & 0 & 0 & 2 & 4 \\ 0 & 0 & 0 & 7 & 5 \\ 0 & 2I & 0 & 0 & 6 \\ 0 & 0 & I & 0 & 7 \\ 9 & 1 & 3 & 2I & 8 \end{bmatrix} \right\}.$$

$= [6\ 5\ 5] \cup [4\ 2I\ 4\ 3\ 4\ 2\ I] \cup [2I\ 2I\ 2I\ 0\ 2I\ 2I] \cup [6\ 6\ 6\ 6] \cup [2I\ 2I\ I\ 3\ 2]$

$= T_1 \cup T_2 \cup T_3 \cup T_4 \cup T_5$

$= T$

where T is a special neutrosophic mixed row vector. Now we can find max min (T, V) and so on.

Now we define special $\underline{\max}$ $\underline{\min}$ function for special neutrosophic mixed matrix. Let $P = P_1 \cup P_2 \cup \ldots \cup P_n$ ($n \geq 2$) where $P_i$ is a $m_i \times t_i$ ($m_i \neq t_i$) are neutrosophic rectangular matrix and $P_j$'s are $p_j \times p_j$ neutrosophic matrix $i \neq j$; $1 \leq i, j \leq n$.

Suppose

$X = X_1 \cup X_2 \cup \ldots \cup X_n$

($n \geq 2$) where $X_i$'s are $1 \times m_i$ neutrosophic mixed row vector and $X_j$'s are $1 \times p_j$ neutrosophic row vector $1 \leq i, j \leq n$.

We define

$\underline{\max}\ \underline{\min}\ (X, P)$

$= \underline{\max}\ \underline{\min}\ \{(X_1 \cup X_2 \cup \ldots \cup X_n), (P_1 \cup P_2 \cup \ldots \cup P_n)\}$

$= \underline{\max}\ \underline{\min}\ \{X_1, P_1\} \cup \underline{\max}\ \underline{\min}\ \{X_2, P_2\} \cup \ldots \cup \underline{\max}\ \underline{\min}\ \{X_n, P_n\}$

$= Y_1 \cup Y_2 \cup \ldots \cup Y_n$

$= Y$

where Y is a neutrosophic mixed row vector. We find at the next step $\underline{max}\ \underline{min}\ \{Y, P^{ST}\}$ (we have defined in pages 85-6 the notion of special transpose).

Now

$$
\begin{aligned}
\underline{max}\ \underline{min}\ \{Y, P^{ST}\} &= \underline{max}\ \underline{min}\ \{(Y_1 \cup Y_2 \cup \ldots \cup Y_n), \left(P_1^{ST} \cup P_2^{ST} \cup \ldots \cup P_n^{ST}\right)\} \\
&= \underline{max}\ \underline{min}\ (Y_1,\ P_1^{ST}) \cup \underline{max}\ \underline{min}\ (Y_2,\ P_2^{ST}) \cup \ldots \cup \underline{max}\ \underline{min}\ (Y_n,\ P_n^{ST}) \\
&= Z_1 \cup Z_2 \cup \ldots \cup Z_n \\
&= Z
\end{aligned}
$$

where Z is again a neutrosophic mixed row vector. We can find $\underline{max}\ \underline{min}\ (Z, P)$ and so on.

We illustrate this situation by the following example.

***Example 1.3.26:*** Let $P = P_1 \cup P_2 \cup P_3 \cup P_4 \cup P_5 \cup P_6$

$$
= \begin{bmatrix} 3 & 0 & I \\ 1 & 2 & 3 \\ 0 & 4I & 0 \\ 5 & 6 & 7 \\ 0 & 0 & 8I \\ 9 & 8 & 7 \\ 4I & 0 & 0 \end{bmatrix} \cup \begin{bmatrix} 5 & 0 & 1 & 2 & I \\ 0 & 3I & 0 & 0 & 4 \\ 0 & 0 & 7 & I & 6 \\ I & 8 & 0 & 0 & 0 \\ 0 & 0 & 0 & 4 & 0 \end{bmatrix} \cup \begin{bmatrix} 3 & 0 & 4 & 2 \\ I & 0 & 0 & 8 \\ 0 & 9 & 7 & 0 \\ 6 & I & 2 & 0 \\ 0 & 0 & 5I & 9 \\ I & 0 & 0 & 0 \\ 5 & 3 & 0 & 5I \end{bmatrix} \cup
$$

$$
\begin{bmatrix} 0 & 0 & 6I & 2 & 0 & 1 & 0 \\ 2 & 0 & 0 & 0 & 8 & 2 & 9 \\ 0 & I & 0 & 0 & 3I & 3 & 0 \\ 0 & 0 & 8 & I & 0 & 7 & 0 \end{bmatrix} \cup
$$

$$\begin{bmatrix} 0 & 4 & 2I & 0 \\ 8 & 0 & 0 & I \\ 3 & I & 0 & 8 \\ 0 & 0 & 5I & 0 \end{bmatrix} \cup \begin{bmatrix} 0 & 9 & 0 \\ 8 & 0 & 0 \\ 0 & 0 & 7 \\ I & 0 & 0 \\ 0 & 2I & 0 \\ 0 & 0 & 3I \end{bmatrix}$$

be a neutrosophic mixed matrix.

Suppose

$X = X_1 \cup X_2 \cup \ldots \cup X_6$

$= [1\ 0\ 0\ 0\ 0\ 4\ I] \cup [0\ I\ 7\ 0\ 0] \cup [0\ 2\ 0\ 0\ 0\ I\ 0] \cup [0\ 3\ 0\ I] \cup [0\ 3\ 4I\ 0] \cup [9\ 0\ 0\ 0\ 1\ 0]$

be the neutrosophic mixed row vector. To find

$\max \min (X, P)$

$= \max \min \{(X_1 \cup X_2 \cup \ldots \cup X_6), (V_1 \cup V_2 \cup \ldots \cup V_6)\}$

$= \max \min (X_1, V_1) \cup \max \min (X_2, V_2) \cup \ldots \cup \max \min (X_6, V_6)$

$$= \max \min \left\{ \begin{bmatrix} 1 & 0 & 0 & 0 & 0 & 4 & I \end{bmatrix}, \begin{bmatrix} 3 & 0 & I \\ 1 & 2 & 3 \\ 0 & 4I & 0 \\ 5 & 6 & 7 \\ 0 & 0 & 8I \\ 9 & 8 & 7 \\ 4I & 0 & 0 \end{bmatrix} \right\} \cup$$

$$\max \min \left\{ \begin{bmatrix} 0 & I & 7 & 0 & 0 \end{bmatrix}, \begin{bmatrix} 5 & 0 & 1 & 2 & I \\ 0 & 3I & 0 & 0 & 4 \\ 0 & 0 & 7 & I & 6 \\ I & 8 & 0 & 0 & 0 \\ 0 & 0 & 0 & 4 & 0 \end{bmatrix} \right\} \cup$$

$$\max\min\left\{[0\ \ 2\ \ 0\ \ 0\ \ 0\ \ I\ \ 0],\ \begin{bmatrix} 3 & 0 & 4 & 2 \\ I & 0 & 0 & 8 \\ 0 & 9 & 7 & 0 \\ 6 & I & 2 & 0 \\ 0 & 0 & 5I & 9 \\ I & 0 & 0 & 0 \\ 5 & 3 & 0 & 5I \end{bmatrix}\right\} \cup$$

$$\max\min\left\{[0\ \ 3\ \ 0\ \ I],\ \begin{bmatrix} 0 & 0 & 6I & 2 & 0 & 1 & 0 \\ 2 & 0 & 0 & 0 & 8 & 2 & 9 \\ 0 & I & 0 & 0 & 3I & 3 & 0 \\ 0 & 0 & 8 & I & 0 & 7 & 0 \end{bmatrix}\right\} \cup$$

$$\max\min\left\{[0\ \ 3\ \ 4I\ \ 0],\ \begin{bmatrix} 0 & 4 & 2I & 0 \\ 8 & 0 & 0 & I \\ 3 & I & 0 & 8 \\ 0 & 0 & 5I & 0 \end{bmatrix}\right\} \cup$$

$$\max\min\left\{[9\ \ 0\ \ 0\ \ 0\ \ 1\ \ 0],\ \begin{bmatrix} 0 & 9 & 0 \\ 8 & 0 & 0 \\ 0 & 0 & 7 \\ I & 0 & 0 \\ 0 & 2I & 0 \\ 0 & 0 & 3I \end{bmatrix}\right\}$$

= [4 4 4] ∪ [0 I 7 I 6] ∪ [I 0 0 2] ∪ [2 0 I I 3 2 3] ∪ [3 I 0 4I] ∪ [0 9 0]

= $Y_1 \cup Y_2 \cup \ldots \cup Y_6$

= Y

is again a special neutrosophic mixed row vector. Using Y we find out

$\max\min \{Y, P^{ST}\}$

$= \max\min \{(Y_1 \cup Y_2 \cup \ldots \cup Y_6) \cup ((P_1^T \cup P_2 \cup P_3^T \cup P_4^T \cup P_5 \cup P_6^T)\}$

$= \max\min \{Y_1, P_1^T\} \cup \max\min \{Y_2, P_2\} \cup \max\min \{Y_3, P_3^T\} \cup \max\min \{Y_4, P_4^T\} \cup \max\min \{Y_5, P_5\} \cup \max\min \{Y_6, P_6^T\}$

=

$$\max\min \left\{ \begin{bmatrix} 4 & 4 & 4 \end{bmatrix}, \begin{bmatrix} 3 & 1 & 0 & 5 & 0 & 9 & 4I \\ 0 & 2 & 4I & 6 & 0 & 8 & 0 \\ I & 3 & 0 & 7 & 8I & 7 & 0 \end{bmatrix} \right\} \cup$$

$$\max\min \left\{ \begin{bmatrix} 0 & I & 7 & I & 6 \end{bmatrix}, \begin{bmatrix} 5 & 0 & 1 & 2 & I \\ 0 & 3I & 0 & 0 & 4 \\ 0 & 0 & 7 & I & 6 \\ I & 8 & 0 & 0 & 0 \\ 0 & 0 & 0 & 4 & 0 \end{bmatrix} \right\} \cup$$

$$\max\min \left\{ \begin{bmatrix} I & 0 & 0 & 2 \end{bmatrix}, \begin{bmatrix} 3 & I & 0 & 6 & 0 & I & 5 \\ 0 & 0 & 9 & I & 0 & 0 & 3 \\ 4 & 0 & 7 & 2 & 5I & 0 & 0 \\ 2 & 8 & 0 & 0 & 9 & 0 & 5I \end{bmatrix} \right\} \cup$$

$$\max\min \left\{ \begin{bmatrix} 2 & 0 & I & I & 3 & 2 & 3 \end{bmatrix}, \begin{bmatrix} 0 & 2 & 0 & 0 \\ 0 & 0 & I & 0 \\ 6I & 0 & 0 & 8 \\ 2 & 0 & 0 & I \\ 0 & 8 & 3I & 0 \\ 1 & 2 & 3 & 7 \\ 0 & 9 & 0 & 0 \end{bmatrix} \right\} \cup$$

$$\max \min \left\{ \begin{bmatrix} 3 & I & 0 & 4I \end{bmatrix}, \begin{bmatrix} 0 & 4 & 2I & 0 \\ 8 & 0 & 0 & I \\ 3 & I & 0 & 8 \\ 0 & 0 & 5I & 0 \end{bmatrix} \right\} \cup$$

$$\max \min \left\{ \begin{bmatrix} 0 & 9 & 0 \end{bmatrix}, \begin{bmatrix} 0 & 8 & 0 & I & 0 & 0 \\ 9 & 0 & 0 & 0 & 2I & 0 \\ 0 & 0 & 7 & 0 & 0 & 3I \end{bmatrix} \right\}$$

$= [3\ 3\ 4I\ 4\ 4\ 4\ 4I] \cup [I\ I\ 7\ 4\ 6] \cup [2\ 2\ 0\ I\ 2\ I\ 2] \cup [I\ 3\ 3I\ 2] \cup [I\ 3\ 4I\ I] \cup [9\ 0\ 0\ 0\ 2I\ 0]$

$= Z_1 \cup Z_2 \cup \ldots \cup Z_6$

$= Z.$

We see Z is a special neutrosophic mixed vector / matrix. Now using Z and P we can find max min {Z, P} and so on.

Now using a special neutrosophic mixed matrix we can find using the min max operator for any special neutrosophic mixed row vector.

We illustrate this by the following example.

***Example 1.3.27:*** Let $T = T_1 \cup T_2 \cup T_3 \cup T_4 \cup T_5$ be a special neutrosophic mixed matrix, where

$$T = \begin{bmatrix} 0 & 5 & 0 \\ 2I & 0 & 0 \\ 0 & 0 & 7 \\ 4 & I & 0 \\ 8 & 0 & 6I \\ 0 & 7 & 3I \\ 2 & 9 & 7 \end{bmatrix} \cup \begin{bmatrix} 0 & 7 & 9 & 0 \\ 4 & 0 & 0 & 2I \\ 0 & 8I & 0 & 8 \\ 3I & 0 & 5I & 0 \end{bmatrix} \cup$$

$$\begin{bmatrix} 0 & 0 & 4 & 8 & 0 & 2 & 0 \\ 1 & 7 & 0 & 0 & 9 & 9 & 4 \\ 2 & I & 0 & I & 0 & 0 & 0 \\ 0 & 0 & I & 0 & 3I & 6 & 2I \end{bmatrix}$$

$$\cup \begin{bmatrix} 8 & 0 & 1 \\ 4 & 5I & 0 \\ 0 & 0 & 6I \end{bmatrix} \cup \begin{bmatrix} 9 & 0 & 2I & 0 & 1 & 6 & 0 & 0 \\ 0 & 8 & 0 & I & 0 & 0 & 7I & 1 \\ 0 & 0 & 6 & 0 & 2I & 0 & 0 & 0 \\ 0 & 6I & 0 & 0 & 0 & 8I & 0 & 0 \\ 0 & 0 & 0 & 3 & 9 & 0 & 8 & 2I \end{bmatrix}$$

be a special neutrosophic mixed matrix.
Suppose

$$\begin{aligned} I &= I_1 \cup I_2 \cup I_3 \cup I_4 \cup I_5 \\ &= [0\ 9\ 0\ 0\ 1\ 0\ 1] \cup [0\ 2\ 0\ 0] \cup [0\ 0\ 0\ 6] \cup [6\ 0\ 0] \cup \\ &\quad [8\ 0\ 0\ 0\ 0] \end{aligned}$$

be a special neutrosophic mixed row vector/matrix. To find the effect of I on T using $m\underline{a}x\ m\underline{i}n$ operator.

$m\underline{a}x\ m\underline{i}n\ (I, T)$

$$\begin{aligned} &= m\underline{a}x\ m\underline{i}n\ \{(I_1 \cup I_2 \cup \ldots \cup I_5), (T_1 \cup T_2 \cup \ldots \cup T_5)\} \\ &= m\underline{a}x\ m\underline{i}n\ (I_1, T_1) \cup m\underline{a}x\ m\underline{i}n\ (I_2, T_2) \cup \ldots \cup m\underline{a}x \\ &\quad m\underline{i}n\ (I_5, T_5) \end{aligned}$$

$$= m\underline{a}x\ m\underline{i}n \left\{ \begin{bmatrix} 0 & 9 & 0 & 0 & 1 & 0 & 1 \end{bmatrix}, \begin{bmatrix} 0 & 5 & 0 \\ 2I & 0 & 0 \\ 0 & 0 & 7 \\ 4 & I & 0 \\ 8 & 0 & 6I \\ 0 & 7 & 3I \\ 2 & 9 & 7 \end{bmatrix} \right\} \cup$$

$$\underline{\max}\,\underline{\min}\left\{\begin{bmatrix}0 & 2 & 0 & 0\end{bmatrix}, \begin{bmatrix}0 & 7 & 9 & 0\\ 4 & 0 & 0 & 2I\\ 0 & 8I & 0 & 8\\ 3I & 0 & 5I & 0\end{bmatrix}\right\}$$

$$\underline{\max}\,\underline{\min}\left\{\begin{bmatrix}0 & 0 & 0 & 6\end{bmatrix}, \begin{bmatrix}0 & 0 & 4 & 8 & 0 & 2 & 0\\ 1 & 7 & 0 & 0 & 9 & 9 & 4\\ 2 & I & 0 & I & 0 & 0 & 0\\ 0 & 0 & I & 0 & 3I & 6 & 2I\end{bmatrix}\right\}$$

$$\underline{\max}\,\underline{\min}\left\{\begin{bmatrix}6 & 0 & 0\end{bmatrix}, \begin{bmatrix}8 & 0 & 1\\ 4 & 5I & 0\\ 0 & 0 & 6I\end{bmatrix}\right\} \cup$$

$$\underline{\max}\,\underline{\min}\left\{\begin{bmatrix}8 & 0 & 0 & 0 & 0\end{bmatrix}, \begin{bmatrix}9 & 0 & 2I & 0 & 1 & 6 & 0 & 0\\ 0 & 8 & 0 & I & 0 & 0 & 7I & 1\\ 0 & 0 & 6 & 0 & 2I & 0 & 0 & 0\\ 0 & 6I & 0 & 0 & 0 & 8I & 0 & 0\\ 0 & 0 & 0 & 3 & 9 & 0 & 8 & 2I\end{bmatrix}\right\}$$

$= [2I\ 1\ 1] \cup [2\ 0\ 0\ 2I] \cup [0\ 0\ I\ 0\ 3I\ 6\ 2I] \cup [6\ 0\ 1] \cup [8\ 0\ 2I\ 0\ I\ 6\ 0\ 0]$

$= P_1 \cup P_2 \cup P_3 \cup P_4 \cup P_5$

$= P$

where P is a special neutrosophic mixed row vector. We now calculate the resultant of P

$\underline{\max}\,\underline{\min}\,\{P, T^{ST})$

$= \underline{\max}\,\underline{\min}\,\{(P_1 \cup P_2 \cup \ldots \cup P_5), \left(\left(T_1^{ST} \cup T_2^{ST} \cup \ldots \cup T_5^{ST}\right)\right\}$

$= \underline{\max}\,\underline{\min}\,(P_1, T_1^{ST}) \cup \underline{\max}\,\underline{\min}\,(P_2, T_2^{ST}) \cup \ldots \cup \underline{\max}\,\underline{\min}\,(P_5, T_5^{ST})$

$$= \max \min \left\{ \begin{bmatrix} 2I & 1 & 1 \end{bmatrix}, \begin{bmatrix} 0 & 2I & 0 & 4 & 8 & 0 & 2 \\ 5 & 0 & 0 & I & 0 & 7 & 9 \\ 0 & 0 & 7 & 0 & 6I & 3I & 7 \end{bmatrix} \right\} \cup$$

$$\max \min \left\{ \begin{bmatrix} 2 & 0 & 0 & 2I \end{bmatrix}, \begin{bmatrix} 0 & 7 & 9 & 0 \\ 4 & 0 & 0 & 2I \\ 0 & 8I & 0 & 8 \\ 3I & 0 & 5I & 0 \end{bmatrix} \right\} \cup$$

$$\max \min \left\{ \begin{bmatrix} 0 & 0 & I & 0 & 3I & 6 & 2I \end{bmatrix}, \begin{bmatrix} 0 & 1 & 2 & 0 \\ 0 & 7 & I & 0 \\ 4 & 0 & 0 & I \\ 8 & 0 & I & 0 \\ 0 & 9 & 0 & 3I \\ 2 & 9 & 0 & 6 \\ 0 & 4 & 0 & 2I \end{bmatrix} \right\} \cup$$

$$\max \min \left\{ \begin{bmatrix} 6 & 0 & 1 \end{bmatrix}, \begin{bmatrix} 8 & 0 & 1 \\ 4 & 5I & 0 \\ 0 & 0 & 6I \end{bmatrix} \right\} \cup$$

$$\max \min \left\{ \begin{bmatrix} 8 & 0 & 2I & 0 & I & 6 & 0 & 0 \end{bmatrix}, \begin{bmatrix} 9 & 0 & 0 & 0 & 0 \\ 0 & 8 & 0 & 6I & 0 \\ 2I & 0 & 6 & 0 & 0 \\ 0 & I & 0 & 0 & 3 \\ 1 & 0 & 2I & 0 & 9 \\ 6 & 0 & 0 & 8I & 0 \\ 0 & 7I & 0 & 0 & 8 \\ 0 & 1 & 0 & 0 & 2I \end{bmatrix} \right\}$$

$$= [1\ 2I\ 1\ 2I\ 2I\ 1\ 2I] \cup [2I\ 2\ 2I\ 0] \cup [2\ 6\ 0\ 6] \cup [6\ 0\ 1] \cup [8\ 0\ I\ 6\ I]$$
$$= Z_1 \cup Z_2 \cup \ldots \cup Z_5$$
$$= Z;$$

where Z is a special neutrosophic mixed row vector.

Z can now work on T as follows.

m$\underline{a}$x m$\underline{i}$n {Z, T}

$$= \text{m}\underline{\text{a}}\text{x m}\underline{\text{i}}\text{n } \{(Z_1 \cup Z_2 \cup \ldots \cup Z_5), (T_1 \cup T_2 \cup \ldots \cup T_5)\}$$
$$= \text{m}\underline{\text{a}}\text{x m}\underline{\text{i}}\text{n } \{Z_1, T\} \cup \text{m}\underline{\text{a}}\text{x m}\underline{\text{i}}\text{n } \{Z_2, T_2\} \cup \ldots \cup \text{m}\underline{\text{a}}\text{x m}\underline{\text{i}}\text{n } \{Z_5, T_5\}$$
$$= R_1 \cup R_2 \cup \ldots \cup R_5$$
$$= R,$$

where R is a special neutrosophic mixed row vector.

We can proceed on to find m$\underline{a}$x m$\underline{i}$n $\{R, T^{ST}\}$ and so on.

Other types of operation on the special fuzzy neutrosophic matrices will be described in the following.

Recall a special neutrosophic matrix will be known as the special fuzzy neutrosophic matrix if its entries are from $N_f = \{[0\ I] \cup [0\ 1]\}$ where $N_f$ is the fuzzy neutrosophic interval.

Any element from $N_f$ will be of the form y = a + bI where a and b $\in$ [0, 1]. Thus any typical element will be 0.7 + 0.5I. If x $\in N_f$ is of the form x = 0.7 I or I then we call x to be an absolute or pure neutrosophic fuzzy number or pure fuzzy neutrosophic number. y will be known as a fuzzy neutrosophic number or neutrosophic fuzzy number if and only if a $\neq$ 0 and b $\neq$ 0. Thus in a pure fuzzy neutrosophic number is one in which a = 0. If in y = a + bI, b = 0 then we call y to be a fuzzy number. Thus all special fuzzy neutrosophic matrices are special neutrosophic matrices but we see all special neutrosophic matrices need not be special fuzzy neutrosophic matrices.

***Example 1.3.28:*** Let

$$P = \begin{bmatrix} 0.8 & 7 & 5I & 2 & 0.7I \end{bmatrix} \cup \begin{bmatrix} 0.8I & 7 & I \\ 0 & 0 & 8 \\ I & 0 & 0.3 \\ 4 & 18 & 0.4I \end{bmatrix} \cup \begin{bmatrix} 9 \\ 0 \\ 0.8I \\ 0.6 \\ 6 \end{bmatrix}$$

is a special neutrosophic mixed rectangular matrix. Clearly P is not a special neutrosophic fuzzy mixed rectangular matrix.

Now as in case of special neutrosophic matrices we can define operations on special fuzzy neutrosophic matrices also. We can apply threshold and updating operation at each stage and obtain the resultant. This sort of operation will be used in special fuzzy models which will be defined in chapter one of this book.

Now we define the special operation on the special fuzzy row vector with special square matrix T.

Let $T = T_1 \cup T_2 \cup \ldots \cup T_n$ be a special fuzzy neutrosophic square matrix. Let each $T_i$ be a $m \times n$ fuzzy neutrosophic matrix. $X = X_1 \cup X_2 \cup \ldots \cup X_n$ be a special fuzzy neutrosophic row vector where each $X_i$ is a fuzzy neutrosophic row vector taking its value from the set {0, 1, I}; i = 1, 2, …, n.

Now we define a special operation using X and T.

$$\begin{aligned} X \text{ o } T &= (X_1 \cup \ldots X_n) \text{ o } (T_1 \cup \ldots \cup T_n) \\ &= X_1 \text{ o } T_1 \cup X_2 \text{ o } T_2 \cup \ldots \cup X_n \text{ o } T_n \end{aligned}$$

where the operation $X_i$ o $T_i$ $1 \leq i \leq$ is described in the following.

Let

$$X_i = [0\ I\ 0\ 0\ 0\ 1]$$

and

$$T_i = \begin{bmatrix} 0 & 0 & 0 & 0 & 0 & 1 \\ 1 & 0 & 0 & 0 & I & 0 \\ 0 & 0 & 0 & 1 & 0 & I \\ I & 0 & 0 & 0 & 1 & 0 \\ 0 & 0 & 0 & 0 & 0 & 1 \\ 0 & 0 & 0 & 1 & 0 & 0 \end{bmatrix}$$

be the fuzzy neutrosophic matrix. We find

$$X_i \text{ o } T_i = \begin{bmatrix} 0 & I & 0 & 0 & 0 & 1 \end{bmatrix} \text{ o } \begin{bmatrix} 0 & 0 & 0 & 0 & 0 & 1 \\ 1 & 0 & 0 & 0 & I & 0 \\ 0 & 0 & 0 & 1 & 0 & I \\ I & 0 & 0 & 0 & 1 & 0 \\ 0 & 0 & 0 & 0 & 0 & 1 \\ 0 & 0 & 0 & 1 & 0 & 0 \end{bmatrix}$$

= [I 0 0 1 I 0] (we threshold any resultant vector in the following way. Suppose

$$X_j \text{ o } T_j = [a_1\ a_2 \ldots a_n]$$

then if $a_i$ is a any real number than

$a_i = 0$ if $a_i \leq 0$
$a_i = 1$ if $a_i > 0$

if $a_i$ is a pure neutrosophic number of the form say $m_iI$ and
if $m_i \leq 0$ then $a_i = 0$
if $m_i > 0$ then $a_i = I$.

If $a_i$ is not a pure neutrosophic number and is of the form $a_i = t_i + s_iI$ and if $t_i > s_i$ and if $t_i \leq 0$ then $a_i = 0$. If $t_i > 0$ then $a_i = 1$. Suppose $s_i > t_i$ and if $s_i \leq 0$ than $a_i = 0$. If $s_i > 0$ then $a_i = 1$. If $s_i = t_i$ then $a_i = I$. We using this type of thresholding technique work with the resultant vector. The updating is done as in case of fuzzy vectors if the $t^{th}$ coordinate is 1 in the starting and if in the resultant it becomes 0 or I we replace the t^th^ coordinate once again by 1).
(So

$$X \text{ o } T_i = [I\ 0\ 0\ 1\ I\ 0] = Y'_i$$

after updating and thresholding $Y'_i$ becomes equal to $Y_i$ = [I I 0 1 I 1]. Now we find

$$Y_i \text{ o } T_i = \begin{bmatrix} I & I & 0 & 1 & I & 1 \end{bmatrix} \text{ o } \begin{bmatrix} 0 & 0 & 0 & 0 & 0 & 1 \\ 1 & 0 & 0 & 0 & I & 0 \\ 0 & 0 & 0 & 1 & 0 & I \\ I & 0 & 0 & 0 & 1 & 0 \\ 0 & 0 & 0 & 0 & 0 & 1 \\ 0 & 0 & 0 & 1 & 0 & 0 \end{bmatrix}$$

$$= \quad [2I\ 0\ 0\ 1\ 1 + I\ 2I]$$

after updating and thresholding we get Z = [I I 0 1 I 1].

$$Z \text{ o } T_i = \begin{bmatrix} I & I & 0 & 1 & I & 1 \end{bmatrix} \text{ o } \begin{bmatrix} 0 & 0 & 0 & 0 & 0 & 1 \\ 1 & 0 & 0 & 0 & I & 0 \\ 0 & 0 & 0 & 1 & 0 & I \\ I & 0 & 0 & 0 & 1 & 0 \\ 0 & 0 & 0 & 0 & 0 & 1 \\ 0 & 0 & 0 & 1 & 0 & 0 \end{bmatrix}$$

$$= \quad [2I\ 0\ 0\ 1\ I + 1\ 2I]$$

after updating and thresholding we get S = [I I 0 1 I 1] which is a fixed point).
Now we find

X o T

$$= \quad (X_1 \cup X_2 \cup \ldots \cup X_n) \text{ o } (T_1 \cup T_2 \cup \ldots \cup T_n)$$
$$= \quad X_1 \text{ o } T_1 \cup X_2 \cup T_2 \cup \ldots \cup X_n \text{ o } T_n$$

(where the operations are carried out as described above)

$= Y'_1 \cup Y'_2 \cup \ldots \cup Y'_n$

$= Y'$.

Now Y' is updated and thresholded to $Y = Y_1 \cup Y_2 \cup \ldots \cup Y_n$. We find

Y o T

$= (Y_1 \cup Y_2 \cup \ldots \cup Y_n) \text{ o } (T_1 \cup T_2 \cup \ldots \cup T_n)$

$= Y_1 \text{ o } T_1 \cup Y_2 \text{ o } T_2 \cup \ldots \cup Y_n \text{ o } T_n$

$= Z'_1 \cup Z'_2 \cup \ldots \cup Z'_n$

$= Z'$.

Now Z' is updated and thresholded so that we obtain Z and we can find Z o T and so on.

Now we illustrate the situation by the following example.

***Example 1.3.29:*** Let $T = T_1 \cup T_2 \cup T_3 \cup T_4 \cup T_5$

$$= \begin{bmatrix} 0 & 1 & 0 & 0 & I \\ 1 & 0 & 0 & 1 & 0 \\ 0 & 0 & 0 & I & 1 \\ I & 1 & 0 & 0 & 0 \\ 0 & 0 & 1 & 0 & 0 \end{bmatrix} \cup \begin{bmatrix} 0 & 1 & 0 & 0 & 0 \\ 1 & 0 & 0 & 0 & 1 \\ 0 & I & 0 & 0 & -1 \\ 0 & 0 & I & 0 & 0 \\ 0 & 0 & 1 & 0 & 0 \end{bmatrix} \cup$$

$$\begin{bmatrix} 0 & I & 0 & 0 & 0 \\ 1 & 0 & 0 & 0 & I \\ 0 & 1 & 0 & I & 0 \\ 0 & 0 & 0 & 0 & 1 \\ 1 & I & 0 & 0 & 0 \end{bmatrix} \cup \begin{bmatrix} 0 & 1 & 0 & 0 & 1 \\ 0 & 0 & 1 & I & 0 \\ 1 & I & 0 & 0 & 0 \\ 0 & 1 & I & 0 & 0 \\ 0 & 1 & 0 & 0 & 0 \end{bmatrix} \cup \begin{bmatrix} 0 & I & 0 & 0 & 0 \\ 1 & 0 & 1 & 0 & 0 \\ 0 & 0 & 0 & 1 & -1 \\ I & 0 & 0 & 0 & 0 \\ 0 & 0 & 0 & 1 & 0 \end{bmatrix}$$

be a special fuzzy neutrosophic square matrix. Each $T_i$ is a $5 \times 5$ square fuzzy neutrosophic matrix i = 1, 2, …, 5.

Let

$$X = X_1 \cup X_2 \cup \ldots \cup X_5$$
$$= [0\ 1\ 0\ 0\ 0] \cup [0\ 0\ 0\ I\ 0] \cup [1\ 0\ 0\ 0\ 0] \cup [0\ 0\ 0\ 0\ 1] \cup [0\ I\ 1\ 0\ 0]$$

be a special neutrosophic fuzzy row vector. To find the effect of X on T, i.e. to find

$$X \text{ o } T = (X_1 \cup X_2 \cup \ldots \cup X_5) \text{ o } (T_1 \cup T_2 \cup \ldots \cup T_5)$$
$$= X_1 \text{ o } T_1 \cup X_2 \text{ o } T_2 \cup \ldots \cup X_5 \text{ o } T_5$$

$$= \begin{bmatrix} 0 & 1 & 0 & 0 & 0 \end{bmatrix} \text{ o } \begin{bmatrix} 0 & 1 & 0 & 0 & I \\ 1 & 0 & 0 & 1 & 0 \\ 0 & 0 & 0 & I & 1 \\ I & 1 & 0 & 0 & 0 \\ 0 & 0 & 1 & 0 & 0 \end{bmatrix} \cup$$

$$\begin{bmatrix} 0 & 0 & 0 & I & 0 \end{bmatrix} \text{ o } \begin{bmatrix} 0 & 1 & 0 & 0 & 0 \\ 1 & 0 & 0 & 0 & 1 \\ 0 & I & 0 & 0 & -1 \\ 0 & 0 & I & 0 & 0 \\ 0 & 0 & 1 & 0 & 0 \end{bmatrix} \cup$$

$$\begin{bmatrix} 1 & 0 & 0 & 0 & 0 \end{bmatrix} \text{ o } \begin{bmatrix} 0 & I & 0 & 0 & 0 \\ 1 & 0 & 0 & 0 & I \\ 0 & 1 & 0 & I & 0 \\ 0 & 0 & 0 & 0 & 1 \\ 1 & I & 0 & 0 & 0 \end{bmatrix} \cup$$

$$\begin{bmatrix} 0 & 0 & 0 & 0 & 1 \end{bmatrix} \text{ o } \begin{bmatrix} 0 & 1 & 0 & 0 & 1 \\ 0 & 0 & 1 & I & 0 \\ 1 & I & 0 & 0 & 0 \\ 0 & 1 & I & 0 & 0 \\ 0 & 1 & 0 & 0 & 0 \end{bmatrix} \cup$$

$$[0 \quad I \quad 1 \quad 0 \quad 0] \text{ o } \begin{bmatrix} 0 & I & 0 & 0 & 0 \\ 1 & 0 & 1 & 0 & 0 \\ 0 & 0 & 0 & 1 & -1 \\ I & 0 & 0 & 0 & 0 \\ 0 & 0 & 0 & 1 & 0 \end{bmatrix}$$

$$\begin{aligned} &= [1\ 0\ 0\ 1\ 0] \cup [0\ 0\ I\ 0\ 0] \cup [0\ I\ 0\ 0\ 0] \cup [0\ 1\ 0\ 0\ 0] \\ &\quad \cup [I\ 0\ I\ 1\ {-1}] \\ &= Y'_1 \cup Y'_2 \cup Y'_3 \cup Y'_4 \cup Y'_5 \\ &= Y'; \end{aligned}$$

$$\begin{aligned} Y &= Y_1 \cup Y_2 \cup Y_3 \cup Y_4 \cup Y_5 \\ &= [1\ 1\ 0\ 1\ 0] \cup [0\ 0\ I\ I\ 0] \cup [1\ I\ 0\ 0\ 0] \cup [0\ 1\ 0\ 0\ 1] \\ &\quad \cup [I\ I\ 1\ 1\ 0] \end{aligned}$$

is the special neutrosophic fuzzy row vector obtained after updating Y′. Now we find

$$\begin{aligned} Y \text{ o } T &= (Y_1 \cup Y_2 \cup \ldots \cup Y_5) \text{ o } (T_1 \cup T_2 \cup \ldots \cup T_5) \\ &= Y_1 \text{ o } T_1 \cup Y_2 \text{ o } T_2 \cup \ldots \cup Y_5 \text{ o } T_5 \end{aligned}$$

$$= [1 \quad 1 \quad 0 \quad 1 \quad 0] \text{ o } \begin{bmatrix} 0 & 1 & 0 & 0 & I \\ 1 & 0 & 0 & 1 & 0 \\ 0 & 0 & 0 & I & 1 \\ I & 1 & 0 & 0 & 0 \\ 0 & 0 & 1 & 0 & 0 \end{bmatrix} \cup$$

$$[0 \quad 0 \quad I \quad I \quad 0] \text{ o } \begin{bmatrix} 0 & 1 & 0 & 0 & 0 \\ 1 & 0 & 0 & 0 & 1 \\ 0 & I & 0 & 0 & -1 \\ 0 & 0 & I & 0 & 0 \\ 0 & 0 & 1 & 0 & 0 \end{bmatrix} \cup$$

$$\begin{bmatrix}1 & I & 0 & 0 & 0\end{bmatrix} \circ \begin{bmatrix} 0 & I & 0 & 0 & 0 \\ 1 & 0 & 0 & 0 & I \\ 0 & 1 & 0 & I & 0 \\ 0 & 0 & 0 & 0 & 1 \\ 1 & I & 0 & 0 & 0 \end{bmatrix} \cup$$

$$\begin{bmatrix}0 & I & 0 & 0 & 1\end{bmatrix} \circ \begin{bmatrix} 0 & 1 & 0 & 0 & 1 \\ 0 & 0 & 1 & I & 0 \\ 1 & I & 0 & 0 & 0 \\ 0 & 1 & I & 0 & 0 \\ 0 & 1 & 0 & 0 & 0 \end{bmatrix} \cup$$

$$\begin{bmatrix}I & I & 1 & 1 & 0\end{bmatrix} \circ \begin{bmatrix} 0 & I & 0 & 0 & 0 \\ 1 & 0 & 1 & 0 & 0 \\ 0 & 0 & 0 & 1 & -1 \\ I & 0 & 0 & 0 & 0 \\ 0 & 0 & 0 & 1 & 0 \end{bmatrix}$$

$$\begin{aligned} &= [1+I\ 2\ 0\ 1\ I] \cup [0\ I\ I\ 0\ -I] \cup [I\ I\ 0\ 0\ I] \cup [0\ 1\ I\ I\ 0] \\ &\quad \cup [2I\ \ I\ \ I\ \ 1\ \ -1] \\ &= Z'_1 \cup Z'_2 \cup Z'_3 \cup \ldots \cup Z'_5 \\ &= Z'. \end{aligned}$$

$$\begin{aligned} Z &= Z_1 \cup Z_2 \cup Z_3 \cup Z_4 \cup Z_5 \\ &= [I\ 1\ 0\ 1\ I] \cup [0\ I\ I\ I\ 0] \cup [1\ I\ 0\ 0\ I] \cup [0\ 1\ 1\ I\ 1] \cup \\ &\quad [I\ I\ 1\ 1\ 0] \end{aligned}$$

is the special fuzzy neutrosophic row vector obtained by updating and thresholding Z′.

Now we can find Z o T and so on till we arrive at a fixed point or a limit cycle.

Now we illustrate the same circle operator 'o' using special fuzzy neutrosophic mixed square matrix. Let $S = S_1 \cup S_2 \cup \ldots \cup S_n$ ($n \geq 2$) be a special fuzzy neutrosophic mixed square matrix; where $S_i$ is a $t_i \times t_i$ neutrosophic fuzzy matrix i = 1, 2, …, n ($t_i \neq t_j$, $i \neq j$, $1 \leq i, j \leq n$). Let $X = X_1 \cup X_2 \cup \ldots \cup X_n$ be a special fuzzy neutrosophic mixed row vector to find

$$
\begin{aligned}
X \text{ o } S &= (X_1 \cup X_2 \cup \ldots \cup X_n) \text{ o } (S_1 \cup S_2 \cup \ldots \cup S_n) \\
&= X_1 \text{ o } S_1 \cup X_2 \text{ o } S_2 \cup \ldots \cup X_n \text{ o } S_n \\
&= Y'_1 \cup Y'_2 \cup \ldots \cup Y'_n \\
&= Y'.
\end{aligned}
$$

Now we update and threshold Y' to Y as Y' may not in general be a special neutrosophic fuzzy mixed row vector. Now let $Y = Y_1 \cup Y_2 \cup \ldots \cup Y_n$ be the special neutrosophic fuzzy mixed row vector. We find $Y \text{ o } S = Z'_1 \cup Z'_2 \cup \ldots \cup Z'_n = Z'$. We update and threshold Z' to obtain Z a special fuzzy neutrosophic mixed row vector. We find Z o S = T', if T the thresholded and updated special neutrosophic fuzzy mixed row vector of T' is a fixed point or a limit cycle we stop the 'o' operation otherwise we proceed on to find T o S and so on.

Now we illustrate this by the following example.

***Example 1.3.30:*** Let $S = S_1 \cup S_2 \cup S_3 \cup S_4 \cup S_5 =$

$$
\begin{bmatrix} 0 & 0 & 1 & 0 \\ I & 0 & 0 & 0 \\ 0 & 1 & 0 & I \\ 1 & 1 & I & 0 \end{bmatrix} \cup \begin{bmatrix} 0 & 1 & I \\ 0 & 0 & 1 \\ 1 & I & 0 \end{bmatrix} \cup \begin{bmatrix} 0 & 1 & 0 & 0 & I & 0 \\ 1 & 0 & I & 0 & 0 & 0 \\ 0 & 0 & 0 & 1 & I & 0 \\ 0 & 0 & 0 & 0 & 1 & I \\ 1 & I & 0 & 0 & 0 & 0 \\ 0 & 0 & 1 & I & 0 & 0 \end{bmatrix} \cup
$$

$$\begin{bmatrix} 0 & 1 & 0 & 0 & 0 \\ I & 0 & 1 & 0 & 0 \\ 0 & I & 0 & 1 & 0 \\ 0 & 0 & I & 0 & 1 \\ 1 & 1 & 0 & 0 & 0 \end{bmatrix} \cup \begin{bmatrix} 0 & 0 & 0 & 1 & 1 & 0 \\ I & 0 & 0 & 0 & 1 & 0 \\ 0 & 0 & 0 & I & 0 & 0 \\ 1 & 0 & 0 & 0 & I & 1 \\ 0 & 0 & 0 & 1 & 0 & I \\ 0 & 1 & I & 0 & 0 & 0 \end{bmatrix}$$

be the given special fuzzy neutrosophic mixed square matrix. Suppose

$$\begin{aligned} X &= X_1 \cup X_2 \cup \ldots \cup X_5 \\ &= [1\ 0\ 0\ 0] \cup [0\ I\ 0] \cup [1\ 0\ 0\ 0\ 0\ I] \cup [0\ 1\ 0\ 0\ 0] \cup [1\ I\ 0\ 0\ 0\ 0] \end{aligned}$$

be the special fuzzy neutrosophic mixed row vector. Now we find

$$\begin{aligned} X \text{ o } S &= [X_1 \cup X_2 \cup \ldots \cup X_5] \text{ o } [S_1 \cup S_2 \cup \ldots \cup S_5] \\ &= X_1 \text{ o } S_1 \cup X_2 \text{ o } S_2 \cup \ldots \cup X_5 \text{ o } S_5 \end{aligned}$$

$$= \begin{bmatrix} 1 & 0 & 0 & 0 \end{bmatrix} \text{ o } \begin{bmatrix} 0 & 0 & 1 & 0 \\ I & 0 & 0 & 0 \\ 0 & 1 & 0 & I \\ 1 & 1 & I & 0 \end{bmatrix} \cup$$

$$\begin{bmatrix} 0 & I & 0 \end{bmatrix} \text{ o } \begin{bmatrix} 0 & 1 & I \\ 0 & 0 & 1 \\ 1 & I & 0 \end{bmatrix} \cup$$

$$\begin{bmatrix} 1 & 0 & 0 & 0 & 0 & I \end{bmatrix} \text{ o } \begin{bmatrix} 0 & 1 & 0 & 0 & I & 0 \\ 1 & 0 & I & 0 & 0 & 0 \\ 0 & 0 & 0 & 1 & I & 0 \\ 0 & 0 & 0 & 0 & 1 & I \\ 1 & I & 0 & 0 & 0 & 0 \\ 0 & 0 & 1 & I & 0 & 0 \end{bmatrix} \cup$$

$$\begin{bmatrix}0 & 1 & 0 & 0 & 0\end{bmatrix} \text{ o } \begin{bmatrix}0 & 1 & 0 & 0 & 0\\ I & 0 & 1 & 0 & 0\\ 0 & I & 0 & 1 & 0\\ 0 & 0 & I & 0 & 1\\ 1 & 1 & 0 & 0 & 0\end{bmatrix} \cup$$

$$\begin{bmatrix}1 & I & 0 & 0 & 0 & 0\end{bmatrix} \text{ o } \begin{bmatrix}0 & 0 & 0 & 1 & 1 & 0\\ I & 0 & 0 & 0 & 1 & 0\\ 0 & 0 & 0 & I & 0 & 0\\ 1 & 0 & 0 & 0 & I & 1\\ 0 & 0 & 0 & 1 & 0 & I\\ 0 & 1 & I & 0 & 0 & 0\end{bmatrix}$$

$= [0\ 0\ 1\ 0] \cup [0\ 0\ I] \cup [0\ 1\ I\ I\ I\ 0] \cup [I\ 0\ 1\ 0\ 0] \cup [I\ 0\ 0\ 1\ I{+}1\ 0]$
$= Y'_1 \cup Y'_2 \cup Y'_3 \cup Y'_4 \cup Y'_5$
$= Y'.$

After updating Y' to

$Y = Y_1 \cup Y_2 \cup Y_3 \cup Y_4 \cup Y_5$
$= [1\ 0\ 1\ 0] \cup [0\ I\ I] \cup [1\ 1\ I\ I\ I\ I] \cup [I\ 1\ 1\ 0\ 0] \cup [1\ I\ 0\ 1\ I\ 0]$

which is again special neutrosophic fuzzy mixed row vector.

$Y \text{ o } S = (Y_1 \cup Y_2 \cup \ldots \cup Y_5) \text{ o } (S_1 \cup S_2 \cup \ldots \cup S_5)$
$= Y_1 \text{ o } S_1 \cup Y_2 \text{ o } S_2 \cup \ldots \cup Y_5 \text{ o } S_5$

$$= \begin{bmatrix}1 & 0 & 1 & 0\end{bmatrix} \text{ o } \begin{bmatrix}0 & 0 & 1 & 0\\ I & 0 & 0 & 0\\ 0 & 1 & 0 & I\\ 1 & 1 & I & 0\end{bmatrix} \cup$$

$$[0 \quad I \quad I] \circ \begin{bmatrix} 0 & 1 & I \\ 0 & 0 & 1 \\ 1 & I & 0 \end{bmatrix} \cup$$

$$[1 \quad 1 \quad I \quad I \quad I \quad I] \circ \begin{bmatrix} 0 & 1 & 0 & 0 & I & 0 \\ 1 & 0 & I & 0 & 0 & 0 \\ 0 & 0 & 0 & 1 & I & 0 \\ 0 & 0 & 0 & 0 & 1 & I \\ 1 & I & 0 & 0 & 0 & 0 \\ 0 & 0 & 1 & I & 0 & 0 \end{bmatrix} \cup$$

$$[I \quad 1 \quad 1 \quad 0 \quad 0] \circ \begin{bmatrix} 0 & 1 & 0 & 0 & 0 \\ I & 0 & 1 & 0 & 0 \\ 0 & I & 0 & 1 & 0 \\ 0 & 0 & I & 0 & 1 \\ 1 & 1 & 0 & 0 & 0 \end{bmatrix} \cup$$

$$[1 \quad I \quad 0 \quad 1 \quad I \quad 0] \circ \begin{bmatrix} 0 & 0 & 0 & 1 & 1 & 0 \\ I & 0 & 0 & 0 & 1 & 0 \\ 0 & 0 & 0 & I & 0 & 0 \\ 1 & 0 & 0 & 0 & I & 1 \\ 0 & 0 & 0 & 1 & 0 & I \\ 0 & 1 & 0 & 0 & 0 & 0 \end{bmatrix}$$

= $[0\ 1\ 1\ I] \cup [I\ I\ I] \cup [1{+}I\ \ 1{+}I\ \ 2I\ \ 2I\ \ 3I\ \ I] \cup [I\ \ 2I\ \ 1\ \ 1\ \ 0] \cup [I{+}1\ \ 0\ \ 0\ \ 1\ \ 1{+}2I\ \ 1{+}I]$

= $Z'_1 \cup Z'_2 \cup Z'_3 \cup Z'_4 \cup Z'_5$

= $Z'$.

Now we update and threshold Z' to get

$$Z = Z_1 \cup Z_2 \cup Z_3 \cup Z_4 \cup Z_5$$
$$= [1\ 1\ I\ I] \cup [I\ I\ I] \cup [1\ I\ I\ I\ I\ I] \cup [I\ 1\ 1\ 1\ 0] \cup [1\ I\ 0\ 1\ I\ I];$$

Z is the special neutrosophic fuzzy mixed row vector. Now

$$Z \text{ o } S = (Z_1 \cup Z_2 \cup \ldots \cup Z_5) \text{ o } (S_1 \cup S_2 \cup \ldots \cup S_5)$$
$$= Z_1 \text{ o } S_1 \cup Z_2 \text{ o } S_2 \cup \ldots \cup Z_5 \text{ o } S_5$$

$$= \begin{bmatrix} 1 & 1 & I & I \end{bmatrix} \text{ o } \begin{bmatrix} 0 & 0 & 1 & 0 \\ I & 0 & 0 & 0 \\ 0 & 1 & 0 & I \\ 1 & 1 & I & 0 \end{bmatrix} \cup$$

$$\begin{bmatrix} I & I & I \end{bmatrix} \text{ o } \begin{bmatrix} 0 & 1 & I \\ 0 & 0 & 1 \\ 1 & I & 0 \end{bmatrix} \cup$$

$$\begin{bmatrix} 1 & I & I & I & I & I \end{bmatrix} \text{ o } \begin{bmatrix} 0 & 1 & 0 & 0 & I & 0 \\ 1 & 0 & I & 0 & 0 & 0 \\ 0 & 0 & 0 & 1 & I & 0 \\ 0 & 0 & 0 & 0 & 1 & I \\ 1 & I & 0 & 0 & 0 & 0 \\ 0 & 0 & 1 & I & 0 & 0 \end{bmatrix} \cup$$

$$\begin{bmatrix} I & 1 & 1 & 1 & 0 \end{bmatrix} \text{ o } \begin{bmatrix} 0 & 1 & 0 & 0 & 0 \\ I & 0 & 1 & 0 & 0 \\ 0 & I & 0 & 1 & 0 \\ 0 & 0 & I & 0 & 1 \\ 1 & 1 & 0 & 0 & 0 \end{bmatrix} \cup$$

$$[1 \quad I \quad 0 \quad 1 \quad I \quad I] \text{ o } \begin{bmatrix} 0 & 0 & 0 & 1 & 1 & 0 \\ I & 0 & 0 & 0 & 1 & 0 \\ 0 & 0 & 0 & I & 0 & 0 \\ 1 & 0 & 0 & 0 & I & 1 \\ 0 & 0 & 0 & 1 & 0 & I \\ 0 & 1 & I & 0 & 0 & 0 \end{bmatrix}$$

= [2I 2I 1+I I] ∪ [I 2I 2I] ∪ [2I 1+I 2I 2I, 3I I] ∪ [I 2I 1+I 1 1] ∪ [I+1 I I 1+2I 1+2I 1+I]

= $P'_1 \cup P'_2 \cup P'_3 \cup P'_4 \cup P'_5$
= P'.

Now P' is not a special fuzzy neutrosophic mixed row vector we now update and threshold P' to P where

$P = P_1 \cup P_2 \cup P_3 \cup P_4 \cup P_5$
= [1 I I I] ∪ [I I I] ∪ [1 I I I I I] ∪ [I 1 I 1 1] ∪ [1 I I I I I]

where P is a special fuzzy neutrosophic mixed row vector.

We can find P o S and continue till we arrive at a fixed point or a limit cycle.

Now we proceed on to show how the operations are performed in case of special fuzzy neutrosophic rectangular matrix.

Let

$T = T_1 \cup T_2 \cup \ldots \cup T_m$,

(m ≥ 2) be a special fuzzy neutrosophic rectangular matrix where $T_i$ is a s × t (s ≠ t) neutrosophic fuzzy rectangular matrix for i = 1, 2, …, m.

Suppose

$X = X_1 \cup X_2 \cup \ldots \cup X_m$

be a special fuzzy neutrosophic row vector where each $X_i$ is a 1 × s neutrosophic fuzzy row vector i = 1, 2, …, m.

Now we find

$$\begin{aligned} X \text{ o } T &= (X_1 \cup X_2 \cup \ldots \cup X_m) \text{ o } (T_1 \cup T_2 \cup \ldots \cup T_m) \\ &= X_1 \text{ o } T_1 \cup X_2 \text{ o } T_2 \cup \ldots \cup X_m \text{ o } T_m \\ &= Y'_1 \cup Y'_2 \cup \ldots \cup Y'_m \\ &= Y', \end{aligned}$$

Y' may not be even a special fuzzy neutrosophic row vector we update and threshold Y' to

$$Y = Y_1 \cup Y_2 \cup \ldots \cup Y_m$$

where each $Y_i$ is a 1 × t fuzzy neutrosophic row vector for i = 1, 2, …, m.

$$\begin{aligned} Y \text{ o } T^T &= (Y_1 \cup Y_2 \cup \ldots \cup Y_m) \text{ o } (T_1^t \cup T_2^t \cup \ldots \cup T_m^t) \\ &= Y_1 \text{ o } T_1^t \cup Y_2 \text{ o } T_2^t \cup \ldots \cup Y_m \text{ o } T_m^t \\ &= Z'_1 \cup Z'_2 \cup \ldots \cup Z'_m \\ &= Z'. \end{aligned}$$

Z' is updated and thresholded to get $Z = Z_1 \cup Z_2 \cup \ldots Z_m$. Now using Z we find Z o T and so on until we arrive at a special limit cycle or a special fixed point.

In case of special fuzzy neutrosophic rectangular matrices we get the resultant as a pair of fuzzy neutrosophic row vector be it a fixed point or a limit cycle which we shall call as special binary pair.

Now we illustrate this by the following example.

***Example 1.3.31:*** Let $T = T_1 \cup T_2 \cup \ldots \cup T_6$ be a special neutrosophic fuzzy rectangular matrix where each $T_i$ is a 7 × 3 special neutrosophic fuzzy rectangular matrix; i = 1, 2, …, 6.

Suppose

$$X = X_1 \cup X_2 \cup \ldots \cup X_6$$
$$= [0\ 1\ 0\ 0\ 0\ 0\ 1] \cup [0\ 0\ 0\ 0\ 1\ 0\ 0] \cup [0\ 0\ 0\ 0\ 0\ 0\ I] \cup [1\ 0\ 0\ 0\ 0\ 0\ 0] \cup [0\ 0\ 0\ 1\ I\ 0\ 0] \cup [0\ 0\ 0\ 0\ 1\ 0\ 0].$$

be a 1 × 7 fuzzy neutrosophic row vector. Let

$$T = \begin{bmatrix} 0 & 1 & 0 \\ 1 & 0 & 0 \\ 0 & 1 & 1 \\ 1 & 1 & 0 \\ 0 & I & 0 \\ 0 & 0 & 1 \\ 0 & 1 & 1 \end{bmatrix} \cup \begin{bmatrix} 0 & 0 & 1 \\ I & 0 & 0 \\ 0 & 1 & 0 \\ 1 & 0 & 0 \\ 1 & 1 & 0 \\ 1 & 0 & 1 \\ 0 & 1 & 1 \end{bmatrix} \cup$$

$$\begin{bmatrix} I & 0 & 0 \\ 0 & 1 & 0 \\ 0 & 0 & 1 \\ 1 & 0 & 0 \\ 1 & 0 & 1 \\ 0 & 1 & 1 \\ 1 & 1 & 0 \end{bmatrix} \cup \begin{bmatrix} 1 & 0 & 0 \\ 0 & 1 & 0 \\ 0 & 0 & 1 \\ 1 & 1 & 0 \\ 0 & 0 & I \\ 1 & 0 & 1 \\ 0 & 1 & 1 \end{bmatrix} \cup$$

$$\begin{bmatrix} 1 & 1 & 0 \\ 0 & 1 & 1 \\ 0 & 0 & 1 \\ 0 & 1 & 0 \\ 1 & 0 & 0 \\ 0 & 0 & 1 \\ 0 & 0 & I \end{bmatrix} \cup \begin{bmatrix} 1 & 0 & 1 \\ 0 & I & 0 \\ 0 & 1 & 1 \\ 1 & 1 & 0 \\ 1 & 0 & 0 \\ 0 & 0 & 1 \\ 1 & 0 & 0 \end{bmatrix}.$$

be the given special neutrosophic fuzzy rectangular matrix.

We find

$$X \text{ o } T = (X_1 \cup X_2 \cup \ldots \cup X_6) \text{ o } (T_1 \cup T_2 \cup \ldots \cup T_6)$$
$$= X_1 \text{ o } T_1 \cup X_2 \text{ o } T_2 \cup \ldots \cup X_6 \text{ o } T_6$$

$$= \begin{bmatrix} 0 & 1 & 0 & 0 & 0 & 0 & 1 \end{bmatrix} \text{ o } \begin{bmatrix} 0 & 1 & 0 \\ 1 & 0 & 0 \\ 0 & 1 & 1 \\ 1 & 1 & 0 \\ 0 & I & 0 \\ 0 & 0 & 1 \\ 0 & 1 & 1 \end{bmatrix} \cup$$

$$\begin{bmatrix} 0 & 0 & 0 & 0 & 1 & 0 & 0 \end{bmatrix} \text{ o } \begin{bmatrix} 0 & 0 & 1 \\ I & 0 & 0 \\ 0 & 1 & 0 \\ 1 & 0 & 0 \\ 1 & 1 & 0 \\ 1 & 0 & 1 \\ 0 & 1 & 1 \end{bmatrix} \cup$$

$$\begin{bmatrix} 0 & 0 & 0 & 0 & 0 & 0 & I \end{bmatrix} \text{ o } \begin{bmatrix} I & 0 & 0 \\ 0 & 1 & 0 \\ 0 & 0 & 1 \\ 1 & 0 & 0 \\ 1 & 0 & 1 \\ 0 & 1 & 1 \\ 1 & 1 & 0 \end{bmatrix} \cup$$

$$[1\ \ 0\ \ 0\ \ 0\ \ 0\ \ 0\ \ 0] \text{ o } \begin{bmatrix} 1 & 0 & 0 \\ 0 & 1 & 0 \\ 0 & 0 & 1 \\ 1 & 1 & 0 \\ 0 & 0 & I \\ 1 & 0 & 1 \\ 0 & 1 & 1 \end{bmatrix} \cup$$

$$[0\ \ 0\ \ 0\ \ 1\ \ I\ \ 0\ \ 0] \text{ o } \begin{bmatrix} 1 & 1 & 0 \\ 0 & 1 & 1 \\ 0 & 0 & 1 \\ 0 & 1 & 0 \\ 1 & 0 & 0 \\ 0 & 0 & 1 \\ 0 & 0 & I \end{bmatrix} \cup$$

$$[0\ \ 0\ \ 0\ \ 0\ \ 1\ \ 0\ \ 0] \text{ o } \begin{bmatrix} 1 & 0 & 1 \\ 0 & I & 0 \\ 0 & 1 & 1 \\ 1 & 1 & 0 \\ 1 & 0 & 0 \\ 0 & 0 & 1 \\ 1 & 0 & 0 \end{bmatrix}$$

$= \quad [1\ 1\ 1] \cup [1\ 1\ 0] \cup [I\ I\ 0] \cup [1\ 0\ 0] \cup [I\ 1\ 0] \cup [1\ 0\ 0]$
$= \quad Y_1 \cup Y_2 \cup \ldots \cup Y_6$
$= \quad Y.$

Now calculate

$YoT^t \quad = \quad (Y_1 \cup Y_2 \cup \ldots \cup Y_6) \text{ o } (T_1^t \cup T_2^t \cup \ldots \cup T_6^t)$

$\qquad\quad = \quad Y_1 \text{ o } T_1^t \cup Y_2 \text{ o } T_2^t \cup \ldots \cup Y_6 \text{ o } T_6^t$

$$= \begin{bmatrix}1 & 1 & 1\end{bmatrix} \text{ o } \begin{bmatrix}0 & 1 & 0 & 1 & 0 & 0 & 0\\ 1 & 0 & 1 & 1 & I & 0 & 1\\ 0 & 0 & 1 & 0 & 0 & 1 & 1\end{bmatrix} \cup$$

$$\begin{bmatrix}1 & 1 & 0\end{bmatrix} \text{ o } \begin{bmatrix}0 & I & 0 & 1 & 1 & 1 & 0\\ 0 & 0 & 1 & 0 & 1 & 0 & 1\\ 1 & 0 & 0 & 0 & 0 & 1 & 1\end{bmatrix} \cup$$

$$\begin{bmatrix}I & I & 0\end{bmatrix} \text{ o } \begin{bmatrix}I & 0 & 0 & 1 & 1 & 0 & 1\\ 0 & 1 & 0 & 0 & 0 & 1 & 1\\ 0 & 0 & 1 & 0 & 1 & 1 & 0\end{bmatrix} \cup$$

$$\begin{bmatrix}1 & 0 & 0\end{bmatrix} \text{ o } \begin{bmatrix}1 & 0 & 0 & 1 & 0 & 1 & 0\\ 0 & 1 & 0 & 1 & 0 & 0 & 1\\ 0 & 0 & 1 & 0 & I & 1 & 1\end{bmatrix} \cup$$

$$\begin{bmatrix}I & 1 & 0\end{bmatrix} \text{ o } \begin{bmatrix}1 & 0 & 0 & 0 & 1 & 0 & 0\\ 1 & 1 & 0 & 1 & 0 & 0 & 0\\ 0 & 1 & 1 & 0 & 0 & 1 & I\end{bmatrix} \cup$$

$$\begin{bmatrix}1 & 0 & 0\end{bmatrix} \text{ o } \begin{bmatrix}1 & 0 & 0 & 1 & 1 & 0 & 1\\ 0 & I & 1 & 1 & 0 & 0 & 0\\ 1 & 0 & 1 & 0 & 0 & 1 & 0\end{bmatrix}$$

= [1 1 2 2 I 1 2] ∪ [0 I 1 1 2 1 1] ∪ [I I 0 I I I 2I] ∪ [1 0 0 1 0 1 0] ∪ [1+I 1 0 1 I 0 0] ∪ [1 0 0 1 1 0 1]

= $S'_1 \cup S'_2 \cup \ldots \cup S'_6$

= S',

we update and threshold S' to obtain S as follows:

S = [1 1 1 1 I 1 1] ∪ [0 I 1 1 1 1 1] ∪ [I I 0 I I I I] ∪ [1 0 0 1 0 1 0] ∪ [I 1 0 1 I 0 0] ∪ [1 0 0 1 1 0 1].

Now we can find

$$S \text{ o } T = (S_1 \cup S_2 \cup \ldots \cup S_6) \text{ o } (T_1 \cup T_2 \cup \ldots \cup T_6)$$
$$= S_1 \text{ o } T_1 \cup \ldots \cup S_6 \text{ o } T_6$$

and so on. Thus we proceed on till we obtain a special pair of fixed points or a special pair of limit cycle or a fixed point and a limit cycle.

Now we illustrate how the operation works on special fuzzy neutrosophic mixed rectangular matrix. Let $W = W_1 \cup W_2 \cup \ldots \cup W_n$ $(n \geq 2)$ be the given special fuzzy neutrosophic mixed rectangular matrix where $W_i$'s are $s_i \times t_i$ $(s_i \neq t_i)$ neutrosophic fuzzy rectangular matrix and $W_j$ $(i \neq j)$ are $p_j \times q_j$ $(p_j \neq q_j$ and $p_j \neq s_i$ or $q_j \neq t_i)$ rectangular neutrosophic fuzzy matrix, $1 \leq i, j \leq n$.

Let $X = X_1 \cup X_2 \cup \ldots \cup X_n$ where $X_i$ are $1 \times s_i$ fuzzy neutrosophic row vector and $X_j$ are $1 \times p_j$ fuzzy neutrosophic row vectors be the special fuzzy neutrosophic mixed row vector $1 \leq i, j \leq n$.

Now we calculate

$$X \text{ o } W = (X_1 \cup X_2 \cup \ldots \cup X_n) \text{ o } (W_1 \cup W_2 \cup \ldots \cup W_n)$$
$$= X_1 \text{ o } W_1 \cup X_2 \text{ o } W_2 \cup \ldots \cup X_n \text{ o } W_n$$
$$= Y'_1 \cup Y'_2 \cup \ldots \cup Y'_n$$
$$= Y'.$$

Y' may or may not be a special fuzzy neutrosophic mixed row vector so we threshold Y' to Y and obtain the special fuzzy neutrosophic mixed row vector i.e. $Y = Y_1 \cup Y_2 \cup \ldots \cup Y_n$ is the special fuzzy neutrosophic mixed row vector. We now find

$$Y \text{ o } W^T = (Y_1 \cup Y_2 \cup \ldots \cup Y_n) \text{ o } \left(W_1^T \cup W_2^T \cup \ldots \cup W_n^T\right)$$
$$= Y_1 \text{ o } W_1^T \cup Y_2 \text{ o } W_2^T \cup \ldots \cup Y_n \text{ o } W_n^T$$
$$= P'_1 \cup P'_2 \cup \ldots \cup P'_n$$
$$= P';$$

P' may or may not be a special fuzzy neutrosophic mixed row vector, we update and threshold P' to P and obtain $P = P_1 \cup P_2 \cup \ldots \cup P_n$ to be the special fuzzy neutrosophic mixed row vector. Now using P we find

$$\begin{aligned} P \text{ o } W &= (P_1 \cup P_2 \cup \ldots \cup P_n) \text{ o } (W_1 \cup W_2 \cup \ldots \cup W_n) \\ &= P_1 \text{ o } W_1 \cup P_2 \text{ o } W_2 \cup \ldots \cup P_n \text{ o } W_n \\ &= R'_1 \cup R'_2 \cup \ldots \cup R'_n \\ &= R'. \end{aligned}$$

R' may or may not be a special fuzzy neutrosophic mixed row vector. We threshold R' to R and obtain the special fuzzy neutrosophic mixed row vector. Now we can find R o $W^T$ and so on until we arrive at a special pair of fixed point or a special pair of limit cycle or a limit cycle and a fixed point.

We illustrate this situation by the following example.

***Example 1.3.32:*** Let $W = W_1 \cup W_2 \cup W_3 \cup W_4 \cup W_5 =$

$$\begin{bmatrix} 0 & 1 & 0 \\ I & 0 & 1 \\ 0 & 0 & 1 \\ 1 & 1 & 0 \\ 0 & I & 0 \\ 1 & 1 & 1 \\ 1 & 0 & 1 \end{bmatrix} \cup \begin{bmatrix} 1 & 0 & 0 & 0 & I & 0 & 0 & 0 \\ 0 & 1 & 0 & 0 & 0 & 1 & 0 & 0 \\ 0 & 0 & 1 & 0 & 0 & 0 & 1 & 0 \\ 1 & 0 & 0 & I & 1 & 0 & 1 & 1 \end{bmatrix} \cup$$

$$\begin{bmatrix} 1 & 0 & 0 & 0 & -1 & 1 & 0 \\ 0 & 0 & 1 & 0 & 0 & 0 & 0 \\ 0 & 0 & 0 & I & 0 & 0 & 0 \\ 0 & 0 & 0 & 0 & 1 & 0 & 1 \\ 1 & 0 & 0 & 1 & 0 & 0 & 1 \\ 0 & 0 & 1 & 0 & 0 & 1 & 0 \end{bmatrix} \cup$$

$$\begin{bmatrix} 0 & I & 0 & 1 & 0 & 0 & 1 & 0 & 1 \\ 1 & 0 & 0 & 0 & I & 0 & 0 & 0 & 1 \\ 0 & 0 & 0 & 1 & 1 & 1 & 0 & 0 & 0 \\ 1 & 0 & 1 & 0 & 0 & 0 & 0 & 1 & 0 \\ 0 & 0 & 1 & 1 & 0 & 0 & 0 & 1 & 0 \end{bmatrix}$$

$$\cup \begin{bmatrix} 0 & 1 & 0 & 0 & 1 & 0 \\ 0 & 0 & I & 0 & 0 & 1 \\ I & 0 & 0 & 0 & 1 & 1 \\ 0 & 0 & 0 & 1 & 0 & 0 \\ 0 & 1 & 1 & 1 & 0 & 0 \end{bmatrix}$$

be a special fuzzy neutrosophic mixed rectangular matrix.
Let

$$\begin{aligned} X &= X_1 \cup X_2 \cup \ldots \cup X_5 \\ &= [0\ 0\ 0\ 0\ 0\ 0\ 1] \cup [I\ 0\ 0\ 0] \cup [0\ 1\ 0\ 0\ 1\ 0] \cup [1\ 0\ 0\ 0\ 1] \cup [0\ 0\ 1\ 0\ 0] \end{aligned}$$

be a special fuzzy neutrosophic mixed row vector.
To find

$$\begin{aligned} X \text{ o } W &= [X_1 \cup X_2 \cup \ldots \cup X_5] \text{ o } [W_1 \cup W_2 \cup \ldots \cup W_5] \\ &= X_1 \text{ o } W_1 \cup X_2 \text{ o } W_2 \cup \ldots \cup X_5 \text{ o } W_5 \end{aligned}$$

$$= \begin{bmatrix} 0 & 0 & 0 & 0 & 0 & 0 & 1 \end{bmatrix} \text{ o } \begin{bmatrix} 0 & 1 & 0 \\ I & 0 & 1 \\ 0 & 0 & 1 \\ 1 & 1 & 0 \\ 0 & I & 0 \\ 1 & 1 & 1 \\ 1 & 0 & 1 \end{bmatrix} \cup$$

$$[I \quad 0 \quad 0 \quad 0] \circ \begin{bmatrix} 1 & 0 & 0 & 0 & I & 0 & 0 & 0 \\ 0 & 1 & 0 & 0 & 0 & 1 & 0 & 0 \\ 0 & 0 & 1 & 0 & 0 & 0 & 1 & 0 \\ 1 & 0 & 0 & I & 1 & 0 & 1 & 1 \end{bmatrix} \cup$$

$$[0 \quad 1 \quad 0 \quad 0 \quad 1 \quad 0] \circ \begin{bmatrix} 1 & 0 & 0 & 0 & -1 & 1 & 0 \\ 0 & 0 & 1 & 0 & 0 & 0 & 0 \\ 0 & 0 & 0 & I & 0 & 0 & 0 \\ 0 & 0 & 0 & 0 & 1 & 0 & 1 \\ 1 & 0 & 0 & 1 & 0 & 0 & 1 \\ 0 & 0 & 1 & 0 & 0 & 1 & 0 \end{bmatrix} \cup$$

$$[1 \quad 0 \quad 0 \quad 0 \quad 1] \circ \begin{bmatrix} 0 & I & 0 & 1 & 0 & 0 & 1 & 0 & 1 \\ 1 & 0 & 0 & 0 & I & 0 & 0 & 0 & 1 \\ 0 & 0 & 0 & 1 & 1 & 1 & 0 & 0 & 0 \\ 1 & 0 & 1 & 0 & 0 & 0 & 0 & 1 & 0 \\ 0 & 0 & 1 & 1 & 0 & 0 & 0 & 1 & 0 \end{bmatrix} \cup$$

$$[0 \quad 0 \quad 1 \quad 0 \quad 0] \circ \begin{bmatrix} 0 & 1 & 0 & 0 & 1 & 0 \\ 0 & 0 & I & 0 & 0 & 1 \\ I & 0 & 0 & 0 & 1 & 1 \\ 0 & 0 & 0 & 1 & 0 & 0 \\ 0 & 1 & 1 & 1 & 0 & 0 \end{bmatrix}$$

= [1 0 1] ∪ [I 0 0 0 I 0 0 0] ∪ [1 0 1 1 0 0 1] ∪ [0 I 1 2 0 0 1 1 1] ∪ [I 0 0 0 1 1]

= $Y'_1 \cup Y'_2 \cup \ldots \cup Y'_5$

= Y'.

Now by thresholding Y' we get

$$Y = Y_1 \cup Y_2 \cup \ldots \cup Y_5$$
$$= [1\ 0\ 1] \cup [I\ 0\ 0\ 0\ I\ 0\ 0\ 0] \cup [1\ 0\ 1\ 1\ 0\ 0\ 1] \cup [0\ I\ 1\ 1\ 0\ 0\ 1\ 1\ 1] \cup [I\ 0\ 0\ 0\ 1\ 1].$$

To find

$$Y \text{ o } W^T = [Y_1 \cup \ldots \cup Y_5] \text{ o } [W_1^T \cup W_2^T \cup \ldots \cup W_5^T]$$
$$= Y_1 \text{ o } W_1^T \cup Y_2 \text{ o } W_2^T \cup \ldots \cup Y_5 \text{ o } W_5^T$$

$$= \begin{bmatrix} 1 & 0 & 1 \end{bmatrix} \text{ o } \begin{bmatrix} 0 & I & 0 & 1 & 0 & 1 & 1 \\ 1 & 0 & 0 & 1 & I & 1 & 0 \\ 0 & 1 & 1 & 0 & 0 & 1 & 1 \end{bmatrix} \cup$$

$$\begin{bmatrix} I & 0 & 0 & 0 & I & 0 & 0 & 0 \end{bmatrix} \text{ o } \begin{bmatrix} 1 & 0 & 0 & 1 \\ 0 & 1 & 0 & 0 \\ 0 & 0 & 1 & 0 \\ 0 & 0 & 0 & I \\ I & 0 & 0 & 1 \\ 0 & 1 & 0 & 0 \\ 0 & 0 & 1 & 1 \\ 0 & 0 & 0 & 1 \end{bmatrix} \cup$$

$$\begin{bmatrix} 1 & 0 & 1 & 1 & 0 & 0 & 1 \end{bmatrix} \text{ o } \begin{bmatrix} 1 & 0 & 0 & 0 & 1 & 0 \\ 0 & 0 & 0 & 0 & 0 & 0 \\ 0 & 1 & 0 & 0 & 0 & 1 \\ 0 & 0 & I & 0 & 1 & 0 \\ -1 & 0 & 0 & 1 & 0 & 0 \\ 1 & 0 & 0 & 0 & 0 & 1 \\ 0 & 0 & 0 & 1 & 1 & 0 \end{bmatrix} \cup$$

$$[0 \quad I \quad 1 \quad 1 \quad 0 \quad 0 \quad 1 \quad 1 \quad 1] \text{ o } \begin{bmatrix} 0 & 1 & 0 & 1 & 0 \\ I & 0 & 0 & 0 & 0 \\ 0 & 0 & 0 & 1 & 1 \\ 1 & 0 & 1 & 0 & 1 \\ 0 & I & 1 & 0 & 0 \\ 0 & 0 & 1 & 0 & 0 \\ 1 & 0 & 0 & 0 & 0 \\ 0 & 0 & 0 & 1 & 1 \\ 1 & 1 & 0 & 0 & 0 \end{bmatrix} \cup$$

$$[I \quad 0 \quad 0 \quad 0 \quad 1 \quad 1] \text{ o } \begin{bmatrix} 0 & 0 & I & 0 & 0 \\ 1 & 0 & 0 & 0 & 1 \\ 0 & I & 0 & 0 & 1 \\ 0 & 0 & 0 & 1 & 1 \\ 1 & 0 & 1 & 0 & 0 \\ 0 & 1 & 1 & 0 & 0 \end{bmatrix}$$

= [0 1+I 1 1 0 2 2] $\cup$ [2I 0 0 2I] $\cup$ [1 1 I 1 3 1] $\cup$ [3+I 1 1 2 3] $\cup$ [1 1 2+I 0 0]
= $R'_1 \cup R'_2 \cup \ldots \cup R'_5$
= R';

we now update and threshold this R′ to obtain R as R′ is not a special fuzzy neutrosophic mixed row vector.

R = [0 I 1 1 0 1 1] $\cup$ [I 0 0 I] $\cup$ [1 1 I 1 1 1] $\cup$ [1 1 1 1 1] $\cup$ [1 1 1 0 0].

We now proceed on to find

$$\begin{aligned} R \text{ o } W &= (R_1 \cup R_2 \cup \ldots \cup R_5) \text{ o } (W_1 \cup W_2 \cup \ldots \cup W_5) \\ &= R_1 \text{ o } W_1 \cup R_2 \text{ o } W_2 \cup \ldots \cup R_5 \text{ o } W_5 \end{aligned}$$

$$= \begin{bmatrix} 0 & I & 1 & 1 & 0 & 1 & 1 \end{bmatrix} \text{ o } \begin{bmatrix} 0 & 1 & 0 \\ I & 0 & 1 \\ 0 & 0 & 1 \\ 1 & 1 & 0 \\ 0 & I & 0 \\ 1 & 1 & 1 \\ 1 & 0 & 1 \end{bmatrix} \cup$$

$$\begin{bmatrix} I & 0 & 0 & I \end{bmatrix} \text{ o } \begin{bmatrix} 1 & 0 & 0 & 0 & I & 0 & 0 & 0 \\ 0 & 1 & 0 & 0 & 0 & 1 & 0 & 0 \\ 0 & 0 & 1 & 0 & 0 & 0 & 1 & 0 \\ 1 & 0 & 0 & I & 1 & 0 & 1 & 1 \end{bmatrix} \cup$$

$$\begin{bmatrix} 1 & 1 & I & 1 & 1 & 1 \end{bmatrix} \text{ o } \begin{bmatrix} 1 & 0 & 0 & 0 & -1 & 1 & 0 \\ 0 & 0 & 1 & 0 & 0 & 0 & 0 \\ 0 & 0 & 0 & I & 0 & 0 & 0 \\ 0 & 0 & 0 & 0 & 1 & 0 & 1 \\ 1 & 0 & 0 & 1 & 0 & 0 & 1 \\ 0 & 0 & 1 & 0 & 0 & 1 & 0 \end{bmatrix} \cup$$

$$\begin{bmatrix} 1 & 1 & 1 & 1 & 1 \end{bmatrix} \text{ o } \begin{bmatrix} 0 & I & 0 & 1 & 0 & 0 & 1 & 0 & 1 \\ 1 & 0 & 0 & 0 & I & 0 & 0 & 0 & 1 \\ 0 & 0 & 0 & 1 & 1 & 1 & 0 & 0 & 0 \\ 1 & 0 & 1 & 0 & 0 & 0 & 0 & 1 & 0 \\ 0 & 0 & 1 & 1 & 0 & 0 & 0 & 1 & 0 \end{bmatrix} \cup$$

$$\begin{bmatrix} 1 & 1 & 1 & 0 & 0 \end{bmatrix} \text{ o } \begin{bmatrix} 0 & 1 & 0 & 0 & 1 & 0 \\ 0 & 0 & I & 0 & 0 & 1 \\ I & 0 & 0 & 0 & 1 & 1 \\ 0 & 0 & 0 & 1 & 0 & 0 \\ 0 & 1 & 1 & 1 & 0 & 0 \end{bmatrix}$$

$= [3+I\ 2\ 3+I] \cup [2I\ 0\ 0\ I\ 2I\ 0\ I\ I] \cup [2\ 0\ 2\ 1+I\ 0\ 2\ 2]$
$\cup [2\ I\ 2\ 3\ 1+I\ 1\ 1\ 2\ 2] \cup [I\ 1\ I\ 0\ 2\ 2]$
$= S'_1 \cup S'_2 \cup \ldots \cup S'_5$
$= S'.$

We threshold S′ and obtain

$= [1\ 1\ 1] \cup [I\ 0\ 0\ I\ I\ 0\ I\ I] \cup [1\ 0\ 1\ I\ 1\ 1\ 1] \cup [1\ I\ 1\ 1$
$I\ 1\ 1\ 1\ 1] \cup [I\ 1\ I\ 0\ 1\ 1]$
$= S_1 \cup S_2 \cup S_3 \cup S_4 \cup S_5$
$= S.$

We can calculate S o $W^T$ and so on until we find the special fixed point and a special limit cycle or a special pair of fixed point or a special pair of limit cycle.

Now we use the same operation on M and find the value of the resultant vector where M is a special fuzzy neutrosophic mixed matrix. Let $M = M_1 \cup M_2 \cup \ldots \cup M_s$ ($s \geq 2$) where $M_i$ are $n_i \times n_i$ square fuzzy neutrosophic matrix and $m_j$'s are $t_j \times p_j$ ($t_j \neq p_j$) fuzzy neutrosophic rectangular matrices $1 \leq i, j \leq s$.

Suppose $X = X_1 \cup X_2 \cup \ldots \cup X_s$ where $X_i$'s are $1\times n_i$ neutrosophic fuzzy row vector and $X_j$'s are $1 \times t_j$ neutrosophic fuzzy row vector $1 \leq i, j \leq s$. To find

$X \text{ o } M = (X_1 \cup X_2 \cup \ldots \cup X_s) \text{ o } (M_1 \cup M_2 \cup \ldots \cup M_s)$
$= X_1 \text{ o } M_1 \cup X_2 \text{ o } M_2 \cup \ldots \cup X_s \text{ o } M_s.$
$= Y'_1 \cup Y'_2 \cup \ldots \cup Y'_s$
$= Y';$

Y' may or may not be a special neutrosophic fuzzy mixed row vector. We threshold Y' to $Y = Y_1 \cup Y_2 \cup \ldots \cup Y_s$ to become a special fuzzy neutrosophic mixed row vector. Now each $Y_i$ is a $1 \times n_i$ fuzzy neutrosophic row vector are $Y_j$'s are $1 \times p_j$ fuzzy neutrosophic row vector. Thus we find now Y o $M^{ST}$ where $M^{ST}$ is the special transpose of M. Let

$Y \text{ o } M^{ST} = Z'_1 \cup Z'_2 \cup \ldots \cup Z'_s$

$= \quad Z',$

we threshold and update Z' to $Z = Z_1 \cup Z_2 \cup \ldots \cup Z_s$ and find Z o M and so on until we arrive at a special fixed binary pair or a special pair of limit cycle or a fixed point and a limit cycle.

We illustrate this by the following example.

***Example 1.3.33:*** Let $M = M_1 \cup M_2 \cup M_3 \cup M_4 \cup M_5 =$

$$\begin{bmatrix} 0 & 1 & 0 & 0 & 0 \\ 0 & 0 & I & 1 & 0 \\ 1 & I & 0 & 1 & 0 \\ 0 & 0 & 0 & 1 & I \\ I & 0 & 1 & 0 & 1 \end{bmatrix} \cup \begin{bmatrix} 0 & 1 & 0 & 0 & 0 & 1 & 0 & 0 & 0 \\ I & 0 & 1 & 0 & 0 & 0 & 0 & 1 & 1 \\ 0 & 0 & 0 & 0 & 1 & 0 & I & 0 & 1 \\ 1 & I & 0 & I & 0 & 0 & 0 & 0 & 1 \end{bmatrix} \cup$$

$$\begin{bmatrix} 0 & 1 & 0 & 0 & 0 \\ I & 0 & 0 & 0 & 1 \\ 0 & 0 & 1 & I & 0 \\ 1 & 0 & I & 1 & 0 \\ 0 & I & 0 & 0 & I \\ 1 & 1 & 0 & 1 & 0 \\ 0 & 1 & 1 & 0 & 1 \end{bmatrix} \cup \begin{bmatrix} 0 & 1 & 0 & 0 \\ I & 0 & 0 & 0 \\ 0 & 0 & 1 & 0 \\ 0 & 0 & 0 & 1 \end{bmatrix} \cup \begin{bmatrix} 0 & 1 & 0 \\ 1 & 0 & 0 \\ 0 & 0 & I \\ 1 & 1 & 0 \\ 0 & 0 & 1 \end{bmatrix}$$

be the given special fuzzy neutrosophic mixed matrix.
Suppose

$$X = [1\ 0\ 0\ 0\ 0] \cup [0\ 0\ 0\ 1] \cup [0\ 0\ I\ 0\ 0\ 1\ 0] \cup [1\ 0\ 0\ 1] \cup [0\ 1\ 0\ 0\ 0]$$

be the special fuzzy neutrosophic mixed row vector. To find

$$X \text{ o } M = (X_1 \cup X_2 \cup \ldots \cup X_5) \text{ o } (M_1 \cup M_2 \cup \ldots \cup M_5)$$
$$= X_1 \text{ o } M_1 \cup X_2 \text{ o } M_2 \cup \ldots \cup X_5 \text{ o } M_5$$

$$= \begin{bmatrix}1 & 0 & 0 & 0 & 0\end{bmatrix} \text{ o } \begin{bmatrix} 0 & 1 & 0 & 0 & 0 \\ 0 & 0 & I & 1 & 0 \\ 1 & I & 0 & 1 & 0 \\ 0 & 0 & 0 & 1 & I \\ I & 0 & 1 & 0 & 1 \end{bmatrix} \cup$$

$$\begin{bmatrix}0 & 0 & 0 & 1\end{bmatrix} \text{ o } \begin{bmatrix} 0 & 1 & 0 & 0 & 0 & 1 & 0 & 0 & 0 \\ I & 0 & 1 & 0 & 0 & 0 & 0 & 1 & 1 \\ 0 & 0 & 0 & 0 & 1 & 0 & I & 0 & 1 \\ 1 & I & 0 & I & 0 & 0 & 0 & 0 & 1 \end{bmatrix} \cup$$

$$\begin{bmatrix}0 & 0 & I & 0 & 0 & 1 & 0\end{bmatrix} \text{ o } \begin{bmatrix} 0 & 1 & 0 & 0 & 0 \\ I & 0 & 0 & 0 & 1 \\ 0 & 0 & 1 & I & 0 \\ 1 & 0 & I & 1 & 0 \\ 0 & I & 0 & 0 & I \\ 1 & 1 & 0 & 1 & 0 \\ 0 & 1 & 1 & 0 & 1 \end{bmatrix} \cup$$

$$\begin{bmatrix}1 & 0 & 0 & 1\end{bmatrix} \text{ o } \begin{bmatrix} 0 & 1 & 0 & 0 \\ I & 0 & 0 & 0 \\ 0 & 0 & 1 & 0 \\ 0 & 0 & 0 & 1 \end{bmatrix} \cup$$

$$\begin{bmatrix}0 & 1 & 0 & 0 & 0\end{bmatrix} \text{ o } \begin{bmatrix} 0 & 1 & 0 \\ 1 & 0 & 0 \\ 0 & 0 & I \\ 1 & 1 & 0 \\ 0 & 0 & 1 \end{bmatrix}$$

$$= [0\ 1\ 0\ 0\ 0] \cup [1\ I\ 0\ I\ 0\ 0\ 0\ 0\ 1] \cup [1\ 1\ I\ 1{+}I\ 0] \cup [0\ 1\ 0\ 1] \cup [1\ 0\ 0]$$
$$= Z'_1 \cup Z'_2 \cup Z'_3 \cup Z'_4 \cup Z'_5$$
$$= Z'.$$

We update and threshold Z' to Z and

$$Z = Z_1 \cup Z_2 \cup Z_3 \cup Z_4 \cup Z_5$$
$$= [1\ 1\ 0\ 0\ 0] \cup [1\ I\ 0\ I\ 0\ 0\ 0\ 0\ 1] \cup [1\ 1\ I\ I\ 0] \cup [1\ 1\ 0\ 1] \cup [1\ 0\ 0].$$

Clearly Z is again a special fuzzy neutrosophic mixed row vector. Now we find

$$Z \text{ o } M^{ST} = [Z_1 \cup Z_2 \cup \ldots \cup Z_5] \text{ o } [M_1 \cup M_2^{ST} \cup M_3^{ST} \cup M_4 \cup M_5^{ST}]$$
$$= Z_1 \text{ o } M_1 \cup Z_2 \text{ o } M_2^t \cup Z_3 \text{ o } M_3^t \cup Z_4 \text{ o } M_4 \cup Z_5 \text{ o } M_5^t$$

$$= \begin{bmatrix} 1 & 1 & 0 & 0 & 0 \end{bmatrix} \text{ o } \begin{bmatrix} 0 & 1 & 0 & 0 & 0 \\ 0 & 0 & I & 1 & 0 \\ 1 & I & 0 & 1 & 0 \\ 0 & 0 & 0 & 1 & I \\ I & 0 & 1 & 0 & 1 \end{bmatrix} \cup$$

$$\begin{bmatrix} 1 & I & 0 & I & 0 & 0 & 0 & 0 & 1 \end{bmatrix} \text{ o } \begin{bmatrix} 0 & I & 0 & 1 \\ 1 & 0 & 0 & I \\ 0 & 1 & 0 & 0 \\ 0 & 0 & 0 & I \\ 0 & 0 & 1 & 0 \\ 1 & 0 & 0 & 0 \\ 0 & 0 & I & 0 \\ 0 & 1 & 0 & 0 \\ 0 & 1 & 1 & 1 \end{bmatrix} \cup$$

$$[1\ 1\ I\ I\ 0] \text{ o } \begin{bmatrix} 0 & I & 0 & 1 & 0 & 1 & 0 \\ 1 & 0 & 0 & 0 & I & 1 & 1 \\ 0 & 0 & 1 & I & 0 & 0 & 1 \\ 0 & 0 & I & 1 & 0 & 1 & 0 \\ 0 & 1 & 0 & 0 & I & 0 & 1 \end{bmatrix} \cup$$

$$[1\ 1\ 0\ 1] \text{ o } \begin{bmatrix} 0 & 1 & 0 & 0 \\ I & 0 & 0 & 0 \\ 0 & 0 & 1 & 0 \\ 0 & 0 & 0 & 1 \end{bmatrix} \cup$$

$$[1\ 0\ 0] \text{ o } \begin{bmatrix} 0 & 1 & 0 & 1 & 0 \\ 1 & 0 & 0 & 1 & 0 \\ 0 & 0 & I & 0 & 1 \end{bmatrix}$$

= [0 1 I 1 0] ∪ [I I+1 1 2 + 2I] ∪ [1 I 2I 1+2I I 2+I 1+I] ∪ [I 1 0 1] ∪ [0 1 0 1 0]
= $P'_1 \cup P'_2 \cup P'_3 \cup P'_4 \cup P'_5$
= P';

we see P' is not even a special fuzzy neutrosophic mixed row vector so we update and threshold P' to

P = $P_1 \cup P_2 \cup \ldots \cup P_5$
= [1 1 I 1 0] ∪ [I I 1 1] ∪ [1 I I I I 1 I] ∪ [1 1 0 1] ∪ [1 1 0 1 0].

P is clearly a special fuzzy neutrosophic mixed row vector.

Now we calculate

P o M = $(P_1 \cup P_2 \cup \ldots \cup P_5)$ o $(M_1 \cup M_2 \cup \ldots \cup M_5)$
= $P_1$ o $M_1 \cup P_2$ o $M_2 \cup \ldots \cup P_5$ o $M_5$

$$= \begin{bmatrix}1 & 1 & I & 1 & 0\end{bmatrix} \circ \begin{bmatrix}0 & 1 & 0 & 0 & 0\\ 0 & 0 & I & 1 & 0\\ 1 & I & 0 & 1 & 0\\ 0 & 0 & 0 & 1 & I\\ I & 0 & 1 & 0 & 1\end{bmatrix} \cup$$

$$\begin{bmatrix}I & I & 1 & 1\end{bmatrix} \circ \begin{bmatrix}0 & 1 & 0 & 0 & 0 & 1 & 0 & 0 & 0\\ I & 0 & 1 & 0 & 0 & 0 & 0 & 1 & 1\\ 0 & 0 & 0 & 0 & 1 & 0 & I & 0 & 1\\ 1 & I & 0 & I & 0 & 0 & 0 & 0 & 1\end{bmatrix} \cup$$

$$\begin{bmatrix}1 & I & I & I & I & 1 & I\end{bmatrix} \circ \begin{bmatrix}0 & 1 & 0 & 0 & 0\\ I & 0 & 0 & 0 & 1\\ 0 & 0 & 1 & I & 0\\ 1 & 0 & I & 1 & 0\\ 0 & I & 0 & 0 & I\\ 1 & 1 & 0 & 1 & 0\\ 0 & 1 & 1 & 0 & 1\end{bmatrix} \cup$$

$$\begin{bmatrix}1 & 1 & 0 & 1\end{bmatrix} \circ \begin{bmatrix}0 & 1 & 0 & 0\\ I & 0 & 0 & 0\\ 0 & 0 & 1 & 0\\ 0 & 0 & 0 & 1\end{bmatrix} \cup$$

$$\begin{bmatrix}1 & 1 & 0 & 1 & 0\end{bmatrix} \circ \begin{bmatrix}0 & 1 & 0\\ 1 & 0 & 0\\ 0 & 0 & I\\ 1 & 1 & 0\\ 0 & 0 & 1\end{bmatrix}$$

$$= [I\ 1+I\ I\ 2+I\ I] \cup [1+I\ 2I\ I\ I\ 1\ I\ I\ I\ 2+I] \cup [2I+1\ 2+2I\ 3I\ 2I+1\ 3I] \cup [I\ 1\ 0\ 1] \cup [2\ 2\ 0]$$

After updating and thresholding we get

$$T = [1\ I\ I\ I\ I] \cup [I\ I\ I\ I\ 1\ I\ I\ I\ I] \cup [I\ I\ I\ I\ I] \cup [1\ 1\ 0\ 1] \cup [1\ 1\ 0]$$
$$= T_1 \cup T_2 \cup T_3 \cup T_4 \cup T_5;$$

T is a special fuzzy neutrosophic mixed row vector. We now calculate T o $M^{ST}$ and so on.

## Chapter Two

# SPECIAL FUZZY MODELS AND SPECIAL NEUTROSOPHIC MODELS AND THEIR GENERALIZATIONS

In this chapter for the first time authors introduce seventeen new types of special fuzzy models and special neutrosophic models.

This chapter has six sections. In the first section we just recall the notion of three fuzzy models; Fuzzy Cognitive Maps (FCMs), Fuzzy Relational Maps (FRMs) and Fuzzy Relational Equations (FREs). For more information refer [232, 239].

Second section introduces the basic notion of neutrosophy from [187-190] and three neutrosophic models analogous to FCMs, FRMs and FREs namely NCMs, NRMs and NREs. Please refer [231-2]. The third section introduces five types of special fuzzy cognitive models and special Neutrosophic cognitive models.

The forth section gives yet another eight new special fuzzy and neutrosophic models which are multi expert models. Yet another set of five new special fuzzy and neutrosophic models are introduced in section five using FRE and NRE. A special super model using all these six fuzzy and neutrosophic multi expert model is also introduced in this section. The final section proposes some simple programming problems for these new special models.

## 2.1 Basic Description of Fuzzy Models

This section has three subsections. In the first subsection the fuzzy cognitive maps model is described. In the second subsection the notion of fuzzy relational maps which are a particular generalization of FCMs that too when the number of attributes are very large and can be divided into two classes are recalled for more about these concepts refer [108, 112]. In the final subsection we briefly recall the definition of FRE. Several types of FRE are introduced in [106, 232].

### 2.1.1 Definition of Fuzzy Cognitive Maps

In this section we recall the notion of Fuzzy Cognitive Maps (FCMs), which was introduced by Bart Kosko [108, 112] in the year 1986. We also give several of its interrelated definitions. FCMs have a major role to play mainly when the data concerned is an unsupervised one. Further this method is most simple and an effective one as it can analyse the data by directed graphs and connection matrices.

**DEFINITION 2.1.1.1:** *An FCM is a directed graph with concepts like policies, events etc. as nodes and causalities as edges. It represents causal relationship between concepts.*

***Example 2.1.1.1****:* In Tamil Nadu (a southern state in India) in the last decade several new engineering colleges have been approved and started. The resultant increase in the production of engineering graduates in these years is disproportionate with the need of engineering graduates. This has resulted in thousands of unemployed and underemployed graduate engineers. Using an expert's opinion we study the effect of such unemployed people on the society. An expert spells out the five major concepts relating to the unemployed graduated engineers as

$E_1$ – Frustration
$E_2$ – Unemployment
$E_3$ – Increase of educated criminals
$E_4$ – Under employment
$E_5$ – Taking up drugs etc.

The directed graph where $E_1$, …, $E_5$ are taken as the nodes and causalities as edges as given by an expert is given in the following Figure 2.1.1.1:

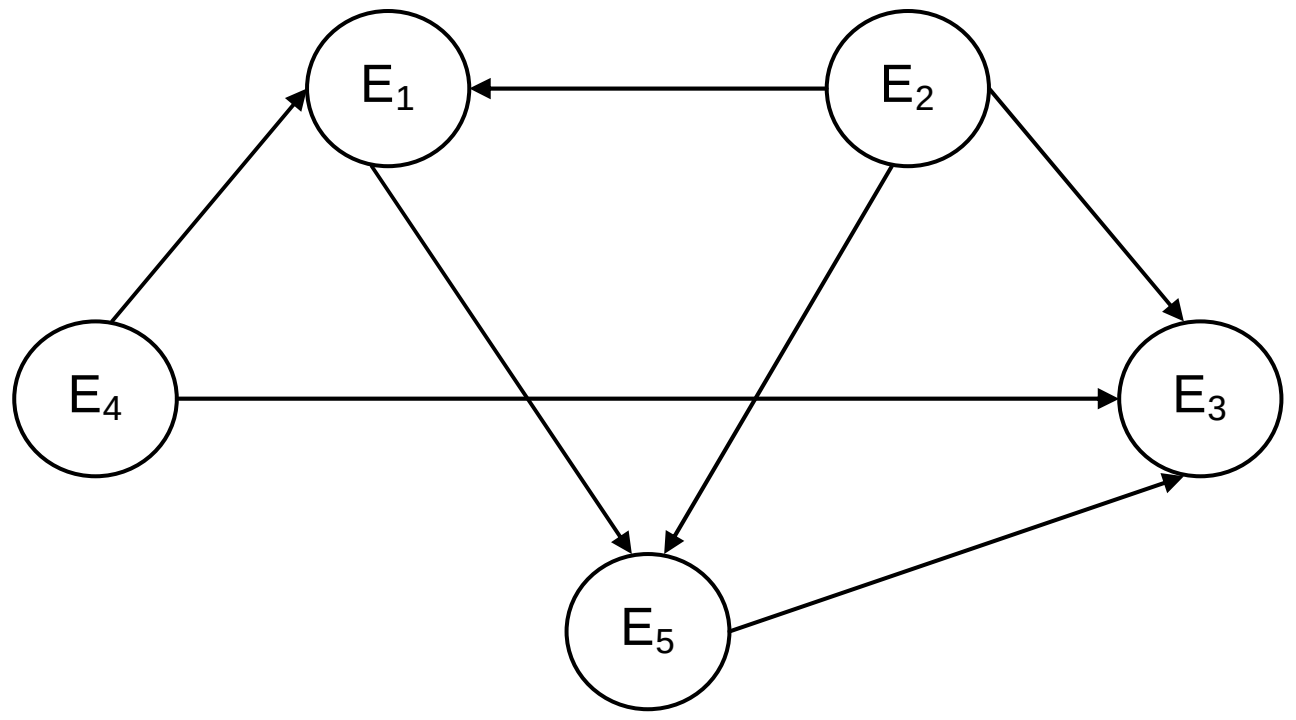

**FIGURE: 2.1.1.1**

According to this expert, increase in unemployment increases frustration. Increase in unemployment, increases the educated criminals. Frustration increases the graduates to take up to evils like drugs etc. Unemployment also leads to the increase in number of persons who take up to drugs, drinks etc. to forget their worries and unoccupied time. Under-employment forces then to do criminal acts like theft (leading to murder) for want of more money and so on. Thus one cannot actually get data for this but can use the expert's opinion for this unsupervised data to obtain some idea about the real plight of the situation. This is just an illustration to show how FCM is described by a directed graph.

{If increase (or decrease) in one concept leads to increase (or decrease) in another, then we give the value 1. If there exists no relation between two concepts the value 0 is given. If increase (or decrease) in one concept decreases (or increases) another, then we give the value –1. Thus FCMs are described in this way.}

**DEFINITION 2.1.1.2:** *When the nodes of the FCM are fuzzy sets then they are called as fuzzy nodes.*

**DEFINITION 2.1.1.3:** *FCMs with edge weights or causalities from the set {–1, 0, 1} are called simple FCMs.*

**DEFINITION 2.1.1.4:** *Consider the nodes / concepts $C_1, \ldots, C_n$ of the FCM. Suppose the directed graph is drawn using edge weight $e_{ij} \in \{0, 1, -1\}$. The matrix E be defined by $E = (e_{ij})$ where $e_{ij}$ is the weight of the directed edge $C_i C_j$. E is called the adjacency matrix of the FCM, also known as the connection matrix of the FCM.*

It is important to note that all matrices associated with an FCM are always square matrices with diagonal entries as zero.

**DEFINITION 2.1.1.5:** *Let $C_1, C_2, \ldots, C_n$ be the nodes of an FCM. $A = (a_1, a_2, \ldots, a_n)$ where $a_i \in \{0, 1\}$. A is called the instantaneous state vector and it denotes the on-off position of the node at an instant;*

$$a_i = 0 \text{ if } a_i \text{ is off and}$$
$$a_i = 1 \text{ if } a_i \text{ is on}$$

*for $i = 1, 2, \ldots, n$.*

**DEFINITION 2.1.1.6:** *Let $C_1, C_2, \ldots, C_n$ be the nodes of an FCM. Let $\overrightarrow{C_1C_2}, \overrightarrow{C_2C_3}, \overrightarrow{C_3C_4}, \ldots, \overrightarrow{C_iC_j}$ be the edges of the FCM ($i \neq j$). Then the edges form a directed cycle. An FCM is said to be cyclic if it possesses a directed cycle. An FCM is said to be acyclic if it does not possess any directed cycle.*

**DEFINITION 2.1.1.7:** *An FCM with cycles is said to have a feedback.*

**DEFINITION 2.1.1.8:** *When there is a feedback in an FCM, i.e., when the causal relations flow through a cycle in a revolutionary way, the FCM is called a dynamical system.*

**DEFINITION 2.1.1.9:** *Let* $\overline{C_1C_2}, \overline{C_2C_3}, \ldots, \overrightarrow{C_{n-1}C_n}$ *be a cycle. When* $C_i$ *is switched on and if the causality flows through the edges of a cycle and if it again causes* $C_i$ *, we say that the dynamical system goes round and round. This is true for any node* $C_i$ *, for* $i = 1, 2, \ldots, n$*. The equilibrium state for this dynamical system is called the hidden pattern.*

**DEFINITION 2.1.1.10:** *If the equilibrium state of a dynamical system is a unique state vector, then it is called a fixed point.*

***Example 2.1.1.2:*** Consider a FCM with $C_1, C_2, \ldots, C_n$ as nodes. For example let us start the dynamical system by switching on $C_1$. Let us assume that the FCM settles down with $C_1$ and $C_n$ on i.e. the state vector remains as (1, 0, 0, …, 0, 1) this state vector (1, 0, 0, …, 0, 1) is called the fixed point. (It is denoted also as [1 0 0 0 … 0 1] as the resultant is a row matrix).

**DEFINITION 2.1.1.11:** *If the FCM settles down with a state vector repeating in the form*

$$A_1 \to A_2 \to \ldots \to A_i \to A_1$$

*then this equilibrium is called a limit cycle.*

**DEFINITION 2.1.1.12:** *Finite number of FCMs can be combined together to produce the joint effect of all the FCMs. Let* $E_1, E_2, \ldots, E_p$ *be the adjacency matrices of the FCMs with nodes* $C_1, C_2, \ldots, C_n$ *then the combined FCM is got by adding all the adjacency matrices* $E_1, E_2, \ldots, E_p$.

*We denote the combined FCM adjacency matrix by* $E = E_1 + E_2 + \ldots + E_p$.

**NOTATION:** Suppose A = ($a_1$, … , $a_n$) is a vector which is passed into a dynamical system E. Then AE = ($a'_1$, … , $a'_n$) after thresholding and updating the vector suppose we get ($b_1$, … , $b_n$) we denote that by

$$(a'_1, a'_2, \ldots , a'_n) \to (b_1, b_2, \ldots , b_n).$$

Thus the symbol '$\to$' means the resultant vector has been thresholded and updated.

FCMs have several advantages as well as some disadvantages. The main advantage of this method it is simple. It functions on expert's opinion. When the data happens to be an unsupervised one the FCM comes handy. This is the only known fuzzy technique that gives the hidden pattern of the situation. As we have a very well known theory, which states that the strength of the data depends on, the number of experts' opinion we can use combined FCMs with several experts' opinions.

At the same time the disadvantage of the combined FCM is when the weightages are 1 and –1 for the same $C_i$ $C_j$, we have the sum adding to zero thus at all times the connection matrices $E_1$, … , $E_k$ may not be conformable for addition.

Combined conflicting opinions tend to cancel out and assisted by the strong law of large numbers, a consensus emerges as the sample opinion approximates the underlying population opinion. This problem will be easily overcome if the FCM entries are only 0 and 1.

We have just briefly recalled the definitions. For more about FCMs please refer Kosko [108, 112].

### 2.1.2 Definition and Illustration of Fuzzy Relational Maps (FRMs)

In this section, we introduce the notion of Fuzzy Relational Maps (FRMs); they are constructed analogous to FCMs described and discussed in the earlier sections. In FCMs we promote the correlations between causal associations among concurrently active units. But in FRMs we divide the very causal associations into two disjoint units, for example, the relation between a teacher and a student or relation between an employee or employer or a relation between doctor and patient and so on. Thus for us to define a FRM we need a domain space and a range space which are disjoint in the sense of concepts. We further assume no intermediate relation exists within the domain elements or node and the range spaces elements. The number of elements in the range space need not in general be equal to the number of elements in the domain space.

Thus throughout this section we assume the elements of the domain space are taken from the real vector space of dimension n and that of the range space are real vectors from the vector space of dimension m (m in general need not be equal to n). We denote by R the set of nodes $R_1,\ldots, R_m$ of the range space, where $R = \{(x_1,\ldots, x_m) \mid x_j = 0 \text{ or } 1 \}$ for j = 1, 2,…, m. If $x_i = 1$ it means that the node $R_i$ is in the on state and if $x_i = 0$ it means that the node $R_i$ is in the off state. Similarly D denotes the nodes $D_1, D_2,\ldots, D_n$ of the domain space where $D = \{(x_1,\ldots, x_n) \mid x_j = 0 \text{ or } 1\}$ for i = 1, 2,…, n. If $x_i = 1$ it means that the node $D_i$ is in the on state and if $x_i = 0$ it means that the node $D_i$ is in the off state.

Now we proceed on to define a FRM.

**DEFINITION 2.1.2.1:** *A FRM is a directed graph or a map from D to R with concepts like policies or events etc, as nodes and causalities as edges. It represents causal relations between spaces D and R .*

*Let $D_i$ and $R_j$ denote that the two nodes of an FRM. The directed edge from $D_i$ to $R_j$ denotes the causality of $D_i$ on $R_j$ called relations. Every edge in the FRM is weighted with a number in the set $\{0, \pm 1\}$. Let $e_{ij}$ be the weight of the edge $D_iR_j$, $e_{ij} \in \{0, \pm 1\}$. The weight of the edge $D_i R_j$ is positive if increase in $D_i$ implies increase in $R_j$ or decrease in $D_i$ implies decrease in $R_j$ ie causality of $D_i$ on $R_j$ is 1. If $e_{ij} = 0$, then $D_i$ does not have any effect on $R_j$ . We do not discuss the cases when increase in $D_i$ implies decrease in $R_j$ or decrease in $D_i$ implies increase in $R_j$ .*

**DEFINITION 2.1.2.2:** *When the nodes of the FRM are fuzzy sets then they are called fuzzy nodes. FRMs with edge weights {0, ±1} are called simple FRMs.*

**DEFINITION 2.1.2.3:** *Let $D_1, \ldots, D_n$ be the nodes of the domain space D of an FRM and $R_1, \ldots, R_m$ be the nodes of the range space R of an FRM. Let the matrix E be defined as $E = (e_{ij})$ where $e_{ij}$ is the weight of the directed edge $D_iR_j$ (or $R_jD_i$), E is called the relational matrix of the FRM.*

*Note:* It is pertinent to mention here that unlike the FCMs the FRMs can be a rectangular matrix with rows corresponding to the domain space and columns corresponding to the range space. This is one of the marked difference between FRMs and FCMs. For more about FRMs refer [241, 250]

**DEFINITION 2.1.2.4:** *Let $D_1, \ldots, D_n$ and $R_1,\ldots, R_m$ denote the nodes of the FRM. Let $A = (a_1,\ldots,a_n)$, $a_i \in \{0, 1\}$. A is called the instantaneous state vector of the domain space and it denotes the on-off position of the nodes at any instant. Similarly let $B = (b_1,\ldots, b_m)$ $b_i \in \{0, 1\}$. B is called instantaneous state vector of the range space and it denotes the on-off position of the nodes at any instant $a_i = 0$ if $a_i$ is off and $a_i = 1$ if $a_i$ is on for $i= 1, 2,\ldots, n$ Similarly, $b_i = 0$ if $b_i$ is off and $b_i = 1$ if $b_i$ is on, for $i= 1, 2,\ldots, m$.*

**DEFINITION 2.1.2.5:** *Let $D_1, \ldots, D_n$ and $R_1,\ldots, R_m$ be the nodes of an FRM. Let $D_iR_j$ (or $R_j D_i$) be the edges of an FRM, $j = 1, 2,\ldots, m$ and $i= 1, 2,\ldots, n$. Let the edges form a directed cycle. An FRM is said to be a cycle if it posses a directed cycle. An FRM is said to be acyclic if it does not posses any directed cycle.*

**DEFINITION 2.1.2.6:** *An FRM with cycles is said to be an FRM with feedback.*

**DEFINITION 2.1.2.7:** *When there is a feedback in the FRM, i.e. when the causal relations flow through a cycle in a revolutionary manner, the FRM is called a dynamical system.*

**DEFINITION 2.1.2.8:** *Let $D_i$ $R_j$ (or $R_j$ $D_i$), $1 \leq j \leq m$, $1 \leq i \leq n$. When $R_i$ (or $D_j$) is switched on and if causality flows through edges of the cycle and if it again causes $R_i$ (or$D_j$), we say that the dynamical system goes round and round. This is true for any node $R_j$ (or $D_i$) for $1 \leq i \leq n$, (or $1 \leq j \leq m$). The equilibrium state of this dynamical system is called the hidden pattern.*

**DEFINITION 2.1.2.9:** *If the equilibrium state of a dynamical system is a unique state vector, then it is called a fixed point. Consider an FRM with $R_1, R_2,\ldots, R_m$ and $D_1, D_2,\ldots, D_n$ as nodes. For example, let us start the dynamical system by switching on $R_1$ (or $D_1$). Let us assume that the FRM settles down with $R_1$ and $R_m$ (or $D_1$ and $D_n$) on, i.e. the state vector remains as (1, 0, …, 0, 1) in R (or 1, 0, 0, … , 0, 1) in D), This state vector is called the fixed point.*

**DEFINITION 2.1.2.10:** *If the FRM settles down with a state vector repeating in the form*

$$A_1 \to A_2 \to A_3 \to \ldots \to A_i \to A_1 \text{ (or } B_1 \to B_2 \to \ldots \to B_i \to B_1)$$

*then this equilibrium is called a limit cycle.*

**METHODS OF DETERMINING THE HIDDEN PATTERN**

Let $R_1, R_2,\ldots, R_m$ and $D_1, D_2,\ldots, D_n$ be the nodes of a FRM with feedback. Let E be the relational matrix. Let us find a hidden pattern when $D_1$ is switched on i.e. when an input is given as vector $A_1 = (1, 0, \ldots, 0)$ in $D_1$, the data should pass through the relational matrix E. This is done by multiplying $A_1$ with the relational matrix E. Let $A_1E = (r_1, r_2,\ldots, r_m)$, after thresholding and updating the resultant vector we get $A_1 E \in R$. Now let $B = A_1E$ we pass on B into $E^T$ and obtain $BE^T$. We update and threshold the vector $BE^T$ so that $BE^T \in D$. This procedure is repeated till we get a limit cycle or a fixed point.

**DEFINITION 2.1.2.11:** *Finite number of FRMs can be combined together to produce the joint effect of all the FRMs. Let $E_1,\ldots, E_p$ be the relational matrices of the FRMs with nodes $R_1, R_2,\ldots,$*

*$R_m$ and $D_1$, $D_2$,…, $D_n$, then the combined FRM is represented by the relational matrix $E = E_1 + \ldots + E_p$.*

## 2.1.3 Properties of Fuzzy Relations and FREs

In this section we just recollect the properties of fuzzy relations like, fuzzy equivalence relation, fuzzy compatibility relations, fuzzy ordering relations, fuzzy morphisms and sup-i-compositions of fuzzy relation. For more about these concepts please refer [231, 240].

Now we proceed on to define fuzzy equivalence relation. A crisp binary relation R(X, X) that is reflexive, symmetric and transitive is called an equivalence relation. For each element x in X, we can define a crisp set $A_x$, which contains all the elements of X that are related to x, by the equivalence relation.

$$A_x = \{y \mid (x, y) \in R(X, X)\}$$

$A_x$ is clearly a subset of X. The element x is itself contained in $A_x$ due to the reflexivity of R, because R is transitive and symmetric each member of $A_x$, is related to all the other members of $A_x$. Further no member of $A_x$, is related to any element of X not included in $A_x$. This set $A_x$ is referred to an as equivalence class of R (X, X) with respect to x. The members of each equivalence class can be considered equivalent to each other and only to each other under the relation R. The family of all such equivalence classes defined by the relation which is usually denoted by X / R, forms a partition on X.

A fuzzy binary relation that is reflexive, symmetric and transitive is known as a fuzzy equivalence relation or similarity relation. In the rest of this section let us use the latter term. While the max-min form of transitivity is assumed, in the following discussion on concepts; can be generalized to the alternative definition of fuzzy transitivity.

While an equivalence relation clearly groups elements that are equivalent under the relation into disjoint classes, the interpretation of a similarity relation can be approached in two different ways. First it can be considered to effectively group elements into crisp sets whose members are similar to each other to some specified degree. Obviously when this degree is equal to 1, the grouping is an equivalence class. Alternatively however we may wish to consider the degree of similarity that the elements of X have to some specified element $x \in X$. Thus for each $x \in X$, a similarity class can be defined as a fuzzy set in which the membership grade of any particular element represents the similarity of that element to the element x. If all the elements in the class are similar to x to the degree of 1 and similar to all elements outside the set to the degree of 0 then the grouping again becomes an equivalence class. We know every fuzzy relation R can be uniquely represented in terms of its $\alpha$-cuts by the formula

$$R = \bigcup_{\alpha \in (0,1]} \alpha . {}^{\alpha}R .$$

It is easily verified that if R is a similarity relation then each $\alpha$-cut, ${}^{\alpha}R$ is a crisp equivalence relation. Thus we may use any similarity relation R and by taking an $\alpha$ - cut ${}^{\alpha}R$ for any value $\alpha \in (0, 1]$, create a crisp equivalence relation that represents the presence of similarity between the elements to the degree $\alpha$. Each of these equivalence relations form a partition of X. Let $\pi({}^{\alpha}R)$ denote the partition corresponding to the equivalence relation ${}^{\alpha}R$. Clearly any two elements x and y belong to the same block of this partition if and only if $R(x, y) \geq \alpha$. Each similarity relation is associated with the set $\pi(R) = \{\pi({}^{\alpha}R) \mid \alpha \in (0,1]\}$ of partition of X. These partitions are nested in the sense that $\pi({}^{\alpha}R)$ is a refinement of $\pi({}^{\beta}R)$ if and only $\alpha \geq \beta$.

The equivalence classes formed by the levels of refinement of a similarity relation can be interpreted as grouping elements that are similar to each other and only to each other to a degree not less than $\alpha$.

Just as equivalences classes are defined by an equivalence relation, similarity classes are defined by a similarity relation. For a given similarity relation R(X, X) the similarity class for each $x \in X$ is a fuzzy set in which the membership grade of each element y $\in$ X is simply the strength of that elements relation to x or R(x, y). Thus the similarity class for an element x represents the degree to which all the other members of X are similar to x. Expect in the restricted case of equivalence classes themselves, similarity classes are fuzzy and therefore not generally disjoint.

Similarity relations are conveniently represented by membership matrices. Given a similarity relation R, the similarity class for each element is defined by the row of the membership matrix of R that corresponds to that element.

Fuzzy equivalence is a cutworthy property of binary relation R(X, X) since it is preserved in the classical sense in each $\alpha$-cut of R. This implies that the properties of fuzzy reflexivity, symmetry and max-min transitivity are also cutworthy. Binary relations are symmetric and transitive but not reflexive are usually referred to as quasi equivalence relations.

The notion of fuzzy equations is associated with the concept of compositions of binary relations. The composition of two fuzzy binary relations P (X, Y) and Q (Y, Z) can be defined, in general in terms of an operation on the membership matrices of P and Q that resembles matrix multiplication. This operation involves exactly the same combinations of matrix entries as in the regular matrix multiplication. However the multiplication and addition that are applied to these combinations in the matrix multiplication are replaced with other operations, these alternative operations represent in each given context the appropriate operations of fuzzy set intersections and union respectively. In the max-min composition for example, the multiplication and addition are replaced with the min and max operations respectively.

We shall give the notational conventions. Consider three fuzzy binary relations P (X, Y), Q (Y, Z) and R (X, Z) which are defined on the sets

$$X = \{x_i \mid i \in I\}$$
$$Y = \{y_j \mid j \in J\} \text{ and}$$
$$Z = \{z_k \mid k \in K\}$$

where we assume that $I = N_n$ $J = N_m$ and $K = N_s$. Let the membership matrices of P, Q and R be denoted by $P = [p_{ij}]$, $Q = [q_{ij}]$, $R = [r_{ik}]$ respectively, where $p_{ij} = P(x_i, y_j)$, $q_{jk} = Q(y_j, z_k)$, $r_{ij} = R(x_i, z_k)$ for all $i \in I$ $(=N_n)$, $j \in J = (N_m)$ and $k \in K$ $(= N_s)$. This clearly implies that all entries in the matrices P, Q, and R are real numbers from the unit interval [0, 1]. Assume now that the three relations constrain each other in such a way that $P \circ Q = R$ where $\circ$ denotes max-min composition. This means that $\max_{j \in J} \min (p_{ij}, q_{jk}) = r_{ik}$ for all $i \in I$ and $k \in K$. That is the matrix equation $P \circ Q = R$ encompasses $n \times s$ simultaneous equations of the form $\max_{j \in J} \min (p_{ij}, q_{jk}) = r_{ik}$. When two of the components in each of the equations are given and one is unknown these equations are referred to as fuzzy relation equations.

When matrices P and Q are given the matrix R is to determined using $P \circ Q = R$. The problem is trivial. It is solved simply by performing the max-min multiplication – like operation on P and Q as defined by $\max_{j \in J} \min (p_{ij}, q_{jk}) = r_{ik}$. Clearly the solution in this case exists and is unique. The problem becomes far from trivial when one of the two matrices on the left hand side of $P \circ Q = R$ is unknown. In this case the solution is guaranteed neither to exist nor to be unique.

Since R in $P \circ Q = R$ is obtained by composing P and Q it is suggestive to view the problem of determining P (or alternatively Q ) from R to Q (or alternatively R and P) as a decomposition of R with respect to Q (or alternatively with respect to P). Since many problems in various contexts can be formulated as problems of decomposition, the utility of any method for solving $P \circ Q = R$ is quite high. The use of fuzzy relation equations in some applications is illustrated. Assume that we have a method for solving $P \circ Q = R$ only for the first decomposition problem (given Q and R).

Then we can directly utilize this method for solving the second decomposition problem as well. We simply write $P \circ Q = R$ in the form $Q^{-1} \text{ o } P^{-1} = R^{-1}$ employing transposed matrices.

We can solve $Q^{-1}$ o $P^{-1}$ = $R^{-1}$ for $Q^{-1}$ by our method and then obtain the solution of P ∘ Q = R by $(Q^{-1})^{-1}$ = Q.

We study the problem of partitioning the equations P ∘ Q = R. We assume that a specific pair of matrices R and Q in the equations P ∘ Q = R is given. Let each particular matrix P that satisfies P ∘ Q = R is called its solution and let S (Q, R) = {P | P ∘ Q = R} denote the set of all solutions (the solution set).

It is easy to see this problem can be partitioned, without loss of generality into a set of simpler problems expressed by the matrix equations $p_i$ o Q = $r_i$ for all i∈I where

$$P_i = [p_{ij} \mid j \in J] \text{ and}$$
$$r_i = [r_{ik} \mid k \in K].$$

Indeed each of the equation in $\max_{j \in J} \min (p_{ij} q_{jk}) = r_{ik}$ contains unknown $p_{ij}$ identified only by one particular value of the index i, that is, the unknown $p_{ij}$ distinguished by different values of i do not appear together in any of the individual equations. Observe that $p_i$, Q, and $r_i$ in $p_i$ ∘ Q = $r_i$ represent respectively, a fuzzy set on Y, a fuzzy relation on Y × Z and a fuzzy set on Z. Let $S_i$ (Q, $r_i$) = [$p_i$ | $p_i$ o Q = $r_i$] denote, for each i ∈ I, the solution set of one of the simpler problem expressed by $p_i$ ∘ Q = $r_i$.

Thus the matrices P in S (Q, R) = [P | P ∘ Q = R ] can be viewed as one column matrix

$$P = \begin{bmatrix} p_1 \\ p_2 \\ \vdots \\ p_n \end{bmatrix}$$

where $p_i \in S_i$ (Q, $r_i$) for all i ∈ I = (=$N_n$). It follows immediately from $\max_{j \in J} \min (p_{ij}, q_{jk}) = r_{ik}$. That if $\max_{j \in J} q_{jk} < r_{ik}$ for some i ∈ I and some k ∈ K, then no values $p_{ij}$ ∈ [0, 1] exists (j ∈ J) that satisfy P ∘ Q = R, therefore no matrix P exists that satisfies the matrix equation.

This proposition can be stated more concisely as follows if

$$\max_{j \in J} q_{jk} < \max_{j \in J} r_{ik}$$

for some $k \in K$ then S (Q, R) = $\phi$. This proposition allows us in certain cases to determine quickly that P $\circ$ Q = R has no solutions its negation however is only a necessary not sufficient condition for the existence of a solution of P $\circ$ Q = R that is for S (Q, R) $\neq \phi$. Since P $\circ$ Q = R can be partitioned without loss of generality into a set of equations of the form $p_i \circ Q = r_i$ we need only methods for solving equations of the later form in order to arrive at a solution.

We may therefore restrict our further discussion of matrix equations of the form P $\circ$ Q = R to matrix equation of the simpler form P $\circ$ Q = r, where $p = [p_j \mid j \in J]$, $Q = [q_{jk} \mid j \in J, k \in K]$ and $r = \{r_k \mid k \in K]$.

We just recall the solution method as discussed by [43]. For the sake of consistency with our previous discussion, let us again assume that p, Q and r represent respectively a fuzzy set on Y, a fuzzy relation on Y $\times$ Z and a fuzzy set on Z. Moreover let $J = N_m$ and $K = N_s$ and let S (Q, r) = {p | p $\circ$ Q = r} denote the solution set of

$$p \circ Q = r.$$

In order to describe a method of solving p $\circ$ Q = r we need to introduce some additional concepts and convenient notation. First let $\wp$ denote the set of all possible vectors.

$$p = \{p_j \mid j \in J\}$$

such that $p_j \in [0, 1]$ for all $j \in J$, and let a partial ordering on $\wp$ be defined as follows for any pair $p^1, p^2 \in \wp$ $p^1 \leq p^2$ if and only if $p_i^2 \leq p_j^2$ for all $j \in J$. Given an arbitrary pair $p^1, p^2 \in \wp$ such that $p^1 \leq p^2$ let $[p^1, p^2] = \{p \in \wp \mid p^1 \leq p < p^2\}$. For any pair $p^1, p^2 \in \wp$ $(\{p^1, p^2\} \leq \}$ is a lattice.

Now we recall some of the properties of the solution set S (Q, r). Employing the partial ordering on $\wp$, let an element $\hat{p}$ of S (Q, r) be called a maximal solution of p ° Q = r if for all p ∈ S (Q, r), p ≥ $\hat{p}$ implies p = $\hat{p}$ if for all p ∈ S (Q, r) p < $\tilde{p}$ then that is the maximum solution. Similar discussion can be made on the minimal solution of p ° Q = r. The minimal solution is unique if p ≥ $\hat{p}$ (i.e. $\hat{p}$ is unique).

It is well known when ever the solution set S (Q, r) is not empty it always contains a unique maximum solution $\hat{p}$ and it may contain several minimal solution. Let $\breve{S}$ (Q, r) denote the set of all minimal solutions.

It is known that the solution set S (Q, r) is fully characterized by the maximum and minimal solution in the sense that it consists exactly of the maximum solution $\hat{p}$ all the minimal solutions and all elements of $\wp$ that are between $\hat{p}$ and the numeral solution.

Thus S (Q, r) = $\bigcup_{p} \left[\tilde{p}, \hat{p}\right]$ where the union is taken for all $\tilde{p} \in \breve{S}$ (Q, r). When S (Q, r) ≠ ϕ, the maximum solution.

$$\hat{p} = [\, \hat{p}_j \mid j \in J] \text{ of } p \circ Q = r$$

is determined as follows:

$$\hat{p}_j = \min_{k \in K} \sigma\,(q_{ik}, r_k) \text{ where } \sigma\,(q_{jk}, r_k) = \begin{cases} r_k & \text{if } q_{jk} > r_k \\ 1 & \text{otherwise} \end{cases}$$

when $\hat{p}$ determined in this way does not satisfy p ° Q = r then S(Q, r) = ϕ. That is the existence of the maximum solution $\hat{p}$ as determined by $\hat{p}_j = \min_{k \in K} \sigma\,(q_{ik}, r_k)$ is a necessary and sufficient condition for S (Q, r) ≠ ϕ. Once $\hat{p}$ is determined by

$$\hat{p}_j = \min_{k \in K} \sigma\,(q_{ik}, r_k),$$

we must check to see if it satisfies the given matrix equations p ° Q = r. If it does not then the equation has no solution (S (Q, r) = ϕ), otherwise $\hat{p}$ in the maximum solution of the equation and we next determine the set $\tilde{S}$ (Q, r) of its minimal solutions.

## 2.2 Neutrosophy and Neutrosophic models

This section has five subsections. In the first subsection a very brief introduction to neutrosophy is given for more refer [187-190]. In the second subsection some basic neutrosophic structures needed to make the book a self contained one are introduced. Sub section three briefly describes Neutrosophic Cognitive Maps (NCMs). Neutrosophic Relational Maps (NRMs) are recollected in the subsection four. The final subsection gives a brief description of binary neutrosophic relation and their properties.

### 2.2.1 An Introduction to Neutrosophy

In this section we introduce the notion of neutrosophic logic created by Florentine Smarandache [187-190], which is an extension / combination of the fuzzy logic in which indeterminacy is included. It has become very essential that the notion of neutrosophic logic play a vital role in several of the real world problems like law, medicine, industry, finance, IT, stocks and share etc. Use of neutrosophic notions will be illustrated/ applied in the later sections of this chapter. Fuzzy theory only measures the grade of membership or the non-existence of a membership in the revolutionary way but fuzzy theory has failed to attribute the concept when the relations between notions or nodes or concepts in problems are indeterminate. In fact one can say the inclusion of the concept of indeterminate situation with fuzzy concepts will form the neutrosophic logic. As in this book the concept of only fuzzy cognitive maps are dealt which mainly deals with the relation / non-relation between two nodes or concepts but it fails to deal the relation between two conceptual nodes when the relation is an indeterminate one. Neutrosophic logic is the only tool known to us, which deals with the notions of indeterminacy, and here we give a brief description of it. For more about Neutrosophic logic please refer Smarandache [187-190].

**DEFINITION 2.2.1.1:** *In the neutrosophic logic every logical variable x is described by an ordered triple x = (T, I, F) where T is the degree of truth, F is the degree of false and I the level of indeterminacy.*

(A). To maintain consistency with the classical and fuzzy logics and with probability there is the special case where T + I + F = 1.

(B). But to refer to intuitionistic logic, which means incomplete information on a variable proposition or event one has T + I + F < 1.

(C). Analogically referring to Paraconsistent logic, which means contradictory sources of information about a same logical variable, proposition or event one has T + I + F > 1.

Thus the advantage of using Neutrosophic logic is that this logic distinguishes between relative truth that is a truth is one or a few worlds only noted by 1 and absolute truth denoted by $1^+$. Likewise neutrosophic logic distinguishes between relative falsehood, noted by 0 and absolute falsehood noted by $^-0$.

It has several applications. One such given by [187-190] is as follows:

***Example 2.2.1.1:*** From a pool of refugees, waiting in a political refugee camp in Turkey to get the American visa, a% have the chance to be accepted – where a varies in the set A, r% to be rejected – where r varies in the set R, and p% to be in pending (not yet decided) – where p varies in P.

Say, for example, that the chance of someone Popescu in the pool to emigrate to USA is (between) 40-60% (considering different criteria of emigration one gets different percentages, we have to take care of all of them), the chance of being rejected is 20-25% or 30-35%, and the chance of being in pending is 10% or 20% or 30%. Then the neutrosophic probability that Popescu emigrates to the Unites States is

NP (Popescu) = ((40-60) (20-25) ∪ (30-35), {10,20,30}), closer to the life.

This is a better approach than the classical probability, where 40 P(Popescu) 60, because from the pending chance – which will be converted to acceptance or rejection – Popescu might get extra percentage in his will to emigrating and also the superior limit of the subsets sum

$$60 + 35 + 30 > 100$$

and in other cases one may have the inferior sum $< 0$, while in the classical fuzzy set theory the superior sum should be 100 and the inferior sum μ 0. In a similar way, we could say about the element Popescu that Popescu ((40-60), (20-25) ∪ (30-35), {10, 20, 30}) belongs to the set of accepted refugees.

***Example 2.2.1.2:*** The probability that candidate C will win an election is say 25-30% true (percent of people voting for him), 35% false (percent of people voting against him), and 40% or 41% indeterminate (percent of people not coming to the ballot box, or giving a blank vote – not selecting any one or giving a negative vote cutting all candidate on the list). Dialectic and dualism don't work in this case anymore.

***Example 2.2.1.3:*** Another example, the probability that tomorrow it will rain is say 50-54% true according to meteorologists who have investigated the past years weather, 30 or 34-35% false according to today's very sunny and droughty summer, and 10 or 20% undecided (indeterminate).

***Example 2.2.1.4:*** The probability that Yankees will win tomorrow versus Cowboys is 60% true (according to their confrontation's history giving Yankees' satisfaction), 30-32% false (supposing Cowboys are actually up to the mark, while Yankees are declining), and 10 or 11 or 12% indeterminate (left to the hazard: sickness of players, referee's mistakes, atmospheric conditions during the game). These parameters act on players' psychology.

As in this book we use mainly the notion of neutrosophic logic with regard to the indeterminacy of any relation in cognitive maps we are restraining ourselves from dealing with several interesting concepts about neutrosophic logic. As FCMs deals with unsupervised data and the existence or non-existence of cognitive relation, we do not in this book elaborately describe the notion of neutrosophic concepts.

However we just state, suppose in a legal issue the jury or the judge cannot always prove the evidence in a case, in several places we may not be able to derive any conclusions from the existing facts because of which we cannot make a conclusion that no relation exists or otherwise. But existing relation is an indeterminate. So in the case when the concept of indeterminacy exists the judgment ought to be very carefully analyzed be it a civil case or a criminal case. FCMs are deployed only where the existence or non-existence is dealt with but however in our Neutrosophic Cognitive Maps we will deal with the notion of indeterminacy of the evidence also. Thus legal side has lot of Neutrosophic (NCMs) applications. Also we will show how NCMs can be used to study factors as varied as stock markets, medical diagnosis, etc.

### 2.2.2 Some Basic Neutrosophic Structures

In this section we define some new neutrosophic algebraic structures like neutrosophic fields, neutrosophic spaces and neutrosophic matrices and illustrate them with examples. For these notions are used in the definition of neutrosophic cognitive maps which is dealt in the later sections of this chapter.

Throughout this book by 'I' we denote the indeterminacy of any notion/ concept/ relation. That is when we are not in a position to associate a relation between any two concepts then we denote it as indeterminacy.

Further in this book we assume all fields to be real fields of characteristic 0 and all vector spaces are taken as real spaces over reals and we denote the indeterminacy by 'I' as i will make a confusion as i denotes the imaginary value viz $i^2 = -1$ that is $\sqrt{-1} = i$.

**DEFINITION 2.2.2.1:** *Let K be the field of reals. We call the field generated by K $\cup$ I to be the neutrosophic field for it involves the indeterminacy factor in it. We define $I^2 = I$, $I + I = 2I$ i.e., $I + \ldots + I = nI$, and if $k \in K$ then $k.I = kI$, $0I = 0$. We denote the neutrosophic field by K(I) which is generated by K $\cup$ I that is $K(I) = \langle K \cup I \rangle$.*

***Example 2.2.2.1:*** Let R be the field of reals. The neutrosophic field is generated by $\langle R \cup I \rangle$ i.e. R(I) clearly $R \subset \langle R \cup I \rangle$.

***Example 2.2.2.2:*** Let Q be the field of rationals. The neutrosophic field is generated by Q and I i.e. $\langle Q \cup I \rangle$ denoted by Q(I).

**DEFINITION 2.2.2.2:** *Let K(I) be a neutrosophic field we say K(I) is a prime neutrosophic field if K(I) has no proper subfield which is a neutrosophic field.*

***Example 2.2.2.3:*** Q(I) is a prime neutrosophic field where as R(I) is not a prime neutrosophic field for $Q(I) \subset R(I)$.

It is very important to note that all neutrosophic fields are of characteristic zero. Likewise we can define neutrosophic subfield.

**DEFINITION 2.2.2.3:** *Let K(I) be a neutrosophic field, $P \subset K(I)$ is a neutrosophic subfield of P if P itself is a neutrosophic field. K(I) will also be called as the extension neutrosophic field of the neutrosophic field P.*

Now we proceed on to define neutrosophic vector spaces, which are only defined over neutrosophic fields. We can define two types of neutrosophic vector spaces one when it is a neutrosophic vector space over ordinary field other being neutrosophic vector space over neutrosophic fields. To this end we have to define neutrosophic group under addition.

**DEFINITION 2.2.2.4:** *We know Z is the abelian group under addition. Z(I) denote the additive abelian group generated by the set Z and I, Z(I) is called the neutrosophic abelian group under '+'.*

Thus to define basically a neutrosophic group under addition we need a group under addition. So we proceed on to define neutrosophic abelian group under addition. *Suppose G is an additive abelian group under '+'. G(I) = $\langle G \cup I \rangle$, additive group generated by G and I, G(I) is called the neutrosophic abelian group under '+'.*

***Example 2.2.2.4:*** Let Q be the group under '+'; Q (I) = $\langle Q \cup I \rangle$ is the neutrosophic abelian group under addition; '+'.

***Example 2.2.2.5:*** R be the additive group of reals, R(I) = $\langle R \cup I \rangle$ is the neutrosophic group under addition.

***Example 2.2.2.6:*** $M_{n \times m}(I) = \{(a_{ij}) \mid a_{ij} \in Z(I)\}$ be the collection of all $n \times m$ matrices under '+' $M_{n \times m}(I)$ is a neutrosophic group under '+'.

Now we proceed on to define neutrosophic subgroup.

**DEFINITION 2.2.2.5:** *Let G(I) be the neutrosophic group under addition. P $\subset$ G(I) be a proper subset of G(I). P is said to be neutrosophic subgroup of G(I) if P itself is a neutrosophic group i.e. P = $\langle P_1 \cup I \rangle$ where $P_1$ is an additive subgroup of G.*

***Example 2.2.2.7:*** Let Z(I) = $\langle Z \cup I \rangle$ be a neutrosophic group under '+'. $\langle 2Z \cup I \rangle$ = 2Z(I) is the neutrosophic subgroup of Z (I).

In fact Z(I) has several neutrosophic subgroups.

Now we proceed on to define the notion of neutrosophic quotient group.

**DEFINITION 2.2.2.6:** *Let G (I) = $\langle G \cup I \rangle$ be a neutrosophic group under '+', suppose P (I) be a neutrosophic subgroup of G (I) then the neutrosophic quotient group*

$$\frac{G(I)}{P(I)} = \{a + P(I) \mid a \in G\ (I)\}.$$

***Example 2.2.2.8:*** Let Z (I) be a neutrosophic group under addition, Z the group of integers under addition P = 2Z(I) is a neutrosophic subgroup of Z(I) the neutrosophic subgroup of Z(I), the neutrosophic quotient group

$$\frac{Z(I)}{P} = \{a + 2Z(I) \mid a \in Z(I)\} = \{(2n+1) + (2n+1)\ I \mid n \in Z\}.$$

Clearly $\frac{Z(I)}{P}$ is a group. For P = 2Z (I) serves as the additive identity. Take a, b $\in \frac{Z(I)}{P}$. If a, b $\in$ Z(I) \ P then two possibilities occur.

a + b is odd times I or a + b is odd or a + b is even times I or even if a + b is even or even times I then a + b $\in$ P. if a + b is odd or odd times I a + b $\in \frac{Z(I)}{P = 2Z(I)}$.

It is easily verified that P acts as the identity and every element in

$$a + 2Z\ (I) \in \frac{Z(I)}{2Z(I)}$$

has inverse. Hence the claim.

Now we proceed on to define the notion of neutrosophic vector spaces over fields and then we define neutrosophic vector spaces over neutrosophic fields.

**DEFINITION 2.2.2.7:** *Let G(I) by an additive abelian neutrosophic group. K any field. If G(I) is a vector space over K then we call G(I) a neutrosophic vector space over K.*

Now we give the notion of strong neutrosophic vector space.

**DEFINITION 2.2.2.8:** *Let G(I) be a neutrosophic abelian group. K(I) be a neutrosophic field. If G(I) is a vector space over K(I) then we call G(I) the strong neutrosophic vector space.*

**THEOREM 2.2.2.1:** *All strong neutrosophic vector space over K(I) are a neutrosophic vector space over K; as $K \subset K(I)$.*

*Proof:* Follows directly by the very definitions.

Thus when we speak of neutrosophic spaces we mean either a neutrosophic vector space over K or a strong neutrosophic vector space over the neutrosophic field K(I). By basis we mean a linearly independent set which spans the neutrosophic space.

Now we illustrate with an example.

***Example 2.2.2.9:*** Let R(I) × R(I) = V be an additive abelian neutrosophic group over the neutrosophic field R(I). Clearly V is a strong neutrosophic vector space over R(I). The basis of V are {(0,1), (1,0)}.

***Example 2.2.2.10:*** Let V = R(I) × R(I) be a neutrosophic abelian group under addition. V is a neutrosophic vector space over R. The neutrosophic basis of V are {(1,0), (0,1), (I,0), (0,I)}, which is a basis of the vector space V over R.

A study of these basis and its relations happens to be an interesting form of research.

**DEFINITION 2.2.2.9:** *Let G(I) be a neutrosophic vector space over the field K. The number of elements in the neutrosophic basis is called the neutrosophic dimension of G(I).*

**DEFINITION 2.2.2.10:** *Let G(I) be a strong neutrosophic vector space over the neutrosophic field K(I). The number of elements in the strong neutrosophic basis is called the strong neutrosophic dimension of G(I).*

We denote the neutrosophic dimension of G(I) over K by $N_k$ (dim) of G (I) and that the strong neutrosophic dimension of G (I) by $SN_{K(I)}$ (dim) of G(I).

Now we define the notion of neutrosophic matrices.

**DEFINITION 2.2.2.11:** *Let $M_{n\times m} = \{(a_{ij}) \;/\; a_{ij} \in K(I)\}$, where K (I), is a neutrosophic field. We call $M_{n\times m}$ to be the neutrosophic matrix.*

***Example 2.2.2.11:*** Let Q(I) = ⟨Q ∪ I⟩ be the neutrosophic field.

$$M_{4\times 3} = \begin{bmatrix} 0 & 1 & I \\ -2 & 4I & 0 \\ 1 & -I & 2 \\ 3I & 1 & 0 \end{bmatrix}$$

is the neutrosophic matrix, with entries from rationals and the indeterminacy I. We define product of two neutrosophic matrices whenever the production is defined as follows:

Let

$$A = \begin{bmatrix} -1 & 2 & -I \\ 3 & I & 0 \end{bmatrix}_{2\times 3}$$

and

$$B = \begin{bmatrix} -I & 1 & 2 & 4 \\ 1 & I & 0 & 2 \\ 5 & -2 & 3I & -I \end{bmatrix}_{3\times 4}$$

$$AB = \begin{bmatrix} -6I+2 & -1+4I & -2-3I & I \\ -2I & 3+I & 6 & 12+2I \end{bmatrix}_{2\times 4}$$

(we use the fact $I^2 = I$).

To define neutrosophic cognitive maps we direly need the notion of neutrosophic matrices. We use square neutrosophic matrices for Neutrosophic Cognitive Maps (NCMs) and use rectangular neutrosophic matrices for Neutrosophic Relational Maps (NRMs).

### 2.2.3 On Neutrosophic Cognitive Maps

The notion of Fuzzy Cognitive Maps (FCMs) which are fuzzy signed directed graphs with feedback are discussed and described in section 2.1.1 of this book. The directed edge $e_{ij}$ from causal concept $C_i$ to concept $C_j$ measures how much $C_i$ causes $C_j$. The time varying concept function $C_i(t)$ measures the non negative occurrence of some fuzzy event, perhaps the strength of a political sentiment, historical trend or opinion about some topics like child labor or school dropouts etc. FCMs model the world as a collection of classes and causal relations between them.

The edge $e_{ij}$ takes values in the fuzzy causal interval [–1, 1] ($e_{ij}$ = 0 indicates no causality, $e_{ij}$ > 0 indicates causal increase; that $C_j$ increases as $C_i$ increases and $C_j$ decreases as $C_i$ decreases, $e_{ij}$ < 0 indicates causal decrease or negative causality $C_j$ decreases as $C_i$ increases or $C_j$, increases as $C_i$ decreases. Simple FCMs have edge value in {-1, 0, 1}. Thus if causality occurs it occurs to maximal positive or negative degree.

It is important to note that $e_{ij}$ measures only absence or presence of influence of the node $C_i$ on $C_j$ but till now any researcher has not contemplated the indeterminacy of any relation between two nodes $C_i$ and $C_j$. When we deal with unsupervised data, there are situations when no relation can be determined between some two nodes. So in this section we try to introduce the indeterminacy in FCMs, and we choose to call this generalized structure as Neutrosophic Cognitive Maps (NCMs). In our view this will certainly give a more appropriate result and also caution the user about the risk of indeterminacy.

Now we proceed on to define the concepts about NCMs.

**DEFINITION 2.2.3.1:** *A Neutrosophic Cognitive Map (NCM) is a neutrosophic directed graph with concepts like policies, events etc. as nodes and causalities or indeterminates as edges. It represents the causal relationship between concepts.*

Let $C_1, C_2, \ldots, C_n$ denote n nodes, further we assume each node is a neutrosophic vector from neutrosophic vector space V. So a node $C_i$ will be represented by $(x_1, \ldots, x_n)$ where $x_k$'s are zero or one or I (I is the indeterminate introduced in Sections 2.2 and 2.3 of the chapter 2) and $x_k = 1$ means that the node $C_k$ is in the on state and $x_k = 0$ means the node is in the off state and $x_k = I$ means the nodes state is an indeterminate at that time or in that situation.

Let $C_i$ and $C_j$ denote the two nodes of the NCM. The directed edge from $C_i$ to $C_j$ denotes the causality of $C_i$ on $C_j$ called connections. Every edge in the NCM is weighted with a number in the set {-1, 0, 1, I}. Let $e_{ij}$ be the weight of the directed edge $C_iC_j$, $e_{ij} \in \{-1, 0, 1, I\}$. $e_{ij} = 0$ if $C_i$ does not have any effect on $C_j$, $e_{ij} = 1$ if increase (or decrease) in $C_i$ causes increase (or decreases) in $C_j$, $e_{ij} = -1$ if increase (or decrease) in $C_i$ causes decrease (or increase) in $C_j$ . $e_{ij} = I$ if the relation or effect of $C_i$ on $C_j$ is an indeterminate.

**DEFINITION 2.2.3.2:** *NCMs with edge weight from {-1, 0, 1, I} are called simple NCMs.*

**DEFINITION 2.2.3.3:** *Let $C_1, C_2, \ldots, C_n$ be nodes of a NCM. Let the neutrosophic matrix N(E) be defined as N(E) = ($e_{ij}$) where $e_{ij}$ is the weight of the directed edge $C_i C_j$, where $e_{ij} \in \{0, 1, -1, I\}$. N(E) is called the neutrosophic adjacency matrix of the NCM.*

**DEFINITION 2.2.3.4:** *Let $C_1, C_2, \ldots, C_n$ be the nodes of the NCM. Let $A = (a_1, a_2, \ldots, a_n)$ where $a_i \in \{0, 1, I\}$. A is called the* *instantaneous state neutrosophic vector and it denotes the on – off – indeterminate state position of the node at an instant*

$a_i$ = 0 if $a_i$ is off (no effect)
$a_i$ = 1 if $a_i$ is on (has effect)
$a_i$ = I if $a_i$ is indeterminate(effect cannot be determined)

*for i = 1, 2,…, n.*

**DEFINITION 2.2.3.5:** *Let $C_1, C_2, \ldots, C_n$ be the nodes of the FCM. Let $\overrightarrow{C_1C_2}$, $\overrightarrow{C_2C_3}$, $\overrightarrow{C_3C_4}$, … , $\overrightarrow{C_iC_j}$ be the edges of the NCM. Then the edges form a directed cycle. An NCM is said to be cyclic if it possesses a directed cyclic. An NCM is said to be acyclic if it does not possess any directed cycle.*

**DEFINITION 2.2.3.6:** *An NCM with cycles is said to have a feedback. When there is a feedback in the NCM i.e. when the causal relations flow through a cycle in a revolutionary manner the NCM is called a dynamical system.*

**DEFINITION 2.2.3.7:** *Let $\overrightarrow{C_1C_2}, \overrightarrow{C_2C_3}, \cdots, \overrightarrow{C_{n-1}C_n}$ be cycle, when $C_i$ is switched on and if the causality flow through the edges of a cycle and if it again causes $C_i$, we say that the dynamical system goes round and round. This is true for any node $C_i$, for i = 1, 2,…, n. the equilibrium state for this dynamical system is called the hidden pattern.*

**DEFINITION 2.2.3.8:** *If the equilibrium state of a dynamical system is a unique state vector, then it is called a fixed point. Consider the NCM with $C_1, C_2, \ldots, C_n$ as nodes. For example let us start the dynamical system by switching on $C_1$. Let us assume that the NCM settles down with $C_1$ and $C_n$ on, i.e. the state vector remain as (1, 0, …, 1) this neutrosophic state vector (1,0, …, 0, 1) is called the fixed point.*

**DEFINITION 2.2.3.9:** *If the NCM settles with a neutrosophic state vector repeating in the form*

$$A_1 \to A_2 \to \ldots \to A_i \to A_1,$$

*then this equilibrium is called a limit cycle of the NCM.*

**METHODS OF DETERMINING THE HIDDEN PATTERN:**

Let $C_1$, $C_2$, …, $C_n$ be the nodes of an NCM, with feedback. Let E be the associated adjacency matrix. Let us find the hidden pattern when $C_1$ is switched on when an input is given as the vector $A_1$ = (1, 0, 0, …, 0), the data should pass through the neutrosophic matrix N(E), this is done by multiplying $A_1$ by the matrix N(E). Let $A_1N(E) = (a_1, a_2, \ldots, a_n)$ with the threshold operation that is by replacing $a_i$ by 1 if $a_i > k$ and $a_i$ by 0 if $a_i < k$ (k – a suitable positive integer) and $a_i$ by I if $a_i$ is not a integer. We update the resulting concept, the concept $C_1$ is included in the updated vector by making the first coordinate as 1 in the resulting vector. Suppose $A_1N(E) \to A_2$ then consider $A_2N(E)$ and repeat the same procedure. This procedure is repeated till we get a limit cycle or a fixed point.

**DEFINITION 2.2.3.10:** *Finite number of NCMs can be combined together to produce the joint effect of all NCMs. If $N(E_1)$, $N(E_2)$, …, $N(E_p)$ be the neutrosophic adjacency matrices of a NCM with nodes $C_1$, $C_2$, …, $C_n$ then the combined NCM is got by adding all the neutrosophic adjacency matrices $N(E_1)$, …, $N(E_p)$. We denote the combined NCMs adjacency neutrosophic matrix by $N(E) = N(E_1) + N(E_2) + \ldots + N(E_p)$.*

**NOTATION:** Let $(a_1, a_2, \ldots, a_n)$ and $(a'_1, a'_2, \ldots, a'_n)$ be two neutrosophic vectors. We say $(a_1, a_2, \ldots, a_n)$ is equivalent to $(a'_1, a'_2, \ldots, a'_n)$ denoted by $((a_1, a_2, \ldots, a_n) \sim (a'_1, a'_2, \ldots, a'_n)$ if $(a'_1, a'_2, \ldots, a'_n)$ is got after thresholding and updating the vector $(a_1, a_2, \ldots, a_n)$ after passing through the neutrosophic adjacency matrix N(E).

The following are very important:

***Note 1:*** The nodes $C_1$, $C_2$, …, $C_n$ are nodes are not indeterminate nodes because they indicate the concepts which are well known. But the edges connecting $C_i$ and $C_j$ may be indeterminate i.e. an expert may not be in a position to say that $C_i$ has some causality on $C_j$ either will he be in a position to state that $C_i$ has no relation with $C_j$ in such cases the relation between $C_i$ and $C_j$ which is indeterminate is denoted by I.

***Note 2:*** The nodes when sent will have only ones and zeros i.e. on and off states, but after the state vector passes through the neutrosophic adjacency matrix the resultant vector will have entries from {0, 1, I} i.e. they can be neutrosophic vectors.

The presence of I in any of the coordinate implies the expert cannot say the presence of that node i.e. on state of it after passing through N(E) nor can we say the absence of the node i.e. off state of it the effect on the node after passing through the dynamical system is indeterminate so only it is represented by I. Thus only in case of NCMs we can say the effect of any node on other nodes can also be indeterminates. Such possibilities and analysis is totally absent in the case of FCMs.

***Note 3:*** In the neutrosophic matrix N(E), the presence of I in the $a_{ij}$ the place shows, that the causality between the two nodes i.e. the effect of $C_i$ on $C_j$ is indeterminate. Such chances of being indeterminate is very possible in case of unsupervised data and that too in the study of FCMs which are derived from the directed graphs.

Thus only NCMs helps in such analysis.

Now we shall represent a few examples to show how in this set up NCMs is preferred to FCMs. At the outset before we proceed to give examples we make it clear that all unsupervised data need not have NCMs to be applied to it. Only data which have the relation between two nodes to be an indeterminate need to be modeled with NCMs if the data has no indeterminacy factor between any pair of nodes one need not go for NCMs; FCMs will do the best job.

### 2.2.4 Neutrosophic Relational Maps

Neutrosophic Cognitive Maps (NCMs) promote the causal relationships between concurrently active units or decides the absence of any relation between two units or the indeterminacy of any relation between any two units. But in Neutrosophic Relational Maps (NRMs) we divide the very causal nodes into two disjoint units.

Thus for the modeling of a NRM we need a domain space and a range space which are disjoint in the sense of concepts. We further assume no intermediate relations exist within the domain and the range spaces. The number of elements or nodes in the range space need not be equal to the number of elements or nodes in the domain space.

Throughout this section we assume the elements of a domain space are taken from the neutrosophic vector space of dimension n and that of the range space are neutrosophic vector space of dimension m. (m in general need not be equal to n). We denote by R the set of nodes $R_1, \ldots, R_m$ of the range space, where $R = \{(x_1, \ldots, x_m) \mid x_j = 0 \text{ or } 1 \text{ for } j = 1, 2, \ldots, m\}$.

If $x_i = 1$ it means that node $R_i$ is in the on state and if $x_i = 0$ it means that the node $R_i$ is in the off state and if $x_i = I$ in the resultant vector it means the effect of the node $x_i$ is indeterminate or whether it will be off or on cannot be predicted by the neutrosophic dynamical system.

It is very important to note that when we send the state vectors they are always taken as the real state vectors for we know the node or the concept is in the on state or in the off state but when the state vector passes through the Neutrosophic dynamical system some other node may become indeterminate i.e. due to the presence of a node we may not be able to predict the presence or the absence of the other node i.e., it is indeterminate, denoted by the symbol I, thus the resultant vector can be a neutrosophic vector.

**DEFINITION 2.2.4.1:** *A Neutrosophic Relational Map (NRM) is a Neutrosophic directed graph or a map from D to R with*

*concepts like policies or events etc. as nodes and causalities as edges. (Here by causalities we mean or include the indeterminate causalities also). It represents Neutrosophic Relations and Causal Relations between spaces D and R .*

*Let $D_i$ and $R_j$ denote the nodes of an NRM. The directed edge from $D_i$ to $R_j$ denotes the causality of $D_i$ on $R_j$ called relations. Every edge in the NRM is weighted with a number in the set {0, +1, –1, I}.*

*Let $e_{ij}$ be the weight of the edge $D_i R_j$, $e_{ij} \in \{0, 1, -1, I\}$. The weight of the edge $D_i R_j$ is positive if increase in $D_i$ implies increase in $R_j$ or decrease in $D_i$ implies decrease in $R_j$ i.e. causality of $D_i$ on $R_j$ is 1. If $e_{ij} = -1$ then increase (or decrease) in $D_i$ implies decrease (or increase) in $R_j$. If $e_{ij} = 0$ then $D_i$ does not have any effect on $R_j$. If $e_{ij} = I$ it implies we are not in a position to determine the effect of $D_i$ on $R_j$ i.e. the effect of $D_i$ on $R_j$ is an indeterminate so we denote it by I.*

**DEFINITION 2.2.4.2:** *When the nodes of the NRM take edge values from {0, 1, –1, I} we say the NRMs are simple NRMs.*

**DEFINITION 2.2.4.3:** *Let $D_1, \ldots, D_n$ be the nodes of the domain space D of an NRM and let $R_1, R_2, \ldots, R_m$ be the nodes of the range space R of the same NRM. Let the matrix N(E) be defined as $N(E) = (e_{ij})$ where $e_{ij}$ is the weight of the directed edge $D_i R_j$ (or $R_j D_i$) and $e_{ij} \in \{0, 1, -1, I\}$. N(E) is called the Neutrosophic Relational Matrix of the NRM.*

The following remark is important and interesting to find its mention in this book.

**Remark**: Unlike NCMs, NRMs can also be rectangular matrices with rows corresponding to the domain space and columns corresponding to the range space. This is one of the marked difference between NRMs and NCMs. Further the number of entries for a particular model which can be treated as disjoint sets when dealt as a NRM has very much less entries than when the same model is treated as a NCM.

Thus in many cases when the unsupervised data under study or consideration can be spilt as disjoint sets of nodes or concepts; certainly NRMs are a better tool than the NCMs.

**DEFINITION 2.2.4.4:** *Let* $D_1, \ldots, D_n$ *and* $R_1, \ldots, R_m$ *denote the nodes of a NRM. Let* $A = (a_1, \ldots, a_n)$, $a_i \in \{0, 1\}$ *is called the neutrosophic instantaneous state vector of the domain space and it denotes the on-off position of the nodes at any instant. Similarly let* $B = (b_1, \ldots, b_n)$ $b_i \in \{0, 1\}$, *B is called instantaneous state vector of the range space and it denotes the on-off position of the nodes at any instant,* $a_i = 0$ *if* $a_i$ *is off and* $a_i = 1$ *if* $a_i$ *is on for* $i = 1, 2, \ldots, n$. *Similarly,* $b_i = 0$ *if* $b_i$ *is off and* $b_i = 1$ *if* $b_i$ *is on for* $i = 1, 2, \ldots, m$.

**DEFINITION 2.2.4.5:** *Let* $D_1, \ldots, D_n$ *and* $R_1, R_2, \ldots, R_m$ *be the nodes of a NRM. Let* $D_i R_j$ *(or* $R_j D_i$ *) be the edges of an NRM,* $j = 1, 2, \ldots, m$ *and* $i = 1, 2, \ldots, n$. *The edges form a directed cycle. An NRM is said to be a cycle if it possess a directed cycle. An NRM is said to be acyclic if it does not possess any directed cycle.*

**DEFINITION 2.2.4.6:** *A NRM with cycles is said to be a NRM with feedback.*

**DEFINITION 2.2.4.7:** *When there is a feedback in the NRM i.e. when the causal relations flow through a cycle in a revolutionary manner the NRM is called a Neutrosophic dynamical system.*

**DEFINITION 2.2.4.8:** *Let* $D_i R_j$ *(or* $R_j D_i$*)* $1 \leq j \leq m$, $1 \leq i \leq n$, *when* $R_j$ *(or* $D_i$ *) is switched on and if causality flows through edges of a cycle and if it again causes* $R_j$ *(or* $D_i$ *) we say that the Neutrosophical dynamical system goes round and round. This is true for any node* $R_j$ *( or* $D_i$ *) for* $1 \leq j \leq m$ *(or* $1 \leq i \leq n$*). The equilibrium state of this Neutrosophical dynamical system is called the Neutrosophic hidden pattern.*

**DEFINITION 2.2.4.9:** *If the equilibrium state of a Neutrosophical dynamical system is a unique Neutrosophic*

*state vector, then it is called the fixed point. Consider an NRM with $R_1$, $R_2$, …, $R_m$ and $D_1$, $D_2$,…, $D_n$ as nodes. For example let us start the dynamical system by switching on $R_1$ (or $D_1$). Let us assume that the NRM settles down with $R_1$ and $R_m$ (or $D_1$ and $D_n$) on, or indeterminate on, i.e. the Neutrosophic state vector remains as (1, 0, 0,…, 1) or (1, 0, 0,…I) (or (1, 0, 0,…1) or (1, 0, 0,…I) in D), this state vector is called the fixed point.*

**DEFINITION 2.2.4.10:** *If the NRM settles down with a state vector repeating in the form $A_1 \to A_2 \to A_3 \to \ldots \to A_i \to A_1$ (or $B_1 \to B_2 \to \ldots \to B_i \to B_1$) then this equilibrium is called a limit cycle.*

**METHODS OF DETERMINING THE HIDDEN PATTERN IN A NRM**

Let $R_1$, $R_2$,…, $R_m$ and $D_1$, $D_2$,…, $D_n$ be the nodes of a NRM with feedback. Let N(E) be the Neutrosophic Relational Matrix. Let us find the hidden pattern when $D_1$ is switched on i.e. when an input is given as a vector; $A_1$ = (1, 0, …, 0) in D; the data should pass through the relational matrix N(E).

This is done by multiplying $A_1$ with the Neutrosophic relational matrix N(E). Let $A_1N(E) = (r_1, r_2, \ldots, r_m)$ after thresholding and updating the resultant vector we get $A_1E \in R$, Now let $B = A_1E$ we pass on B into the system $(N(E))^T$ and obtain $B(N(E))^T$. We update and threshold the vector $B(N(E))^T$ so that $B(N(E))^T \in D$.

This procedure is repeated till we get a limit cycle or a fixed point.

**DEFINITION 2.2.4.11:** *Finite number of NRMs can be combined together to produce the joint effect of all NRMs. Let $N(E_1)$, $N(E_2)$,…, $N(E_r)$ be the Neutrosophic relational matrices of the NRMs with nodes $R_1$,…, $R_m$ and $D_1$,…,$D_n$, then the combined NRM is represented by the neutrosophic relational matrix $N(E) = N(E_1) + N(E_2) + \ldots + N(E_r)$.*

## 2.2.5 Binary Neutrosophic Relation and their Properties

In this section we introduce the notion of neutrosophic relational equations and fuzzy neutrosophic relational equations and analyze and apply them to real-world problems, which are abundant with the concept of indeterminacy. We also mention that most of the unsupervised data also involve at least to certain degrees the notion of indeterminacy.

Throughout this section by a neutrosophic matrix we mean a matrix whose entries are from the set N = [0, 1] $\cup$ *I* and by a fuzzy neutrosophic matrix we mean a matrix whose entries are from N' = [0, 1] $\cup$ {n*I* / n $\in$ (0, 1]}.

Now we proceed on to define binary neutrosophic relations and binary neutrosophic fuzzy relation.

A binary neutrosophic relation $R_N(x, y)$ may assign to each element of X two or more elements of Y or the indeterminate *I*. Some basic operations on functions such as the inverse and composition are applicable to binary relations as well. Given a neutrosophic relation $R_N(X, Y)$ its domain is a neutrosophic set on X $\cup$ *I* domain R whose membership function is defined by domR(x) = $\max_{y \in X \cup I} R_N(x, y)$ for each x $\in$ X $\cup$ *I*.

That is each element of set X $\cup$ *I* belongs to the domain of R to the degree equal to the strength of its strongest relation to any member of set Y $\cup$ *I*. The degree may be an indeterminate *I* also. Thus this is one of the marked difference between the binary fuzzy relation and the binary neutrosophic relation. The range of $R_N(X, Y)$ is a neutrosophic relation on Y, ran R whose membership is defined by ran R(y) = $\max_{x \in X} R_N(x, y)$ for each y $\in$ Y, that is the strength of the strongest relation that each element of Y has to an element of X is equal to the degree of that element's membership in the range of R or it can be an indeterminate *I*.

The height of a neutrosophic relation $R_N(x, y)$ is a number h(R) or an indeterminate *I* defined by $h_N(R) = \max_{y \in Y \cup I} \max_{x \in X \cup I} R_N(x, y)$. That is $h_N(R)$ is the largest membership grade attained by any pair (x, y) in R or the indeterminate *I*.

A convenient representation of the neutrosophic binary relation $R_N(X, Y)$ are membership matrices $R = [\gamma_{xy}]$ where $\gamma_{xy} \in R_N(x, y)$. Another useful representation of a binary neutrosophic relation is a neutrosophic sagittal diagram. Each of the sets X, Y represented by a set of nodes in the diagram, nodes corresponding to one set are clearly distinguished from nodes representing the other set. Elements of X' × Y' with non-zero membership grades in $R_N(X, Y)$ are represented in the diagram by lines connecting the respective nodes. These lines are labeled with the values of the membership grades.

An example of the neutrosophic sagittal diagram is a binary neutrosophic relation $R_N(X, Y)$ together with the membership

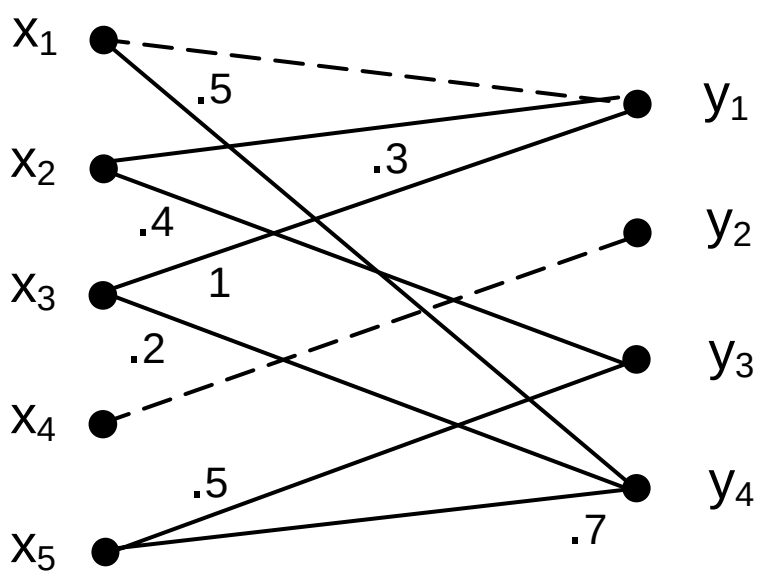

FIGURE: 2.2.5.1

neutrosophic matrix is given below.

$$\begin{array}{c} \\ x_1 \\ x_2 \\ x_3 \\ x_4 \\ x_5 \end{array}\begin{array}{c} \begin{array}{cccc} y_1 & y_2 & y_3 & y_4 \end{array} \\ \begin{bmatrix} I & 0 & 0 & 0.5 \\ 0.3 & 0 & 0.4 & 0 \\ 1 & 0 & 0 & 0.2 \\ 0 & I & 0 & 0 \\ 0 & 0 & 0.5 & 0.7 \end{bmatrix} \end{array}.$$

The inverse of a neutrosophic relation $R_N(X, Y) = R(x, y)$ for all $x \in X$ and all $y \in Y$. A neutrosophic membership matrix $R^{-1} = [r^{-1}_{yx}]$ representing $R^{-1}_N(Y, X)$ is the transpose of the matrix R for $R_N(X, Y)$ which means that the rows of $R^{-1}$ equal the columns of R and the columns of $R^{-1}$ equal rows of R. Clearly $(R^{-1})^{-1} = R$ for any binary neutrosophic relation.

Consider any two binary neutrosophic relation $P_N(X, Y)$ and $Q_N(Y, Z)$ with a common set Y. The standard composition of these relations which is denoted by $P_N(X, Y) \bullet Q_N(Y, Z)$ produces a binary neutrosophic relation $R_N(X, Z)$ on $X \times Z$ defined by $R_N(x, z) = [P \bullet Q]_N(x, z) = \max_{y \in Y} \min[P_N(x, y), Q_N(x, y)]$ for all $x \in X$ and all $z \in Z$.

This composition which is based on the standard $t_N$-norm and $t_N$-co-norm, is often referred to as the max-min composition. It can be easily verified that even in the case of binary neutrosophic relations $[P_N(X, Y) \bullet Q_N(Y, Z)]^{-1} = Q_N^{-1}(Z, Y) \bullet P_N^{-1}(Y, X)$. $[P_N(X, Y) \bullet Q_N(Y, Z)] \bullet R_N(Z, W) = P_N(X, Y) \bullet [Q_N(Y, Z) \bullet R_N(Z, W)]$, that is, the standard (or max-min) composition is associative and its inverse is equal to the reverse composition of the inverse relation. However, the standard composition is not commutative, because $Q_N(Y, Z) \bullet P_N(X, Y)$ is not well defined when $X \neq Z$. Even if X = Z and $Q_N (Y, Z) \circ P_N (X, Y)$ are well defined still we can have $P_N (X, Y) \circ Q (Y, Z) \neq Q (Y, Z) \circ P (X, Y)$.

Compositions of binary neutrosophic relation can the performed conveniently in terms of membership matrices of the relations. Let $P = [p_{ik}]$, $Q = [q_{kj}]$ and $R = [r_{ij}]$ be membership matrices of binary relations such that $R = P \circ Q$. We write this using matrix notation

$$[r_{ij}] = [p_{ik}] \text{ o } [q_{kj}]$$

where $r_{ij} = \max_{k} \min (p_{ik}, q_{kj})$.

A similar operation on two binary relations, which differs from the composition in that it yields triples instead of pairs, is known as the relational join. For neutrosophic relation $P_N (X, Y)$ and $Q_N (Y, Z)$ the relational join P * Q corresponding to the neutrosophic standard max-min composition is a ternary relation $R_N (X, Y, Z)$ defined by $R_N (x, y, z) = [P * Q]_N (x, y, z) = \min [P_N (x, y), Q_N (y, z)]$ for each $x \in X$, $y \in Y$ and $z \in Z$.

This is illustrated by the following Figure 2.2.5.2.

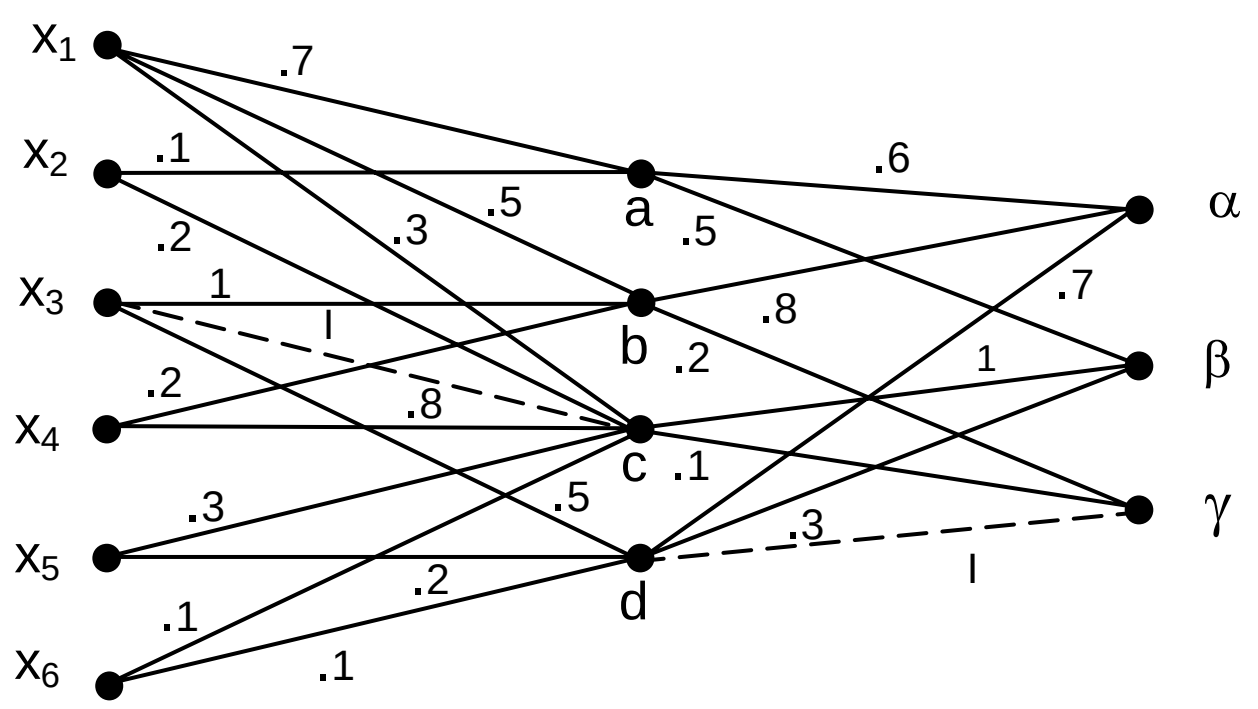

FIGURE: 2.2.5.2

In addition to defining a neutrosophic binary relation there exists between two different sets, it is also possible to define neutrosophic binary relation among the elements of a single set X.

A neutrosophic binary relation of this type is denoted by $R_N(X, X)$ or $R_N(X^2)$ and is a subset of $X \times X = X^2$.

These relations are often referred to as neutrosophic directed graphs or neutrosophic digraphs.

Neutrosophic binary relations $R_N(X, X)$ can be expressed by the same forms as general neutrosophic binary relations. However they can be conveniently expressed in terms of simple diagrams with the following properties:

I. Each element of the set X is represented by a single node in the diagram.

II. Directed connections between nodes indicate pairs of elements of X for which the grade of membership in R is non zero or indeterminate.

III. Each connection in the diagram is labeled by the actual membership grade of the corresponding pair in R or in indeterminacy of the relationship between those pairs.

The neutrosophic membership matrix and the neutrosophic sagittal diagram is as follows for any set X = {a, b, c, d, e}.

$$\begin{bmatrix} 0 & I & 0.3 & 0.2 & 0 \\ 1 & 0 & I & 0 & 0.3 \\ I & 0.2 & 0 & 0 & 0 \\ 0 & 0.6 & 0 & 0.3 & I \\ 0 & 0 & 0 & I & 0.2 \end{bmatrix}.$$

Neutrosophic membership matrix for x is given above and the neutrosophic sagittal diagram is given below.

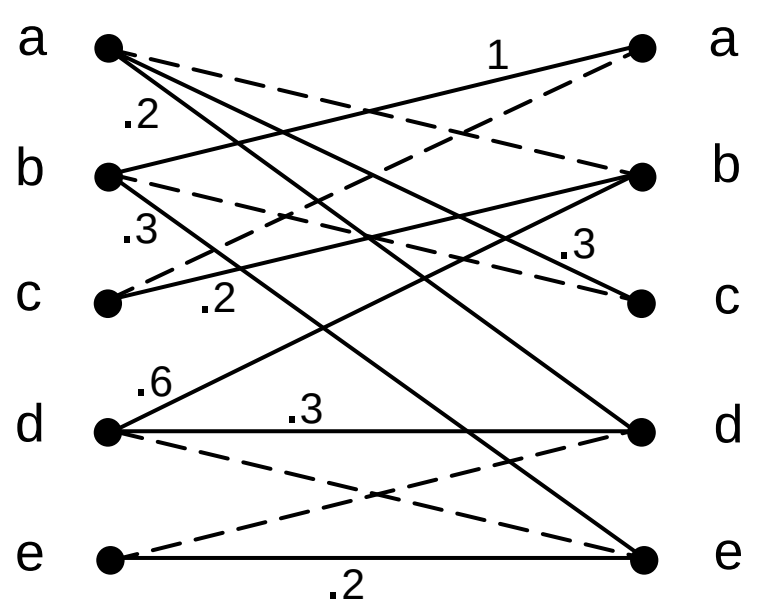

FIGURE 2.2.5.3

Neutrosophic diagram or graph is left for the reader as an exercise.

The notion of reflexivity, symmetry and transitivity can be extended for neutrosophic relations $R_N$ (X, Y) by defining them in terms of the membership functions or indeterminacy relation.

*Thus $R_N$ (X, X) is reflexive if and only if $R_N$ (x, x) = 1 for all x ∈ X.* If *this is not the case for some x ∈ X the relation is irreflexive.*
*A weaker form of reflexivity, if for no x in X; $R_N$(x, x) = 1 then we call the relation to be anti-reflexive referred to as ∈-reflexivity, is sometimes defined by requiring that*

$$R_N(x, x) \geq \in \text{ where } 0 < \in < 1.$$

*A fuzzy relation is symmetric if and only if*

$$R_N(x, y) = R_N(y, x) \text{ for all } x, y, \in X.$$

*Whenever this relation is not true for some x, y $\in$ X the relation is called asymmetric. Furthermore when $R_N(x, y) > 0$ and $R_N(y, x) > 0$ implies that x = y for all x, y $\in$ X the relation R is called anti-symmetric.*

*A fuzzy relation $R_N(X, X)$ is transitive (or more specifically max-min transitive) if*

$$R_N(x, z) \geq \max_{y \in Y} \min [R_N(x, y), R_N(y, z)]$$

*is satisfied for each pair (x, z) $\in X^2$. A relation failing to satisfy the above inequality for some members of X is called non-transitive and if $R_N(x, x) < \max_{y \in Y} \min [RN(x, y), RN(y, z)]$ for all (x, x) $\in X^2$, then the relation is called anti-transitive.*

*Given a relation $R_N(X, X)$ its transitive closure $\overline{R}_{NT}(x, X)$ can be analyzed in the following way.*

The transitive closure on a crisp relation $R_N$ (X, X) is defined as the relation that is transitive, contains

$$R_N(X, X) < \max_{y \in Y} \min [R_N(x, y) R_N(y, z)]$$

for all (x, x) $\in X^2$, then the relation is called anti-transitive. Given a relation $R_N$ (x, x) its transitive closure $\overline{R}_{NT}$ (X, X) can be analyzed in the following way.

The transitive closure on a crisp relation $R_N$ (X, X) is defined as the relation that is transitive, contains $R_N$ and has the fewest possible members. For neutrosophic relations the last requirement is generalized such that the elements of transitive closure have the smallest possible membership grades, that still allow the first two requirements to be met.

Given a relation $R_N$ (X, X) its transitive closure $\overline{R}_{NT}$ (X, X) can be determined by a simple algorithm.

Now we proceed on to define the notion of neutrosophic equivalence relation.

**DEFINITION 2.2.5.1:** *A crisp neutrosophic relation $R_N(X, X)$ that is reflexive, symmetric and transitive is called an neutrosophic equivalence relation. For each element x in X, we can define a crisp neutrosophic set $A_x$ which contains all the elements of X that are related to x by the neutrosophic equivalence relation.*

*Formally $A_x = [y \mid (x, y) \in R_N (X, X)\}$. $A_x$ is clearly a subset of X. The element x is itself contained in $A_x$, due to the reflexivity of R because R is transitive and symmetric each member of $A_x$ is related to all other members of $A_x$. Further no member of $A_x$ is related to any element of X not included in $A_x$. This set $A_x$ is clearly referred to as an neutrosophic equivalence class of $R_N$ (X, x) with respect to x. The members of each neutrosophic equivalence class can be considered neutrosophic equivalent to each other and only to each other under the relation R.*

But here it is pertinent to mention that in some X; (a, b) may not be related at all to be more precise there may be an element $a \in X$ which is such that its relation with several or some elements in $X \setminus \{a\}$ is an indeterminate. The elements which cannot determine its relation with other elements will be put in as separate set.

A neutrosophic binary relation that is reflexive, symmetric and transitive is known as a neutrosophic equivalence relation.

Now we proceed on to define Neutrosophic intersections neutrosophic t – norms ($t_N$ – norms)

Let A and B be any two neutrosophic sets, the intersection of A and B is specified in general by a neutrosophic binary operation on the set $N = [0, 1] \cup I$, that is a function of the form

$$i_N: N \times N \to N.$$

For each element x of the universal set, this function takes as its argument the pair consisting of the elements membership grades in set A and in set B, and yield the membership grade of the element in the set constituting the intersection of A and B. Thus,

$$(A \cap B)(x) = i_N [A(x), B(x)] \text{ for all } x \in X.$$

In order for the function $i_N$ of this form to qualify as a fuzzy intersection, it must possess appropriate properties, which ensure that neutrosophic sets produced by $i_N$ are intuitively acceptable as meaningful fuzzy intersections of any given pair of neutrosophic sets. It turns out that functions known as $t_N$-norms, will be introduced and analyzed in this section. In fact the class of $t_N$- norms is now accepted as equivalent to the class of neutrosophic fuzzy intersection. We will use the terms $t_N$ – norms and neutrosophic intersections inter changeably.

Given a $t_N$ – norm, $i_N$ and neutrosophic sets A and B we have to apply:

$$(A \cap B)(x) = i_N [A(x), B(x)]$$

for each $x \in X$, to determine the intersection of A and B based upon $i_N$.

However the function $i_N$ is totally independent of x, it depends only on the values A (x) and B(x). Thus we may ignore x and assume that the arguments of $i_N$ are arbitrary numbers $a, b \in [0, 1] \cup I = N$ in the following examination of formal properties of $t_N$-norm.

A neutrosophic intersection/ $t_N$-norm $i_N$ is a binary operation on the unit interval that satisfies at least the following axioms for all $a, b, c, d \in N = [0, 1] \cup I$.

$1_N$ $i_N(a, 1) = a$

$2_N$ $i_N(a, I) = I$

$3_N$ $b \leq d$ implies $i_N(a, b) \leq i_N(a, d)$

$4_N$ $i_N(a, b) = i_N(b, a)$

$5_N$ $i_N(a, i_N(b, d)) = i_N(a, b), d)$.

We call the conditions $1_N$ to $5_N$ as the axiomatic skeleton for neutrosophic intersections / $t_N$ – norms. Clearly $i_N$ is a continuous function on $N \setminus I$ and $i_N(a, a) < a \ \forall a \in N \setminus I$

$$i_N(I\, I) = I.$$

If $a_1 < a_2$ and $b_1 < b_2$ implies $i_N(a_1, b_1) < i_N(a_2, b_2)$.

Several properties in this direction can be derived as in case of t-norms.

The following are some examples of $t_N$ –norms

1. $i_N$ (a, b) = min (a, b)
   $i_N$ (a, *I*) = min (a, *I*) = *I* will be called as standard neutrosophic intersection.
2. $i_N$ (a, b) = ab for a, b $\in$ N \ *I* and $i_N$ (a, b) = *I* for a, b $\in$ N where a = *I* or b = *I* will be called as the neutrosophic algebraic product.
3. Bounded neutrosophic difference.
   $i_N$ (a, b) = max (0, a + b – 1) for a, b $\in$ *I*
   $i_N$ (a, *I*) = *I* is yet another example of $t_N$ – norm.

1. Drastic neutrosophic intersection
2.

$$i_N(a, b) = \begin{cases} a \text{ when } b = 1 \\ b \text{ when } a = 1 \\ I \text{ when } a = I \\ \quad \text{or } b = I \\ \quad \text{or } a = b = I \\ 0 \quad \text{otherwise.} \end{cases}$$

As *I* is an indeterminate adjoined in $t_N$ – norms. It is not easy to give then the graphs of neutrosophic intersections. Here also we leave the analysis and study of these $t_N$ – norms (i.e. neutrosophic intersections) to the reader.

The notion of neutrosophic unions closely parallels that of neutrosophic intersections. Like neutrosophic intersection the general neutrosophic union of two neutrosophic sets A and B is specified by a function

$$\mu_N: N \times N \to N \text{ where } N = [0, 1] \cup I.$$

The argument of this function is the pair consisting of the membership grade of some element x in the neutrosophic set A and the membership grade of that some element in the neutrosophic set B, (here by membership grade we mean not only the membership grade in the unit interval [0, 1] but also the indeterminacy of the membership). The function returns the membership grade of the element in the set $A \cup B$.

Thus $(A \cup B)(x) = \mu_N [A(x), B(x)]$ for all $x \in X$. Properties that a function $\mu_N$ must satisfy to be initiatively acceptable as neutrosophic union are exactly the same as properties of functions that are known. Thus neutrosophic union will be called as neutrosophic t-co-norm; denoted by $t_N$ – co-norm.

A neutrosophic union / $t_N$ – co-norm $\mu_N$ is a binary operation on the unit interval that satisfies at least the following conditions for all a, b, c, d $\in N = [0, 1] \cup I$

$C_1$ $\mu_N(a, 0) = a$

$C_2$ $\mu_N(a, I) = I$

$C_3$ $b \leq d$ implies
$\mu_N(a, b) \leq \mu_N(a, d)$

$C_4$ $\mu_N(a, b) = \mu_N(b, a)$

$C_5$ $\mu_N(a, \mu_N(b, d))$
$= \mu_N(\mu_N(a, b), d)$.

Since the above set of conditions are essentially neutrosophic unions we call it the axiomatic skeleton for neutrosophic unions / $t_N$-co-norms.

The addition requirements for neutrosophic unions are

i. $\mu_N$ is a continuous functions on $N \setminus \{I\}$
ii. $\mu_N(a, a) > a$.
iii. $a_1 < a_2$ and $b_1 < b_2$ implies $\mu_N(a_1. b_1) < \mu_N(a_2, b_2)$; $a_1, a_2, b_1, b_2 \in N \setminus \{I\}$.

We give some basic neutrosophic unions.

Let $\mu_N : [0, 1] \times [0, 1] \to [0, 1]$

$$\mu_N(a, b) = \max(a, b)$$

$\mu_N$ (a, *I*) = *I* is called as the standard neutrosophic union.
$\mu_N$ (a, b) = a + b – ab and
$\mu_N$ (a, *I*) = *I* .

This function will be called as the neutrosophic algebraic sum.

$$\mu_N (a, b) = \min (1, a + b) \text{ and } \mu_N (a, I) = I$$

will be called as the neutrosophic bounded sum. We define the notion of neutrosophic drastic unions

$$\mu_N (a, b) = \begin{cases} \text{a when b} = 0 \\ \text{b when a} = 0 \\ \text{I when a} = \text{I} \\ \quad \text{or b} = \text{I} \\ 1 \ \text{ otherwise.} \end{cases}$$

Now we proceed on to define the notion of neutrosophic Aggregation operators. Neutrosophic aggregation operators on neutrosophic sets are operations by which several neutrosophic sets are combined in a desirable way to produce a single neutrosophic set.

Any neutrosophic aggregation operation on n neutrosophic sets ($n \geq 2$) is defined by a function $h_N$: $N^n \to N$ where N = [0, 1] $\cup$ *I* and $N^n = \underbrace{N \times \ldots \times N}_{n-\text{times}}$ when applied to neutrosophic sets $A_1, A_2, \ldots, A_n$ defined on X the function $h_N$ produces an aggregate neutrosophic set A by operating on the membership grades of these sets for each $x \in X$ (Here also by the term membership grades we mean not only the membership grades from the unit interval [0, 1] but also the indeterminacy *I* for some $x \in X$ are included). Thus

$$A_N (x) = h_N (A_1 (x), A_2 (x), \ldots, A_n(x))$$

for each $x \in X$.

## 2.3 Special Fuzzy Cognitive Models and their Neutrosophic Analogue

In this section we define five types of special fuzzy models and their neutrosophic analogue. This is the first time such models are defined. They are different from the usual combined fuzzy models like combined fuzzy cognitive maps, combined fuzzy relational maps and so on. These models helps not only in easy comparison also it gives a equal opportunity to study the opinion of several experts that is why these models can be thought of as multi expert models preserving the individual opinion even after calculations of the resultant. Thus all special fuzzy models are multi expert models as well as a multi models for the same model can at the same time use a maximum of four different models in the study. Thus this is a special feature of these special fuzzy models. We also give the special neutrosophic analogue of them. For in many problem we may not have a clear cut feeling i.e., the expert may not be in a position to give his /her opinion it can also be an indeterminate. When the situation of indeterminacy arises fuzzy models cannot play any role only the neutrosophic models can tackle the situation. We also build a special mixed models which will be both having neutrosophic models as well as fuzzy models in their dynamical system.

In this section we define for the first time a new type of Special Fuzzy Cognitive Models (SFCM) and their neutrosophic analogue which is defined as Special Neutrosophic Cognitive Models (SNCM). Further we define Special Fuzzy and Neutrosophic Cognitive Maps (models) (SFNCM).

We illustrate them with examples from real world problems. It is pertinent to mention here that we give only illustrative examples not always any specific study. Now we proceed on to define the notion of special fuzzy cognitive maps models.

**DEFINITION 2.3.1:** *Suppose we have some m experts working on a problem P and suppose all of them agree to work with the same set of attributes say n attributes using only the Fuzzy Cognitive Maps then this gives a multiexpert model called the Special Fuzzy Cognitive Model, and if* $M_i$ *denotes the fuzzy* *connection matrix given by the $i^{th}$ expert using the set of n attributes i = 1, 2, …, m then we call the special fuzzy square matrix $M = M_1 \cup M_2 \cup \ldots \cup M_m$ to be the special fuzzy connection matrix associated with the Special Fuzzy Cognitive Map (model) (SFCM).*

*Let $C_1^1\, C_2^1 \ldots C_n^1, C_1^2\, C_2^2 \ldots C_n^2, \ldots, C_1^m\, C_2^m \ldots C_n^m$ be the special nodes of the SFCM*

$$A = \left[a_1^1\ a_2^1\ \ldots\ a_n^1\right] \cup \left[a_1^2\ a_2^2\ \ldots\ a_n^2\right] \cup \ldots \cup \left[a_1^m\ a_2^m\ \ldots\ a_n^m\right]$$

*where $a_j^i \in \{0, 1\}$, A is called the instantaneous special state vector and denotes the on-off position of the special node at an instant*

$$a_j^i = 0 \text{ if } a_j^i \text{ is off and}$$
$$= 1 \text{ if } a_j^i \text{ is on } I = 1, 2, \ldots, m \text{ and } j = 1, 2, \ldots, n.$$

*The attributes or nodes of a SFCM will be known as the special nodes or special attributes; context wise one may use just nodes or attributes.*

The main advantage of using this model is that it can work simultaneously at one stretch and give opinion of all the m experts. This is clearly a multi expert model.

Now we shall illustrate the functioning of this model. Suppose $M = M_1 \cup M_2 \cup \ldots \cup M_m$ be a $n \times n$ special fuzzy square matrix associated with the SFCM for the given problem P i.e. M is the multiexpert special fuzzy model. Suppose we want to study the effect of $X = X_1 \cup X_2 \cup \ldots \cup X_m$ where X is a special fuzzy row state vector where each $X_i$ is a $1 \times n$ fuzzy row state vector suggested by the $i^{th}$ expert $1 \leq i \leq m$ with entries from the set {0, 1} i.e. the special fuzzy row vector gives the on or off states of the attributes associated with the problem P that is under study. Now the effect of X on M using the special operator described in pages 20-1 of this book is given by

$$\begin{aligned} X \text{ o } M &= (X_1 \cup X_2 \cup \ldots \cup X_m) \text{ o } (M_1 \cup M_2 \cup \ldots \cup M_m) \\ &= X_1 \text{ o } M_1 \cup X_2 \text{ o } M_2 \cup \ldots \cup X_m \text{ o } M_m \\ &= Y'_1 \cup Y'_2 \cup \ldots \cup Y'_m \\ &= Y'; \end{aligned}$$

Y' may or may not be a special fuzzy row vector with entries from the set {0, 1}. It may also happen that the elements in $Y'_i$; $1 \leq i \leq m$ may not belong to the set {0, 1}. Further the nodes which were in the on state in X may not be in the on state in Y′. To overcome all these problems we update and threshold Y′ to $Y = Y_1 \cup Y_2 \cup \ldots \cup Y_m$. Now each $Y_i$ has its entries from the set {0, 1} and those nodes which were in the on state in X remain to be in the on state in Y also. Now we find

$$\begin{aligned} Y \text{ o } M &= (Y_1 \cup Y_2 \cup \ldots \cup Y_m) \text{ o } (M_1 \cup M_2 \cup \ldots \cup M_m) \\ &= Y_1 \text{ o } M_1 \cup Y_2 \text{ o } M_2 \cup \ldots \cup Y_m \text{ o } M_m \\ &= Z_1' \cup Z_2' \cup \ldots \cup Z_m' \\ &= Z'; \end{aligned}$$

now Z′ may not be a special fuzzy row vector so we threshold and update Z′ to $Z = Z_1 \cup Z_2 \cup \ldots \cup Z_m$ and now proceed on to find Z o M, we continue this process until we arrive at a special fixed point or a special limit cycle or a special fixed points and limit cycles. This final resultant special fuzzy row vector will be known as the special hidden pattern of the SFCM model.

We illustrate this by a simple real world model.

***Example 2.3.1:*** Suppose some 5 experts are interested to study the problem, the nation will face due to the production of more number engineering graduates which is very disproportionate to the job opportunities created by the nation (India) for them.

All the five experts wish to work with the five concepts relating the unemployed engineering graduates.

$E_1$ - Frustration
$E_2$ - Unemployment
$E_3$ - Increase of educated criminals
$E_4$ - Under employment
$E_5$ - Taking up drugs, alcohol etc.

The special fuzzy square connection matrix related with the SFCM model given by the 5 experts be denoted by

$$M = M_1 \cup M_2 \cup M_3 \cup M_4 \cup M_5$$

$$\begin{array}{c} \\ E_1 \\ E_2 \\ E_3 \\ E_4 \\ E_5 \end{array}\begin{array}{c} \begin{array}{ccccc} E_1 & E_2 & E_3 & E_4 & E_5 \end{array} \\ \begin{bmatrix} 0 & 1 & 1 & 1 & 0 \\ 1 & 0 & 1 & 0 & 1 \\ 1 & 1 & 0 & 1 & 0 \\ 1 & 0 & 0 & 0 & 1 \\ 1 & 1 & 0 & 0 & 0 \end{bmatrix} \end{array} \cup \begin{array}{c} \\ E_1 \\ E_2 \\ E_3 \\ E_4 \\ E_5 \end{array}\begin{array}{c} \begin{array}{ccccc} E_1 & E_2 & E_3 & E_4 & E_5 \end{array} \\ \begin{bmatrix} 0 & 1 & 0 & 1 & 0 \\ 1 & 0 & 1 & 0 & 1 \\ 0 & 1 & 0 & 1 & 0 \\ 1 & 0 & 0 & 0 & 1 \\ 0 & 1 & 1 & 0 & 0 \end{bmatrix} \end{array}$$

$$\cup \begin{array}{c} \\ E_1 \\ E_2 \\ E_3 \\ E_4 \\ E_5 \end{array}\begin{array}{c} \begin{array}{ccccc} E_1 & E_2 & E_3 & E_4 & E_5 \end{array} \\ \begin{bmatrix} 0 & 1 & 0 & 0 & 1 \\ 1 & 0 & 1 & 0 & 0 \\ 1 & 0 & 0 & 0 & 1 \\ 0 & 0 & 1 & 0 & 1 \\ 1 & 0 & 0 & 1 & 0 \end{bmatrix} \end{array} \cup$$

$$\begin{array}{c} \\ E_1 \\ E_2 \\ E_3 \\ E_4 \\ E_5 \end{array}\begin{array}{c} \begin{array}{ccccc} E_1 & E_2 & E_3 & E_4 & E_5 \end{array} \\ \begin{bmatrix} 0 & 0 & 1 & 0 & 1 \\ 0 & 0 & 0 & 1 & 1 \\ 0 & 1 & 0 & 0 & 1 \\ 0 & 1 & 1 & 0 & 1 \\ 1 & 0 & 1 & 0 & 0 \end{bmatrix} \end{array} \cup \begin{array}{c} \\ E_1 \\ E_2 \\ E_3 \\ E_4 \\ E_5 \end{array}\begin{array}{c} \begin{array}{ccccc} E_1 & E_2 & E_3 & E_4 & E_5 \end{array} \\ \begin{bmatrix} 0 & 0 & 0 & 1 & 1 \\ 1 & 0 & 0 & 0 & 1 \\ 0 & 0 & 0 & 1 & 1 \\ 1 & 1 & 0 & 0 & 0 \\ 0 & 0 & 1 & 1 & 0 \end{bmatrix} \end{array}.$$

Clearly M is a special fuzzy square matrix. Now we want to study the effect of unemployment alone in the on state by all the 5 experts and all other nodes are in, the off state i.e. X = [0 1 0 0 0] $\cup$ [0 1 0 0 0] $\cup$ [0 1 0 0 0] $\cup$ [0 1 0 0 0] $\cup$ [0 1 0 0 0]. X is a special fuzzy row state vector. To find the special hidden pattern associated with X using the SFCM. Now

$$\begin{aligned} X \text{ o } M &= (X_1 \cup X_2 \cup \ldots \cup X_5) \text{ o } (M_1 \cup M_2 \cup \ldots \cup M_5) \\ &= X_1 \text{ o } M_1 \cup X_2 \text{ o } M_2 \cup \ldots \cup X_5 \text{ o } M_5 \end{aligned}$$

$$= \begin{bmatrix}0 & 1 & 0 & 0 & 0\end{bmatrix} \text{ o } \begin{bmatrix}0 & 1 & 1 & 1 & 0\\ 1 & 0 & 1 & 0 & 1\\ 1 & 1 & 0 & 1 & 0\\ 1 & 0 & 0 & 0 & 1\\ 1 & 1 & 0 & 0 & 0\end{bmatrix} \cup$$

$$\begin{bmatrix}0 & 1 & 0 & 0 & 0\end{bmatrix} \text{ o } \begin{bmatrix}0 & 1 & 0 & 1 & 0\\ 1 & 0 & 1 & 0 & 1\\ 0 & 1 & 0 & 1 & 0\\ 1 & 0 & 0 & 0 & 1\\ 0 & 1 & 1 & 0 & 0\end{bmatrix} \cup$$

$$\begin{bmatrix}0 & 1 & 0 & 0 & 0\end{bmatrix} \text{ o } \begin{bmatrix}0 & 1 & 0 & 0 & 1\\ 1 & 0 & 1 & 0 & 0\\ 1 & 0 & 0 & 0 & 1\\ 0 & 0 & 1 & 0 & 1\\ 1 & 0 & 0 & 1 & 0\end{bmatrix} \cup$$

$$\begin{bmatrix}0 & 1 & 0 & 0 & 0\end{bmatrix} \text{ o } \begin{bmatrix}0 & 0 & 1 & 0 & 1\\ 0 & 0 & 0 & 1 & 1\\ 0 & 1 & 0 & 0 & 1\\ 0 & 1 & 1 & 0 & 1\\ 1 & 0 & 1 & 0 & 0\end{bmatrix} \cup$$

$$\begin{bmatrix}0 & 1 & 0 & 0 & 0\end{bmatrix} \text{ o } \begin{bmatrix}0 & 0 & 0 & 1 & 1\\ 1 & 0 & 0 & 0 & 1\\ 0 & 0 & 0 & 1 & 1\\ 1 & 1 & 0 & 0 & 0\\ 0 & 0 & 1 & 1 & 0\end{bmatrix}$$

$$= [1\ 0\ 1\ 0\ 1] \cup [1\ 0\ 1\ 0\ 1] \cup [1\ 0\ 1\ 0\ 0] \cup [0\ 0\ 0\ 1\ 1]$$

$\cup$ [1 0 0 0 1]

$= Y'_1 \cup Y'_2 \cup Y'_3 \cup Y'_4 \cup Y'_5$

$= Y'$;

we update and threshold Y' to

$$Y = Y_1 \cup Y_2 \cup Y_3 \cup Y_4 \cup Y_5$$
$$= [1\ 1\ 1\ 0\ 1] \cup [1\ 1\ 1\ 0\ 1] \cup [1\ 1\ 1\ 0\ 0] \cup [0\ 1\ 0\ 1\ 1] \cup [1\ 1\ 0\ 0\ 1].$$

Now we find the effect of Y on M i.e.

$$Y \text{ o } M = (Y_1 \cup Y_2 \cup \ldots \cup Y_5) \text{ o } (M_1 \cup M_2 \cup \ldots \cup M_5)$$
$$= Y_1 \text{ o } M_1 \cup Y_2 \text{ o } M_2 \cup \ldots \cup Y_5 \text{ o } M_5$$
$$= [3\ 3\ 2\ 2\ 1] \cup [1\ 3\ 2\ 2\ 1] \cup [2\ 1\ 1\ 0\ 2] \cup [1\ 1\ 2\ 1\ 2] \cup [1\ 0\ 1\ 2\ 2]$$
$$= Z'_1 \cup Z'_2 \cup Z'_3 \cup Z'_4 \cup Z'_5$$
$$= Z';$$

we update and threshold Z' to obtain

$$Z = Z_1 \cup Z_2 \cup Z_3 \cup Z_4 \cup Z_5$$
$$= [1\ 1\ 1\ 1\ 1] \cup [1\ 1\ 1\ 1\ 1] \cup [1\ 1\ 1\ 0\ 1] \cup [1\ 1\ 1\ 1\ 1] \cup [1\ 1\ 1\ 1\ 1].$$

Now we find the effect of Z on M i.e.

$$Z \text{ o } M = (Z_1 \cup Z_2 \cup \ldots \cup Z_5) \text{ o } (M_1 \cup M_2 \cup \ldots \cup M_5)$$
$$= Z_1 \text{ o } M_1 \cup Z_2 \text{ o } M_2 \cup \ldots \cup Z_5 \text{ o } M_5$$
$$= [4\ 3\ 2\ 2\ 2] \cup [2\ 3\ 2\ 2\ 2] \cup [3\ 1\ 1\ 1\ 2] \cup [1\ 2\ 3\ 1\ 4] \cup [2\ 1\ 1\ 3\ 3]$$
$$= P'_1 \cup P'_2 \cup P'_3 \cup P'_4 \cup P'_5$$
$$= P'.$$

Clearly P' is not a special fuzzy row vector; so update and threshold P' to obtain

$$P = P_1 \cup P_2 \cup \ldots \cup P_5$$
$$= [1\ 1\ 1\ 1\ 1] \cup [1\ 1\ 1\ 1\ 1] \cup [1\ 1\ 1\ 1\ 1] \cup [1\ 1\ 1\ 1\ 1] \cup [1\ 1\ 1\ 1\ 1]$$

Thus unemployment according to all experts leads to all other problems. All experts agree upon it. We can for instance find the effect of the special state vector like

$$X = [1\ 0\ 0\ 0\ 0] \cup [0\ 1\ 0\ 0\ 0] \cup [0\ 0\ 1\ 0\ 0] \cup [0\ 0\ 0\ 1\ 0] \cup [0\ 0\ 0\ 0\ 1]$$

where the first experts wishes to study on state of first node frustration, second expert the on state of the second node unemployment, the third expert the on state on the third node viz. increase of educated criminals and the forth expert under employments and taking up drugs alcohol etc. is taken up by the fifth expert.

$$\begin{aligned} X \text{ o } M &= (X_1 \cup X_2 \cup \ldots \cup X_5) \text{ o } (M_1 \cup M_2 \cup \ldots \cup M_5) \\ &= X_1 \text{ o } M \cup X_2 \text{ o } M_2 \cup \ldots \cup X_5 \text{ o } M_5 \\ &= [0\ 1\ 1\ 1\ 0] \cup [1\ 0\ 1\ 0\ 1] \cup [1\ 0\ 0\ 0\ 1] \cup [0\ 1\ 1\ 0\ 1] \cup [0\ 0\ 1\ 1\ 0] \\ &= Y'_1 \cup Y'_2 \cup Y'_3 \cup Y'_4 \cup Y'_5 \\ &= Y' \end{aligned}$$

after updating and thresholding Y' we get

$$\begin{aligned} Y &= Y_1 \cup Y_2 \cup \ldots \cup Y_5 \\ &= [1\ 1\ 1\ 1\ 0] \cup [1\ 1\ 1\ 0\ 1] \cup [1\ 0\ 1\ 0\ 1] \cup [0\ 1\ 1\ 1\ 1] \cup [0\ 0\ 1\ 1\ 1]. \end{aligned}$$

One can find Y o M and so on until one gets the special hidden pattern.

This we have just mentioned for the reader should to know the expert can take any special node or special nodes to be in the on state for the study and this gives a single solution of all the experts. The two main advantages over other multi expert systems are

1. Stage wise comparison is possible for the same state special vector by all experts or stage wise comparison for different state special vector by different experts.
2. The working is not laborious as a C-program will do the job in no time.

Next we proceed on to define the notion of Special Mixed Fuzzy Cognitive maps/models (SMFCMs). This model comes in handy when experts have different sets of attributes and the number of attributes are also varying in number.

**DEFINITION 2.3.2:** *Suppose we have some t experts working on a problem P and suppose they do not agree to work on the same set of attributes or the same attributes, they choose different sets of attributes but they all agree to work with FCM model then we can use the Special Mixed Fuzzy Cognitive Maps (SMFCMs) model and get a new dynamical system which can cater to the opinion of each and every expert simultaneously. Suppose the* $i^{th}$ *expert works with* $n_i$ *attributes and the fuzzy connection matrix associated with the FCM be given by* $M_i$, $i = 1, 2, \ldots, t$, *we use these t experts to get a special fuzzy mixed square matrix* $M = M_1 \cup M_2 \cup \ldots \cup M_t$ *where M corresponds to the dynamical system which gives the opinion of all t experts and M is called the special fuzzy connection matrix and the model related with this M will be known as the Special Mixed Fuzzy Cognitive Maps (SMFCMs) model associated with t experts.*

Thus we see M is also a multi expert system model. It is important to mention here that SMFCM the multi expert FCM model is different from the Combined FCM (CFCM[108, 112]) as well as Special Fuzzy Cognitive Maps (SFCM) model. The salient features of SMFCM model is

1. SMFCM is better than the CFCM model as the CFCM model can function only with the same set of attributes for the problem P.
2. Here in the SMFCM contradictory opinions do not cancel out as in case of CFCMs.
3. In the SMFCM the opinion of each and every expert is under display and see their opinion not as CFCMs which gives only a collective resultant.
4. SMFCM model is better than the SFCM model as SFCM can permit all the expert to work only with the same set of attributes for the given problem.

We give an example of SMFCMs as a model for the reader to understand this new model.

***Example 2.3.2:*** We here give the prediction of electoral winner or how the preference of a particular politician and so on. Suppose we have four experts working on the problem and the first expert wishes to work with the four attributes.

$e_1^1$ - Language of the politician

$e_2^1$ - Community of the politician

$e_3^1$ - Service to people, public figure configuration, personality and nature

$e_4^1$ - Party's strength and the opponents strength.

Let $M_1$ be the fuzzy connection matrix given by the first expert, i.e.,

$$M_1 = \begin{array}{c} \\ e_1^1 \\ e_2^1 \\ e_3^1 \\ e_4^1 \end{array} \begin{array}{c} \begin{array}{cccc} e_1^1 & e_2^1 & e_3^1 & e_4^1 \end{array} \\ \begin{bmatrix} 0 & 1 & 1 & 0 \\ 1 & 0 & 1 & 0 \\ 0 & 0 & 0 & 1 \\ 1 & 1 & 1 & 0 \end{bmatrix} \end{array}.$$

Suppose the second expert wishes to work with 5 attributes given by

$e_1^2$ - Service done by the politician to people

$e_2^2$ - Finance and media accessibility

$e_3^2$ - Party's strength and opponents strength

$e_4^2$ - Working member of the party

$e_5^2$ - His community and the locals community.

Let $M_2$ be the fuzzy connection matrix associated with his opinion.

$$M_2 = \begin{array}{c} \\ e_1^2 \\ e_2^2 \\ e_3^2 \\ e_4^2 \\ e_5^2 \end{array} \begin{array}{c} \begin{array}{ccccc} e_1^2 & e_2^2 & e_3^2 & e_4^2 & e_5^2 \end{array} \\ \begin{bmatrix} 0 & 1 & 1 & 0 & 0 \\ 1 & 0 & 1 & 0 & 0 \\ 0 & 0 & 0 & 1 & 1 \\ 1 & 1 & 1 & 0 & 1 \\ 0 & 0 & 1 & 1 & 0 \end{bmatrix} \end{array}.$$

Suppose the third expert wishes to work with only four attributes given by

- $e_1^3$ - Community and nativity of the politician
- $e_2^3$ - His interest for working with people and their opinion on his personality
- $e_3^3$ - Amount of money he can spend on propaganda in media
- $e_4^3$ - Working members of the party i.e., public mass support.

Let $M_3$ be the connection matrix given by the third expert using FCM model.

$$M_3 = \begin{array}{c} \\ e_1^3 \\ e_2^3 \\ e_3^3 \\ e_4^3 \end{array} \begin{array}{c} \begin{array}{cccc} e_1^3 & e_2^3 & e_3^3 & e_4^3 \end{array} \\ \begin{bmatrix} 0 & 1 & 1 & 0 \\ 1 & 0 & 1 & 0 \\ 1 & 0 & 0 & 1 \\ 1 & 1 & 0 & 0 \end{bmatrix} \end{array}.$$

Let $e_1^4, e_2^4, \ldots, e_6^4$ be the attributes given by the fourth expert.

- $e_1^4$ - Community, nativity and gender of the politician
- $e_2^4$ - The constructive and progressive work done by him in his locality

| | | |
|---|---|---|
| $e_3^4$ | - | His social and economic status |
| $e_4^4$ | - | Support of locals in propagating for his winning the election |
| $e_5^4$ | - | Money he spends on propaganda for election |
| $e_6^4$ | - | The strength of the party's campaigning in favour of him. |

Let $M_4$ denote connection matrix given by the fourth expert.

$$M_4 = \begin{array}{c} \\ e_1^4 \\ e_2^4 \\ e_3^4 \\ e_4^4 \\ e_5^4 \\ e_6^4 \end{array} \begin{array}{c} \begin{array}{cccccc} e_1^4 & e_2^4 & e_3^4 & e_4^4 & e_5^4 & e_6^4 \end{array} \\ \begin{bmatrix} 0 & 1 & 1 & 1 & 0 & 1 \\ 1 & 0 & 0 & 1 & 0 & 1 \\ 0 & 0 & 0 & 0 & 1 & 1 \\ 1 & 1 & 0 & 0 & 0 & 1 \\ 0 & 0 & 1 & 1 & 0 & 1 \\ 0 & 1 & 1 & 1 & 1 & 0 \end{bmatrix} \end{array}.$$

Now let $M = M_1 \cup M_2 \cup M_3 \cup M_4$ give the SMFCM which is the multi expert opinion of the four experts. Now suppose the four experts wish to work on the special state vector

$$X = [1\ 0\ 0\ 0] \cup [0\ 1\ 0\ 0\ 0] \cup [1\ 0\ 0\ 0] \cup [1\ 0\ 0\ 0\ 0\ 0]$$

where the first expert wants to work with language of the politician alone to be in the on state and the second expert wants to work with the node finance and media accessibility alone in the on state and all the other nodes to be in the off state; the third expert wants to work with the node community and nativity alone in the on state and all other nodes to be in the off state and the forth expert wants to work with the node community, nativity and the gender of the politician alone in the on state. To find the special hidden pattern of X on the special dynamical system,

$$M \quad = \quad M_1 \cup M_2 \cup M_3 \cup M_4$$

$$= \begin{bmatrix} 0 & 1 & 1 & 0 \\ 1 & 0 & 1 & 0 \\ 0 & 0 & 0 & 1 \\ 1 & 1 & 1 & 0 \end{bmatrix} \cup \begin{bmatrix} 0 & 1 & 1 & 0 & 0 \\ 1 & 0 & 1 & 0 & 0 \\ 0 & 0 & 0 & 1 & 1 \\ 1 & 1 & 1 & 0 & 1 \\ 0 & 0 & 1 & 1 & 0 \end{bmatrix} \cup$$

$$\begin{bmatrix} 0 & 1 & 1 & 0 \\ 1 & 0 & 1 & 1 \\ 1 & 0 & 0 & 1 \\ 1 & 1 & 0 & 0 \end{bmatrix} \cup \begin{bmatrix} 0 & 1 & 1 & 1 & 0 & 1 \\ 1 & 0 & 0 & 1 & 0 & 1 \\ 0 & 0 & 0 & 0 & 1 & 1 \\ 1 & 1 & 0 & 0 & 0 & 1 \\ 0 & 0 & 1 & 1 & 0 & 1 \\ 0 & 1 & 1 & 1 & 1 & 0 \end{bmatrix}.$$

Now

$X \text{ o } M = X_1 \text{ o } M_1 \cup X_2 \text{ o } M_2 \cup X_3 \text{ o } M_3 \cup X_4 \text{ o } M_4$

$$= \begin{bmatrix} 1 & 0 & 0 & 0 \end{bmatrix} \circ \begin{bmatrix} 0 & 1 & 1 & 0 \\ 1 & 0 & 1 & 0 \\ 0 & 0 & 0 & 1 \\ 1 & 1 & 1 & 0 \end{bmatrix}$$

$$\cup \begin{bmatrix} 0 & 1 & 0 & 0 & 0 \end{bmatrix} \circ \begin{bmatrix} 0 & 1 & 1 & 0 & 0 \\ 1 & 0 & 1 & 0 & 0 \\ 0 & 0 & 0 & 1 & 1 \\ 1 & 1 & 1 & 0 & 1 \\ 0 & 0 & 1 & 1 & 0 \end{bmatrix}$$

$$\cup \begin{bmatrix} 1 & 0 & 0 & 0 \end{bmatrix} \circ \begin{bmatrix} 0 & 1 & 1 & 0 \\ 1 & 0 & 1 & 0 \\ 1 & 0 & 0 & 1 \\ 1 & 1 & 0 & 0 \end{bmatrix}$$

$$\cup \begin{bmatrix}1 & 0 & 0 & 0 & 0 & 0\end{bmatrix} \circ \begin{bmatrix} 0 & 1 & 1 & 1 & 0 & 1 \\ 1 & 0 & 0 & 1 & 0 & 1 \\ 0 & 0 & 0 & 0 & 1 & 1 \\ 1 & 1 & 0 & 0 & 0 & 1 \\ 0 & 0 & 1 & 1 & 0 & 1 \\ 0 & 1 & 1 & 1 & 1 & 0 \end{bmatrix}$$

$$\begin{aligned} &= [0\ 1\ 1\ 0] \cup [1\ 0\ 1\ 0\ 0] \cup [0\ 1\ 1\ 0] \cup [0\ 1\ 1\ 1\ 0\ 1] \\ &= Y'_1 \cup Y'_2 \cup Y'_3 \cup Y'_4 \\ &= Y'. \end{aligned}$$

After updating and thresholding Y' we get

$$\begin{aligned} Y &= Y_1 \cup Y_2 \cup Y_3 \cup Y_4 \\ &= [1\ 1\ 1\ 0] \cup [1\ 1\ 1\ 0\ 0] \cup [1\ 1\ 1\ 0] \cup [1\ 1\ 1\ 1\ 0\ 1] \\ &= Y. \end{aligned}$$

Clearly Y is a special fuzzy row vector.

Now we find the effect of Y on the special dynamical system M and get Z'.

$$\begin{aligned} Y \text{ o } M &= [1\ 1\ 2\ 1] \cup [1\ 1\ 2\ 1\ 1] \cup [2\ 1\ 2\ 1] \cup [2\ 3\ 2\ 3\ 2\ 4] \\ &= Z'_1 \cup Z'_2 \cup Z'_3 \cup Z'_4 \cup Z'_5 \\ &= Z'. \end{aligned}$$

Let Z be the resultant vector got by updating and thresholding Z'. i.e.,

$$Z = [1\ 1\ 1\ 1] \cup [1\ 1\ 1\ 1\ 1] \cup [1\ 1\ 1\ 1] \cup [1\ 1\ 1\ 1\ 1\ 1].$$

It is left for the reader to verify the special hidden pattern in a fixed point given by the special fuzzy row vector

$$Z = [1\ 1\ 1\ 1] \cup [1\ 1\ 1\ 1\ 1] \cup [1\ 1\ 1\ 1] \cup [1\ 1\ 1\ 1\ 1\ 1].$$

Now we proceed on to define the new special neutrosophic fuzzy cognitive model in which multi experts opinion are given.

**DEFINITION 2.3.3:** *Suppose we have some n experts interested in working with a problem P and all of them agree to work with the same m number of attributes related with the problem and further all of them agree to work with the neutrosophic cognitive maps(NCMs) model; i.e., all experts feel that some of the interrelation between the attributes is an indeterminate. Suppose let the* $i^{th}$ *expert using the NCMs give the neutrosophic connection matrix to be the* $N_i$, *i = 1, 2, …, n. Thus each* $N_i$ *is a* $m \times m$ *neutrosophic matrix and the multiexpert Special Neutrosophic Cognitive Maps (SNCMs) model is denoted by N, which is a special fuzzy neutrosophic square matrix and is given by* $N = N_1 \cup N_2 \cup \ldots \cup N_n$. *Thus N is a special neutrosophic connection square matrix associated with the SNCMs model.*

*We just mention how it works. Suppose all the n experts give their preference of node which is to be in the on state then the special fuzzy row vector given by them collectively be denoted by* $X = X_1 \cup X_2 \cup \ldots \cup X_n$ *where* $X_i$ *is the special state vector given by the* $i^{th}$ *expert; i = 1, 2, …, n. Now the effect of X on the special fuzzy dynamical system N is given by*

$$\begin{aligned} X \text{ o } N &= (X_1 \cup X_2 \cup \ldots \cup X_n) \text{ o } (N_1 \cup N_2 \cup \ldots \cup N_n) \\ &= X_1 \text{ o } N_1 \cup X_2 \text{ o } N_2 \cup \ldots \cup X_n \text{ o } N_n \\ &= Y'_1 \cup Y'_2 \cup \ldots \cup Y'_n \\ &= Y'. \end{aligned}$$

*The working is described in page 139-141 of chapter 1. Now after updating and thresholding Y' we get* $Y = Y_1 \cup Y_2 \cup \ldots \cup Y_n$. *Now we find the effect of Y on N i.e.,*

$$\begin{aligned} Y \text{ o } N &= (Y_1 \cup Y_2 \cup \ldots \cup Y_n) \text{ o } (N_1 \cup N_2 \cup \ldots \cup N_n) \\ &= Y_1 \text{ o } N_1 \cup Y_2 \text{ o } N_2 \cup \ldots \cup Y_n \text{ o } N_n \\ &= Z'_1 \cup Z'_2 \cup \ldots Z'_n \\ &= Z'. \end{aligned}$$

*We see Z' may or may not be a special fuzzy neutrosophic row vector and further the nodes which we started to work with may not be in the on state so we update and threshold Z' to Z and let* $Z = Z_1 \cup Z_2 \cup \ldots \cup Z_n$. *Clearly Z is a special fuzzy neutrosophic row vector. We now find Z o N and so on till we arrive at a fixed point or a limit cycle (say) T. This T will be known as the*

*special hidden pattern of the special dynamical neutrosophic system N.*

Now we see this new Special Neutrosophic Cognitive Maps (SNCMs) model has the following advantages.

1. When all the experts want to make use of the element of indeterminacy it serves the purpose.
2. This is also a multi expert model.

Now the reader can construct real world problems multi expert SNCM model to study the results as the working is analogous to the SFCM model described in definition 2.3.1 of this chapter.

Now we proceed on to define yet another new multi expert special neutrosophic model.

**DEFINITION 2.3.4:** *Let us take n experts who wishes to work with the same problem P. But the experts have different sets of attributes which they want to work with, but all of them agree to work with the neutrosophic cognitive maps model i.e., they all uniformly agree upon the fact that they should choose their own attributes which involve some degrees of indeterminacy. Now we give a multi expert model to solve this problem. Let the $i^{th}$ expert choose to work with $n_i$ attributes using a NCM model and let $N_i$ be the neutrosophic connection matrix given by him; this is true i = 1, 2, …, m. Let us take $N = N_1 \cup N_2 \cup \ldots \cup N_m$. Clearly N is a special fuzzy neutrosophic mixed square matrix. This N will be known as the special dynamical system associated with the problem and this special dynamical system is defined as the Special Mixed Neutrosophic Cognitive Maps (SMNCMs ) model.*

*We must give the functioning of the SMNCMs model. Let N be the special fuzzy neutrosophic connection matrix associated with the SMNCMs model. Let each of the expert(say $i^{th}$ expert) give $X_i$ to be the state row vector with which he wishes to work; for i = 1, 2, …, m. Then the special fuzzy mixed row vector denoted by $X = X_1 \cup X_2 \cup \ldots \cup X_m$ will be the integrated special fuzzy mixed row vector of all the experts whose effect we*

*want to study and find its resultant on the special neutrosophic dynamical system N. Now*

$$X \text{ o } N = (X_1 \cup X_2 \cup \ldots \cup X_m) \text{ o } (N_1 \cup N_2 \cup \ldots \cup N_m)$$
$$= X_1 \text{ o } N_1 \cup X_2 \text{ o } N_2 \cup \ldots \cup X_m \text{ o } N_m.$$

*(This operation is defined in chapter one, pages 146-151 of this book). Let*

$$X \text{ o } N = Y'_1 \cup Y'_2 \cup \ldots \cup Y'_m$$
$$= Y'.$$

*$Y'$ may or may not be a special fuzzy neutrosophic mixed row vector so we update and threshold $Y'$ and obtain $Y = Y_1 \cup Y_2 \cup \ldots \cup Y_n$. Now we find the effect of Y on N.*

$$X \text{ o } N = (Y_1 \cup Y_2 \cup \ldots \cup Y_m) \text{ o } (N_1 \cup N_2 \cup \ldots \cup N_m)$$
$$= Y_1 \text{ o } N_1 \cup Y_2 \text{ o } N_2 \cup \ldots \cup Y_m \text{ o } N_m$$
$$= Z'_1 \cup Z'_2 \cup \ldots \cup Z'_m$$
$$= Z' \text{ (say);}$$

*$Z'$ may or may not be a special fuzzy neutrosophic mixed row vector, we update and threshold $Z'$ to obtain $Z = Z_1 \cup Z_2 \cup \ldots \cup Z_n$. We find Z o N and so on till we obtain a fixed point or a limit cycle. This fixed point or the limit cycle will be known as the special fixed point or the special limit cycle and is defined as the special hidden pattern of the special dynamical system SMNCM.*

The main advantage of this multi expert model SMNCMs are :

1. One can obtain simultaneously every expert opinion and compare them stage by stage also.
2. SMNCMs is the dynamical system which gives the special hidden pattern.
3. SMNCMs is best suited when the relations between attributes involves indeterminacy.
4. SMNCMs gives the liberty to the experts to choose any desired number of attributes which he/she chooses to work with.

5. SMNCMs are better than SNCMs, SFCMs and SMFCMs when the problem involves indeterminacy and gives liberty to the experts to use any desired number of attributes to work with the same problem.

Now the reader is expected to apply SMNCM model in any of the multi expert problem and find the special hidden pattern. The working is exactly analogous to the example 2.3.2 given in this book.

Now we proceed on to define a new multi expert model which is different from these four models known as Special Fuzzy Neutrosophic Cognitive Maps model (SFNCMs-model).

**DEFINITION 2.3.5:** *Suppose we have a problem P for which we want to obtain a multi experts opinion. Suppose we have two sets of experts and one set of experts are interested in only using Fuzzy Cognitive Maps models and another set of experts keen on using Neutrosophic Cognitive Maps models for the same problem P and experts from both the sets demand for different sets of attributes to be fixed by them for the problem P. Suppose we have $t_1$ experts who want to use FCM model with different sets of attributes then, we use SMFCMs to model them say $M^1$ gives the special connection matrix associated with SMFCMs of the $t_1$ experts, the $t_2$ experts who choose to work with NCM model with varying attributes is modelled using the SMNCMs model. Let $M^2$ give the special connection matrix of the SMNCMs of the $t_2$ experts. Now we define $M^1 \cup M^2$ to be the combined special connection matrix of both the models and*

$$M^1 \cup M^2 = (M_1^1 \cup M_2^1 \cup \ldots \cup M_{t_1}^1) \cup (M_1^2 \cup M_2^2 \cup \ldots \cup M_{t_2}^2)$$

*to be the associated special connection matrix of the new Special Fuzzy Neutrosophic Cognitive Maps (SFNCMs) model. $M^1 \cup M^2$ will be known as the special combined fuzzy neutrosophic mixed square matrix.*

*We just give the working of the SFNCMs model. Suppose $t_1$ experts give their on state of their nodes with which they wish to work as $X^1 = X_1^1 \cup X_2^1 \cup \ldots \cup X_{t_1}^1$ and $t_2$ experts give their preference nodes as $X^2 = X_1^2 \cup X_2^2 \cup \ldots \cup X_{t_2}^2$ then the* *special state vector for which we have to find the special hidden pattern is given by*

$X^1 \cup X^2 \quad = \quad (X_1^2 \cup X_2^2 \cup \ldots \cup X_{t_1}^2)$

$\qquad\qquad (X_1^2 \cup X_2^2 \cup \ldots \cup X_{t_2}^2)$

*using the SFNCMs dynamical system* $M^1 \cup M^2$*. The effect of* $X^1 \cup X^2$ *on* $M^1 \cup M^2$ *is given by*

$(X^1 \cup X^2)\ o\ (M^1 \cup M^2)$

$= \quad (X^1 \cup X^2)\ o\ (M^1 \cup M^2)$

$= \quad X^1\ o\ M^1 \cup X^2\ o\ M^2$

$= \quad (X_2^1 \cup X_2^1 \cup \ldots \cup X_{t_1}^1)\ o\ (M_1^1 \cup M_2^1 \cup \ldots \cup M_{t_1}^1) \cup$

$\quad (X_1^2 \cup X_2^2 \cup \ldots \cup X_{t_2}^2)\ o\ (M_1^2 \cup M_2^2 \cup \ldots \cup M_{t_2}^2)$

$= \quad X_1^1\ o\ M_1^1 \cup X_2^1\ o\ M_2^1 \cup X_3^1\ o\ M_3^1 \cup \ldots \cup$

$\quad X_{t_1}^1\ o\ M_{t_1}^1 \cup X_1^2\ o\ M_1^2 \cup X_2^2\ o\ M_2^2 \cup X_3^2\ o\ M_3^2$

$\quad \cup \ldots \cup X_{t_2}^2\ o\ M_{t_2}^2$

$= \quad (Z_1^2 \cup Z_2^2 \cup \ldots \cup Z_{t_1}^2) \cup (T_1^2 \cup T_2^2 \cup \ldots \cup T_{t_2}^2)$

$= \quad Z' \cup T'.$

1. $Z' \cup T'$ *may or may not be a special fuzzy neutrosophic mixed row vector*
2. $Z' \cup T'$ *may have the nodes with which started to work with to be in the off state so* $Z' \cup T'$ *is thresholded and updated to (say)* $P^1 \cup Q^2 = (P_1^1 \cup P_2^1 \cup \ldots \cup P_{t_1}^1) \cup (Q_1^2 \cup Q_2^2 \cup \ldots \cup Q_{t_2}^2)$.

*Now we find the effect of* $P^1 \cup Q^2$ *on the special dynamical system* $M^1 \cup M^2$ *i.e.,*

$(P^1 \cup Q^2)\ o\ (M^1 \cup M^2) = \quad P^1\ o\ M^1 \cup Q^2\ o\ M^2$

$= \quad (P_1^1 \cup P_2^1 \cup \ldots \cup P_{t_1}^1)\ o\ (M_1^1 \cup M_2^1 \cup \ldots \cup M_{t_1}^1) \cup$

$\quad (P_1^2 \cup P_2^2 \cup \ldots \cup P_{t_2}^2)\ o\ (M_1^2 \cup M_2^2 \cup \ldots \cup M_{t_2}^2)$

$= \quad (V_1^1 \cup V_2^1 \cup \ldots \cup V_{t_1}^1) \cup (W_1^2 \cup W_2^2 \cup \ldots \cup W_{t_2}^2)$

$= \quad V_2 \cup W_2$.

*This resultant may not be even a special fuzzy neutrosophic mixed row vector so we update and threshold $V_2 \cup W_2$ to $S^1 \cup S^2 = (S_1^1 \cup S_2^1 \cup \ldots \cup S_{t_1}^1) \cup (S_1^2 \cup S_2^2 \cup \ldots \cup S_{t_2}^2)$ we proceed on to find $(S^1 \cup S^2)$ o $(M^1 \cup M^2)$ so on until we arrive at a special fixed point or a special limit cycle. This resultant vector will correspond to the special hidden pattern of the special dynamical system SFNFCM.*

Now we will illustrate this situation by an explicit model.

***Example 2.3.3:*** Let the four experts work with the problem of finding prediction of electoral winner using SFNFCM. Suppose the first expert works with four attributes [for the attributes refer example 2.3.2] and the fuzzy connection matrix given by him is

$$M_1 = \begin{bmatrix} 0 & 1 & 1 & 0 \\ 1 & 0 & 1 & 0 \\ 0 & 0 & 0 & 1 \\ 1 & 1 & 1 & 0 \end{bmatrix}.$$

The second expert uses the NCM and the connection matrix given by him using the 5 attributes given in the example 2.3.2 is given below

$$M_2 = \begin{bmatrix} 0 & 1 & 1 & 0 & 0 \\ 1 & 0 & 1 & 1 & 0 \\ I & 0 & 0 & 0 & 1 \\ 1 & 0 & I & 0 & 0 \\ 0 & 0 & 1 & 1 & 0 \end{bmatrix},$$

which is neutrosophic model of the NCM. The third expert works with four attributes and gives the following NCM model which is a neutrosophic 4 × 4 matrix given by $M_3$; i.e.,

$$M_3 = \begin{bmatrix} 0 & 0 & 1 & 0 \\ 1 & 1 & I & 0 \\ 0 & 0 & 0 & 1 \\ I & 1 & 0 & 0 \end{bmatrix}$$

and the fourth expert using the six attributes given in the example 2.3.2 gives the following neutrosophic connection matrix

$$M_4 = \begin{bmatrix} 0 & 1 & 1 & 0 & 0 & 0 \\ 1 & 0 & 0 & 0 & 0 & 0 \\ 0 & 0 & 0 & 0 & 1 & 1 \\ 1 & 1 & 0 & 0 & 0 & 0 \\ 0 & 0 & 1 & 1 & 0 & 0 \\ 0 & 0 & 0 & 1 & 1 & 0 \end{bmatrix}.$$

The special fuzzy neutrosophic mixed square matrix associated with the SFNFCM model is given by $M_1 \cup M_2 \cup M_3 \cup M_4$

$$= \begin{bmatrix} 0 & 1 & 1 & 0 \\ 1 & 0 & 1 & 0 \\ 0 & 0 & 0 & 1 \\ 1 & 1 & 1 & 0 \end{bmatrix} \cup \begin{bmatrix} 0 & 1 & 1 & 0 & 0 \\ 1 & 0 & 1 & 1 & 0 \\ I & 0 & 0 & 0 & 1 \\ 1 & 0 & I & 0 & 0 \\ 0 & 0 & 1 & 1 & 0 \end{bmatrix} \cup$$

$$\begin{bmatrix} 0 & 0 & 1 & 0 \\ 1 & 1 & I & 0 \\ 0 & 0 & 0 & 1 \\ I & 1 & 0 & 0 \end{bmatrix} \cup \begin{bmatrix} 0 & 1 & 1 & 0 & 0 & 0 \\ 1 & 0 & 0 & 0 & 0 & 0 \\ 0 & 0 & 0 & 0 & 1 & 1 \\ 1 & 1 & 0 & 0 & 0 & 0 \\ 0 & 0 & 1 & 1 & 0 & 0 \\ 0 & 0 & 0 & 1 & 1 & 0 \end{bmatrix}.$$

Suppose the experts wish to work with special fuzzy mixed row input vector

$$X = X_1 \cup X_2 \cup X_3 \cup X_4$$
$$= [1\ 0\ 0\ 0] \cup [0\ 0\ 0\ 1\ 0] \cup [0\ 1\ 0\ 0] \cup [0\ 0\ 0\ 0\ 0\ 0\ 1].$$

Now the effect of the special input vector X on M is given by

$$X \text{ o } M = (X_1 \cup X_2 \cup X_3 \cup X_4) \text{ o } (M_1 \cup M_2 \cup M_3 \cup M_4)$$
$$= X_1 \text{ o } M_1 \cup X_2 \text{ o } M_2 \cup X_3 \text{ o } M_3 \cup X_4 \text{ o } M_4$$

$$= \begin{bmatrix} 1 & 0 & 0 & 0 \end{bmatrix} \text{ o } \begin{bmatrix} 0 & 1 & 1 & 0 \\ 1 & 0 & 1 & 0 \\ 0 & 0 & 0 & 1 \\ 1 & 1 & 1 & 0 \end{bmatrix} \cup$$

$$\begin{bmatrix} 0 & 0 & 0 & 1 & 0 \end{bmatrix} \text{ o } \begin{bmatrix} 0 & 1 & 1 & 0 & 0 \\ 1 & 0 & 1 & 1 & 0 \\ I & 0 & 0 & 0 & 1 \\ 1 & 0 & I & 0 & 0 \\ 0 & 0 & 1 & 1 & 0 \end{bmatrix} \cup$$

$$\begin{bmatrix} 0 & 1 & 0 & 0 \end{bmatrix} \text{ o } \begin{bmatrix} 0 & 0 & 1 & 0 \\ 1 & 0 & I & 0 \\ 0 & 0 & 0 & 1 \\ I & 1 & 0 & 0 \end{bmatrix} \cup$$

$$\begin{bmatrix} 0 & 0 & 0 & 0 & 0 & 1 \end{bmatrix} \text{ o } \begin{bmatrix} 0 & 1 & 1 & 0 & 0 & 0 \\ 1 & 0 & 0 & 0 & 0 & 0 \\ 0 & 0 & 0 & 0 & 1 & 1 \\ 1 & 1 & 0 & 0 & 0 & 0 \\ 0 & 0 & 1 & 1 & 0 & 0 \\ 0 & 0 & 0 & 1 & 1 & 0 \end{bmatrix}$$

$= [0\ 1\ 1\ 0] \cup [1\ 0\ I\ 0\ 0\ ] \cup [1\ 0\ I\ 0] \cup [0\ 0\ 0\ 1\ 1\ 0]$ .

This is updated and thresholded and we obtain

$$Y = Y_1 \cup Y_2 \cup Y_3 \cup Y_4$$
$$= [1\ 1\ 1\ 0] \cup [1\ 0\ I\ 1\ 0\ ] \cup [1\ 1\ I\ 0] \cup [0\ 0\ 0\ 1\ 1\ 0] .$$

We can find the effect of Y on the dynamical system M. The reader is requested to find Y o M and the special hidden pattern.

## 2.4 Special FRMs and NRMs and their generalizations

In this section we define for the first time the new notion of Special Fuzzy Relational Maps (SFRM) and illustrate it by examples. Further we define the new multi expert model known as special mixed model known as special mixed fuzzy relational maps, we also define the neutrosophic analogue for these two models. Finally we define a new Special Fuzzy Neutrosophic Relational Maps (SFNRMs). SFNRMS can be thought of as a generalization of FCM, SFCMs, SNCMs, SMFCMs, SMNCMs and SFNCMs. These models can be applied when the attributes related with the problem can be divided into two different disjoint classes. This is more time saving and space economic, so we seek to define this model.

**DEFINITION 2.4.1:** *Suppose n experts want to work with a problem P and if all of them wish to work using the fuzzy relational maps model with m attributes in the domain space and with t attributes in the range space The* $i^{th}$ *expert gives his opinion which is rechanged into the connection matrix of the FRM as* $M_i$*; i = 1, 2, …, n is a m × t fuzzy matrix with entries of* $M_i$ *taken from the set S = {–1, 0, 1}. Let the collective views of the n experts in the form of FRMs be given by* $M = M_1 \cup M_2 \cup \ldots \cup M_n$ *which is a special fuzzy rectangular matrix known as the fuzzy relational connection matrix of the Special Fuzzy Relational Maps (SFRMs) model.*

*We shall just describe how it functions. Given* $M = M_1 \cup M_2 \cup \ldots \cup M_n$ *to be the associated matrix of the SFRM, each of the FRMs in the SFRM is formed or obtained using the m attributes as the domain space and the t-attributes as the range space. Now the domain space of the SFRM will be known as the special domain space* $D_s$ *and it contains special fuzzy row vectors each of the fuzzy row vectors in this special fuzzy row vector is a* $1 \times m$ *matrix. Thus* $X = X_1 \cup X_2 \cup \ldots \cup X_n$ *then each* $X_i = \left[x_1^i, x_2^i, \ldots, x_m^i\right]$; $1 \le i \le n$ *and* $x_j^i \in \{0, 1\}$; $1 \le j \le m$. *X is a special fuzzy row vector from the special domain space* $D_s$.

*Similarly the range space of SFRM will be known as the special range space denoted by* $R_s$ *and every element* $R_s$ *is a special fuzzy row vector of the form* $Y = Y_1 \cup Y_2 \cup \ldots \cup Y_n$ *where each* $Y_i = \left[y_1^i, y_2^i, \ldots, y_t^i\right]$ *with* $y_j^i \in \{0, 1\}$, $i = 1, 2, \ldots, n$; $1 \le j \le t$. *These special fuzzy row vectors form the special range space or special domain space, only indicate the ON or OFF state of the node/attribute. Any input vector of the dynamical system of the SFRM would be a special fuzzy row vector. Suppose* $X = X_1 \cup X_2 \cup \ldots \cup X_n$ *from* $D_s$ *be the input special vector given by all the n experts. To find the effect X on the SFRM special dynamical system M.*

$$
\begin{aligned}
X \circ M &= (X_1 \cup X_2 \cup \ldots \cup X_n) \circ (M_1 \cup M_2 \cup \ldots \cup M_n) \\
&= X_1 \circ M_1 \cup X_2 \circ M_2 \cup \ldots \cup X_n \circ M_n \\
&= Y'_1 \cup Y'_2 \cup \ldots \cup Y'_n \\
&= Y'.
\end{aligned}
$$

*Y' may or may not be a special fuzzy row vector. We threshold Y' to Y and find* $Y \circ M^t$ *as* $Y \circ M$ *is not defined as in the FRMs, we are endorsed to use the transpose also. Each* $Y_i$ *is a* $1 \times t$ *fuzzy row vector* $(1 \le i \le n)$. *We find*

$$
\begin{aligned}
Y \circ M^t &= Z'_1 \cup Z'_2 \cup \ldots \cup Z'_m \\
&= Z'.
\end{aligned}
$$

*Z' may or may not be a special fuzzy row vector so we first threshold it; now it may so happen the coordinate with which we started in the on state would have become off so we update*

*Z' to obtain Z. Let* $Z = Z_1 \cup Z_2 \cup \ldots \cup Z_n$; $Z \in D_s$ *and Z is a special fuzzy row vector with each* $Z_i$ *a* $1 \times m$ *fuzzy row vector* $1 \le i \le n$.

*Now we find the effect of Z on the dynamical system M, i.e., Z o M and so on till we arrive at a fixed point or a limit cycle or a combination of both which will be known as the special fixed point or special limit cycle. When the system reaches a special fixed point or a special limit cycle we call that the special hidden pattern of the SFRM for the special input vector* $X = X_1 \cup X_2 \cup \ldots \cup X_n \in D_s$. *It is important to note that there is no error if one takes a input special vector* $Y = Y_1 \cup Y_2 \cup \ldots \cup Y_n \in R_s$ *and work with the dynamical system of the SFRM. We find*

$$\begin{aligned} Y \text{ o } M^t &= (Y_1 \cup Y_2 \cup \ldots \cup Y_n) \text{ o } (M_1 \cup M_2 \cup \ldots \cup M_n)^t \\ &= (Y_1 \cup Y_2 \cup \ldots \cup Y_n) \text{ o } (M_1^t \cup M_2^t \cup \ldots \cup M_n^t) \\ &= Y_1 \text{ o } M_1^t \cup Y_2 \text{ o } M_2^t \cup \ldots \cup Y_n \text{ o } M_n^t \\ &= Z'_1 \cup Z'_2 \cup \ldots \cup Z'_n \\ &= Z', \end{aligned}$$

*(we threshold Z' and obtain* $Z = Z_1 \cup \ldots \cup Z_n$*). Now one can find the effect of Z on M as*

$$\begin{aligned} (Z_1 \cup Z_2 \cup \ldots \cup Z_n) \text{ o } (M_1 \cup M_2 \cup \ldots \cup M_n) &= Z_1 \text{ o } M_1 \cup Z_2 \text{ o } M_2 \cup \ldots \cup Z_n \text{ o } M_n \\ &= P'_1 \cup P'_2 \cup \ldots \cup P'_n \, ; \end{aligned}$$

*we threshold P' to find* $P = P_1 \cup P_2 \cup \ldots \cup P_n$. *Now we find P o* $M^t$ *and so on until we arrive at a special hidden pattern which may be a special fixed point or a special limit cycle. It is important and interesting to note that in case of SFRM we get a pair of special fuzzy row vectors (S, T) where* $S \in D_s$ *and* $T \in R_s$ *to be the special hidden pattern. We see in case of SFCMs we get only one special fuzzy row vector.*

Now we illustrate this situation by an example.

***Example 2.4.1:*** Suppose three experts want to work with the problem of finding the relationship between the industries profit and the type of salaries to be given to the employee by the industries.

The nodes given by the three experts related to the pay techniques for the employee and their performance is taken as the nodes of the domain spaces which is as follows.

$D_1$ – Pay with allowances and bonus to workers
$D_2$ – Only pay
$D_3$ – Pay with allowance(or bonus)
$D_4$ – Best performance by the employee
$D_5$ – Average performance by the employee
$D_6$ – Poor performance by the employee
$D_7$ – Employee works for more number of hours
$D_8$ – Employee works for less number of hours.

Attributes related with the industry which is taken as the nodes of the range space.

$R_1$ – Maximum profit to the industry
$R_2$ – Only profit to the industry
$R_3$ – Neither profit nor loss to the industry
$R_4$ – Loss to the industry
$R_5$ – Heavy loss to the industry.

The connection matrix $M_1$ of the FRM given by the first expert related to this problem.

$$M_1 = \begin{bmatrix} 0 & 0 & 0 & 0 & 1 \\ 1 & 0 & 0 & 0 & 0 \\ 0 & 0 & 1 & 0 & 0 \\ 1 & 0 & 0 & 0 & 0 \\ 0 & 1 & 0 & 0 & 0 \\ 0 & 0 & 0 & 0 & 1 \\ 1 & 0 & 0 & 0 & 0 \\ 0 & 0 & 0 & 1 & 0 \end{bmatrix}$$

where $D_i$'s are taken along the rows and $R_j$'s along the column of the relational matrix, $(1 \le i \le 8$ and $1 \le j \le 5)$.

The relational matrix given by the second expert for the same problem with the same set of nodes using the FRM model is given by $M_2$.

$$M_2 = \begin{bmatrix} 0 & 0 & 0 & 1 & 1 \\ 0 & 0 & 1 & 0 & 0 \\ 0 & 0 & 1 & 0 & 0 \\ 1 & 1 & 0 & 0 & 0 \\ 0 & 0 & 1 & 0 & 0 \\ 0 & 0 & 0 & 1 & 0 \\ 0 & 0 & 0 & 0 & 1 \\ 0 & 0 & 0 & 1 & 1 \end{bmatrix}.$$

The relational matrix given by the third expert related to the problem of pay-profit in industries using the FRM model is given by $M_3$.

$$M_3 = \begin{bmatrix} 0 & 0 & 0 & 1 & 0 \\ 0 & 0 & 1 & 0 & 0 \\ 0 & 1 & 0 & 0 & 0 \\ 1 & 0 & 0 & 0 & 0 \\ 0 & 0 & 1 & 0 & 0 \\ 0 & 0 & 0 & 0 & 1 \\ 0 & 1 & 0 & 0 & 0 \\ 0 & 0 & 0 & 0 & 1 \end{bmatrix}.$$

Thus $M = M_1 \cup M_2 \cup M_3$ is the special fuzzy rectangular matrix related with the SFRM model i.e.,

$$M = \begin{bmatrix} 0 & 0 & 0 & 0 & 1 \\ 1 & 0 & 0 & 0 & 0 \\ 0 & 0 & 1 & 0 & 0 \\ 1 & 0 & 0 & 0 & 0 \\ 0 & 1 & 0 & 0 & 0 \\ 0 & 0 & 0 & 0 & 1 \\ 1 & 0 & 0 & 0 & 0 \\ 0 & 0 & 0 & 1 & 0 \end{bmatrix} \cup \begin{bmatrix} 0 & 0 & 0 & 1 & 1 \\ 0 & 0 & 1 & 0 & 0 \\ 0 & 0 & 1 & 0 & 0 \\ 1 & 1 & 0 & 0 & 0 \\ 0 & 0 & 1 & 0 & 0 \\ 0 & 0 & 0 & 1 & 0 \\ 0 & 0 & 0 & 0 & 1 \\ 0 & 0 & 0 & 1 & 1 \end{bmatrix} \cup$$

$$\begin{bmatrix} 0 & 0 & 0 & 1 & 0 \\ 0 & 0 & 1 & 0 & 0 \\ 0 & 1 & 0 & 0 & 0 \\ 1 & 0 & 0 & 0 & 0 \\ 0 & 0 & 1 & 0 & 0 \\ 0 & 0 & 0 & 0 & 1 \\ 0 & 1 & 0 & 0 & 0 \\ 0 & 0 & 0 & 0 & 1 \end{bmatrix}.$$

Suppose the experts want to work with the special input vector

$$\begin{aligned} X &= X_1 \cup X_2 \cup X_3 \\ &= [1\ 0\ 0\ 0\ 0\ 0\ 0\ 0] \cup [1\ 0\ 0\ 0\ 0\ 0\ 0\ 0] \cup [1\ 0\ 0\ 0\ 0\ 0\ 0\ 0]. \end{aligned}$$

Now the effect of X on the dynamical system M is given by

$$\begin{aligned} X \text{ o } M &= (X_1 \cup X_2 \cup X_3) \text{ o } (M_1 \cup M_2 \cup M_3) \\ &= X_1 \text{ o } M_1 \cup X_2 \text{ o } M_2 \cup X_3 \text{ o } M_3 \end{aligned}$$

$$= [1\ \ 0\ \ 0\ \ 0\ \ 0\ \ 0\ \ 0\ \ 0] \text{ o } \begin{bmatrix} 0 & 0 & 0 & 0 & 1 \\ 1 & 0 & 0 & 0 & 0 \\ 0 & 0 & 1 & 0 & 0 \\ 1 & 0 & 0 & 0 & 0 \\ 0 & 1 & 0 & 0 & 0 \\ 0 & 0 & 0 & 0 & 1 \\ 1 & 0 & 0 & 0 & 0 \\ 0 & 0 & 0 & 1 & 0 \end{bmatrix} \cup$$

$$\begin{bmatrix}1 & 0 & 0 & 0 & 0 & 0 & 0 & 0\end{bmatrix} \text{ o } \begin{bmatrix}0 & 0 & 0 & 1 & 1 \\ 0 & 0 & 1 & 0 & 0 \\ 0 & 0 & 1 & 0 & 0 \\ 1 & 1 & 0 & 0 & 0 \\ 0 & 0 & 1 & 0 & 0 \\ 0 & 0 & 0 & 1 & 0 \\ 0 & 0 & 0 & 0 & 1 \\ 0 & 0 & 0 & 1 & 1\end{bmatrix} \cup$$

$$\begin{bmatrix}1 & 0 & 0 & 0 & 0 & 0 & 0 & 0\end{bmatrix} \text{ o } \begin{bmatrix}0 & 0 & 0 & 1 & 0 \\ 0 & 0 & 1 & 0 & 0 \\ 0 & 1 & 0 & 0 & 0 \\ 1 & 0 & 0 & 0 & 0 \\ 0 & 0 & 1 & 0 & 0 \\ 0 & 0 & 0 & 0 & 1 \\ 0 & 1 & 0 & 0 & 0 \\ 0 & 0 & 0 & 0 & 1\end{bmatrix}$$

= [0 0 0 0 1] ∪ [0 0 0 1 1] ∪ [0 0 0 1 0]
= $Y_1 \cup Y_2 \cup Y_3$
= $Y \in R_s$.

Now we find

Y o M = [1 0 0 0 0 1 0 0] ∪ [2 0 0 0 0 1 1 2] ∪ [1 0 0 0 0 0 0 0]
= Z' .

Z' is thresholded and updated and Z is obtained. Now

Z = [1 0 0 0 0 1 0 0] ∪ [1 0 0 0 0 1 1 1] ∪ [1 0 0 0 0 0 0 0].

We find

$$Z \text{ o } M = [0\ 0\ 0\ 0\ 2] \cup [0\ 0\ 0\ 3\ 3] \cup [0\ 0\ 0\ 1\ 0] = T'$$

after thresholding we get the resultant of T' as

$$T = [0\ 0\ 0\ 0\ 1] \cup [0\ 0\ 0\ 1\ 1] \cup [0\ 0\ 0\ 1\ 0] = Y \in R_s.$$

Thus the special hidden pattern of the special dynamical system is given by the pair of special fuzzy row vectors or the special binary pair {[1 0 0 0 0 0 1 0 0] ∪ [1 0 0 0 0 1 1 1] ∪ [1 0 0 0 0 0 0 0], [0 0 0 0 1] ∪ [0 0 0 1 1] ∪ [0 0 0 1 0]}.

The reader can work with any special input state vector and find the special hidden pattern as the binary pair which would be a pair of special fuzzy row vector.

Now we proceed on to define yet another new special model viz special mixed FRM model.

**DEFINITION 2.4.2:** *Let P be a problem which is analyzed by a set of m experts. The only common feature for all of the experts is that all of them agree to work with FRM model on the same problem but they demand their liberty to choose any set of attributes as per their desire. Let the $i^{th}$ expert choose to work with $t_i$ number domain attributes of the FRM which forms the rows of the relational matrix and with $s_i$ number of range attributes which form the columns of the relational matrix. Let $M_i$ denote the connection relational matrix which is a $t_i \times s_i$ fuzzy matrix with entries from the set {-1, 0, 1} for i = 1, 2, …, m. (i.e., the same procedure is adopted by all the m experts). Let us define the Special Mixed Fuzzy Relational Maps (SMFRMs) model to be the model represented by the special fuzzy mixed rectangular matrix $M = M_1 \cup M_2 \cup \ldots \cup M_m$. The special domain space $D_s$ associated with this model is represented by the set of all special fuzzy mixed row vectors $X = X_1 \cup X_2 \cup \ldots \cup X_m$ where each $X_i$ is a $1 \times t_i$ fuzzy row vector with $X_i = \left[x_1^i\, x_2^i\, \ldots\, x_{t_i}^i\right]$ and $x_j^i \in \{0, 1\}$, $1 \le j \le t_i$, $1 \le i \le m$. The* *special range space $R_s$ associated with the SMFRM consists of all special fuzzy mixed row vectors of the form $Y = Y_1 \cup Y_2 \cup \ldots \cup Y_m$ where each $Y_i$ is a $1 \times s_i$ fuzzy row vector with $Y_i = \left[ y_1^i \; y_2^i \ldots y_{s_i}^i \right]$; and $y_j^i \in \{0, 1\}$; $j = 1, 2, \ldots, s_i$ and $1 \leq i \leq m$. The special dynamical system of the SMFRM comprises of $\{M = M_1 \cup M_2 \cup \ldots \cup M_n, D_s$ and $R_s\}$.*

*Now we will show how the special dynamical system of the SMFRM functions. Let $X = X_1 \cup X_2 \cup \ldots \cup X_m$ be the special fuzzy mixed row vector given as the input vector by all the m experts; $X \in D_s$ to find the effect of X on the special dynamical system M.*

$$
\begin{aligned}
X \text{ o } M &= (X_1 \cup X_2 \cup \ldots \cup X_n) \text{ o } (M_1 \cup M_2 \cup \ldots \cup M_m) \\
&= X_1 \text{ o } M_1 \cup X_2 \text{ o } M_2 \cup \ldots \cup X_m \text{ o } M_m \\
&= Y_1^{'} \cup Y_2^{'} \cup \ldots \cup Y_m^{'} \\
&= Y^{'}.
\end{aligned}
$$

*We threshold $Y'$ to $Y = Y_1 \cup Y_2 \cup \ldots \cup Y_m$, $Y \in R_S$. Now we calculate*

$$
\begin{aligned}
Y \text{ o } M^t &= (Y_1 \cup Y_2 \cup \ldots \cup Y_m) \text{ o } (M_1^t \cup M_2^t \cup \ldots \cup M_m^t) \\
&= Y_1 \text{ o } M_1^t \cup Y_2 \text{ o } M_2^t \cup \ldots \cup Y_m \text{ o } M_m^t \\
&= Z'_1 \cup Z'_2 \cup \ldots \cup Z'_m \\
&= Z'.
\end{aligned}
$$

*$Z'$ may or may not belong to $D_s$; we update and threshold Z' to Z and if $Z = Z_1 \cup Z_2 \cup \ldots \cup Z_m$. We find Z o M and so on, until we arrive at a special fixed point or a special limit cycle which is a binary pair (T, P) with T and P special fuzzy mixed row vectors belonging to $D_s$ and $R_s$ respectively.*

Now we proceed on to define the new notion of special neutrosophic relational maps (SNRM) model.

**DEFINITION 2.4.3:** *Suppose we have a problem P at hand and let t experts want to work with it using Neutrosophic Relational*

*Maps and all of them agree upon to work with the same set of domain attributes.*

*Now if in the SFRM dynamical system the special fuzzy rectangular matrix is replaced by the special fuzzy neutrosophic rectangular matrix and in the special range space and in the special domain space, the special fuzzy row vectors are replaced by special fuzzy neutrosophic row vectors then the resulting dynamical system will be known as the Special Neutrosophic Relational Maps (SNRMs) model. if all the t experts agree to work on the problem P using the same fixed number say m attributes which forms the domain space and n attributes forming the range space then if we denote the fuzzy neutrosophic relational matrix given by the* $i^{th}$ *expert by* $M_i$ *where* $M_i$ *is a* $m \times n$ *fuzzy neutrosophic rectangular matrix true for* $i = 1, 2, \ldots, t$, *then* $M = M_1 \cup M_2 \cup \ldots \cup M_t$ *together with the special range and special domain space is the dynamical system related with all the t experts.*

*Now the special domain space of this system would be given by* $X = X_1 \cup X_2 \cup \ldots \cup X_t$ *where each* $X_i = \left[x_1^i \; x_2^i \; \ldots \; x_m^i\right]$; $x_j^i \in \{0, I, 1\}$ $i = 1, 2, \ldots, t$ *and* $1 \leq j \leq m$ *and* $X \in D_s$, $D_s$ *denotes the set of all X with entries from {0, I, 1}. Clearly X is a special neutrosophic row vector*.

*Note:* It is important to note that the special input vectors in SNRMs will always be from the set {0, 1} i.e., it denotes the on or off state of the node; however the resultant special row vectors given after finding X o M can be a special neutrosophic row vector, of the special domain space or special range space. That is why we in the definition take the special domain space and special range space to be a special neutrosophic row vectors.

Now the special range space $Y = Y_1 \cup Y_2 \cup \ldots \cup Y_t$ where each $Y_i = \left(y_1^i \; y_2^i \; \ldots \; y_n^i\right)$; $y_j^i \in \{0, I, 1\}$; $i = 1, 2, \ldots, t$ and $1 \leq j \leq n$ and $Y \in R_s$, $R_s$ denotes the set of all Y with entries from {0, 1}. Y is also a special neutrosophic input row vector. Now the triple $\{M = M_1 \cup M_2 \cup \ldots \cup M_n, D_s, R_s\}$ forms the special dynamical neutrosophic system and the model would be known as the Special Neutrosophic Relational Maps (SNRM-model).

Next we proceed on to give a brief description of another special neutrosophic model.

**DEFINITION 2.4.4:** *Suppose in the definition 2.4.2 of the SMFRM model if the FRMs are replaced by NRM we get the new Special Mixed Neutrosophic Relational Maps (SMNRMs) model i.e., the special dynamical system of the SMNRM model is given by the triple $\{M = M_1 \cup M_2 \cup \ldots \cup M_n, D_s, R_s\}$ where $M_i$'s are fuzzy neutrosophic rectangular $t_i \times s_i$ matrices with entries from the set {-I, -1, 0, I, 1}; i.e., M is a special fuzzy neutrosophic mixed rectangular matrices; where $D_s$ and $R_s$ are special fuzzy neutrosophic mixed row vectors.*

*Here it is pertinent to mention that however when the experts give any input special vector $X = X = X_1 \cup X_2 \cup \ldots \cup X_n$ then $X_i$'s can take only entries from the set {0, 1}; i.e., the dynamical set accepts only the on or off state of the nodes for no expert can claim the existence or non existence or indeterminacy of the node only while working with the problem of the dynamical system with certain nodes in the on state either from the range attributes or the domain attributes it may so happen that the state of the node may become indeterminate with some input vectors.*

Now we proceed on to define yet another special new model viz., Special Fuzzy Neutrosophic Relational Maps (SFNRM) model.

**DEFINITION 2.4.5:** *Suppose we have a problem P with which some t experts wish to work. It so happens that out of the t experts some m experts wants to work using only the Fuzzy Relational Maps (FRM) model and rest (t – m) experts wants to work using the Neutrosophic Relational Maps (NRM) model. Here also they i.e., the m experts may work with the same set of attributes or different sets of for their domain and range spaces* *of the FRMs, like wise the t-m experts may work with same set of attributes or different sets of attributes for their domain and range spaces of the NRMs. Now this new model given by these t experts for the problem will consists of a special dynamical system, the special matrix which will be a mixture of rectangular fuzzy matrices as well as rectangular neutrosophic matrices given by* $M = \{M_1 \cup M_2 \cup \ldots \cup M_m\} \cup \{M_{m+1} \cup \ldots \cup M_t\}$ *where* $M_1 \cup M_2 \cup \ldots \cup M_m$ *is a special fuzzy mixed rectangular matrix and* $\{M_{m+1} \cup \ldots \cup M_t\}$ *is a special neutrosophic rectangular matrix. Thus* $\{M_1 \cup M_2 \cup \ldots \cup M_m\} \cup \{M_{m+1}\ M_{m+1} \cup \ldots \cup M_t\}$ *will be known as the special mixed fuzzy and neutrosophic mixed rectangular matrix. Now the special domain space* $X = \{X_1 \cup X_2 \cup \ldots \cup X_m\} \cup \{X_{m+1} \cup X_{m+2} \cup \ldots \cup X_t\}$ *will be such that* $\{X_1 \cup X_2 \cup \ldots \cup X_m\}$ *is a special fuzzy mixed row vectors and* $\{X_{m+1} \cup X_{m+2} \cup \ldots \cup X_t\}$ *will be a special fuzzy neutrosophic mixed row vectors; thus* $X = \{X_1 \cup X_2 \cup \ldots \cup X_m\} \cup \{X_{m+1} \cup X_{m+2} \cup \ldots \cup X_t\}$ *will be a special mixed fuzzy and neutrosophic mixed row vectors. This will be known as the Special Domain Space of the new model. Similarly* $Y = \{Y_1 \cup Y_2 \cup \ldots \cup Y_m\} \cup \{Y_{m+1} \cup Y_{m+2} \cup \ldots \cup Y_t\}$ *will be a special mixed fuzzy and neutrosophic mixed row vector known as the Special Range Space of the new model.*

*So the Special Fuzzy and Neutrosophic Relational Maps (SFNRM) model will consist of Special Fuzzy Neutrosophic mixed rectangular matrices, special domain space* $D_s$ *and special range* $R_s$ *where* $D_s$ *will be the set of all special fuzzy neutrosophic mixed row vectors of the form* $X = \{X_1 \cup X_2 \cup \ldots \cup X_m\} \cup \{X_{m+1} \cup X_{m+2} \cup \ldots \cup X_t\}$ *where* $X_1 \cup X_2 \cup \ldots \cup X_m$ *are special fuzzy mixed row vectors with entries from the set {0, 1} and* $\{X_{m+1} \cup X_{m+2} \cup \ldots \cup X_t\}$ *are special neutrosophic mixed row vectors with entries from the set {I, 0 1}. Clearly even in this model any input vector from the t experts say* $X = \{X_1 \cup X_2 \cup \ldots \cup X_m\} \cup \{X_{m+1} \cup X_{m+2} \cup \ldots \cup X_t\}$*will take entries only from the set {0, 1} i.e., on or off state of the nodes or attributes.*

*Thus* $M = \{M_1 \cup M_2 \cup \ldots \cup M_m\} \cup \{M_{m+1} \cup \ldots \cup M_t, D_s, R_s\}$ *will form the special dynamical system of SFNRM model.*

The main advantage of this model is that the set of experts can choose to work either using FRM or NRM. For all the while, the models which we have given in this section can either work in FRM or in NRM. So in this manner this model is more powerful or advantageous than other models.

Secondly in this model one can use any set of attributes for the domain space and the range space, it is left for the choice of the expert, thus the special domain space and the special range space contain special fuzzy and neutrosophic mixed row vectors.

This model can be realized as a generalization of SFRMs, SMFRMs, SNRMs and SMNRMs models.

Now we proceed on to define yet another new model called the Special Mixed Fuzzy Cognitive and Fuzzy Relational Maps (SMFCMs) model.

**DEFINITION 2.4.6:** *Suppose we have a problem P at hand and some t experts want to work with this problem. Here some of the experts want to work using FCMs and some of them prefer to work using FRMs and they may opt for same set of attributes or for different sets of attributes according to their choice. In this case we use the Special Mixed Fuzzy Cognitive and Fuzzy Relational Maps Model. This model can be thought of as a mixture of SMFCM and SMFRM model. Let m experts choose to work with FCM and t-m of them choose to work with FRM then these m experts choose to work with SMFCM model and t-m of them work with SMFRM model. Suppose* $M = \{M_1 \cup M_2 \cup \ldots \cup M_m\}$ *denotes the SMFCMs associated special fuzzy mixed square matrices and* $T = \{M_{m+1} \cup M_{m+2} \cup \ldots \cup M_t\}$ *denotes the SMFRMs associated special mixed fuzzy rectangular matrices then* $M \cup T = \{M_1 \cup M_2 \cup \ldots \cup M_m\} \cup \{M_{m+1} \cup M_{m+2} \cup \ldots \cup M_t\}$ *denotes the special fuzzy mixed matrices associated with the SMFCRM-model and the special domain space is given by* $\{X_1 \cup X_2 \cup \ldots \cup X_m\} \cup \{X_{m+1} \cup X_{m+2} \cup \ldots \cup X_t\}$ *the set of all special fuzzy mixed row vectors with entries from the set {0, 1}* *and the special range space* $Y = \{Y_1 \cup Y_2 \cup \ldots \cup Y_m\} \cup \{Y_{m+1} \cup Y_{m+2} \cup \ldots \cup Y_t\}$ *denotes the set of all special fuzzy mixed row matrices with entries from the set {0, 1}. Thus* $M = \{M_1 \cup M_2 \cup \ldots \cup M_m\} \cup \{M_{m+1}\ M_{m+2} \cup \ldots \cup M_t\}\ D_s, R_s\}$ *denotes the special fuzzy dynamical system of the SMFCFRM.*

The main advantage of this special model is that for the same problem some experts can choose to work with FCM and the rest with FRM so this multi expert model need not model it in two stages or as two models for now in this new model they function as a single model with double purpose.

Now we proceed on to describe the next model viz. the Special Mixed Neutrosophic Cognitive and Neutrosophic Relational Maps Model (SMNCNRMs model).

**DEFINITION 2.4.7:** *Suppose we have a problem P at hand and all the experts agree to work using only the neutrosophic models. One set of experts out of t experts say m experts wish to work with the NCM model with different sets of attributes i.e., the SMNCM model and the t-m experts wish to work with the NRM model they too with different sets of domain and range space i.e., SMNRM model. Then as in case of the above definition given in 2.4.6 we form the Special Mixed Neutrosophic Cognitive and Neutrosophic Relational Map Model i.e., in the definition 2.4.6 we can replace the FCM by NCM and FRM by NRM and obtain the SMNCNRM-model. The dynamical system associated with this model is given by* $M = \{M_1 \cup M_2 \cup \ldots \cup M_m\} \cup \{M_{m+1}\ M_{m+2} \cup \ldots \cup M_t, D_s, R_s\}$ *where* $D_s$ *is the special domain space consisting of all special neutrosophic mixed row vectors of the form* $X = \{X_1 \cup X_2 \ldots \cup X_m\} \cup \{X_{m+1} \cup X_{m+2} \cup \ldots \cup X_t\}$*, where* $X_i$ *takes its values from the set {I, 0, 1}, i = 1, 2, …, t and the special range space* $R_s$ *consists of all special neutrosophic mixed row vectors of the form* $Y = \{Y_1 \cup Y_2 \ldots \cup Y_m\} \cup \{Y_{m+1} \cup Y_{m+2} \cup \ldots \cup Y_t\}$ *where* $Y_i$ *takes their entries from the set {I, 0, 1}, i = 1, 2, …, t. It is pertinent to mention here that in case of special input vectors from the t experts would be given by* $P = \{P_1 \cup P_2 \ldots \cup P_m\} \cup$

*$\{P_{m+1} \cup P_{m+2} \cup \ldots \cup P_t\}$ we see the entries of each $P_i$; $1 \leq i \leq t$ are only are from the set $\{0, 1\}$ i.e., the on or off state of the attribute in the special range space or the special domain space.*

The main advantage of this model over the other models is that this model can function with two types of models at a time like the SMFCFRMs and secondly when every expert agrees to work with neutrosophic models i.e., they feel the notion of indeterminacy is present in the problem.

Next we proceed on to define yet a new model which can be thought of as the generalization of all the models so far introduced.

**DEFINITION 2.4.8:** *Let P be the problem which is to be investigated. Suppose t experts wish to work with the problem, but some i.e., $t_1$ experts wish to use FCM model, a few of them say $t_2$ experts use FRM model, some $t_3$ of them want to use NCM model as they think the problem involves some degree of indeterminacy and the rest of experts $t - (t_1 + t_2 + t_3)$ of them wish to work with NRM model. The model should be so connected that it should work as a single unit and give the results simultaneously. Now this new special model can be got as the union of SMFCM $\cup$ SMFRM $\cup$ SMNCM $\cup$ SMNRM so that this can be represented by the union of four sets of special matrices given by these four special models. Suppose M denotes the associated matrix of the new model them*

$$\begin{aligned} M &= M_1 \cup M_2 \cup M_3 \cup M_4 \\ &= (M_1 \cup M_2 \cup \ldots \cup M_{t_1}) \cup \{M_{t_1+1} \cup M_{t_{1+2}} \cup \ldots \cup M_{t_2}\} \\ &\quad \cup \{M_{t_2+1} \cup M_{t_{1+2}} \cup \ldots \cup M_{t_3}\} \cup \{M_{t_{3+1}} \cup M_{t_{3+2}} \cup \ldots \cup M_t\} \end{aligned}$$

*where $M_1$ is a special fuzzy mixed square matrices, $M_2$ is a special fuzzy mixed rectangular matrices, $M_3$ is a special fuzzy neutrosophic mixed square matrices and $M_4$ is a special fuzzy neutrosophic mixed rectangular matrices related with the models SMFCM, SMFRM, SMNCM and SMNRM respectively. Now this special dynamical system will have three components given by the triple $\{M = M_1 \cup M_2 \cup M_3 \cup M_4 = \{M_1 \cup M_2 \cup$*

$\ldots \cup M_{t_1}$ *} ∪{* $M_{t_1+1} \cup M_{t_{1+2}} \cup \ldots \cup M_{t_2}$ *} ∪ {* $M_{t_{2+1}} \cup M_{t_{2+2}} \cup \ldots \cup M_{t_3}$ *} ∪{* $M_{t_{3+1}} \cup M_{t_{3+2}} \cup \ldots \cup M_t$ *}, $D_s$, $R_s$}*

*where $D_s$ is the special fuzzy neutrosophic mixed row vectors of the form {$X_1 \cup X_2 \cup \ldots \cup X_{t_1}$ } ∪{ $X_{t_1+1} \cup X_{t_{1+2}} \cup \ldots \cup X_{t_2}$ } ∪{ $X_{t_{2+1}} \cup X_{t_{2+2}} \cup \ldots \cup X_{t_3}$ } ∪{ $X_{t_{3+1}} \cup X_{t_{3+2}} \cup \ldots \cup X_t$} = $X_1 \cup X_2 \cup X_3 \cup X_4$ where each $X_i = \left[x_1^i \; x_2^i \; \ldots \; x_l^i\right]$ where $x_j^i \in$ {0, 1} with $1 \leq l < t_1, t_2$ and $X_i = \left[x_1^k \; x_2^k \ldots x_p^k\right]$($1 \leq p \leq t_3$, $t-(t_1+ t_2+ t_3)$ and $x_r^k \in$ {0, 1, I}. Similarly $R_s$ is also a special fuzzy neutrosophic mixed row vector like $D_s$.*

Now X o M = $(X_1 \cup X_2 \cup X_3 \cup X_4)$ o $(M_1 \cup M_2 \cup M_3 \cup M_4)$ = $X_1$ o $M_1 \cup X_2$ o $M_2 \cup X_3$ o $M_3 \cup X_4$ o $M_4$ where each $X_i$ o $M_1$ works like the dynamical system SMFCM, $X_2$ o $M_2$ works like the dynamical system SMFRM, $X_3$ o $M_3$ works like the dynamical system SMNCM and $X_4$ o $M_4$ works like the dynamical system of SMNRM respectively.

However the new SMFCRNCRM model can be realized as generalized model of every other models so far described in sections 2.1 to 2.4 of this book. Further this model combines the fuzzy and indeterminacy together.

Now we will show the working of the model by an explicit example.

***Example 2.4.2:*** Suppose one is interested in studying the Employee-Employers model so that the industry should run in profit as well as the workers should be satisfied and the performance of the industry must be good in both quality and production. There are six experts who wish to work with the model and each of them choose to work with a distinct model. The set of attributes collectively put by all of them is as follows:

| | | |
|---|---|---|
| $A_1$ | - | Employee gets pay, allowance and bonus |
| $A_2$ | - | He gets only pay |
| $A_3$ | - | Pay with either allowances or bonus alone |
| $A_4$ | - | Best performance by the employee |

$A_5$ - Maximum profit for the employer
$A_6$ - Average performance by the employee
$A_7$ - Poor performance by the employee
$A_8$ - Only profit to the employer
$A_9$ - Employee works for more number of hours
$A_{10}$ - Neither profit nor loss to the employer
$A_{12}$ - Heavy loss to the employer
$A_{13}$ - Employee works for less number of hours.

The first expert wishes to work with 8 attributes using the FCM model. The connection matrix $M_1$ given by him is as follows:

$$M_1 = \begin{array}{c} \\ A_1 \\ A_2 \\ A_3 \\ A_4 \\ A_5 \\ A_9 \\ A_{11} \\ A_{13} \end{array} \begin{array}{c} \begin{array}{cccccccc} A_1 & A_2 & A_3 & A_4 & A_5 & A_9 & A_{11} & A_{13} \end{array} \\ \begin{bmatrix} 0 & 0 & 0 & 1 & 1 & 0 & 0 & 0 \\ 0 & 0 & 0 & 0 & 0 & 0 & 1 & 1 \\ 0 & 0 & 0 & 0 & 0 & 0 & 1 & 0 \\ 1 & 0 & 0 & 0 & 1 & 0 & 0 & 0 \\ 1 & 0 & 0 & 0 & 0 & 1 & 0 & 0 \\ 1 & 0 & 0 & 0 & 0 & 1 & 0 & 0 \\ 0 & 1 & 0 & 0 & 0 & 0 & 1 & 1 \\ 0 & 1 & 0 & 0 & 0 & 0 & 1 & 0 \end{bmatrix} \end{array}.$$

The second expert also wants to work with FCM but with different sets of attributes and more so only with six attributes the connection matrix $M_2$ given by him is a $6 \times 6$ matrix, $M_2$ given by:

$$M_2 = \begin{array}{c} \\ A_1 \\ A_2 \\ A_3 \\ A_4 \\ A_8 \\ A_7 \end{array} \begin{array}{c} \begin{array}{cccccc} A_1 & A_2 & A_3 & A_4 & A_8 & A_7 \end{array} \\ \begin{bmatrix} 0 & 0 & 0 & 1 & 0 & 0 \\ 0 & 0 & 0 & 0 & 0 & 1 \\ 0 & 0 & 0 & 0 & 1 & 0 \\ 1 & 0 & 0 & 0 & 1 & 0 \\ 1 & 0 & 1 & 0 & 0 & 0 \\ 0 & 1 & 0 & 0 & 0 & 0 \end{bmatrix} \end{array}.$$

The third expert wants to work with a NCM model just having only 7 attributes the connection matrix $M_3$ given by him is as follows. However he uses some set of attributes from the 13 attributes.

$$M_3 = \begin{array}{c} \\ A_2 \\ A_3 \\ A_6 \\ A_7 \\ A_8 \\ A_9 \\ A_{10} \end{array} \begin{array}{c} \begin{array}{ccccccc} A_2 & A_3 & A_6 & A_7 & A_8 & A_9 & A_{10} \end{array} \\ \begin{bmatrix} 0 & 0 & 1 & 0 & 0 & I & 1 \\ 0 & 0 & 0 & 1 & 0 & 0 & 0 \\ 0 & 0 & 0 & 0 & 0 & 0 & 1 \\ 1 & 0 & 0 & 0 & 0 & 0 & 0 \\ 1 & 0 & 0 & 0 & 0 & 0 & 0 \\ 0 & 1 & 0 & 0 & 1 & 0 & I \\ 0 & 0 & I & 0 & 0 & 0 & 0 \end{bmatrix} \end{array}.$$

The fourth expert works with the FRM taking all the attributes and forms the relational matrix $M_4$ which is given by him is as follows:

$$M_4 = \begin{array}{c} \\ A_1 \\ A_2 \\ A_3 \\ A_4 \\ A_6 \\ A_7 \\ A_9 \\ A_{13} \end{array} \begin{array}{c} \begin{array}{ccccc} A_5 & A_8 & A_{10} & A_{11} & A_{12} \end{array} \\ \begin{bmatrix} 0 & 0 & 0 & 0 & 1 \\ 1 & 0 & 0 & 0 & 0 \\ 0 & 0 & 1 & 0 & 0 \\ 1 & 0 & 0 & 0 & 0 \\ 0 & 1 & 0 & 0 & 0 \\ 0 & 0 & 0 & 0 & 1 \\ 1 & 0 & 0 & 0 & 0 \\ 0 & 0 & 0 & 1 & 0 \end{bmatrix} \end{array}.$$

The fifth expert gives a relational matrix $M_5$ using FRM model he does not use all the attributes but uses only a fewer number of attributes.

$$M_5 = \begin{array}{c} \\ A_2 \\ A_3 \\ A_4 \\ A_6 \\ A_7 \\ A_{13} \end{array} \begin{array}{c} \begin{array}{cccc} A_8 & A_{10} & A_{11} & A_{12} \end{array} \\ \left[\begin{array}{cccc} 0 & 0 & 1 & 1 \\ 0 & 1 & 0 & 0 \\ 0 & 1 & 0 & 0 \\ 1 & 0 & 0 & 0 \\ 0 & 1 & 0 & 0 \\ 0 & 0 & 1 & 1 \end{array}\right] \end{array}.$$

The sixth expert however using yet a different set of attributes and uses NRM and the related matrix $M_6$ given by him is as follows:

$$M_6 = \begin{array}{c} \\ A_1 \\ A_2 \\ A_3 \\ A_4 \\ A_6 \\ A_7 \\ A_{13} \end{array} \begin{array}{c} \begin{array}{cccc} A_8 & A_{10} & A_{11} & A_{12} \end{array} \\ \left[\begin{array}{cccc} 0 & 0 & 0 & 1 \\ 0 & 0 & 1 & 0 \\ 0 & 0 & 1 & 0 \\ 1 & 1 & 0 & 0 \\ 0 & 0 & 1 & 0 \\ 0 & 0 & 0 & I \\ 0 & 0 & 0 & 1 \end{array}\right] \end{array}.$$

The six experts choose to work with the SMFCRNCRM model given by the
$(M_1 \cup M_2 \cup M_3 \cup M_4 \cup M_5 \cup M_6)$

$$\begin{bmatrix} 0 & 0 & 0 & 1 & 1 & 0 & 0 & 0 \\ 0 & 0 & 0 & 0 & 0 & 0 & 1 & 1 \\ 0 & 0 & 0 & 0 & 0 & 0 & 1 & 0 \\ 1 & 0 & 0 & 0 & 1 & 0 & 0 & 0 \\ 1 & 0 & 0 & 0 & 0 & 1 & 0 & 0 \\ 1 & 0 & 0 & 0 & 0 & 1 & 0 & 0 \\ 0 & 1 & 0 & 0 & 0 & 0 & 1 & 1 \\ 0 & 1 & 0 & 0 & 0 & 0 & 1 & 0 \end{bmatrix} \cup \begin{bmatrix} 0 & 0 & 0 & 1 & 0 & 0 \\ 0 & 0 & 0 & 0 & 0 & 1 \\ 0 & 0 & 0 & 0 & 1 & 0 \\ 1 & 0 & 0 & 0 & 1 & 0 \\ 1 & 0 & 1 & 0 & 0 & 0 \\ 0 & 1 & 0 & 0 & 0 & 0 \end{bmatrix} \cup$$

$$\begin{bmatrix} 0 & 0 & 1 & 0 & 0 & I & 1 \\ 0 & 0 & 0 & 1 & 0 & 0 & 0 \\ 0 & 0 & 0 & 0 & 0 & 0 & 1 \\ 1 & 0 & 0 & 0 & 0 & 0 & 0 \\ 1 & 0 & 0 & 0 & 0 & 0 & 0 \\ 0 & 1 & 0 & 0 & 1 & 0 & I \\ 0 & 0 & I & 0 & 0 & 0 & 0 \end{bmatrix} \cup \begin{bmatrix} 0 & 0 & 0 & 0 & 1 \\ 1 & 0 & 0 & 0 & 0 \\ 0 & 0 & 1 & 0 & 0 \\ 1 & 0 & 0 & 0 & 0 \\ 0 & 1 & 0 & 0 & 0 \\ 0 & 0 & 0 & 0 & 1 \\ 1 & 0 & 0 & 0 & 0 \\ 0 & 0 & 0 & 1 & 0 \end{bmatrix} \cup$$

$$\begin{bmatrix} 0 & 0 & 1 & 1 \\ 0 & 1 & 0 & 0 \\ 0 & 1 & 0 & 0 \\ 1 & 0 & 0 & 0 \\ 0 & 1 & 0 & 0 \\ 0 & 0 & 1 & 1 \end{bmatrix} \cup \begin{bmatrix} 0 & 0 & 0 & 1 \\ 0 & 0 & 1 & 0 \\ 0 & 0 & 1 & 0 \\ 1 & 1 & 0 & 0 \\ 0 & 0 & 1 & 0 \\ 0 & 0 & 0 & I \\ 0 & 0 & 0 & 1 \end{bmatrix};$$

is a special fuzzy-neutrosophic mixed matrix. Now the set of special attributes or special initial domain state vector with which the model works is

$$\begin{aligned} X = & \left\{ \left[x_1^1\ x_2^1 \ldots\ x_8^1\right] \cup \left[x_1^2\ x_2^2 \ldots\ x_6^2\right] \cup \left[x_1^3\ x_2^3 \ldots\ x_7^3\right] \cup \right. \\ & \left. \left[x_1^4\ x_2^4 \ldots x_8^4\right] \cup \left[x_1^5\ x_2^5 \ldots x_6^5\right] \cup \left[x_1^6\ x_2^6\ \ldots\ x_7^6\right] \right\} \\ = & \ X_1 \cup X_2 \cup X_3 \cup \ldots \cup X_6 \end{aligned}$$

where X is a special fuzzy mixed row vector $x_j^i \in \{0, 1\}$ or $\{0, 1, I\}$; $1 \leq i \leq 6$ and $1 \leq j \leq 8$. The special initial range space state vector

$$Y = Y_1 \cup Y_2 \ldots \cup Y_6$$

$$= \left\{\left[y_1^1\ y_2^1 \dots\ y_8^1\right] \cup \left[y_1^2\ y_2^2 \dots\ y_6^2\right] \cup \left[y_1^3\ y_2^3 \dots\ y_7^3\right] \cup \left[y_1^4\ y_2^4 \dots y_5^4\right] \cup \left[y_1^5\ y_2^5 \dots y_4^5\right] \cup \left[y_1^6\ y_2^6\ \dots\ y_4^6\right]\right\}$$

where $y_t^j \in \{0, 1, I\}$ but if this Y is started as a initial vector then $y_t^i \in \{0, 1\}$only i.e., on or off state.

Let us suppose the experts agree and give the special initial state vector

$$X = \{[1\ 0\ 0\ 0\ 0\ 0\ 0\ 0] \cup [0\ 0\ 0\ 0\ 1\ 0] \cup [0\ 1\ 0\ 0\ 0\ 0\ 0] \cup [0\ 0\ 0\ 0\ 0\ 0\ 0\ 1] \cup [0\ 0\ 1\ 0\ 0\ 0] \cup [0\ 0\ 0\ 1\ 0\ 0\ 0]\}$$

Now

$$X \text{ o}_m M = [X_1 \cup X_2 \cup \dots \cup X_6] \text{ o } [M_1 \cup M_2 \cup \dots \cup M_6]$$
$$= X_1 \text{ o } M_1 \cup X_2 \text{ o } M_2 \cup \dots \cup X_6 \text{ o } M_6$$

$$= [0\ 0\ 0\ 1\ 1\ 0\ 0\ 0] \cup [1\ 0\ 1\ 0\ 0\ 0] \cup [0\ 0\ 0\ 1\ 0\ 0\ 0] \cup [0\ 0\ 0\ 1\ 0] \cup [0\ 1\ 0\ 0] \cup [1\ 1\ 0\ 0]$$
$$= Y'_1 \cup Y'_2 \cup Y'_3 \cup Y'_4 \cup Y'_5 \cup Y'_6$$
$$= Y'$$

after updating Y' we get

$$Y = Y_1 \cup Y_2 \cup Y_3 \cup Y_4 \cup Y_5 \cup Y_6$$
$$= [1\ 0\ 0\ 1\ 1\ 0\ 0\ 0] \cup [1\ 0\ 1\ 0\ 1\ 0] \cup [0\ 1\ 0\ 1\ 0\ 0\ 0] \cup [0\ 0\ 0\ 1\ 0] \cup [0\ 1\ 0\ 0] \cup [1\ 1\ 0\ 0]\ .$$

Now Y o M is not defined and M is a mixture of both square and rectangular matrices we have to find only the special transpose of the model and determine

$$Y \text{ o}_m M = Y_1 \text{ o } M_1 \cup Y_2 \text{ o } M_2 \cup Y_3 \text{ o } M_3 \cup Y_4 \text{ o } M_4^t \cup Y_5 \text{ o } M_5^t \cup Y_6 \text{ o } M_6^t$$

$$= [2\ 0\ 0\ 1\ 2\ 1\ 0\ 0] \cup [1\ 0\ 1\ 1\ 1\ 0] \cup [1\ 0\ 0\ 1\ 0\ 0\ 0] \cup [0\ 0\ 0\ 0\ 0\ 0\ 0\ 1] \cup [0\ 1\ 1\ 0\ 1\ 0] \cup [0\ 0\ 0\ 2\ 0\ 0\ 0]$$
$$= Z'_1 \cup Z'_2 \cup Z'_3 \cup Z'_4 \cup Z'_5 \cup Z'_6$$

$= Z'$

after thresholding $Z'$ we get

$$Z = Z_1 \cup Z_2 \cup \ldots \cup Z_6$$
$$= [1\ 0\ 0\ 1\ 1\ 1\ 0\ 0] \cup [1\ 0\ 1\ 1\ 1\ 0] \cup [1\ 1\ 0\ 1\ 0\ 0\ 0] \cup [0\ 0\ 0\ 0\ 0\ 0\ 0\ 1] \cup [0\ 1\ 1\ 0\ 1\ 0] \cup [0\ 0\ 0\ 1\ 0\ 0\ 0] .$$

Now using the special dynamical system M we find

$$Z \text{ o } M = Z_1 \text{ o } M_1 \cup Z_2 \text{ o } M_2 \cup Z_3 \text{ o } M_3 \cup Z_4 \text{ o } M_4 \cup Z_5 \text{ o } M_5 \cup Z_6 \text{ o } M_6$$
$$= [3\ 0\ 0\ 1\ 2\ 2\ 0\ 0] \cup [2\ 0\ 1\ 1\ 2\ 0] \cup [1\ 0\ 1\ 1\ 0\ I\ 1] \cup [0\ 0\ 0\ 1\ 0] \cup [0\ 3\ 0\ \ 0] \cup [1\ 1\ 0\ 0]$$
$$= T'_1 \cup T'_2 \cup \ldots \cup T'_6.$$

After updating and thresholding T' we get

$$T = T_1 \cup T_2 \cup \ldots \cup T_6$$
$$= [1\ 0\ 0\ 1\ 1\ 1\ 0\ 0] \cup [1\ 0\ 1\ 1\ 1\ 0] \cup [1\ 1\ 1\ 1\ 0\ I\ 0] \cup [0\ 0\ 0\ 1\ 0] \cup [0\ 1\ 0\ 0] \cup [1\ 1\ 0\ 0].$$

We see $T_1$, $T_2$, $T_4$, $T_5$ and $T_6$ remain as fixed points. We need to work only with $M_3$ keeping this in mind we just work with $M_3$ i.e., finding $T_3 \circ M_3$ one gets the special limit cycle or a special fixed point.

Now having seen the working procedure of the model we now define yet another new special fuzzy and special neutron neutrosophic model and generalize them in the following section.

## 2.5 Special FRE models and NRE models

In this section for the first time we introduce the notion of Special Fuzzy relational equations and their generalizations. Any Fuzzy relational equations (FRE) model has only the capacity to tackle a single experts opinion on any problem at a time. A brief description of this model has been described in page 180-6 of this book. For more about FRE please refer [106, 232]. The FRE model makes use of min max or max min operators. This new model which we will be calling as special Fuzzy Relational Equations Model can deal with any desired number of experts at a time. This is the first time such types of fuzzy models have been introduced.

In this section we introduce 5 types of Special fuzzy relational equations models.

**DEFINITION 2.5.1:** *Suppose we have a problem P at hand and at least t experts wish to work with the same problem using the same set of attributes both for the domain space as well as range space then how to model this problem with multi expert models so that the views or solution is given simultaneously and collectively and not individually and separately. This model we would be defining will be known as Special Fuzzy Relational Equations (SFRE) model.*

*Suppose we have some t experts giving their opinion with the same set of m domain and n-range attributes. Let $T_i$ denote the FRE model given by the $i^{th}$ experts using these m domain attributes and n-range attributes, then $T_i$ is a $m \times n$ fuzzy membership matrix i.e., entries of $T_i$ are from the unit interval [0, 1]. Let $T_1, T_2, \ldots, T_t$ be the FRE matrices given by the t experts. Let $T = T_1 \cup T_2 \cup \ldots \cup T_t$ then T would be known as the Special Fuzzy Relational Equation(SFRE) model i.e., all the t experts use the same set of m row attributes and same set of n column attributes(m need not always be equal to n).*

*Thus the membership matrix of the $i^{th}$ expert is a $m \times n$ fuzzy matrix. The special fuzzy rectangular matrix associated with the SFRE will be known as the special fuzzy membership function matrix of the SFRE.*

*For more about FREs please refer [106, 232]. Depending on the type of FRE and the problem we can adopt any tool say neural net works or matrix equation or any suitable method to solve the problem.*

*Now if in the definition 2.5.1 the t experts choose to work with different sets of column attributes say $m_i$ column attributes is chosen by the $i^{th}$ expert and different sets of row attributes say $n_i$ attributes then the membership matrix would be a $m_i \times n_i$ fuzzy matrix, true for i =1, 2, …, t. Now the special membership* *matrix associated with this new model will be a special fuzzy mixed rectangular matrix and the associated special membership matrix will be given by $M_1 \cup M_2 \cup \ldots \cup M_t$ where each $m_i \times n_i$ fuzzy rectangular matrix. We call this model as Special Mixed Fuzzy Relational Equations(SMFRE) model. Depending on the type of FRE and the problem one can work for the solution or the solution close to the expected solution.*

The next special model is using the Neutrosophic Relational Equations (NRE).

**DEFINITION 2.5.2:** *Suppose we have p experts interested in using the NRE model for the problem P(say). Further suppose all of them agree upon to work with the same set of m row attributes and the same set of n column attributes. Let the neutrosophic membership matrix given by the $i^{th}$ expert be a $m \times n$ neutrosophic matrix say $M_i$; $i = 1, 2, \ldots, p$. Then we have $M = M_1 \cup M_2 \cup \ldots \cup M_p$ where M is a special fuzzy neutrosophic rectangular matrix.*

*The model associated with this M will be known as the Special Neutrosophic Relational Equations(SNRE) model. In other words SNRE model is just the SFRE model in which instead of working with the fuzzy membership matrices we have to work with the fuzzy neutrosophic membership matrices. So by replacing in SFRE model all fuzzy membership matrices by the appropriate fuzzy neutrosophic membership matrices we get the SNRE model. Now as in case of SMFRE if we replace each of the $m_i \times n_i$ fuzzy membership matrix by a fuzzy neutrosophic membership matrix we would get the new Special Mixed Neutrosophic Relational Equations Model i.e., the SMNRE model.*

**DEFINITION 2.5.3:** *Now suppose we have a set of m experts working on a problem P where $m_1$ of them expect to work using the SFRE model put together or equivalently the FRE model individually. A set $m_2$ experts choose to work using SMFRE model for the same problem P, $m_3$ of the experts choose to use NRE model i.e., they collectively work using the SNRE model* *and the rest of the experts collectively use the SMNRE model thus all the m experts works either on NRE and FRE models.*

*We call this model as the special mixed Fuzzy Neutrosophic Relational Equations model (SMFNRE-model). This model is a combination of all the four special models viz SMFRE, SFRE, SNRE and SMNRE.*

The main advantages of this model is:

1. This is a new multi expert special model which can simultaneously cater to any number of experts.
2. This model allows them to work with FRE and NRE.
3. This model enables them to use any desired number of row attributes and column attributes.

Since the very choice of the type of FRE or NRE is highly dependent on the problem and the expected solution for the instance it can be a multi objective optimization problem with FRE (NRE) constraints or it may be a non linear optimization problem with FRE(NRE) constraints or by using max min or min max operations of FRE or use of neutral net work and so on. For several types refer[106, 232].

Lastly we define yet a very generalized special model for solving problems of the following type. The problem can be thought of as the multiexpert multi model problem for we not only have multi experts but we also have the experts wishing to work with very many models. Then can we have a single dynamical system to tackle the situation.

The answer is yes and we define a new Special Fuzzy-Neutrosophic-Cognitive-Relational- Equations S-F-N-C-R-E model.

**DEFINITION 2.5.4:** *Suppose we have some m experts working on a problem P. There are six sets of experts, first set of experts wish to work with FCM, each one with the attributes of their own choice, second set of experts wish to work with NCM each one with the set of attributes of their own choice. The third set of experts wants to work with FRM with a set of attributes different for each expert. The fourth set of experts want to work* *with NRM with a set of attributes varying from expert to expert. The fifth set of experts want to work using FRE and with their own set of rows and columns attributes give their fuzzy membership matrix.*

*The last set of experts work using NRE. Now a model which can cater to the needs of all these 6 sets of experts with different types of models is termed as Super Special Hexagonal Model denoted by SSHM = SMFCM ∪ SMNCM ∪ SMFRE ∪ SMNRE ∪ SMFRM ∪ SMNRM where '∪' is just notationally describing all the six models are present and the union of the six model Special fuzzy and Special fuzzy neutrosophic mixed square/rectangular matrices will be the dynamical system associated with this super model.*

Further the special operation described in chapter one will give the desired method of working.

*Note:* When we use only 5 models and five sets of experts we call it as Super Special Quintuple model, if we use only four models with four sets of experts we call it as Super Special Quadruple model, with three models and three sets of experts will be known Super Special Triple model.

The main advantages of this model over the other models are as follows :

a. This multiexpert model can simultaneously deal with both multi experts and multi models.
b. With the advent of computer finding a program is not very difficult for these SSHM models. When programmed finding the solution is very easy and simple with any number of experts and any number of models.
c. Simultaneous solution is possible using special operations. This model can be both fuzzy as well as neutrosophic simultaneously.
d. The special fuzzy operation mentioned in – helps in finding the special resultant vector given any special input vector.

## 2.6 Some Programming Problems for Computer Experts

In this section we suggest some simple programming problems for a interested computer scientist or a programmer.

a. Give a C or C++ program to find the special hidden pattern of the following special models.

i. SFCM
ii. SMFCM
iii. SNCM
iv. SMNCM
v. SFRM
vi. SMFRM
vii. SNRM
viii. SMNRM
ix. SMFNRM
x. SMFNCRM
xi. SSHM

*Hint*: In case of single experts we have both Java and C program for FCM and FRM [228, 231].

b. Find a program for calculating the min max or max min operation for the given fuzzy row vector X = [$x_1$ $x_2$ … $x_n$]; $x_i \in [0, 1]$ and A = ($a_{ij}$); A a n × m fuzzy matrix with $a_{ij} \in [0, 1]$; $1 \leq i \leq n$ and $1 \leq i \leq m$; i.e., finding

max min {X,A}

$$= \max\min\left\{ \left[x_1 x_2 ... x_n\right], \begin{bmatrix} a_{11} & a_{12} & \cdots & a_{1m} \\ a_{21} & a_{22} & \cdots & a_{2m} \\ \vdots & \vdots & & \vdots \\ a_{n1} & a_{n2} & \cdots & a_{nm} \end{bmatrix} \right\}$$

$$= \Big[\max\left\{\min\left(x_1, a_{11}\right), \min\left(x_2, a_{22}\right), \ldots, \min\left(x_n, a_{n1}\right)\right\}$$

$$\max\left\{\min\left(x_1, a_{12}\right), \min\left(x_2, a_{22}\right), \ldots, \min\left(x_n, a_{n2}\right)\right\}$$

$$\max\{\min(x_1, a_{1m}), \min(x_2, a_{2m}), \ldots, \min(x_n, a_{nm})\}]$$

$$= \quad (b_1\ b_2\ \ldots b_m); \quad b_i \in [0,1].$$

c. Given M = $M_1 \cup M_2 \cup \ldots \cup M_n$, M is a super fuzzy neutrosophic mixed matrix. X = $X_1 \cup X_2 \cup \ldots \cup X_n$ is a super fuzzy mixed row vector. Using the special operator $\circ^{sp}$ find a C or C++ program for X $\circ^{sp}$ M = Y'; after thresholding Y' we get Y find Y $\circ^{sp}$ M = Z'; after updating and thresholding Z' we get Z find Z $\circ^{sp}$ M and so on until we get T a special fixed point or a special limit cycle.

*Hint*: Refer [228, 231] for a similar procedure.

# FURTHER READING

We have given a long list for further reading so that the interested reader can make use of them:

# INDEX

## A

## B

## C

## D

## E

## F

## H

## I

## L

## M

## N

## O

## P

## R

## S

## T

# U

# Z

# ABOUT THE AUTHORS

**Dr.W.B.Vasantha Kandasamy** is an Associate Professor in the Department of Mathematics, Indian Institute of Technology Madras, Chennai. In the past decade she has guided 11 Ph.D. scholars in the different fields of non-associative algebras, algebraic coding theory, transportation theory, fuzzy groups, and applications of fuzzy theory of the problems faced in chemical industries and cement industries. Currently, four Ph.D. scholars are working under her guidance.

She has to her credit 638 research papers. She has guided over 57 M.Sc. and M.Tech. projects. She has worked in collaboration projects with the Indian Space Research Organization and with the Tamil Nadu State AIDS Control Society. This is her 31st book.

On India's 60th Independence Day, Dr.Vasantha was conferred the Kalpana Chawla Award for Courage and Daring Enterprise by the State Government of Tamil Nadu in recognition of her sustained fight for social justice in the Indian Institute of Technology (IIT) Madras and for her contribution to mathematics. (The award, instituted in the memory of Indian-American astronaut Kalpana Chawla who died aboard Space Shuttle Columbia). The award carried a cash prize of five lakh rupees (the highest prize-money for any Indian award) and a gold medal.

She can be contacted at vasanthakandasamy@gmail.com
You can visit her on the web at: http://mat.iitm.ac.in/~wbv or: http://www.vasantha.net

---

**Dr. Florentin Smarandache** is an Associate Professor of Mathematics at the University of New Mexico in USA. He published over 75 books and 100 articles and notes in mathematics, physics, philosophy, psychology, literature, rebus. In mathematics his research is in number theory, non-Euclidean geometry, synthetic geometry, algebraic structures, statistics, neutrosophic logic and set (generalizations of fuzzy logic and set respectively), neutrosophic probability (generalization of classical and imprecise probability). Also, small contributions to nuclear and particle physics, information fusion, neutrosophy (a generalization of dialectics), law of sensations and stimuli, etc.
He can be contacted at smarand@unm.edu

---

**K. Ilanthenral** is the editor of The Maths Tiger, Quarterly Journal of Maths. She can be contacted at ilanthenral@gmail.com